\input amstex
\documentstyle{amsppt}
\loadmsbm

\nologo

\TagsOnRight

\NoBlackBoxes

\input epsf
\input supp-pdf

%
%
%
%

%
%
%

\hsize = 5.5 true in
\vsize = 8.5 true in

\define\Arg{\operatorname{Arg}}

\define\supp{\operatorname{supp}}

\define\diam{\operatorname{diam}}
\define\dist{\operatorname{dist}}

\def\floor{\mathbin{\hbox{\vrule height1.2ex width0.8pt depth0pt
        \kern-0.8pt \vrule height0.8pt width1.2ex depth0pt}}}

\font\bigletter=cmss10 scaled 1500

\font\letter=cmss10 

\font\normal=cmss10 scaled 700

\define\fraction{\operatorname{frac}}

\define\integer{\operatorname{int}}

\define\real{\operatorname{Re}}

\define\imag{\operatorname{Im}}

\define\res{\operatorname{res}}


\define\ima{\operatorname{\dot{\text{\rm \i\! \i}}}}

\define\ave{\text{\normal ave}}

\define\crit{\text{\normal crit}}

\define\sym{\text{\normal sym}}

\define\Box{\text{\normal B}}

\define\Mellin{\text{\letter M}}

\accentedsymbol\hati{\widehat{\bold i}}

\define\Vol{\text{\letter V}}

\define\Scale{\text{\letter S}}

\define\Min{\text{\normal M}}

\define\Haus{\text{\normal H}}

\define\Renyi{\text{\normal R}}

\define\inner{\text{\normal inner}}

\topmatter
\title
Multifractal tubes
\endtitle
\endtopmatter

\centerline{\smc Lars Olsen}
\centerline{Department of Mathematics}
\centerline{University of St\. Andrews}
\centerline{St\. Andrews, Fife KY16 9SS, Scotland}
\centerline{e-mail: {\tt lo\@st-and.ac.uk}}

\bigskip
\bigskip
\bigskip
\bigskip

\topmatter
\abstract
{
Tube formulas
refer
to the
study of 
volumes of $r$
neighbourhoods
of sets.
For sets satisfying
some (possible very weak) convexity conditions,
this 
has a long history.
However,
within the past 20 years
Lapidus has initiated and pioneered a
systematic study
of tube formulas for fractal sets.
Following this,
it is natural to ask to what extend it is possible to
develop
a theory of multifractal tube formulas for multifractal measures.
In this paper we
propose and develop
 a framework for
such a  theory.
Firstly,
we 
define
multifractal
tube formulas and, more generally,
multifractal
tube measures
for general multifractal measures.
Secondly,
we
introduce and develop two
approaches 
for analysing these concepts
for
self-similar multifractal measures, namely:

\medskip

\noindent
(1)
{\it Multifractal tubes of self-similar measures and renewal theory.}
Using techniques from renewal theory
we 
give a complete description
of
the 
asymptotic behaviour of the 
multifractal tube formulas and tube measures
 of self-similar measures satisfying the Open Set Condition.
  
\medskip

\noindent
(2)
{\it Multifractal tubes of self-similar measures and zeta-functions.}
Unfortunately,
renewal theory techniques
 do not yield
 \lq\lq explicit"
expressions
 for the functions
 describing the
 asymptotic behaviour of the 
multifractal tube formulas and tube measures
 of self-similar measures.
 This is clearly undesirable.
 For this reason,
we
introduce and develop a second framework
for studying
multifractal tube formulas 
 of self-similar measures.
 This approach is based on 
multifractal zeta-functions
and 
allow us
obtain
\lq\lq explicit"
expressions 
for the
multifractal tube formulas
of self-similar measures, namely,
using the Mellin transform and the residue theorem, we are able to express the
multifractal
tube formulas as sums involving the residues of the zeta-function.

}
\endabstract
\endtopmatter

\footnote""
{
\!\!\!\!\!\!\!\!
}

\rightheadtext{Multifractal tubes}

\rightheadtext{L\. Olsen}

\newpage

${}$

\bigskip
\bigskip
\bigskip
\bigskip
\bigskip
\bigskip
\bigskip
\bigskip

\centerline{\bigletter Contents}


\bigskip
\bigskip
\bigskip
\bigskip

\noindent
{\bigletter Part 1: Statements of Results}
\roster
\item"1." Fractal tubes
\item"2." Multifractals
\item"3." Multifractal tubes
\item"4." Multifractal tube measures
\item"5." Symbolic multifractal tubes of self-similar measures:
Multifractal zeta-functions
and
explicit formulas
\endroster

\bigskip

\noindent
{\bigletter Part 2: Proofs of the Results from Section 3}
\roster
\item"6." Proving that
$\lambda_{q,m}(r)
 \le
 \sum_{
  |\bold i|=|\bold j|=m,
  \bold i\not=\bold j
  }
 Q_{\bold i,\bold j}^{q}(r)$
\item"7." Proving that
 $Q_{\bold i,\bold j}^{q}(r)
 \le
 \text{\rm constant}
 \,
 Z_{m}^{q}(r)$
\item"8." Proving that
 $Z_{m}^{q}(r)
 \le
 \text{\rm constant}
 \,
 r^{-\gamma(q)}$
\item"9." Proof of Theorem 3.3
\endroster

\bigskip

\noindent
{\bigletter Part 3: Proofs of the Results from Section 4}
\roster
\item"10." Analysis of $\Cal H_{\mu}^{q,\beta(q)}$
\item"11." Proof of Theorem 4.1
\item"12." Proof of Theorem 4.2
\endroster

\bigskip

\noindent
{\bigletter Part 4: Proofs of the Results from Section 5}
\roster
\item"13." Analysis of the poles of $\zeta_{\mu}^{q}$ 
\item"14." Proof of Theorem 5.4
\item"15." Proof of Theorem 5.5
\item"16." Proof of Theorem 5.7
\endroster

\newpage

${}$

\bigskip
\bigskip
\bigskip
\bigskip
\bigskip
\bigskip
\bigskip
\bigskip

\centerline{\bigletter Part 1:}

\bigskip

\centerline{\bigletter Statements of Results}


\bigskip
\bigskip
\bigskip
\bigskip

Tube formulas
refer
to the
study of 
volumes of $r$
neighbourhoods
of sets.
For sets satisfying
some (possible very weak) convexity conditions,
this 
has a long history going back to Steiner in the early 19'th century.
However,
within the past 20 years
Lapidus has initiated and pioneered a
systematic study
of tube formulas for fractal sets.
Following this line of investigation,
it is natural to ask to what extend it is possible to
develop
a theory of multifractal tube formulas for multifractal measures.
The purpose of this paper is 
to propose a framework for
developing such a  theory.
Firstly,
we 
define
multifractal
tube formulas and, more generally,
multifractal
tube measures
for general multifractal measures.
Secondly,
we
introduce and develop two
approaches 
for analysing these concepts
for
self-similar multifractal measures, namely: 

\medskip

\noindent
(1)
{\it Multifractal tubes of self-similar measures and renewal theory.}
Using techniques from renewal theory
we 
give a complete description
of
the 
asymptotic behaviour of the 
multifractal tube formulas and tube measures
 of self-similar measures satisfying the Open Set Condition.
 This 
 is
 presented in
 Section 3 (for tube formulas)
 and
 Section 4 (for tube measures).
 
%
 
\medskip

\noindent
(2)
{\it Multifractal tubes of self-similar measures and zeta-functions.}
While 
renewal theory techniques
 are powerful tools,
 they do not yield
 \lq\lq explicit"
expressions
 for the functions
 describing the
 asymptotic behaviour of the 
multifractal tube formulas and tube measures
 of self-similar measures.
 This is clearly undesirable.
 For this reason,
we
introduce and develop a second framework
for studying
multifractal tube formulas 
 of self-similar measures.
 This approach is based on 
multifractal zeta-functions
and 
allow us
obtain
\lq\lq explicit"
expressions 
for the
multifractal tube formulas
of self-similar measures, namely,
using the Mellin transform and the residue theorem, we are able to express the
multifractal
tube formulas as sums involving the residues of the zeta-function.
This is done in Section 5.

%

\bigskip
\bigskip

\heading{1. Fractal tubes}\endheading

Let $E$ be a subset of $\Bbb R^{d}$ and $r>0$. 
We will write
$B(E,r)$ for the open $r$ neighbourhood
of $E$, i\.e\.
 $$
 B(E,r)
 =
 \Big\{
 x\in\Bbb R^{d}
 \,\Big|\,
 \dist(x,E)<r
 \Big\}\,.
 \tag1.1
 $$
Intuitively
we think of the set $B(E,r)$ as consisting of the $E$
surrounded 
by a 
\lq\lq tube" of width $r$.
Our main interest is 
to compute the 
volume of the 
\lq\lq tube" of width $r$
surrounding $E$ or equivalently
computing the volume of the set $B(E,r)$
and then subtracting the volume of $E$.
To make this formal, 
we define the
Minkowski volume $V_{r}(E)$ of $E$ by
 $$
 V_{r}
 (E)
 =
 \frac{1}{r^{d}}
 \,
 \Cal L^{d}(B(E,r))\,;
 \tag1.2
  $$ 
here and below $\Cal L^{d}$ denotes the $d$-dimensional Lebesgue measure in 
$\Bbb R^{d}$
and  
the normalising factor
$\frac{1}{r^{d}}$ is included to make
the subsequent results
simpler - we note that different authors use different normalising factors.
Tube formulas refers to
formulas for computing the
Minkowski volume $V_{r}(E)$ as a function of the
width $r$ of the \lq\lq tube" surrounding $E$. 
In particular,
one is typically interested in the following two types
of results:
\roster
\item"$\bullet$"
Asymptotic behaviour:
finding a formula for the 
asymptotic behaviour of
$V_{r}(E)$ as $r\searrow 0$;
\item"$\bullet$"
Explicit formulas:
finding an explicit formula
for 
$V_{r}(E)$ valid for all sufficiently small $r$.
\endroster
For convex sets $E$,
this problem has a rich and fascinating history
starting with the work of Steiner in the early 19'th century.
Indeed,
Steiner
showed
that
if $C$ is a bounded
convex subset of $\Bbb R^{d}$,
then
there are 
constants
$\kappa^{0}(C),\kappa^{1}(C),\ldots,\kappa^{d}(C)$
such that
$$
 \Cal L^{d}(B(C,r))
 =
 \sum_{l}\kappa^{l}(C)\,r^{d-l}
 \tag1.3
 $$
for $r>0$.
The coefficients $\kappa^{l}(C)$ are called the Quermassintegrale
or mixed volumes, 
and the polynomial $\sum_{l}\kappa^{l}(C)\,r^{d-l}$ is called the Steiner polynomial.
We also note that
the coefficients
have clear geometric interpretations. 
For example,
$\kappa^{d}(C)$ equals the volume of $C$ and $\kappa^{d-1}(C)$ 
is equal to the surface area of $C$.
Steiner's formula has subsequently been
extended  to 
more general classes of sets.
For example, in the late 1930's Weyl proved  
that a similar 
result
holds for
compact oriented $d$-dimensional Riemannian manifolds $C$
(with or without boundary) isometrically embedded into Euclidean space.
This theory reached its mature form in the 1960's
where
Federer [Fed1,Fed2]
unified the tube formulas of Steiner for convex bodies
and of Weyl for smooth submanifolds, as described
in [BeGo,Gray,We],
and extended these results to sets of positive reach.
Federer's tube formula has since been extended in various
directions by a number of researchers in integral geometry and geometric measure
theory, including 
[Fu1,Fu2,Schn1,Schn2,St,Z\"a1,Z\"a2]
and most recently
(and most generally) in [HuLaWe]. 
The books [Gray,Mo,Schn2] contain extensive
endnotes with further information and many other references.
While the above references 
investigate tube formulas for sets that 
satisfy some (possibly very weak) convexity and/or smoothness conditions,
very recently 
there has been significant interest in developing 
a theory of tube formulas
for fractal sets 
and a number of exciting
works
have appeared.
Indeed,
in the early 1990's 
Lapidus introduced the notion of \lq\lq complex dimensions"
and has during the past 20 years
very successfully pioneered the use of \lq\lq complex dimensions"
to
obtain explicit tube formulas for 
certain classes of fractal subsets of (mainly) the real line;
this exciting theory is described in detail in 
Lapidus \& van Frankenhuysen's
intriguing 
books
[Lap-vF1,Lap-vF2].
In a parallel development,
 and
building on earlier work by Lalley [Lal1,Lal2,Lal3] 
and Gatzouras [Ga] (see also [Fa3]),
Winter [Wi] has initiated the systematic study
of
curvatures of fractal sets 
and applied this theory to study
the asymptotic behaviour of the Minkowski volume 
$V_{r}(E)$
of fractal sets $E$
using methods from renewal theory.
The work in this paper may be viewed
as a natural higher dimensional
multifractal development
of this line of research.

The
Minkowski volume 
$V_{r}(E)$ is 
closely related to 
various notions from Fractal Geometry.
Indeed, using the Minkowski volume $V_{r}(E)$, 
we define the lower and upper Minkowski dimension
of $E$ 
by
 $$
 \aligned
 \underline\dim_{\Min}(E)
&=
 \,
 \liminf_{r\searrow0}
 \,
 \frac{\log V_{r}(E)}{-\log r}\,,\\
 \overline\dim_{\Min}(E)
&=
 \limsup_{r\searrow0}
 \frac{\log V_{r}(E)}{-\log r}\,.
 \endaligned
 \tag1.4
 $$
The link with Fractal Geometry is now explained as follows.
Namely, box dimensions play an important role in Fractal Geometry
and
it is not difficult to see
that
the lower Minkowski dimension equals the lower box dimension 
and that
the upper Minkowski dimension equals the upper box dimension;
for the definition of the box dimensions the reader is referred to Falconer's text book [Fa1].

It is clearly 
also of interest to 
analyse the behaviour
of the Minkowski volume $V_{r}(E)$
itself as $r\searrow0$.
Indeed, if, for example, 
$a_{1},\ldots,a_{d},b_{1},\ldots,b_{d}$ are real
numbers
with $a_{i}\le b_{i}$ for all $i$,
and
$U$
denotes the  rectangle $[a_{1},b_{1}]\times\cdots\times[a_{d},b_{d}]$ in $\Bbb R^{d}$, 
then
it is clear that
$ \frac{1}{r^{-d}}
 V_{r}(U)
 =
 ((b_{1}+r)-(a_{1}-r))\cdots((b_{d}+r)-(a_{d}-r))
\to
(b_{1}-a_{1})\cdots(b_{d}-a_{d})
=
\Cal L^{d}(U)$.
This suggests that if $t$ is a real number, then the limit
$\lim_{r\searrow0}
 \frac{1}{r^{-t}}
 V_{r}(E)$ (if it exists)
may be interpreted as  the $t$-dimensional volume of $E$. 
 Motivated by this, for 
 a real number $t$,
 we define the 
lower and upper $t$-dimensional Minkowski
content of $E$ by
 $$
 \aligned
 \underline M^{t}(E)
&=
 \,
 \liminf_{r\searrow0}
 \,
 \frac{1}{r^{-t}}
 \,
 V_{r}(E)\,,\\
 \overline M^{t}(E)
&=
 \limsup_{r\searrow0}
 \frac{1}{r^{-t}}
 \,
 V_{r}(E)\,.
 \endaligned
 \tag1.5
 $$
If 
$\underline M^{t}(E)
=
\overline M^{t}(E)$,
i\.e\. if the
limit
$ \lim_{r\searrow0}
 \,
 \frac{1}{r^{-t}}
 \,
 V_{r}(E)$
 exists, 
then we say the $E$ is $t$ Minkowski measurable,
and we 
denote the  
common value of
$\underline M^{t}(E)$ and
$\overline M^{t}(E)$
by
$M^{t}(E)$,
i\.e\. we write
 $$
 M^{t}(E)
 =
 \underline M^{t}(E)
 =
 \overline M^{t}(E)\,.
 \tag1.6
 $$
Of course,
a set $E$ may not be Minkowski measurable, i\.e\.
the limit
$ \lim_{r\searrow0}
 \,
 \frac{1}{r^{-t}}
 \,
 V_{r}(E)$
may not exist.
In this case it is natural to study
the limiting behaviour of 
suitably defined
\lq\lq averages" of
$ \frac{1}{r^{-t}}V_{r}(E)$.
We therefore define the lower and upper
average
$t$-dimensional Minkowski content of $E$ by
 $$
 \aligned
 \underline M_{\ave}^{t}(E)
&=
 \,
 \liminf_{r\searrow0}
 \frac{1}{-\log r}
 \int_{r}^{1}
 \frac{1}{s^{-t}}
 \,
 V_{s}(E)
 \,
 \frac{ds}{s}
 \,,\\
 \overline M_{\ave}^{t}(E)
&=
 \,
 \liminf_{r\searrow0}
 \frac{1}{-\log r}
 \int_{r}^{1}
 \frac{1}{s^{-t}}
 \,
 V_{s}(E)
 \,
 \frac{ds}{s}\,.
 \endaligned
 \tag1.7
 $$
If 
$\underline M_{\ave}^{t}(E)
=
\overline M_{\ave}^{t}(E)$,
i\.e\. if the
limit
$ \lim_{r\searrow0}
 \frac{1}{-\log r}
 \int_{r}^{1}
 \frac{1}{s^{-t}}
 \,
 V_{s}(E)
 \,
 \frac{ds}{s}$
 exists, 
then we say the $E$ is $t$ averagely Minkowski measurable,
and we 
denote the  
common value of
$\underline M_{\ave}^{t}(E)$ and
$\overline M_{\ave}^{t}(E)$
by
$M_{\ave}^{t}(E)$,
i\.e\. we write
 $$
 M_{\ave}^{t}(E)
 =
 \underline M_{\ave}^{t}(E)
 =
 \overline M_{\ave}^{t}(E)\,.
 \tag1.8
 $$
While 
the Minkowski dimensions in  many cases can be computed rigorously
relatively easy,
it is a notoriously difficult problem to compute the 
Minkowski content.
In fact, it is only within the past 15 years
that the Minkowski content
of non-trivial examples have been computed.
Indeed, using techniques from complex analysis,
Lapidus and collaborators [Lap1-vF1,Lap1-vF2] 
have computed the Minkowski content of certain self-similar subsets of the real line,
and 
using ideas from the theory of Mercerian 
theorems,
Falconer [Fa3]
have obtained 
similar results.

It is our intention to extend
the notion of 
Minkowski volume 
$V_{r}(E)$ to multifractals
and
investigate
the
asymptotic behaviour
of the corresponding multifractal Minkowski volume 
as $r\searrow 0$
for self-similar multifractals.
In order to motivate our definitions we will now explain
what the term
\lq\lq multifractal analysis" covers.



\bigskip
\bigskip

\heading{2. Multifractals.}\endheading

{\bf 2.1. Multifractal spectra.}
Distributions with widely varying intensity
occur often in the physical sciences, e\.g\.
the spatial-temporal distribution of rainfall,
the spatial distribution of oil and gas in the underground,
the distribution of galaxies in the universe,
the dissipation of energy in a highly turbulent fluid flow, or
the occupation measure on strange attractors.
Such distributions are called multifractals and have recently been
the focus of much attention in the physics literature.

%

%

%

%


%

%

%
For a Borel measure $\mu$ on a $\Bbb R^{d}$ and
a real number $\alpha$,
let us consider
the set $\Delta_{\mu}(\alpha)$ of
those points
$x$ in $\Bbb R^{d}$ for which the measure
$\mu(B(x,r))$ of the ball
$B(x,r)$ with center $x$ and radius $r$ behaves like
$r^{\alpha}$ for small $r$,
i\.e\. the set
 $$
   \Delta_{\mu}(\alpha)
   =
   \Bigg\{x\in\Bbb R^{d}
   \,\Bigg|\,
   \lim_{r\searrow 0}
   \frac{\log\mu(B(x,r))}{\log r}
   =\alpha
   \Bigg\}\,.
 $$
If the intensity of the measure $\mu$ varies very widely, it may
happen that the sets $\Delta_{\mu}(\alpha)$
display a
fractal-like character for a range of values of $\alpha$. If this is the
case, then the measure is called a multifractal measure
or simply a multifractal, and
it is natural to study
the sizes of the sets $\Delta_{\mu}(\alpha)$ as $\alpha$
varies.
We do this by
studying the
Hausdorff dimension of 
$\Delta_{\mu}(\alpha)$. 
More precisely,
we define the
multifractal spectrum $f_{\mu}:\Bbb R\to\Bbb R$
of $\mu$ by
 $$
 f_{\mu}(\alpha)
 =
 \dim\Delta_{\mu}(\alpha)\,,
 \tag2.1
$$
of the sets $\Delta_{\mu}(\alpha)$ as a function of $\alpha$
where $\dim$ denotes the Hausdorff dimension.
The function in (2.1) and similar functions are generically known
as ``the multifractal spectrum of $\mu$'',
``the singularity spectrum
of $\mu$''
or  ``the spectrum of scaling indices'', 
and one of the main problems
in multifractal analysis is to study these and related functions.
The function $f_{\mu}(\alpha)$ was
first explicitly defined by the physicists Halsey et al. in 1986 in
their seminal paper [HaJeKaPrSh].
The concepts underlying the above mentioned multifractal
decompositions go back to two early papers by Mandelbrot [Man1,Man2]
from 1972 and 1974, respectively, where Mandelbrot suggested that
the bulk of intermittent dissipation of energy in a highly turbulent
fluid flow occurs over a set of fractal dimension. The ideas
introduced in [Man1,Man2] were taken up by Frisch \& Parisi [FrPa] in 1985
and finally by
Halsey et al\.
[HaJeKaPrSh]
in 1986.
Of course, for many measures the limit
$\lim_{r\searrow 0}\frac{\log\mu(B(x,r))}{\log r}$ may fail to exist
for all or many $x$, in which case we need to work with lower or 
upper limits as $r$ tends to $0$
and (perhaps) replace
\lq\lq$=\alpha$" in the definition of $\Delta_{\mu}(\alpha)$
with
\lq\lq$\le\alpha$" or \lq\lq$\ge\alpha$".

\bigskip

{\bf 2.2. Renyi dimensions.}
Based on a remarkable insight
together with a clever heuristic argument
Halsey et al\. [HaJeKaPrSh]
suggested
that the multifractal spectrum
$f_{\mu}(\alpha)$ can be computed
using a principle known
as the Multifractal Formalism.
The
Multifractal Formalism
involves the so-called Renyi
dimensions  which we will now define.
Let $\mu$ be a Borel measure on $\Bbb R^{d}$.
For 
$q\in\Bbb R$ and $r>0$, we 
define 
the $q$-th  moment $I^{q}_{\mu,r}(E)$ of
a subset $E$ of $\Bbb R^{d}$ with respect to $\mu$ at scale $r$ 
by
 $$
 I_{\mu,r}^{q}(E)
 =
 \,
 \int
 \limits_{E}
 \mu(B(x,r))^{q-1}
 \,
 d\mu(x)\,.
  \tag2.2
 $$  
Next, the lower and upper Renyi dimensions 
of $E$ with respect to $\mu$ are defined by 
 $$
 \aligned
 \underline\dim_{\Renyi,\mu}^{q}(E)
&=
 \liminf_{r\searrow 0}\frac{\log I_{\mu,r}^{q}(E)}{-\log r}
 \,,\\
  \overline\dim_{\Renyi,\mu}^{q}(E)
 &= 
 \limsup_{r\searrow 0}\frac{\log I_{\mu,r}^{q}(E)}{-\log r}
 \,.
 \endaligned
 \tag2.3
 $$
In particular, the
Renyi dimensions of the support $\supp\mu$  of $\mu$
play an important role in the statement of the 
Multifractal Formalism.
For this reason it is useful 
to denote these dimensions by a separate notation,
and 
we therefore define the 
lower and upper
Renyi spectra
$\underline\tau_{\mu}(q),\overline\tau_{\mu}(q):\Bbb R\to[-\infty,\infty]$
 of $\mu$ by
 $$
 \aligned
 \underline\tau_{\mu}(q)
&=
 \underline\dim_{\Renyi,\mu}^{q}(\supp\mu)\\
&= 
 \,\liminf_{r\searrow 0}\,\frac{\log I_{\mu,r}^{q}(\supp\mu)}{-\log r}
 \,,\\
  \overline\tau_{\mu}(q)
 &= 
 \overline\dim_{\Renyi,\mu}^{q}(\supp\mu)\\
&=  
 \limsup_{r\searrow 0}\frac{\log I_{\mu,r}^{q}(\supp\mu)}{-\log r}
 \,.
 \endaligned
 \tag2.4
 $$

\bigskip

{\bf 2.3. The Multifractal Formalism.}
We can now state the Multifractal Formalism.
Loosely speaking the 
Multifractal Formalism
says that
the multifractal spectrum $f_{\mu}$
and the 
Renyi dimensions $\underline\tau_{\mu}(q)$ and $\overline\tau_{\mu}(q)$
carry the same information.
More precisely,
the 
Multifractal Formalism
says that
the
multifractal spectrum 
equals the 
Legendre transform of the 
Renyi dimensions.
Before stating this formally,
we remind the reader
that if
$\varphi:\Bbb R\to\Bbb R$ is a real valued function,
then the Legendre transform $\varphi^{*}:\Bbb R\to[-\infty,\infty]$
of $\varphi$ is defined by
 $$
 \varphi^{*}(x)
 =
 \inf_{y}(xy+\varphi(y))\,.
 \tag2.5
 $$
We can state the 
Multifractal Formalism.

\bigskip

\proclaim{The Multifractal Formalism -- A Physics Folklore Theorem}
The multifractal spectrum
$f_{\mu}$ 
of $\mu$
equals the 
Legendre transforms,
$\underline\tau_{\mu}^{*}$ and $\overline\tau_{\mu}^{*}$,
of 
the Renyi dimensions, i\.e\.
 $$
 f_{\mu}(\alpha)
 =
 \underline\tau_{\mu}^{*}(\alpha)
 =
 \overline\tau_{\mu}^{*}(\alpha)
 $$
 for all $\alpha\ge 0$.
\endproclaim

\bigskip

The Multifractal Formalism is a truly remarkable result:
it says that the locally defined multifractal spectrum 
$f_{\mu}$ can be computed in terms of the Legendre 
transforms of the globally defined moment scaling functions
$\underline \tau_{\mu}^{*}$ and
$\overline\tau_{\mu}^{*}$.
There is apriori no reason to expect that the Legendre transforms of the 
moment scaling function 
$\underline\tau_{\mu}^{*}$ and
$\overline\tau_{\mu}^{*}$ should provide any information about the 
fractal dimension of the set of points $x$ such that
$\mu(B(x,r))\approx r^{\alpha}$ for $r\approx 0$.
In some sense the 
Multifractal Formalism
is a genuine mystery.

During the past 20 years
there has been an enormous interest 
in
verifying the
Multifractal Formalism 
and
computing the multifractal spectra of measures
in 
the mathematical literature.
In the mid 1990's
Cawley \& Mauldin [CaMa] and Arbeiter \& Patzschke [ArPa]
verified the Multifractal Formalism for self-similar measures 
satisfying the OSC,
and within the last
20 years the multifractal spectra of various classes of measures
in Euclidean space $\Bbb R^{d}$ 
exhibiting some degree of self-similarity have been computed 
rigorously, cf\. 
the textbooks [Fa2,Pes2]
and the references therein.



\bigskip
\bigskip

\heading{3. Multifractal tubes}\endheading

{\bf 3.1. Multifractal tubes of general measures.}
Motivated by
Lapidus \& van Frankenhuysen
investigations [Lap-vF1,Lap-vF2]
of tube formulas
for fractal sets,
it it natural to seek
to develop a theory of multifractal tube formulas
for multifractal measures.
We will now present a framework for developing such a theory
and
as an application illustrating these ideas
we 
give a complete description
of
the 
asymptotic behaviour of the 
multifractal tube formulas for self-similar measures satisfying the Open Set Condition.

Multufractal tube formulas are defined as follows.
First note that if 
$r>0$ and $E$ is a subset of $\Bbb R^{d}$, then
the
Minkowski volume $V_{r}(E)$ is given by
 $$
 \align
 V_{r}(E)
&=
 \frac{1}{r^{d}}
 \,
 \Cal L^{d}(B(E,r))\\
&=
 \frac{1}{r^{d}}
 \int\limits_{B(E,r)}
 \,
 d\Cal L^{d}(x)\,,
 \endalign
 $$
where we have rewritten the
Lebesgue measure 
$\Cal L^{d}(B(E,r))$ of $B(E,r)$
as the integral
$ \int_{B(E,r)}
 d\Cal L^{d}(x)$.
Motivated by 
the Renyi dimensions (i\.e\. (2.2) and (2.3))
and the above expression for  
$V_{r}(E)$, we define the
multifractal Minkowski volume as follows.
Namely,
let $r>0$ and $E$ be a subset of $\Bbb R^{d}$.
For  real number $q$
and
a Borel measure $\mu$ on $\Bbb R^{d}$,
 we define the 
 multifractal $q$ Minkowski  volume
$V_{\mu,r}^{q}(E)$ of $E$ with respect to the measure $\mu$ by
 $$
 V_{\mu,r}^{q}
 (E)
 =
 \frac{1}{r^{d}}
 \int\limits_{B(E,r)}
 \mu(B(x,r))^{q}
 \,
 d\Cal L^{d}(x)\,.
 \tag3.1
 $$
Note, that if $q=0$, then the $q$ multifractal Minkowski volume
$V_{\mu,r}^{q}(E)$ reduces to the usual 
Minkowski volume, i\.e
 $$
 V_{\mu,r}^{0}
 (E)
 =
 V_{r}
 (E)\,.
 $$
The
importance of the Renyi dimensions in multifractal analysis together with the 
formal 
resemblance 
between the 
multifractal Minkowski volume 
$V_{\mu,r}^{q}(E)$
and the 
moments
$I_{\mu,r}^{q}(E)$ used in the definition
the Renyi dimensions
may be seen as a justification 
for
calling
the quantity
$V_{\mu,r}^{q}(E)$
for
the {\it multifractal} Minkowski volume;
a further justification for this terminology will be provided below.

Using the   multifractal Minkowski volume 
we can define 
multifractal Monkowski 
dimensions.
For  real number $q$
and
a Borel measure $\mu$
 on $\Bbb R^{d}$,
 we define the 
 lower and upper
 multifractal $q$ Minkowski  dimension of $E$,
by
 $$
 \aligned
 \underline\dim_{\Min,\mu}^{q}(E)
&=
 \,
 \liminf_{r\searrow0}
 \,
 \frac{\log V_{\mu,r}^{q}(E)}{-\log r}\,,\\
 \overline\dim_{\Min,\mu}^{q}(E)
&=
 \limsup_{r\searrow0}
 \frac{\log V_{\mu,r}^{q}(E)}{-\log r}\,.
 \endaligned
 \tag3.2
 $$
Again we note the close similarity between the 
multifractal Minkowski dimensions and the Renyi dimensions.
Indeed, the next proposition shows that
this similarity is not merely a formal resemblance.
In fact, for $q\ge 0$, the 
multifractal Minkowski dimensions and the Renyi dimensions
coincide.
This clearly provides  a  further justification 
for calling
the quantity
$V_{\mu,r}^{q}(E)$
for the {\it multifractal} Minkowski volume.

\bigskip

\proclaim{Proposition 3.1}
Let $\mu$ be a Borel measure on $\Bbb R^{d}$
and
$E\subseteq\Bbb R^{d}$.
If $q\ge 0$, then
 $$
 \align
 \underline\dim_{\Min,\mu}^{q}(E)
&=
 \underline\dim_{\Renyi,\mu}^{q}(E)\,,\\
 \overline\dim_{\Min,\mu}^{q}(E)
&=
  \overline\dim_{\Renyi,\mu}^{q}(E)\,.
 \endalign
 $$
In particular,  
if $q\ge 0$, then
 $$
 \align
 \underline\dim_{\Min,\mu}^{q}(\supp\mu)
&=
 \underline\tau_{\mu}(q)\,,\\
 \overline\dim_{\Min,\mu}^{q}(\supp\mu)
&=
  \overline\tau_{\mu}(q)\,.
 \endalign
 $$
\endproclaim
\noindent{\it  Proof}\newline
\noindent
This  follows easily from the definitions
and the proof is therefore omitted.
\hfill$\square$

\bigskip

Having defined multifractal Minkowski
dimensions, we 
also define multifractal Minkowski content
and average mutltifractal Minkowski content.
For real numbers $q$ and $t$,
we define the
lower and upper $(q,t)$-dimensional multifractal Minkowski
content of $E$ with respect to $\mu$ by
 $$
 \aligned
 \underline M_{\mu}^{q,t}(E)
&=
 \,
 \liminf_{r\searrow0}
 \,
 \frac{1}{r^{-t}}
 \,
 V_{\mu,r}^{q}(E)\,,\\
 \overline M_{\mu}^{q,t}(E)
&=
 \limsup_{r\searrow0}
 \frac{1}{r^{-t}}
 \,
 V_{\mu,r}^{q}(E)\,.
 \endaligned
 \tag3.3
 $$
If 
$\underline M_{\mu}^{q,t}(E)
=
\overline M_{\mu}^{q,t}(E)$,
i\.e\. if the
limit
$ \lim_{r\searrow0}
 \,
 \frac{1}{r^{-t}}
 \,
 V_{\mu,r}^{q}(E)$
 exists, 
then we say the $E$ is $(q,t)$ multifractal Minkowski measurable
with respect to $\mu$,
and we 
denote the  
common value of
$\underline M_{\mu}^{q,t}(E)$ and
$\overline M_{\mu}^{q,t}(E)$
by
$M_{\mu}^{q,t}(E)$,
i\.e\. we write
 $$
 M_{\mu}^{q,t}(E)
 =
 \underline M_{\mu}^{q,t}(E)
 =
 \overline M_{\mu}^{q,t}(E)\,.
 \tag3.4
 $$
Of course,
sets may not be multifractal Minkowski measurable, and it is therefore useful
to introduce 
a suitable 
averaging procedure 
when computing the 
multifractal Minkowski content.
Motivated by this we define 
 the
lower and upper $(q,t)$-dimensional average multifractal Minkowski
content of $E$ with respect to $\mu$ by
 $$
 \aligned
 \underline M_{\mu,\ave}^{q,t}(E)
&=
 \,
 \liminf_{r\searrow0}
 \,
 \frac{1}{-\log r}
 \int_{r}^{1}
 \frac{1}{s^{-t}}
 \,
 V_{\mu,s}^{q}(E)
 \,
 \frac{ds}{s}
 \,,\\
 \overline M_{\mu,\ave}^{q,t}(E)
&=
 \limsup_{r\searrow0}
 \frac{1}{-\log r}
 \int_{r}^{1}
 \frac{1}{s^{-t}}
 \,
 V_{\mu,s}^{q}(E)
 \,
 \frac{ds}{s}\,.
 \endaligned
 \tag3.5
 $$
If 
$\underline M_{\mu,\ave}^{q,t}(E)
=
\overline M_{\mu,\ave}^{q,t}(E)$,
i\.e\. if the
limit
$ \lim_{r\searrow0}
 \frac{1}{-\log r}
 \int_{r}^{1}
 \frac{1}{s^{-t}}
 \,
 V_{\mu,s}^{q}(E)
 \,
 \frac{ds}{s}$
 exists, 
then we say the $E$ is $(q,t)$ averagely multifractal Minkowski measurable
with respect to $\mu$,
and we 
denote the  
common value of
$\underline M_{\mu,\ave}^{q,t}(E)$ and
$\overline M_{\mu,\ave}^{q,t}(E)$
by
$M_{\mu,\ave}^{q,t}(E)$,
i\.e\. we write
 $$
 M_{\mu,\ave}^{q,t}(E)
 =
 \underline M_{\mu,\ave}^{q,t}(E)
 =
 \overline M_{\mu,\ave}^{q,t}(E)\,.
 \tag3.6
 $$
Note that definitions (3.3), (3.4), (3.5) and (3.6)
reduce to
(1.5), (1.6), (1.7) and (1.8),
respectively,
for $q=0$.
We will now give a complete description of
the multifractal Minkowski 
contents for self-similar measures $\mu$.
%

\bigskip

{\bf 3.2. Multifractal tubes of self-similar measures.}
We
will now compute the 
multifractal Minkowski 
content of self-similar measures.
We begin by recalling the definition of a self-similar measure.
Let $S_{i}:\Bbb R^{d}\to\Bbb R^{d}$ for $i=1,\dots,N$ be contracting 
similarities and 
let $(p_{1},\dots,p_{N})$ be a probability vector.
We denote the Lipschitz constant of $S_{i}$ by
$r_{i}\in(0,1)$.
The self-similar set $K$ and the self-similar measure $\mu$ associated
with 
the list 
$(S_{1},\dots,S_{N},p_{1},\dots,p_{N})$
are defined as follows.
Namely,
$K$
is the  unique
non-empty compact 
subset of $\Bbb R^{d}$ such that
 $$
 K
 =
 \bigcup_{i}S_{i}(K)\,,
 \tag3.7
 $$
and $\mu$
the unique 
Borel probability measure on $\Bbb R^{d}$
such that
 $$
 \mu
 =
 \sum_{i}p_{i}\mu\circ S_{i}^{-1}\,,
 \tag3.8
 $$
cf\. [Hu].  
We note that it is
well-known that
$\supp\mu = K$.

We will frequently assume that the list
$(S_{1},\dots,S_{N})$ satisfies certain \lq\lq disjointness" 
conditions, viz\. the Open Set Condition (OSC) or the 
Strong Separation Condition (SSC) defined below.

\medskip

{\it The Open Set Condition}:
There exists an open non-empty and bounded subset $U$ of $\Bbb R^{d}$ 
with
$\cup_{i}S_{i}U\subseteq U$ and
$S_{i}U\cap S_{j}U=\varnothing$ for all $i,j$ with $i\not=j$.

\medskip

{\it The Strong Separation Condition}:
There exists an open non-empty and bounded subset $U$ of $\Bbb R^{d}$ 
with
$\cup_{i}S_{i}U\subseteq U$ and
$\overline{S_{i}U}\cap\overline{S_{j}U}=\varnothing$ for all 
$i,j$ with $i\not=j$.

\medskip

Multifractal analysis of self-similar measures has attracted an enormous interest 
during the past 20 years.
For example,
using methods from ergodic theory,
Peres \& Solomyak [PeSo] have recently shown that for any 
self-similar measure $\mu$,
the 
Renyi dimensions always exists, 
i\.e\. the limit
$ \lim_{r\searrow 0}\frac{\log I^{q}_{\mu,r}(K)}{-\log r}$
always exists,
regardless of whether or not the OSC is satisfied provided $q\ge 0$.
If in addition the OSC is satisfied, an explicit expression
for the two limits
$\underline\tau_{\mu}(q)= \liminf_{r\searrow 0}\frac{\log I^{q}_{\mu,r}(K)}{-\log r}$
and
$\overline\tau_{\mu}(q)= \limsup_{r\searrow 0}\frac{\log I^{q}_{\mu,r}(K)}{-\log r}$
 can be obtained.
Indeed,
Arbeiter \& Patzschke [ArPa] and Cawley \& Mauldin [CaMa]
proved that if the OSC is satisfied, then
 $$
 \aligned
 \underline\tau_{\mu}(q)
&=
\,
 \liminf_{r\searrow 0}
 \,
 \frac{\log I_{r}^{q}(K)}{-\log r}\\
&=
 \beta(q)\,,\\
 \overline\tau_{\mu}(q)
&= 
 \limsup_{r\searrow 0}
 \frac{\log I_{r}^{q}(K)}{-\log r}\\
&=
 \beta(q)\,,
 \endaligned
 \tag3.9
 $$
for $q\in\Bbb R$,
where $\beta(q)$ is defined by
 $$
 \sum_{i}p_{i}^{q}r_{i}^{\beta(q)}=1\,.
 \tag3.10
 $$
Arbeiter \& Patzschke [ArPa] and Cawley \& Mauldin [CaMa]
also verified
the Multifractal Formalism for self-similar measures satisfying the OSC.
Namely, in
[ArPa,CaMa]
it is proved that
if $\mu$ is a self-similar 
measure satisfying the OSC, then
 $$
 f_{\mu}(\alpha)
 =
 \beta^{*}(\alpha)
 $$
 for all $\alpha\ge 0$;
 recall, that the definition of the Legendre transform
  is given 
 in (2.5).
We will now
compute  the 
multifractal Minkowski dimensions and 
multifractal Minkowski content
of self-similar measures satisfying various separation conditions.
First, we note that the
multifractal Minkowski dimensions
coincide with $\beta(q)$.
This is not a deep fact and 
is included mainly for completeness.

\bigskip

\proclaim{Theorem 3.2}
Let $K$ and $\mu$ be given by (3.7) and (3.8).
Fix $q\in\Bbb R$
and
assume 
that one of the following conditions is satisfied:
\roster
\item"(i)"
The OSC is satisfied and $0\le q$;
\item"(ii)"
The SSC is satisfied.
\endroster
Then we have
 $$
 \underline\dim_{\Min,\mu}^{q}(K)
 =
 \overline\dim_{\Min,\mu}^{q}(K)
 =
 \beta(q)
 $$
for all $q\in\Bbb R$.
\endproclaim
\noindent{\it  Proof}\newline
\noindent
As noted above, this is not a deep fact and 
can be proven directly
from the definitions using standard arguments 
(similar to those in 
[ArPa] or Falconer's textbook [Fa2]).
The result also follows immediately from
from the main Theorem 3.3 below.
\hfill$\square$

\bigskip

\noindent
Next, we give a 
complete description
of
the 
asymptotic behaviour of the 
multifractal tube formulas for self-similar measures satisfying the OSC.
In particular, we prove that
if
the set
$\{\log r_{1}^{-1},\ldots,\log r_{N}^{-1}\}$
is not contained in a discrete additive subgroup of $\Bbb R$,
then 
$K$ is 
$(q,\beta(q))$ multifractal Minkowski measurable 
with respect to $\mu$,
and 
if
the set
$\{\log r_{1}^{-1},\ldots,\allowmathbreak\log r_{N}^{-1}\}$
is  contained in a discrete additive subgroup of $\Bbb R$,
then 
$K$ is 
$(q,\beta(q))$ average multifractal Minkowski measurable 
with respect to $\mu$.
This is the content
of Theorem 3.2.
The proof of Theorem 3.2 is based on renewal theory and will be 
discussed
 after the statement of the theorem.

\bigskip

\proclaim{Theorem 3.3}
Let $K$ and $\mu$ be given by (3.7) and (3.8).
Fix $q\in\Bbb R$
and
assume 
that one of the following conditions is satisfied:
\roster
\item"(i)"
The OSC is satisfied and $0\le q$;
\item"(ii)"
The SSC is satisfied.
\endroster
Define $\lambda_{q}:(0,\infty)\to\Bbb R$ by
 $$
 \lambda_{q}(r)
 =
 V_{\mu,r}^{q}(K)
 -
 \sum_{i}p_{i}^{q}
 \,
 \bold 1_{(0,r_{i}]}(r)
 \,
 V_{\mu,r_{i}^{-1}r}^{q}(K)
 $$ 
Then we have:
\roster
\item"(1)"
The
non-arithmetic case:
If 
the set
$\{\log r_{1}^{-1},\ldots,\log r_{N}^{-1}\}$
is not contained in a discrete additive subgroup of $\Bbb R$,
then
 $$
 \frac{1}{r^{-\beta(q)}}
 \,
 V_{\mu,r}^{q}(K)
 =
 c_{q}
 +
 \varepsilon_{q}(r)
 $$
where
$c_{q}\in\Bbb R$ is the constant given by
 $$
 c_{q}
 =
  \frac{1}{-\sum_{i}p_{i}^{q}r_{i}^{\beta(q)}\log r_{i}	}
 \,
 \int_{0}^{1}
 r^{\beta(q)}
 \lambda_{q}(r)
  \,
 \frac{dr}{r}
 $$ 
and
$\varepsilon_{q}(r)
 \to0$
 as $r\searrow0$.
In addition,
$K$ is 
$(q,\beta(q))$ multifractal Minkowski measurable 
with respect to $\mu$
with
 $$
 \align
  M_{\mu}^{q,\beta(q)}(K)
 &=
  \frac{1}{-\sum_{i}p_{i}^{q}r_{i}^{\beta(q)}\log r_{i}	}
 \,
 \int_{0}^{1}
 r^{\beta(q)}
 \lambda_{q}(r)
  \,
 \frac{dr}{r}\,.\\
 \endalign
 $$
\item"(2)"
The
arithmetic case:
If 
the set
$\{\log r_{1}^{-1},\ldots,\log r_{N}^{-1}\}$
is contained in a discrete additive subgroup of $\Bbb R$
and
$\langle\log r_{1}^{-1},\ldots,\log r_{N}^{-1}\rangle=u\Bbb Z$
with $u>0$,
then
 $$
 \frac{1}{r^{-\beta(q)}}
 \,
 V_{\mu,r}^{q}(K)
 =
 \pi_{q}(r)
 +
 \varepsilon_{q}(r)
 $$
where
$\pi_{q}:(0,\infty)\to\Bbb R$ is the multiplicatively periodic function
with period equal to $e^{u}$
(i\.e\.
$\pi_{q}(e^{u}r)
 =
 \pi_{q}(r)$
for all 
$r\in(0,\infty)$)
given by
 $$
 \pi_{q}(r)
 =
  \frac{1}{-\sum_{i}p_{i}^{q}r_{i}^{\beta(q)}\log r_{i}	}
 \,
 \sum
 \Sb
  n\in\Bbb Z\\
  {}\\
  re^{nu}\le 1
  \endSb
 (re^{un})^{\beta(q)}
 \,
 \lambda_{q}(re^{un})
 \,
 u$$ 
and
$\varepsilon_{q}(r)
 \to0$
 as $r\searrow0$.
In addition,
$K$ is 
$(q,\beta(q))$ averagely multifractal Minkowski measurable
with respect to $\mu$
with
 $$
 \align
 M_{\mu,\ave}^{q,\beta(q)}(K)
&=
  \frac{1}{-\sum_{i}p_{i}^{q}r_{i}^{\beta(q)}\log r_{i}	}
 \,
 \int_{0}^{1}
 r^{\beta(q)}
 \lambda_{q}(r)
  \,
 \frac{dr}{r}\,.
 \endalign
 $$ 
\endroster
\endproclaim

\bigskip

\noindent
It is instructive to consider the special case $q=0$.
Indeed, since
the multifractal $q$ Minkowski volume for $q=0$
equals the
usual Minkowski volume
and
since the 
$(q,t)$-dimensional
multifractal Minkowski content for 
$q=0$
equals the
usual $t$-dimensional Minkowski content, 
the following corollary, 
providing formulas for the asymptotic behaviour of the 
Minkowski volume of self-similar sets,
follows immediately from Theorem 3.3 by putting $q=0$
This result was first obtained by Gatzouras [Ga]
and later by Winter [Wi].

\bigskip

\proclaim{Corollary 3.4 [Ga]}
Let $K$ be given by (3.7).
Assume 
that the OSC is satisfied.
Let 
$t$ denote the common value of the box dimensions and the Hausdorff dimension of $K$, 
i\.e\. $t$ is the unique number such that
$\sum_{i}r_{i}^{t}=1$ (see [Fa2] or [Hu]).
Define $\lambda:(0,\infty)\to\Bbb R$ by
 $$
 \lambda(r)
 =
 V_{r}(K)
 -
 \sum_{i}
 \,
 \bold 1_{(0,r_{i}]}(r)
 \,
 V_{r_{i}^{-1}r}(K)
 $$ 
Then we have:
\roster
\item"(1)"
The
non-arithmetic case:
If 
the set
$\{\log r_{1}^{-1},\ldots,\log r_{N}^{-1}\}$
is not contained in a discrete additive subgroup of $\Bbb R$,
then
 $$
 \frac{1}{r^{-t}}
 \,
 V_{r}(K)
 =
 c
 +
 \varepsilon(r)
 $$
where
$c\in\Bbb R$ is the constant given by
 $$
 c
 =
  \frac{1}{-\sum_{i}r_{i}^{t}\log r_{i}	}
 \,
 \int_{0}^{1}
 r^{t}
 \lambda(r)
  \,
 \frac{dr}{r}
 $$ 
and
$\varepsilon(r)
 \to0$
as $r\searrow0$.
In addition,
$K$ is 
$t$ Minkowski measurable 
with
 $$
 \align
  M^{t}(K)
 &=
  \frac{1}{-\sum_{i}r_{i}^{t}\log r_{i}	}
 \,
 \int_{0}^{1}
 r^{t}
 \lambda(r)
  \,
 \frac{dr}{r}\,.\\
 \endalign
 $$
\item"(2)"
The
arithmetic case:
If 
the set
$\{\log r_{1}^{-1},\ldots,\log r_{N}^{-1}\}$
is contained in a discrete additive subgroup of $\Bbb R$
and
$\langle\log r_{1}^{-1},\ldots,\log r_{N}^{-1}\rangle=u\Bbb Z$
with $u>0$,
then
 $$
 \frac{1}{r^{-t}}
 \,
 V_{r}(K)
 =
 \pi(r)
 +
 \varepsilon(r)
 $$
where
$\pi:(0,\infty)\to\Bbb R$ is the multiplicatively periodic function
with period equal to $e^{u}$
(i\.e\.
$\pi(e^{u}r)
 =
 \pi(r)$
for all 
$r\in(0,\infty)$)
given by
 $$
 \pi(r)
 =
  \frac{1}{-\sum_{i}r_{i}^{t}\log r_{i}	}
 \,
 \sum
 \Sb
  n\in\Bbb Z\\
  {}\\
  re^{nu}\le 1
  \endSb
 (re^{un})^{t}
 \,
 \lambda(re^{un})
 \,
 u$$ 
and
$\varepsilon(r)
 \to0$
as $r\searrow0$.
In addition,
$K$ is 
$t$ averagely Minkowski measurable
with
 $$
 \align
 M_{\ave}^{t}(K)
&=
  \frac{1}{-\sum_{i}r_{i}^{t}\log r_{i}	}
 \,
 \int_{0}^{1}
 r^{t}
 \lambda(r)
  \,
 \frac{dr}{r}\,.
 \endalign
 $$ 
\endroster
\endproclaim
\noindent{\it  Proof}\newline
\noindent
Since
$\beta(0)
=
\underline\dim_{\Box}(K)
=
\overline\dim_{\Box}(K)
=
\dim(K)
=
t$
(see [Fa2] or [Hu])
and
$V_{\mu,r}^{0}(K)
=
V_{r}(K)$,
this follows from Theorem 3.3 by putting $q=0$.
\hfill$\square$

\bigskip

{\bf 3.3. How does one prove Theorem 3.3?
}
How does one prove Theorem 3.3
on the asymptotic behaviour
of multifractal tubes of self-similar measures?
The proof  is based on renewal theory and, in particular, on 
a very recent renewal theorem by 
Levitin \& Vassiliev [LeVa].
Below we state Levitin \& Vassiliev's renewal theorem.

\bigskip

\proclaim{Theorem 3.5.
Levitin \& Vassiliev's renewal theorem [LeVa]}
Let
$t_{1},\dots, t_{N}>0$
and
$p_{1},\ldots,p_{N}>0$ with $\sum_{i}p_{i}=1$.
Define the probability measure $P$ by
 $$
 P
 =
 \sum_{i}p_{i}\delta_{t_{i}}\,.
 $$

\noindent 
Let $\lambda,\Lambda:\Bbb R\to\Bbb R$
be real valued functions 
satisfying the following conditions:
\roster 
\item"(i)"
The function $\lambda$ is piecewise continuous;
\item"(ii)"
There are constants 
$c,k>0$ such that
 $$
 |\lambda(t)|
 \le 
 ce^{-k|t|}
 $$
for all $t\in\Bbb R$;
\item"(iii)"
We have
 $$
 \Lambda(t)
 \to
 0
 \,\,
 \text{as $t\to-\infty$;}
 $$
\item"(iv)"
We have
 $$
 \Lambda(t)
 =
 \int\Lambda(t-s)\,dP(s)
 +
 \lambda(t)
 $$
for all $t\in\Bbb R$.
\endroster

\noindent
Then the following holds:
\roster
\item"(1)"
The
non-arithmetic case:
If $\{t_{1},\ldots,t_{N}\}$
is not contained in a discrete additive subgroup of $\Bbb R$,
then
 $$
 \Lambda(t)
 =
 c
 +
 \varepsilon(t)
 $$
for all $t\in\Bbb R$ where
 $$
 c
 =
 \frac{1}{\int s\,dP(s)}
 \int\lambda(s)\,ds
 $$
and 
$\varepsilon(t) 
 \to
 0$
as $t\to\infty$.
In addition,
 $$
  \frac{1}{T}\int_{0}^{T}\Lambda(t)\,dt
 \to
  c
 =
 \frac{1}{\int s\,dP(s)}
 \int\lambda(s)\,ds
 \,\,\text{as $T\to\infty$.}
 \tag3.11
 $$

\item"(2)"
The
arithmetic case:
If $\{t_{1},\ldots,t_{N}\}$
is contained in a discrete additive subgroup of $\Bbb R$
and
$\langle t_{1},\ldots,t_{N}\rangle=u\Bbb Z$
with $u>0$,
then
 $$
 \Lambda(t)
 =
 \pi(t)
 +
 \varepsilon(t)
 $$
for all $t\in\Bbb R$ where
$\pi:\Bbb R\to\Bbb R$ 
is the periodic function with period equal to $u$,
(i\.e\.
$\pi(t+u)
 =
 \pi(t)$
for all $t\in\Bbb R$) 
given by
 $$
 \pi(t)
 =
 \frac{1}{\int s\,dP(s)}
 u
 \sum_{n\in\Bbb Z}
 \lambda(t+nu)
 $$
and
$\varepsilon(t) 
 \to
 0$
as $t\to\infty$.
In addition
 $$
  \frac{1}{T}\int_{0}^{T}\Lambda(t)\,dt
 \to
  c
 =
 \frac{1}{\int s\,dP(s)}
 \int\lambda(s)\,ds
 \,\,\text{as $T\to\infty$.}
 \tag3.12
 $$
 \endroster
\endproclaim
\noindent{\it  Proof}\newline
\noindent
All statements, except (3.11) and (3.12), follow 
[LeVa],
and statements (3.11) and (3.12) are easily proved
from the remaining parts of the theorem.
\hfill$\square$

\bigskip

\noindent
The key difference between 
Levitin \& Vassiliev's renewal theorem
and the classical renewal theorem 
from Feller's books [Fel1,Fel2] is the conclusion 
in the arithmetic case.
While the assumptions in the classical renewal theorem are weaker, the
conclusion 
in the arithmetic case is also weaker. 
More precisely, in the arithmetic case
Levitin \& Vassiliev's renewal theorem
says that the error-term
$\varepsilon(t)$
tends to $0$ as $t$ tends to infinity, i\.e\.
$$
 \lim_{t\to\infty}\varepsilon(t)=0\,,
 $$
whereas
the
 classical renewal theorem only allows us to conclude that the error-term 
$\varepsilon(t)$
tends to $0$ as $t$ tends to infinity
through \lq\lq steps" of length $u$, 
i\.e\.
 $$
 \lim
 \Sb
 n\in\Bbb N\\
 n\to\infty
 \endSb
 \varepsilon(nu+s)=0
 $$
for all $s\in\Bbb R$.

Using 
Levitin \& Vassiliev's renewal theorem (Theorem 3.5)
we 
can now prove Theorem 3.3.
Below is a sketch of the proof;
the detailed arguments are presented in Sections 6--9.
In order to prove Theorem 3.3,
we will apply
Levitin \& Vassiliev's Renewal Theorem
to the probability measure $P=P_{q}$
and the
functions
$\lambda=\lambda_{q}^{0}$
and
$\Lambda=\Lambda_{q}^{0}$
defined below.
First,
recall that $\lambda_{q}:(0,\infty)\to\Bbb R$ is defined by
 $$
 \align
 \lambda_{q}(r)
&=
V_{\mu,r}^{q}(K)
 -
 \sum_{i}p_{i}^{q}
 \,
 \bold 1_{(0,r_{i}]}(r)
\,
V_{\mu,r_{i}^{-1}r}^{q}(K)\,.
\endalign
$$
Next,
define $\Lambda_{q}:(0,\infty)\to\Bbb R$ by
 $$
 \Lambda_{q}(r)
 =
 V_{\mu,r}^{q}(K)\,.
 \qquad\qquad
 \qquad\qquad
 \qquad\qquad
 \,\,\,\,{}
 $$
We now define the functions 
$\lambda_{q}^{0},\Lambda_{q}^{0}:\Bbb R\to\Bbb R$ as follows.
 Namely, define
$\lambda_{q}^{0}:\Bbb R\to\Bbb R$ by
 $$
 \lambda_{q}^{0}(t)
 =
 \bold 1_{[0,\infty)}(t)
\,
e^{-t\beta(q)}
 \lambda_{q}(e^{-t})\,,
  \qquad\qquad
 \qquad\,\,\,
 {}
 $$
and define $\Lambda_{q}^{0}:\Bbb R\to\Bbb R$ by
 $$
 \Lambda_{q}^{0}(t)
= 
 \bold 1_{[0,\infty)}(t)
\,
e^{-t\beta(q)}
 \Lambda_{q}(e^{-t})\,.
 \qquad\qquad
 \qquad
 {}
 $$
Finally,  define the probability measure $P_{q}$ by
 $$
 P_{q}
 =
 \sum_{i}p_{i}^{q}r_{i}^{\beta(q)}\delta_{\log r_{i}^{-1}}\,.
 $$
The crux of the matter now is to show that
 the
 probability measure $P=P_{q}$
and the
functions
$\lambda=\lambda_{q}^{0}$
and
$\Lambda=\Lambda_{q}^{0}$
satisfy Conditions (i)--(iv)
in 
Levitin \& Vassiliev's renewal theorem. 
Conditions (i), (iii) and (iv) are
not difficult to verify.
The main 
 difficulty is to prove that
 Condition (ii) is satisfied.
 This 
 is
 highly technical and requires a number
very delicate estimates.
Finally,
applying
Levitin \& Vassiliev's renewal theorem 
to
the probability measure $P=P_{q}$
and the
functions
$\lambda=\lambda_{q}^{0}$
and
$\Lambda=\Lambda_{q}^{0}$
yields Theorem 3.3.



\bigskip
\bigskip

\heading{4. Multifractal tubes measures}\endheading

{\bf 4.1. Multifractal tube measures of general measures.}
The statement in Theorem 3.3 is a global one:
it provides information about the limiting behaviour of
the suitably 
normalized multifractal Minkowski
volume
 $$
 \frac{1}{r^{-\beta(q)}}
 \,
 V_{\mu,r}^{q}(K)
 $$
 of the entire support $K$ of $\mu$ as $r\searrow 0$.
However, it is equally
natural
to ask for 
local versions of Theorem 3.3 describing 
 the limiting behaviour 
of the  
normalized multifractal Minkowski
volume
 $$
 \frac{1}{r^{-\beta(q)}}
 \,
 V_{\mu,r}^{q}(E)
 $$
 of (well behaved) subsets $E$ of the support of $\mu$ as $r\searrow 0$.
In order to address this question, we now introduce 
multifractal
tube measures.
A further motivation for introducing 
multifractal
tube measures comes  from convex geometry
and will be discussed below.

The multifractal tube measures are defined as follows.
Fix a Borel measure  $\mu$ on $\Bbb R^{d}$ and $r>0$.
For a real number $q$, we define the 
multifractal Minkowski tube measure
$\Cal I_{\mu,r}^{q}$ by
 $$
 \Cal I_{\mu,r}^{q}(E)
 =
 \frac{1}{r^{d}}
 \int\limits_{E\cap B(\supp\mu,r)}
 \mu(B(x,r))^{q}
 \,
 d\Cal L^{d}(x)
 \tag4.1
 $$
for Borel subsets $E$ of $\Bbb R^{d}$;
recall, that $\supp\mu$
denotes the support of $\mu$.
Of course, 
the measures
$\Cal I_{\mu,r}^{q}$ will, in general, not converge weakly 
as $r\searrow0$
(indeed, it follows immediately from Theorem 3.3 that, in general,
$\Cal I_{\mu,r}^{q}(\Bbb R^{d})=V_{\mu,r}^{q}(K)$
does not converge as $r\searrow0$).
Hence in order to ensure weak convergence of 
$\Cal I_{\mu,r}^{q}$ as $r\searrow0$
it is necessary
to normalize the measures
$\Cal I_{\mu,r}^{q}$.
There are two natural ways to normalized.
Firstly, we can normalize by volume.
More precisely, we define the
volume normalized multifractal tube measure 
$\Cal V_{\mu,r}^{q}$
by 
 $$
 \Cal V_{\mu,r}^{q}
 =
 \frac{1}{\Cal I_{\mu,r}^{q}(\Bbb R^{d})}
 \,\,
 \Cal I_{\mu,r}^{q}\,.
 \tag4.2
 $$
Secondly, we can normalize by scaling.
More precisely, we defined the 
lower  and upper
scaling
normalized multifractal tube measures
$\underline{\Cal S}_{\mu,r}^{q}$
and
$\overline{\Cal S}_{\mu,r}^{q}$
by
 $$
 \aligned
 \underline{\Cal S}_{\mu,r}^{q}
&=
 \frac{1}{r^{-\underline\dim_{\Min,\mu}^{q}(\supp\mu)}}
 \,\,
 \Cal I_{\mu,r}^{q}\,,\\
&{}\\
 \overline{\Cal S}_{\mu,r}^{q}
&=
 \frac{1}{r^{-\overline\dim_{\Min,\mu}^{q}(\supp\mu)}}
 \,\,
 \Cal I_{\mu,r}^{q}\,;
 \endaligned
 \tag4.3
 $$
recall, that
$\underline\dim_{\Min,\mu}^{q}$
and
$\overline\dim_{\Min,\mu}^{q}$
denote the lower and upper multifractal $q$
Minkowski 
dimension, respectively, see (3.2).

It is instructive to 
consider the particular case $q=0$.
To discuss this case we first make the following definition.
Namely, if $U$ is a closed subset of $\Bbb R^{d}$ and $r>0$, 
the
parallel volume measure 
${\Vol}_{U,r}$
of $U$ is defined by
 $$
 {\Vol}_{U,r}(E)
 =
  \frac
  {\Cal L^{d}(E\cap B(U,r))}
   {\Cal L^{d}(B(U,r))}
   \,,
   \tag4.4
   $$
 see, for example, the texts [Gray,Mo,Schn2].
We now 
note that if $q=0$ and $\mu$ is any Borel measure with $\supp\mu=U$, then
the
volume normalized multifractal tube measure 
$\Cal V_{\mu,r}^{q}$ simplifies to
 $$
 \align
 \Cal V_{\mu,r}^{0}(E)
&=
  \frac
  {\Cal L^{d}(E\cap B(\supp\mu,r))}
   {\Cal L^{d}(B(\supp\mu,r))}\\
&=
  \frac
  {\Cal L^{d}(E\cap B(U,r))}
   {\Cal L^{d}(B(U,r))}
\\
 &=
  {\Vol}_{U,r}(E)\,.  
  \tag4.5
   \endalign
   $$
This observation provides a further motivation for
introducing multifractal 
tube measures.
Namely,
the measure
$ \Cal V_{\mu,r}^{0}(E)={\Vol}_{U,r}(E)$
is closely related
to the notion of 
curvature measures
in
convex geometry.
Curvature measures were introduced in the 1950's
and 
are now recognized as a very powerful tool
for analyzing 
geometric properties of convex sets, see
[Gray,Mo,Schn2].
Indeed, if $U$ is a closed convex subset of $\Bbb R^{d}$
with non-empty interior
and
$l=0,1,2,\ldots,d$, then 
the 
$l$-th order curvature 
measure
${\Vol}_{U}^{l}$
associated with $U$
 is defined as 
the
weak limit 
${\Vol}_{U}^{l}=\lim_{r\searrow0}{\Vol}_{U,r}^{l}$ 
of a certain family $(\,{\Vol}_{U,r}^{l}\,)_{r>0}$ of measures.
While
we will not provide 
the reader with
the definition of the measures ${\Vol}_{U,r}^{l}$
for a general integer $l=0,1,2,\dots,d$
(instead the interested reader can find the definition in previously mentioned texts 
[Gray,Mo,Schn2]),
we do note that 
if $l=d$, then
${\Vol}_{U,r}^{d}={\Vol}_{U,r}$.
In particular,
the
$d$-th order curvature measure ${\Vol}_{U}^{d}$ is defined by
 $$
 \align
 {\Vol}_{U}^{d}
&=
 \lim_{r\searrow0}
 {\Vol}_{U,r}^{d}\\
&= 
 \lim_{r\searrow0}
 {\Vol}_{U,r}\\
&=
 \lim_{r\searrow0}
  \Cal V_{\mu,r}^{0}\,,
 \endalign
 $$
where we have used the fact that
$ \Cal V_{\mu,r}^{0}
=
  {\Vol}_{U,r}$
(see (4.5))  
  and 
$\lim$ denotes the limit with respect to the weak topology.
This 
shows 
that the weak limit
 $$
 \lim_{r\searrow0}
 \Cal V_{\mu,r}^{q}
 $$
(if it exists) may be viewed as a $d$-th order
multifractal curvature measure
and the study of 
multifractal tube measures 
can therefore be seen as a first attempt to create
a theory of multifractal curvatures.

It is, of course,
also possible to define versions of the 
parallel volume measure 
analogous 
to $\underline{\Cal S}_{\mu,r}^{q}$ and $\underline{\Cal S}_{\mu,r}^{q}$.
Indeed,
 if $U$ is a closed subset of $\Bbb R^{d}$ and $r>0$, we define the
 lower and upper
scaling
parallel volume measures
$\underline{\Scale}_{U,r}$ and $\overline{\Scale}_{U,r}$
of $U$ by
 $$
 \aligned
 \underline{\Scale}_{U,r}(E)
&=
  \frac
  {1}
  {r^{-\underline\dim_{\Min}(U)+d}}
  \,
  {\Cal L^{d}(E\cap B(U,r))}\,,\\
&{}\\  
 \overline{\Scale}_{U,r}(E)
&=
  \frac
  {1}
  {r^{-\overline\dim_{\Min}(U)+d}}
  \,
  {\Cal L^{d}(E\cap B(U,r))}\,;
  \endaligned
  \tag4.6
     $$
 recall,
 that $\underline\dim_{\Min}$ and $\overline\dim_{\Min}$ denote
 the lower and upper Minkowski dimension, respectively, see (1.4).
As above, we
note that if $q=0$ and $\mu$ is any probability measure with $\supp\mu=U$, then
the
scaling normalized multifractal tube measure 
$\underline{\Cal S}_{\mu,r}^{q}$ and $\underline{\Cal S}_{\mu,r}^{q}$ simplify to
 $$
 \aligned
 \underline{\Cal S}_{\mu,r}^{0}(E)
 &=
  \underline{\Scale}_{U,r}(E)\,,\\
 \overline{\Cal S}_{\mu,r}^{0}(E)
 &=
  \overline{\Scale}_{U,r}(E)\,.
   \endaligned
   \tag4.7
   $$

\bigskip

{\bf 4.2. Multifractal tube measures of self-similar measures.}
For self-similar measures $\mu$ satisfying the OSC,
we will now
investigate 
the existence 
of the weak
limits of the 
multifractal tube measures
$\Cal V_{\mu,r}^{q}$, 
$\underline{\Cal S}_{\mu,r}^{q}$
and
$\overline{\Cal S}_{\mu,r}^{q}$
as $r\searrow 0$.
In fact,
in many cases these limits exist and equal 
 normalized
 the
multifractal Hausdorff measure (defined below) restricted to the support of $\mu$.

The 
multifractal Hausdorff measure is defined as follows.
Namely,
in an attempt to develop a theoretical framework
for studying the multifractal structure of general
Borel measures,
Olsen [Ol1], Pesin [Pes1] and 
Peyri\`ere [Pey]
introduced
a two parameter family 
$\{\Cal H_{\mu}^{q,t}\,|\, q,t\in\Bbb R\}$ 
of measures
based on certain generalizations of the Hausdorff measure.
The measures 
$\Cal H_{\mu}^{q,t}$ 
 have subsequently been investigated further by a 
large number of authors, including
[Col,\allowlinebreak
Da1,\allowlinebreak
Da2,\allowlinebreak
HoRaSt,\allowlinebreak
Ol2,\allowlinebreak
O'N1,\allowlinebreak
O'N2,\allowlinebreak
Sche],
and are defined as follows.
Let $E\subseteq\Bbb R^{d}$ and $\delta>0$. A countable family 
$(B(x_i,r_i))_i$ of closed balls in $\Bbb R^{d}$ is called a centered
$\delta$-covering of $E$ if $E\subseteq\cup_i\,B(x_i,r_i)$, 
$x_i\in E$ and $0<r_i<\delta$ for all
$i$. 
For $E\subseteq\Bbb R^{d}$, $q,t\in\Bbb R$ and $\delta>0$
write
 $$
 \aligned
 \overline{\Cal H}_{\mu,\delta}^{q,t}(E)
&=\,\,
 \inf
   \Bigg\{
     \sum_i\mu(B(x_i,r_i))^{q}(2r_i)^t
       \,\Bigg|\,
      \text{$(B(x_i,r_i))_i$ is a centered $\delta$-covering of $E$}
   \Bigg\}\,,\\
 \overline{\Cal H}_\mu^{q,t}(E)
&=\,\,\sup_{\delta>0}
  \overline{\Cal H}_{\mu,\delta}^{q,t}(E)\,,\\
 \Cal H_\mu^{q,t}(E)
&=\,\,\sup_{F\subseteq E}\overline{\Cal H}_\mu^{q,t}(F)\,.
 \endaligned
 $$
It follows from [Ol1] that $\Cal H_\mu^{q,t}$ is a
measure on the family of Borel subsets of $\Bbb R^{d}$. The measure 
$\Cal H_\mu^{q,t}$ is, of course, a multifractal generalization of the
centered Hausdorff measure.
In fact, it is easily seen that if $t\ge0$, then
$2^{-t}\Cal H_\mu^{0,t}\le\Cal H^t\le\Cal H_\mu^{0,t}$
where $\Cal H^{t}$ denotes the $t$-dimensional Hausdorff measure.
It is also easily seen that the measure $\Cal H_\mu^{q,t}$ 
in the usual way assign a dimension
to each subset $E$ of $\Bbb R^{d}$
(see  [Ol1]):
there exist  a unique number
$\dim_\mu^q(E)\in[-\infty,\infty]$ such
that
 $$
 \align
 \Cal H_\mu^{q,t}(E)
&=
 \cases
 \infty &\quad\text{for $t<\dim_\mu^q(E)$}\\
 0      &\quad\text{for $\dim_\mu^q(E)<t$}
 \endcases\,.
 \endalign
 $$
The number $\dim_\mu^q(E)$ is an obvious multifractal analogue of the
Hausdorff dimension $\dim(E)$ of $E$.
In fact, it follows immediately from the definitions that
$\dim(E)=\,\,\dim_\mu^0(E)$. 
One of the main importances of the multifractal Hausdorff measure
$\Cal H_{\mu}^{q,t}$  is its  relationship  with the multifractal spectrum
of $\mu$.
Indeed, if we define the dimension function
$b_{\mu}:\Bbb R\to[-\infty,\infty]$ by
 $$
 \align 
 b_{\mu}(q)
&=\,\dim_\mu^q(\supp\mu)\,,
 \endalign
 $$
then 
it follows from [Ol1]
that the
multifractal spectra 
$f_{\mu}$ of $\mu$
(recall, that the multifractal spectrum $f_{\mu}$ is defined in (2.1))
is bounded above by the Legendre transform $b_{\mu}^{*}$
of
$b_{\mu}$  i\.e\.
 $$
 f_{\mu}(\alpha)
 \le
 b_{\mu}^{*}(\alpha)
 $$ 
for all $\alpha\ge 0$, see [Ol1];
recall, that the definition of the Legendre transform
$\varphi^{*}$ of a real valued function $\varphi:\Bbb R\to\Bbb R$
is given in section 2.3.
This inequality may be viewed as a rigorous version of the 
Multifractal Formalism.
Furthermore, for many natural families of measures 
we have 
$f_{\mu}(\alpha)
 =b_{\mu}^{*}(\alpha)$ 
for all $\alpha\ge 0$, cf\. 
[Col,\allowlinebreak
Da1,\allowlinebreak
Da2,\allowlinebreak
Ol1,\allowlinebreak
Ol2]

Using the
multifractal Hausdorff measures $\Cal H_{\mu}^{q,t}$,
we will now explicitly identify
the weak
limits of the 
multifractal tube measures
$\Cal V_{\mu,r}^{q}$, 
$\underline{\Cal S}_{\mu,r}^{q}$
and
$\overline{\Cal S}_{\mu,r}^{q}$
as $r\searrow 0$
for self-similar measures $\mu$.
The first result shows that the
weak limit of
$\Cal V_{\mu,r}^{q}$
as $r\searrow 0$
always exists and equals the 
normalized multifractal Hausdorff measure.
In Theorem 4.1 and the subsequent parts of the paper
we use the following
notation.
Namely,
if $\Cal M$
is a Borel masure on $\Bbb R^{d}$
and $E$ is a Borel subset of $\Bbb R^{d}$, 
then 
we denote the restriction of $\Cal M$ to $E$ by
$\Cal M\floor E$, i\.e\.
 $$
 (\Cal M\floor E)(B)
 =
 \Cal M(E\cap B)
 \tag4.8
 $$
for all Borel subsets $B$ of $\Bbb R^{d}$.
We can now state Theorem 4.1.

\bigskip

\proclaim{Theorem 4.1}
Let $K$ and $\mu$ be given by (3.7) and (3.8).
Fix $q\in\Bbb R$
and
assume 
that one of the following conditions is satisfied:
\roster
\item"(i)"
The OSC is satisfied and $0\le q$;
\item"(ii)"
The SSC is satisfied.
\endroster
Then we have
$$
 \Cal V_{\mu,r}^{q}
 \,
 \to
 \,
 \frac{1}{\Cal H_{\mu}^{q,\beta(q)}(K)}
 \,\,
 \Cal H_{\mu}^{q,\beta(q)}\floor K
 \qquad
 \text{weakly.}
 $$ 
\endproclaim

 \bigskip

\noindent 
Next, we study
the 
limiting behaviour of
$\underline{\Cal S}_{\mu,r}^{q}$
and
$\overline{\Cal S}_{\mu,r}^{q}$
as $r\searrow 0$
for self-similar measures $\mu$.
Contrary to
Theorem 4.1, the weak limits of
$\underline{\Cal S}_{\mu,r}^{q}$
and
$\overline{\Cal S}_{\mu,r}^{q}$
as $r\searrow 0$
may not exist.
Indeed,
if 
the set
$\{\log r_{1}^{-1},\ldots,\log r_{N}^{-1}\}$
is contained in a discrete additive subgroup of $\Bbb R$,
then the weak limits
of
$\underline{\Cal S}_{\mu,r}^{q}$
and
$\overline{\Cal S}_{\mu,r}^{q}$
as $r\searrow 0$
do not necessarily
exist;
however
the weak limits of certain averages of 
$\underline{\Cal S}_{\mu,r}^{q}$
and
$\overline{\Cal S}_{\mu,r}^{q}$
exist
and
equal  a multiple
of the 
normalized multifractal Hausdorff measure.
On the other hand,
if 
the set
$\{\log r_{1}^{-1},\ldots,\log r_{N}^{-1}\}$
is not contained in a discrete additive subgroup of $\Bbb R$,
then the weak limits
of
$\underline{\Cal S}_{\mu,r}^{q}$
and
$\overline{\Cal S}_{\mu,r}^{q}$
as $r\searrow 0$
always
exist and
equal  a multiple
of the 
normalized multifractal Hausdorff measure.

\bigskip

\proclaim{Theorem 4.2}
Let $K$ and $\mu$ be given by (3.7) and (3.8).
Fix $q\in\Bbb R$
and
assume 
that one of the following conditions is satisfied:
\roster
\item"(i)"
The OSC is satisfied and $0\le q$;
\item"(ii)"
The SSC is satisfied.
\endroster
Then the following holds.
\roster
\item"(1)"
We have
 $$
 \underline{\Cal S}_{\mu,r}^{q}
 =
 \overline{\Cal S}_{\mu,r}^{q}
 =
 \frac{1}{r^{-\beta(q)}}
 \,\,
 \Cal I_{\mu,r}^{q}\,.
 $$
\endroster
Write
$\Cal S_{\mu,r}^{q}$
for the 
common value of 
$\underline{\Cal S}_{\mu,r}^{q}$
and
$\overline{\Cal S}_{\mu,r}^{q}$, i\.e\. write
 $$
 \Cal S_{\mu,r}^{q}
 =
 \frac{1}{r^{-\beta(q)}}
 \,\,
 \Cal I_{\mu,r}^{q}\,.
 $$
Also, define the average measure 
$\Cal S_{\mu,r,\ave}^{q}$
by
 $$
 {}
 \qquad\qquad
 \quad\,\,
 \Cal S_{\mu,r,\ave}^{q}
 =
 \frac{1}{-\log r}
 \int_{r}^{1}
 \,\,
 \frac{1}{s^{-\beta(q)}}
 \,\,
 \Cal I_{\mu,s}^{q}
 \,
 \frac{ds}{s}\,.
 $$
Then the following holds.
\roster
\item"(2)"
The non-arithmetic case:
If 
the set
$\{\log r_{1}^{-1},\ldots,\log r_{N}^{-1}\}$
is not contained in a discrete additive subgroup of $\Bbb R$,
then
$$
 \align
 \Cal S_{\mu,r}^{q}
 \,
&\to
 \,
 M_{\mu}^{q,\beta(q)}(K)
 \,\,
 \frac{1}{\Cal H_{\mu}^{q,\beta(q)}(K)}
 \,\,
 \Cal H_{\mu}^{q,\beta(q)}\floor K
 \qquad
 \text{weakly,}\\
&{}\\ 
  \Cal S_{\mu,r,\ave}^{q}
 \,
&\to
 \,
 M_{\mu,\ave}^{q,\beta(q)}(K)
 \,\,
 \frac{1}{\Cal H_{\mu}^{q,\beta(q)}(K)}
 \,\,
 \Cal H_{\mu}^{q,\beta(q)}\floor K
 \qquad
 \text{weakly;}
 \endalign
 $$ 
recall, that $K$
is $(q,\beta(q))$ mutifractal Minkowski measurable with respect to $\mu$
and
$(q,\beta(q))$ average multifractal Minkowski measurable
with respect to $\mu$
by Theorem 3.3 and 
the
multifractal Minkowski content $M_{\mu}^{q,\beta(q)}(K)$
and the
average multifractal Minkowski content $M_{\mu,\ave}^{q,\beta(q)}(K)$
are therefore well-defined.
\item"(3)"
The arithmetic case:
If 
the set
$\{\log r_{1}^{-1},\ldots,\log r_{N}^{-1}\}$
is contained in a discrete additive subgroup of $\Bbb R$,
then
$$
 \align 
 \Cal S_{\mu,r,\ave}^{q}
 \,
&\to
 \,
 M_{\mu,\ave}^{q,\beta(q)}(K)
 \,\,
 \frac{1}{\Cal H_{\mu}^{q,\beta(q)}(K)}
 \,\,
 \Cal H_{\mu}^{q,\beta(q)}\floor K
 \qquad
 \text{weakly;}
 \endalign
 $$ 
recall, that $K$
is
$(q,\beta(q))$ average multifractal Minkowski measurable
with respect to $\mu$
by Theorem 3.3 and 
the
average multifractal Minkowski content $M_{\mu,\ave}^{q,\beta(q)}(K)$
is therefore well-defined.
\endroster
\endproclaim

\bigskip

\noindent
As with Theorem 3.3, it is 
instructive to consider the special case $q=0$.
Indeed, 
we first note that 
if
$K$ and $\mu$ are given by (3.7) and (3.8),
then
(see (4.5))
 $$
 \align
 \Cal V_{\mu,r}^{0}(E)
&=
 \frac
 {\Cal L^{d}(E\cap B(K,r))}
 {\Cal L^{d}(B(K,r))}
 \\
&=
 \Vol_{K,r}(E) 
 \endalign
 $$
i\.e\.
$\Cal V_{\mu,r}^{0}$ equals the normalised parallel body measure
$\Vol_{K,r}$.
Next,
writing $t$ for the common value of the the box dimensions 
and Hausdorff dimension
of $K$,
we note (see (4.7)) that
 $$
 \align
 \underline{\Cal S}_{\mu,r}^{0}(E)
 =
 \overline{\Cal S}_{\mu,r}^{0}(E)
&=
 \frac{1}{r^{-t+d}}
 \,
 \Cal L^{d}(E\cap B(K,r))\\
&=
\underline\Scale_{K,r}(E)
 =
 \overline\Scale_{K,r}(E)\,,
 \qquad\qquad
 \qquad\qquad
 \quad
 \endalign
 $$
i\.e\.
$\underline{\Cal S}_{\mu,r}^{0}$
and
$\overline{\Cal S}_{\mu,r}^{0}$
 equal the scaling parallel body measures
$\underline\Scale_{K,r}$
 and
$\overline\Scale_{K,r}$.
 The following corollaries 
therefore follow immediately from Theorem 4.1 and Theorem 4.2 
by putting $q=0$.
These results were first obtained by Winter in his doctoral dissertation [Wi].
Recall, that for $t\ge 0$, we write
$\Cal H^{t}$ for the $t$-dimensional Hausdorff measure.

\bigskip

\proclaim{Corollary 4.3. [Wi]}
Let $K$ be given by (3.7).
Assume 
that the OSC is satisfied.
Let 
$t$ denote the common value of the box dimensions and the Hausdorff dimension of $K$, 
i\.e\. $t$ is the unique number such that
$\sum_{i}r_{i}^{t}=1$.
For $r>0$,
the normalised parallel body measure $\Vol_{K,r}$ is given by
 $$
 \Vol_{K,r}(E)
 =
 \frac
 {1}
 {\Cal L^{d}(B(K,r))}
 \,
 \Cal L^{d}(E\cap B(K,r))\,.
 $$
Then we have
$$
 \Vol_{K,r}
 \,
 \to
 \,
 \frac{1}{\Cal H^{t}(K)}
 \,\,
 \Cal H^{t}\floor K
 \qquad
 \text{weakly.}
 $$ 
\endproclaim
\noindent{\it  Proof}\newline
\noindent
Recall, that $K$ is given by (3.7),
and let
$\mu$ be given by $(3.8)$.
Since
$\Cal V_{\mu,r}^{0}=\Vol_{K,r}$, 
the statement now follows from Theorem 4.1 by putting $q=0$.
\hfill$\square$

\bigskip

\proclaim{Corollary 4.4. [Wi]}
Let $K$ be given by (3.7).
Assume 
that the OSC is satisfied.
Let 
$t$ denote for the common value of the 
box dimensions 
and the
Hausdorff dimension of $K$, 
i\.e\. $t$ is the unique number such that
$\sum_{i}r_{i}^{t}=1$.
\roster
\item"(1)"
We have 
$$
 \underline\Scale_{K,r}(E)
 =
 \overline\Scale_{K,r}(E)
 =
 \frac{1}{r^{-t+d}}
 \,
 \Cal L^{d}(E\cap B(K,r))\,.
 $$
\endroster 
Write 
$\Scale_{K,r}$
for the common value of
$ \underline\Scale_{K,r}$
and 
$\overline\Scale_{K,r}$, i\.e\. write
$$
 \Scale_{K,r}(E)
 =
 \frac{1}{r^{-t+d}}
 \,
 \Cal L^{d}(E\cap B(K,r))\,.
 $$
Also, define the average
measure
$\Scale_{K,r,\ave}$
by
 $$
 {}
 \qquad\qquad
 \Scale_{K,r,\ave}(E)
 =
 \frac{1}{-\log r}
 \int_{r}^{1}
 \,\,
 \frac{1}{s^{-t+d}}
 \,
 \Cal L^{d}(E\cap B(K,s))
 \,
 \frac{ds}{s}\,.
 $$
Then the following holds.
\roster
\item"(2)"
The non-arithmetic case:
If 
the set
$\{\log r_{1}^{-1},\ldots,\log r_{N}^{-1}\}$
is not contained in a discrete additive subgroup of $\Bbb R$,
then
$$
 \align
 \Scale_{K,r}
 \,
&\to
 \,\,
 \,
 M^{t}(K)
 \,
 \,\,
 \frac{1}{\Cal H^{t}(K)}
 \,\,
 \Cal H^{t}\floor K
 \qquad
 \text{weakly,}\\
&{}\\ 
  \Scale_{K,r,\ave}
 \,
&\to
 \,
 M_{\ave}^{t}(K)
 \,\,
 \frac{1}{\Cal H^{t}(K)}
 \,\,
 \Cal H^{t}\floor K
 \qquad
 \text{weakly;}
 \endalign
 $$ 
recall, that $K$
is $t$ Minkowski measurable 
and
$t$ average Minkowski measurable
by Corollary 3.4 and 
the
Minkowski content $M^{t}(K)$
and the
average Minkowski content $M_{\ave}^{t}(K)$
are therefore well-defined.
\item"(3)"
The arithmetic case:
If 
the set
$\{\log r_{1}^{-1},\ldots,\log r_{N}^{-1}\}$
is contained in a discrete additive subgroup of $\Bbb R$
then
$$
 \align 
 \Scale_{K,r,\ave}
 \,
&\to
 \,
 M_{\ave}^{t}(K)
 \,\,
 \frac{1}{\Cal H^{t}(K)}
 \,\,
 \Cal H^{t}\floor K
 \qquad
 \text{weakly;}
 \endalign
 $$ 
recall, that $K$
is
$t$ average Minkowski measurable
by Corollary 3.4 and 
the
average multifractal Minkowski content $M_{\ave}^{t}(K)$
is therefore well-defined.
\endroster
\endproclaim
\noindent{\it  Proof}\newline
\noindent
Recall, that $K$ is given by (3.7),
and let
$\mu$ be given by $(3.8)$.
Since
$\underline{\Cal S}_{\mu,r}^{0}
=
\overline{\Cal S}_{\mu,r}^{0}
=
\Scale_{K,r}$, 
the statement now follows from
 Theorem 4.2 by putting $q=0$.
\hfill$\square$

\bigskip

In Section 4.1
it was suggested
that 
the 
limiting behaviour of
the multifractal tube measures
$\Cal V_{\mu,r}^{q}$
may be viewed
as providing a 
local  version of Theorem 3.3.
Indeed, Theorem 3.3 clearly
describes the 
limiting behaviour of $\frac{1}{r^{-\beta(q)}}V_{\mu,r}^{q}(K)$ as $r\searrow0$
whereas
the main results in Section 4.1 equally clearly provide information about the
the 
limiting behaviour of $\frac{1}{r^{-\beta(q)}}V_{\mu,r}^{q}(E)$ as $r\searrow0$
for 
\lq\lq well-behaved" subsets $E$ of $K$.
The view point is made even more explicit (and precise) in the next corollary.
Below we use the following notation, namely, 
if $E$ is a subset of $\Bbb R^{d}$,  then we  denote 
the boundary of $E$ by
$\partial E$.

\bigskip

\proclaim{Corollary 4.5}
Let $K$ and $\mu$ be given by (3.7) and (3.8).
Fix $q\in\Bbb R$
and
assume 
that one of the following conditions is satisfied:
\roster
\item"(i)"
The OSC is satisfied and $0\le q$;
\item"(ii)"
The SSC is satisfied.
\endroster
Let $E\subseteq \Bbb R^{d}$ be a Borel set with:

\roster
\item"(a)"
$\Cal H_{\mu}^{q,\beta(q)}(E\cap K)>0$.

\item"(b)"
$ \Cal H_{\mu}^{q,\beta(q)}(\partial E\cap K)=0$.

\item"(c)"
$E\cap B(K,r)=B(E\cap K,r)$ for $r$ small enough.

\endroster
(Observe that, for example, the set $E=\Bbb R^{d}$
satisfies the above conditions,
and if
$K=L\cup M$
with
$\dist(L,M)>0$
and
$\Cal H_{\mu}^{q,\beta(q)}(L)>0$
and
$0<\delta<\dist(L,M)$, then
 the set 
 $E=B(L,\delta)$
satisfies the above conditions.)

\medskip

\noindent
Then the following holds.
\roster
\item"(1)"
The non-arithmetic case:
If 
the set
$\{\log r_{1}^{-1},\ldots,\log r_{N}^{-1}\}$
is not contained in a discrete additive subgroup of $\Bbb R$,
then
$E\cap K$
is $(q,\beta(q))$ mutifractal Minkowski measurable with respect to $\mu$
with
 $$
 M_{\mu}^{q,\beta(q)}(E\cap K)
 =
 M_{\mu}^{q,\beta(q)}(K)
 \,\,
 \frac{\Cal H_{\mu}^{q,\beta(q)}(E\cap K)}{\Cal H_{\mu}^{q,\beta(q)}(K)}\,;
 $$
recall, that
$K$ is
$(q,\beta(q))$ average multifractal Minkowski measurable
with respect to $\mu$
by Theorem 3.3 and 
the
multifractal Minkowski content $M_{\mu}^{q,\beta(q)}(K)$
is therefore well-defined.
\item"(2)"
The arithmetic case:
If 
the set
$\{\log r_{1}^{-1},\ldots,\log r_{N}^{-1}\}$
is contained in a discrete additive subgroup of $\Bbb R$
then
$E\cap K$
is $(q,\beta(q))$ average mutifractal Minkowski measurable with respect to $\mu$
with
 $$
 M_{\mu,\ave}^{q,\beta(q)}(E\cap K)
 =
 M_{\mu,\ave}^{q,\beta(q)}(K)
 \,\,
 \frac{\Cal H_{\mu}^{q,\beta(q)}(E\cap K)}{\Cal H_{\mu}^{q,\beta(q)}(K)}\,;
 $$
recall, that $K$
is
$(q,\beta(q))$ average multifractal Minkowski measurable
with respect to $\mu$
by Theorem 3.3 and 
the
average multifractal Minkowski content $M_{\mu,\ave}^{q,\beta(q)}(K)$
is therefore well-defined.
\endroster
\endproclaim
\noindent{\it  Proof}\newline
\noindent
This follows immedately
from Theorem 4.2 since
the condition
$E\cap B(K,r)=B(E\cap K,r)$, implies that
$\Cal I_{\mu,r}^{q}(E)
=
\frac{1}{r^{d}}
\int_{E\cap B(K,r)}\mu(B(x,r))^{q}\,d\Cal L^{d}(x)
=
\frac{1}{r^{d}}
\int_{B(E\cap K,r)}\mu(B(x,r))^{q}\,d\Cal L^{d}(x)
=
V_{\mu,r}^{q}(E\cap K)$.
\hfill$\square$

\bigskip

\noindent
Note that Corollary 4.5 is a
genuine 
extension of Theorem 3.3:
namely, 
if we let $E=K$ in Corollary 4.5, then Corollary 4.5 simplifies to Theorem 3.3.


\bigskip
\bigskip

\centerline{\smc 5. Symbolic multifractal tubes of self-similar measures:}
\centerline{\smc Multifractal zeta-functions}
\centerline{\smc and}
\centerline{\smc explicit formulas}

%

%

\bigskip

Throughout this section we will let 
$K$ and $\mu$ denote the self-similar set and the self-similar
measure given by (3.7) and (3.8), respectively.
While Theorem 3.3 provides complete information about the asymptotic 
behaviour of the 
multifractal Minkowski volume 
$V_{\mu,r}^{q}(K)$
of  $K$, it
does not provide
\lq\lq explicit"
formulas
for the
multifractal Minkowski content 
$M_{\mu}^{q,\beta(q)}(K)=\lim_{r\searrow0}\frac{1}{r^{-\beta(q)}}V_{\mu,r}^{q}(K)$.
Indeed, the formulas in Theorem 3.3
for
multifractal Minkowski content of $K$
involve 
the integral of an auxiliary function $\lambda_{q}$.
Even in very simple cases it is highly unlikely that this integral can be computed
explicitly.
This is clearly unsatisfactory
and it would be desirable
if more explicit expressions could be found.
In fact, 
even in the fractal case,
the problem of
finding 
explicit formulas for the Minkowski content 
is  highly non-trivial.
However, 
despite, or perhaps in spite,
 of the difficulties, 
this problem has recently attracted considerable interest.
In particular, 
Lapidus and collaborators [LapPea1,LapPea2,LapPeaWi,Lap-vF1,Lap-vF2] have
during the past 20 years
and
with spectacular success
pioneered
the use of zeta-functions
to 
obtain 
explicit formulas for the Minkowski content 
of self-similar subsets of the line
and certain self-similar sets in higher dimensions.
It would clearly be desirable if analogous 
 formulas
for the 
multifractal Minkowski content
could
be found.
However, the significant difficulties encountered by Lapidus and collaborators
when computing the 
Minkowski content 
of self-similar subsets of the line suggests that this problem is
exceptionally difficult.
For this reason we
introduce 
\lq\lq symbolic" multifractal Minkowski
volumes.
The 
\lq\lq symbolic" multifractal Minkowski
volumes
are defined in
such a way
that
they are 
\lq\lq compatible" with the 
usual Minkowski volumes 
(see Theorem 5.1 below for 
a precise formulation of this)
and such that
the 
zeta-function technique can be applied to give
explicit formulas for the 
corresponding
\lq\lq symbolic" multifractal Minkowski content.
A multifractal zeta-function is a meromorphic 
function whose residues are closely related to
the asymptotic behaviour of 
the
\lq\lq symbolic" 
multifractal Minkowski volume.
Namely, using the residue theorem it is possible to
relate 
the \lq\lq symbolic"  multifractal Minkowski volume to the 
residues of the zeta-function,
and
a careful analysis of the residues
will then provide explicit formulas for
the
\lq\lq symbolic" 
multifractal Minkowski volume.
The idea
 of using zeta-functions
in order to
obtain
 explicit formulas for the \lq\lq symbolic" multifractal Minkowski content
 has classical origins.
  For example,
 the 
 \lq\lq standard" proofs of the Prime Number Theorem
 is based on applying this technique
 to the Riemann zeta-function, see [Ed,Pat].
The zeta-function technique 
  for 
  is not only restricted to problems in number theory, but
  has also been
  used
  successfully 
  to obtain explicit formulas for 
  \lq\lq counting functions" in many
  other areas in mathematics.
  For example,
  in dynamical systems,
   Parry \& Pollicott [ParPo1,ParPo2]
   obtained
   asymptotic formulas for the number of closed geodesics 
whose length is less than $x$ as $x\to\infty$
by
  applying this technique to Ruelle's zeta-function for Axiom A flows.
  For other applications of this technique in dynamical systems the reader is referred to Ruelle's
text [Rue2].

\bigskip

{\bf 5.1. Symbolic multifractal tubes of self-similar measures.}
We will now define the symbolic multifractal Moinkowski volume.
We first introduce the following notation.
 Let
$\Sigma
 =
 \{1,\ldots,N\}$ and write
 $$
 \aligned
 \Sigma^{m}
&=\{1,\dots,N\}^{m}\,,\\
 \Sigma^{*}
&=\bigcup_{m}\Sigma^{m}\,.
 \endaligned
 \tag5.1
 $$
i\.e\. $\Sigma^{m}$ is the family of all
strings
$\bold i=i_{1}\ldots i_{m}$
of length $m$ 
with $i_{j}\in\{1,\ldots,N\}$
and
$\Sigma^{*}$ is the family of all finite strings
$\bold i=i_{1}\ldots i_{m}$
with $i_{j}\in\{1,\ldots,N\}$.
Also, for  $\bold i=i_{1}\dots i_{m}\in\Sigma^{*}$, we will write
 $$
 \aligned
 r_{\bold i}
&=
 r_{i_{1}}\cdots r_{i_{m}}\,,\\
 p_{\bold i}
&=
 p_{i_{1}}\cdots p_{i_{m}}\,.
 \endaligned
 \tag5.2
 $$
Next, for brevity, put
 $$
 r_{\min}
 =\min_{i=1,\ldots,N}r_{i}\,\,\,
 r_{\max}
 =
 \max_{i=1,\ldots,N}r_{i}\,.
 $$
Finally, if 
$\bold i=i_{1}\dots i_{m}\in\Sigma^{*}$, then
we will write
$\widehat{\bold i}$ for the 
\lq\lq parent" of $\bold i$, i\.e\. we will write
 $$
 \widehat{\bold i}
 =
 i_{1}\dots i_{m-1}\,.
 $$
We can now define the symbolic multifractal Minkowski
volume.
We provide several comments
discussing  the motivation behind  the
definition of
the symbolic multifractal Minkowski
volume
immediately after the statement of the definition.

\bigskip

\proclaim{Definition. Symbolic multifractal Minkowski volume}
Fix $q\in\Bbb R$ and $l=0,1,\ldots,d$.
For brevity write
 $$
 \sigma_{q,l}=\sum_{i=1}^{N}p_{i}^{q}r_{i}^{l-dq}\,,
 \tag5.3
 $$ 
and let
  $$
 C_{\mu,r}^{q,l,\sym}(K)
 =
   \sum
 \Sb
 \bold i\\
 r_{\bold i}<r<r_{\hat{\bold i}}
 \endSb
 p_{\bold i}^{q}r_{\bold i}^{l-dq}
 \,\,
 +
 \,\,
 \tfrac{1+\frac{1}{\sigma_{q,l}}}{2}
 \,\,
   \sum
 \Sb
 \bold i\\
 r=r_{\hat{\bold i}}
 \endSb
 p_{\bold i}^{q}r_{\bold i}^{l-dq}\,.
 \tag5.4
 $$ 
Let
 $ \kappa_{\mu}^{q,0}(K),
 \kappa_{\mu}^{q,1}(K),\ldots,
 \kappa_{\mu}^{q,d}(K)$
 be 
  real numbers
satisfying the following consistency condition
  $$
   \sum_{l}
 \kappa_{\mu}^{q,l}(K)\,(\sigma_{q,l}-1)
 =
 0\,.
 \tag5.5
 $$
 For $r>0$,
we
define the symbolic $q$ multifractal Minkowski volume
$V_{\mu,r}^{q,\sym}(K)$
 of $K$ 
with respect to the measure $\mu$ by
 $$
 \align
 V_{\mu,r}^{q,\sym}(K)
&=
 \frac{1}{r^{d}}
 \sum_{l}
  \kappa_{\mu}^{q,l}(K)
  \,
 C_{\mu,r}^{q,l,\sym}(K)
\,
r^{(d-l)+dq}\\
&=
 \sum_{l}
  \kappa_{\mu}^{q,l}(K)
  \,
 C_{\mu,r}^{q,l,\sym}(K)
\,
r^{-l+dq}\,.
\tag5.6
\endalign
$$
\endproclaim

\bigskip

\noindent
{\bf Comment. Motivating the definition of $V_{\mu,r}^{q,\sym}(K)$.}
We will now make a number of  comments
explaining the motivation behind the definition
of
the
 symbolic $q$ multifractal 
 Minkowski volume
$V_{\mu,r}^{q,\sym}(K)$.

\noindent
{\it (1) Motivating definition (5.6):}
It is clear
that definition (5.6)
is motivated by Steiner's formula (1.3):
 the quantity
 $$
  C_{\mu,r}^{q,l,\sym}(K)
\,
r^{(d-l)+dq}
=
\left(
 \sum
 \Sb
 \bold i\\
 r_{\bold i}<r<r_{\hat{\bold i}}
 \endSb
 p_{\bold i}^{q}r_{\bold i}^{l-dq}
 \,\,
 +
 \,\,
 \tfrac{1+\frac{1}{\sigma_{q,l}}}{2}
 \,\,
   \sum
 \Sb
 \bold i\\
 r=r_{\hat{\bold i}}
 \endSb
 p_{\bold i}^{q}r_{\bold i}^{l-dq}
 \right)
 \,
r^{(d-l)+dq}
$$
clearly
corresponds 
to the 
term $r^{d-l}$ in Steiner's formula,
and the
quantities
$$
 \kappa_{\mu}^{q,0}(K),
 \kappa_{\mu}^{q,1}(K),
 \ldots,
 \kappa_{\mu}^{q,d}(K)
 $$
 correspond to the
  Quermassintegrale 
$ \kappa^{0}(C),
 \kappa^{1}(C),\ldots,
 \kappa^{d}(C)$
in Steiner's formula.

\noindent
{\it (2) Motivating consistency condition (5.5):}
Consistency condition (5.5) is motivated by the following argument.
If $C$ is a bounded convex subset of $\Bbb R^{d}$ with smooth boundary $\partial C$,
 then it follows from  Weyl's extension of
Steiner's formula 
(to
compact oriented $d$-dimensional Riemannian manifolds
 isometrically embedded into Euclidean space)
 applied
 to $\partial C$ 
 that there are
constants 
$k^{0}(C),k^{1}(C),\ldots,k^{d-1}(C)$ 
such that
the  volume of 
$B(\Bbb R^{d}\setminus C,r)\cap C
=
B(\partial C,r)\cap C$
is given by
$$
 \Cal L^{d}(B(\Bbb R^{d}\setminus C,r)\cap C)
 =
 \sum_{l}k^{l}(C)\,r^{d-l}
 \tag5.7
 $$
for all $r$ less than 
the inner radius $r_{\inner}$ of $C$.
Trivially, we also have
$$
 \Cal L^{d}(B(\Bbb R^{d}\setminus C,r)\cap C)
 =
\Cal L^{d}(C)
\qquad\quad\,\,
\tag5.8
 $$
for all $r$ greater than 
the inner radius $r_{\inner}$ of $C$.
Since the volume $ \Cal L^{d}(B(\Bbb R^{d}\setminus C,r)\cap C)$ 
is
  a continuous function of $r$,  
 it now follows from (5.7) and (5.8) that
  the constants
  $k^{0}(C),\allowmathbreak
  k^{1}(C),\allowmathbreak
  \ldots,\allowmathbreak
  k^{d-1}(C)$ must satisfy the following consistency 
  condition, namely,
$\sum_{l<d}k^{l}(C)\,r_{\inner}^{d-l}
=
\Cal L^{d}(C)$.
Writing  $k^{d}(C)=-\Cal L^{d}(C)$, this 
condition can be
rewritten as
 $$
 \sum_{l<d}k^{l}(C)\,r_{\inner}^{d-l}
 =
 -k^{d}(C)\,r_{\inner}^{d-d}
 $$
i\.e\.
 $$
  \sum_{l}k^{l}(C)\,r_{\inner}^{d-l}
 =
 0\,.
 \tag5.9
 $$
However, typically the set $K$ has zero $d$-dimensional volume,
and the symbolic multifractal Minkowski volume 
$ V_{\mu,r}^{q,\sym}(K)$ 
of $B(K,r)$
can therefore also be thought 
of as being equal to the 
symbolic multifractal 
Minkowski volume of 
$B(K,r)\cap(\Bbb R^{d}\setminus K)$.
Comparison with
(5.9) therefore
suggests 
that the 
coefficients
$\kappa_{\mu}^{q,0}(K),
 \kappa_{\mu}^{q,1}(K),
 \ldots,
 \kappa_{\mu}^{q,d}(K)$
 must satisfy the following consistency condition,  namely,
 $$
 \sum_{l}\kappa_{\mu}^{q,l}(K)\,(\sigma_{q,l}-1)=0\,.
 \tag5.10
 $$
 While
consistency condition (5.5) is motivated by the informal discussion above,
we note that it, nevertheless, plays
a crucial  role in the proofs in Section 16.

\noindent
{\it (3) Motivating the weight factor 
$\frac{1+\frac{1}{\sigma_{q,l}}}{2}$
in (5.4):}
Recall definition (5.4) of $C_{\mu,r}^{q,l,\sym}(K)$, namely,
 $$
 C_{\mu,r}^{q,l,\sym}(K)
 =
   \sum
 \Sb
 \bold i\\
 r_{\bold i}<r<r_{\hat{\bold i}}
 \endSb
 p_{\bold i}^{q}r_{\bold i}^{l-dq}
 \,\,
 +
 \,\,
 \tfrac{1+\frac{1}{\sigma_{q,l}}}{2}
 \,\,
   \sum
 \Sb
 \bold i\\
 r=r_{\hat{\bold i}}
 \endSb
 p_{\bold i}^{q}r_{\bold i}^{l-dq}\,.
 $$ 
In this definition,
the 
\lq\lq boundary" term
$  \sum_{
 \bold i\,,\,
 r=r_{\hat{\bold i}}
 }
 p_{\bold i}^{q}r_{\bold i}^{l-dq}$
 has been weighted by the
 factor $\frac{1+\frac{1}{\sigma_{q,l}}}{2}$.
 The motivation behind this is the following.
 Namely, below we intend to apply the Mellin transform 
 to the function
 $r\to C_{\mu,r}^{q,l,\sym}(K)$.
 However, the Mellin transform only applies to 
 piecewise continuous functions $f$
 for which
 $\frac{\lim_{x\nearrow x_{0}}f(x)+\lim_{x\searrow x_{0}}f(x)}{2}=f(x_{0}) $
 at all discontinuity points $x_{0}$.
 Weighting 
 the 
\lq\lq boundary" term
$  \sum_{
 \bold i\,,\,
 r=r_{\hat{\bold i}}
 }
 p_{\bold i}^{q}r_{\bold i}^{l-dq}$
by the
 factor $\frac{1+\frac{1}{\sigma_{q,l}}}{2}$
 ensures that the function $r\to C_{\mu,r}^{q,l,\sym}(K)$
 satisfies this condition.
 A similar practice is 
 also commonly used in number theory 
 where 
 analogous
  \lq\lq boundary"  terms
  are weighted by the factor $\frac{1}{2}$.
  Borrowing terminology from physics 
  where the parameter $q$ is interpreted as the inverse temperature of the physical system
  associated with $\mu$
  (see, for example, [BaPo, pp\. 128--132;
  BeSc, pp. 114--126;
  Ot, pp. 309--910]),
we may therefore, somewhat poetically, say that
  the factor $\frac{1+\frac{1}{\sigma_{q,l}}}{2}$
  represents the usual weight 
  factor $\frac{1}{2}$
 when \lq\lq raised to the temperature $\frac{1}{q}$".

\bigskip

\noindent
{\bf Comment. Comparing $V_{\mu,r}^{q}(K)$ and $V_{\mu,r}^{q,\sym}(K)$.}
The definition
of
the symbolic multifractal Minkowski volume
may be viewed as a natural multifractal analogue
of the usual 
multifractal Mikskowski volume $V_{\mu,r}^{q}(K)$ 
 given by
 $$
  V_{\mu,r}^{q}
 (K)
 =
 \frac{1}{r^{d}}
 \int\limits_{B(K,r)}
 \mu(B(x,r))^{q}
 \,
 d\Cal L^{d}(x)\,.
 $$
 Indeed,
even though
the symbolic 
$q$ multifractal Minkowski volume
$V_{\mu,r}^{q,\sym}(K)$
does not necessarily
equal the usual
$q$ multifractal Minkowski volume
$V_{\mu,r}^{q}(K)$,
it is nevertheless 
\lq\lq compatible"
with $V_{\mu,r}^{q}(K)$.
More precisely,
the usual $q$ multifractal Minkowski volume
and
the
symbolic
$q$ multifractal Minkowski volume
give rise to the same dimensions.
This is the content of Theorem 5.1.
While 
this is not a deep fact 
it may be seen as providing 
further
 justification for the study
the
symbolic $q$ multifractal Minkowski volume.

\bigskip

\proclaim{Theorem 5.1}
Let $q\in\Bbb R$.
Recall that we define 
the 
 lower and upper
 multifractal $q$ Minkowski  dimension of $K$ by
 $$
 \aligned
 \underline\dim_{\Min,\mu}^{q}(K)
&=
 \,
 \liminf_{r\searrow0}
 \,
 \frac{\log V_{\mu,r}^{q}(K)}{-\log r}\,,\\
 \overline\dim_{\Min,\mu}^{q}(K)
&=
 \limsup_{r\searrow0}
 \frac{\log V_{\mu,r}^{q}(K)}{-\log r}\,.
 \endaligned
 $$
Similarly, 
we define 
the symbolic 
 lower and upper
 multifractal $q$ Minkowski  dimension of $K$ by 
 $$
 \aligned
 \underline\dim_{\Min,\mu}^{q,\sym}(K)
&=
 \,
 \liminf_{r\searrow0}
 \,
 \frac{\log V_{\mu,r}^{q,\sym}(K)}{-\log r}\,,\\
 \overline\dim_{\Min,\mu}^{q,\sym}(K)
&=
 \limsup_{r\searrow0}
 \frac{\log V_{\mu,r}^{q,\sym}(K)}{-\log r}\,.
 \endaligned
 $$
Assume 
that one of the following conditions is satisfied:
\roster
\item"(i)"
The OSC is satisfied and $0\le q$;
\item"(ii)"
The SSC is satisfied.
\endroster
Then we have
 $$
 \align
  \underline\dim_{\Min,\mu}^{q}(K)
 &=
  \underline\dim_{\Min,\mu}^{q,\sym}(K)\,,\\
  \overline\dim_{\Min,\mu}^{q}(K)
&=   
  \overline\dim_{\Min,\mu}^{q,\sym}(K)\,.
  \endalign
  $$
\endproclaim
\noindent{\it  Proof}\newline
\noindent
As noted above, this is not a deep fact and follows
from the definitions using standard arguments.
\hfill$\square$

\bigskip

In analogy with the usual multifractal Minkowski content
(see (3.3)--(3.6))
we 
also define symbolic multifractal Minkowski content.
For real numbers $q$ and $t$,
we define the
lower and upper $(q,t)$-dimensional symbolic multifractal Minkowski
content of $K$ with respect to $\mu$ by
 $$
 \aligned
 \underline M_{\mu}^{q,t,\sym}(K)
&=
 \,
 \liminf_{r\searrow0}
 \,
 \frac{1}{r^{-t}}
 \,
 V_{\mu,r}^{q,\sym}(K)\,,\\
 \overline M_{\mu}^{q,t,\sym}(K)
&=
 \limsup_{r\searrow0}
 \frac{1}{r^{-t}}
 \,
 V_{\mu,r}^{q,\sym}(K)\,.
 \endaligned
 \tag5.11
 $$
If 
$\underline M_{\mu}^{q,t,\sym}(K)
=
\overline M_{\mu}^{q,t,\sym}(K)$,
i\.e\. if the
limit
$ \lim_{r\searrow0}
 \,
 \frac{1}{r^{-t}}
 \,
 V_{\mu,r}^{q,\sym}(K)$
 exists, 
then we say that $K$ is $(q,t)$ symbolic multifractal Minkowski measurable
with respect to $\mu$,
and we 
denote the  
common value of
$\underline M_{\mu}^{q,t,\sym}(K)$ and
$\overline M_{\mu}^{q,t,\sym}(K)$
by
$M_{\mu}^{q,t,\sym}(K)$,
i\.e\. we write
 $$
 M_{\mu}^{q,t,\sym}(K)
 =
 \underline M_{\mu}^{q,t,\sym}(K)
 =
 \overline M_{\mu}^{q,t,\sym}(K)\,.
 \tag5.12
 $$
Of course,
$K$ may not be symbolic multifractal Minkowski measurable, and it is therefore useful
to introduce 
a suitable 
averaging procedure 
when computing the 
symbolic
multifractal Minkowski content.
Motivated by this we define 
 the
lower and upper $(q,t)$-dimensional symbolic average multifractal Minkowski
content of $K$ with respect to $\mu$ by
 $$
 \aligned
 \underline M_{\mu,\ave}^{q,t,\sym}(K)
&=
 \,
 \liminf_{r\searrow0}
 \frac{1}{-\log r}
 \int_{r}^{1}
 \frac{1}{s^{-t}}
 \,
 V_{\mu,s}^{q,\sym}(K)
 \,
 \frac{ds}{s}
 \,,\\
 \overline M_{\mu,\ave}^{q,t,\sym}(K)
&=
 \,
 \liminf_{r\searrow0}
 \frac{1}{-\log r}
 \int_{r}^{1}
 \frac{1}{s^{-t}}
 \,
 V_{\mu,s}^{q,\sym}(K)
 \,
 \frac{ds}{s}\,.
 \endaligned
 \tag5.13
 $$
If 
$\underline M_{\mu,\ave}^{q,t,\sym}(K)
=
\overline M_{\mu,\ave}^{q,t,\sym}(K)$,
i\.e\. if the
limit
$ \lim_{r\searrow0}
 \frac{1}{-\log r}
 \int_{r}^{1}
 \frac{1}{s^{-t}}
 \,
 V_{\mu,s}^{q,\sym}(K)
 \,
 \frac{ds}{s}$
 exists, 
then we say the $K$ is $(q,t)$ averagely symbolic multifractal Minkowski measurable
with respect to $\mu$,
and we 
denote the  
common value of
$\underline M_{\mu,\ave}^{q,t,\sym}(K)$ and
$\overline M_{\mu,\ave}^{q,t,\sym}(K)$
by
$M_{\mu,\ave}^{q,t,\sym}(K)$,
i\.e\. we  write
 $$
 M_{\mu,\ave}^{q,t,\sym}(K)
 =
 \underline M_{\mu,\ave}^{q,t,\sym}(K)
 =
 \overline M_{\mu,\ave}^{q,t,\sym}(K)\,.
 \tag5.14
 $$

How does one obtain explicit expressions for
$M_{\mu,}^{q,t,\sym}(K)$
and/or
$M_{\mu,\ave}^{q,t,\sym}(K)$?
The main tool for this is the notion of a
multifractal zeta-function.
A multifractal zeta-function is a certain meromorphic 
function whose residues are closely related to
the asymptotic behaviour of 
$V_{\mu,r}^{q,\sym}(K)$ as $r\searrow0$.
In the next section we define multifractal zeta-functions,
and the subsequent sections 
explain
how multifractal zeta-functions can be used
to analyse 
the symbolic multifratal Minkowski volume.

\bigskip

{\bf 5.2. Multifractal zeta-functions of self-similar measures --
a tool for finding explicit formulas for the symbolic multifratal Minkowski volume
of self-similar measures}.
Informally,
the multifractal zeta-function
$\zeta_{\mu}^{q}$ is defined by
 $$
 \zeta_{\mu}^{q}(s)
 =
 \sum_{\bold i}p_{\bold i}^{q}r_{\bold i}^{s}
 $$
for those complex numbers $s$ for which the series  
$ \sum_{\bold i}p_{\bold i}^{q}r_{\bold i}^{s}$
converges.
Formally, we proceed as described in the definition below.
In the definition below and in the later sections of the paper we use the following notation, 
namely,
if $f:\Bbb C\to\Bbb C$ is a complex valued function on $\Bbb C$, then we write
$Z(f)$ for the zeros of $f$, i\.e\. we write
 $$
 Z(f)
 =
 \Big\{s\in\Bbb C\,\Big|\,
 f(s)=0
 \Big\}\,.
 $$

\bigskip

\proclaim{Definition. The multifractal zeta-function $\zeta_{\mu}^{q}$
and 
the 
modified multifractal zeta-function
$Z_{\mu}^{q}$ }
Fix $q\in\Bbb R$.
For $s\in\Bbb C$ with $\real s>\beta(q)$,
 the series
$\sum_{\bold i}p_{\bold i}^{q}r_{\bold i}^{s}$
is convergent  with
 $$
 \sum_{\bold i}p_{\bold i}^{q}r_{\bold i}^{s}
 =
 \frac{\sum_{i}p_{i}^{q}r_{i}^{s}}{1-\sum_{i}p_{i}^{q}r_{i}^{s}}\,.
 \tag5.15
 $$
We can therefore define the multifractal zeta-function $\zeta_{\mu}^{q}$
by
$$
 \zeta_{\mu}^{q}(s)
 =
 \sum_{\bold i}p_{\bold i}^{q}r_{\bold i}^{s}
 \qquad
 \quad
  \qquad\qquad
 \qquad\qquad
 \quad\,\,\,\,\,
 \text{for $s\in\Big\{w\in\Bbb C\,\Big|\,\real s>\beta(q)\Big\}$.}
 \,\,\,\,\,\,\,\,
 \tag5.16
  $$
It follows from (5.15)
that $\zeta_{\mu}^{q}$ can be extended to 
$\Bbb C\setminus Z(w\to 1-\sum_{i}p_{i}^{q}r_{i}^{w})$ by
$$
{}\,
\,
 \zeta_{\mu}^{q}(s)
 =
 \frac{\sum_{i}p_{i}^{q}r_{i}^{s}}{1-\sum_{i}p_{i}^{q}r_{i}^{s}}
 \quad
 \qquad\qquad
 \qquad\qquad
 \quad\,\,\,\,\,\,
 \text{for $s\in\Bbb C\setminus Z\Big(w\to 1-\sum_{i}p_{i}^{q}r_{i}^{w}\Big)$.}
 \tag5.17
  $$
We define the modified zeta-function
$Z_{\mu}^{q}$ by
 $$
 Z_{\mu}^{q}(s)
 =
  \Bigg(
 \sum_{l}
 \frac
 {\kappa_{\mu}^{q,l}(K)\,(\sigma_{q,l}-1)}
 {s-(l-dq)}
 \Bigg)
 \,
 \zeta_{\mu}^{q}(s)
  \qquad\quad
 \text{for $s\in\Bbb C\setminus Z\Big(w\to 1-\sum_{i}p_{i}^{q}r_{i}^{w}\Big)$.}
 \tag5.18
 $$
\endproclaim
\noindent{\it  Proof}\newline
\noindent
For $s\in\Bbb C$ with $\real s>\beta(q)$,
we have
$|\sum_{i}p_{i}^{q}r_{i}^{s}|
\le
\sum_{i}p_{i}^{q}|r_{i}^{s}|
=
\sum_{i}p_{i}^{q}r_{i}^{\real s}
<1$,
and the series
$\sum_{\bold i}p_{\bold i}^{q}r_{\bold i}^{s}
=
\sum_{n}\sum_{|\bold i|=n} p_{\bold i}^{q}r_{\bold i}^{s}
=
\sum_{n}( \sum_{i}p_{i}^{q}r_{i}^{s})^{n}$
is therefore convergent  with
$\sum_{\bold i}p_{\bold i}^{q}r_{\bold i}^{s}
 =
 \frac{\sum_{i}p_{i}^{q}r_{i}^{s}}{1-\sum_{i}p_{i}^{q}r_{i}^{s}}$.
 \hfill$\square$

\bigskip

Zeta-functions 
similar
to $\zeta_{\mu}^{q}$
have appear 
frequently
in the study of dynamical systems, see, for example, 
[Lap-vF2,ParPo1,ParPo2,Rue1,Rue2] and the references therein.
In addition, we note
that Lapidus and collaborators have
introduced various intriguing multifractal zeta-functions [LapRo,LapLe-VeRo].
However, the multifractal zeta-functions in [LapRo,LapLe-VeRo]
serve very different purposes and are 
significantly different 
from the multifractal zeta-function
$\zeta_{\mu}^{q}$ introduced above.
Indeed,
our motivation for introducing
the  function
$\zeta_{\mu}^{q}$ is that
explicit formulas for 
 $V_{\mu,r}^{q,\sym}(K)$
 involving the residues of $\zeta_{\mu}^{q}$
can be obtained
by first use the Mellin transform to write
$V_{\mu,r}^{q,\sym}(K)$
as a complex contour integral of $\zeta_{\mu}^{q}$
and then
use the residue theorem to
express this integral as a sum involving 
the residues of $\zeta_{\mu}^{q}$.
As with all applications of the residue theorem,
this requires
information about the
 poles and residues of $\zeta_{\mu}^{q}$.
For this reason, 
the next section lists
some of the main properties of 
the poles and residues of $\zeta_{\mu}^{q}$.

\bigskip

{\bf 5.3. An intermezzo: the poles and residues of $\zeta_{\mu}^{q}$ and 
the sequence $(t_{q,n})_{n}$.}
For $q\in\Bbb R$,
 define
$\alpha(q)$ by
 $$
 \alpha(q)
 =
 \inf
 \left\{
 t\in\Bbb R
 \,
 \left|
 \,
 \sum
 \Sb
  i\\
  r_{i}=r_{\min}
 \endSb 
 p_{i}^{q}r_{\min}^{t}
 \le
  1
 +
 \sum
 \Sb
  i\\
  r_{i}>r_{\min}
 \endSb 
 p_{i}^{q}r_{i}^{t}
 \right.
 \right\}\,.
 \tag5.19
 $$
Also, recall that $\beta(q)$ is defined by
 $$
 \sum_{i}p_{i}^{q}r_{i}^{\beta(q)}
=1\,.
$$  
  Using the numbers $\alpha(q)$ and $\beta(q)$ we can now 
  describe the location of the poles of $\zeta_{\mu}^{q}$. 
Recall, that
if
$f$
is 
a meromorphic function, 
then $Z(f)$ denotes the set of zeros of $f$, i\.e\.
 $$
 Z(f)
 =
 \Big\{s\in\Bbb C\,\Big|\,f(s)=0\Big\}\,;
 \qquad\,\,\,\,\,
 $$
in addition, we let
$P(f)$ denotes the set of poles of $f$, i\.e\.
 $$
 P(f)
 =
 \Big\{s\in\Bbb C\,\Big|\,\text{$s$ is a pole of $f$} \Big\}\,.
 $$

\bigskip

\proclaim{Proposition 5.2.
The poles of $\zeta_{\mu}^{q}$}
Fix $q\in\Bbb R$.
\roster
\item"(1)"
We have $-\infty<\alpha(q)\le\beta(q)<\infty$.
\item"(2)"
We have
 $$
 P(\zeta_{\mu}^{q})
 =
 Z
 \Big(
 s\to 1-\sum_{i}p_{i}^{q}r_{i}^{s}
 \Big)\,.
 \qquad\qquad
 \quad\,\,
  $$
\item"(3)"
We have
   $$
 P(\zeta_{\mu}^{q})
\subseteq
 \Big\{
 s\in\Bbb C
 \,\Big|\,
 \alpha(q)
 \le
 \real(s)
 \le
 \beta(q)
 \Big\}\,.
 \,\,\,\,\,\,\,\,\,
 $$
\item"(4.1)"
Poles $\omega$ with
$\real(\omega)=\beta(q)$
in the
non-arithmetic case:
If 
the set
$\{\log r_{1}^{-1},\ldots,\log r_{N}^{-1}\}$
is not contained in a discrete additive subgroup of $\Bbb R$,
then
 $$
  P(\zeta_{\mu}^{q})
\cap \Big\{
 s\in\Bbb C
 \,\Big|\,
 \real(s)
 =
 \beta(q)
 \Big\}
 =
 \big\{
 \beta(q)
 \big\}\,.
 $$
\item"(4.2)"
Poles $\omega$ with
$\real(\omega)=\beta(q)$
in the
arithmetic case:
If 
the set
$\{\log r_{1}^{-1},\ldots,\log r_{N}^{-1}\}$
is contained in a discrete additive subgroup of $\Bbb R$
and
$\langle\log r_{1}^{-1},\ldots,\log r_{N}^{-1}\rangle=u\Bbb Z$
with $u>0$,
then
 $$
 {}\qquad\,\,\,\,\,\,\,\,\,
  P(\zeta_{\mu}^{q})
\cap 
\Big\{
 s\in\Bbb C
 \,\Big|\,
 \real(s)
 =
 \beta(q)
 \Big\}
 =
 \,\,\,\,\beta(q)\,+\tfrac{2\pi}{u}\ima\Bbb Z\,,
 $$
and for each $i$, there is a unique integer $k_{i}$ such that
$\log r_{i}^{-1}=k_{i}u$ 
and, in addition,
 $$
 \align
 {}\qquad\qquad
   \qquad\qquad
   \qquad
 P(\zeta_{\mu}^{q})
&=
 \Big(\,\beta(q)+\tfrac{2\pi}{u}\ima\Bbb Z\,\Big)\\
&\qquad\qquad
 \,\,
 \cup
 \,\,
 \bigcup
 \Sb
 w\in Z(z\to 1-\sum_{i}p_{i}^{q}z^{k_{i}})\\
 {}\\
 w\not=e^{-u\beta(q)}
 \endSb
 \Big(
 \,
 -\tfrac{\log |w|}{u}
 -
 \tfrac{\Arg w}{u}
 \ima
 +
 \tfrac{2\pi}{u}\ima \Bbb Z
 \,
 \Big)
 \endalign
 $$
(where $\Arg z$ denotes the principal argument of $z\in\Bbb C$).
\item"(5)"
Density of poles of $\zeta_{\mu}^{q}$:
Writing $\gamma=-\frac{1}{\pi}\log r_{\min}$, then we have
 $$
 \Big|
 \,
 \big\{\omega\in P(\zeta_{\mu}^{q})
 \,
 \big|
 \,
 |\imag(\omega)|\le t
 \big\}
 \,
 \Big|
 =
 \gamma t+\Cal O(\log t)\,.
 $$
 \endroster
\endproclaim

\bigskip

\noindent
Proposition 5.2 
play an important role in the proofs
of the main results in this section and  
is
proved in Section 13:
statements (1)--(4.2) 
are proved in Proposition 13.1
and the density statement (5) is proved in Theorem 13.8.
Statements (1)--(4.2)
are not deep and the proofs of those statements 
are straight forward.
However, the density statement (5)
requires a more involved and careful proof
based on the Argument Principle (see [Con, p\. 123]) together with
various variants
 of Jensen's formula from complex analysis
(see [Con, p\. 280]), namely Proposition 13.6 and Proposition 13.7.
We note that the intriguing
books 
by
Lapidus \& van Frankenhuysen
[Lap-vF1,Lap-vF2]
prove related density results.
Indeed, 
in the first editor of their book
[Lap-vF1],
Lapidus \& van Frankenhuysen
prove a density result similar to (5)
for the poles of a zeta-function related to $\zeta_{\mu}^{0}$
involving an error term of the form
$\Cal O(\sqrt{t})$.
However, in 
the second 
editor of their book
[Lap-vF2],
Lapidus \& van Frankenhuysen
present
a density
result for a larger class of zeta-functions
with an improved error term
of the form $\Cal O(1)$.
While we have not been able to 
prove the density statement (5) with an error term of the form $\Cal O(1)$,
we note that 
the error term 
in (5), namely $\Cal O(\log t)$,
in consistent with the results in
[JoLaGo, Theorem 6.2 and pp\. 58--59]
where related results 
(also involving error terms of the form $\Cal O(\log t)$) 
 are considered
in very general and abstract settings.

It follows from Proposition 5.2
that there is a critical strip
 $$
 \Cal S_{\crit}^{q}
 =
 \Big\{
 s\in\Bbb C
 \,\Big|\,
 \alpha(q)\le \real(s)\le\beta(q)
 \Big\}
 $$
and a 
critical
line
 $$
 \Cal L_{\crit}^{q}
 =
 \Big\{
 s\in\Bbb C
 \,\Big|\,
 \real(s)=\beta(q)
 \Big\}
 \qquad\quad
 $$
such that all poles of $\zeta_{\mu}^{q}$ belong to 
$\Cal S_{\crit}^{q}$
and lie
to the left of the 
line 
$\Cal L_{\crit}^{q}$.
We also observe 
that the nature of the poles
on the critical line is  determined 
by the algebraic properties
of (the logarithms of) the contracting ratios $r_{1},\ldots, r_{N}$.

Below we 
first use 
the Mellin transform
to obtain an explicit
formula for $V_{\mu,r}^{q,\sym}(K)$, namely  Theorem 5.4, 
expressing $V_{\mu,r}^{q,\sym}(K)$
as 
a complex contour integral
of $\zeta_{\mu}^{q}$.
Next,
we use the residue theorem
to \lq\lq compute"
the contour integral from Theorem 5.4
thus obtaining
the second explicit formula for $V_{\mu,r}^{q,\sym}(K)$,
namely Theorem 5.5,
expressing $V_{\mu,r}^{q,\sym}(K)$ as 
a sum of residues
of $\zeta_{\mu}^{q}$.
However, since all the 
 all the poles of $\zeta_{\mu}^{q}$  lie in the critical strip
$\Cal S_{\crit}^{q}$,
any application of the residue theorem to  $\zeta_{\mu}^{q}$
is likely to
involve integrating  $\zeta_{\mu}^{q}$
over line segments 
crossing the critical strip $\Cal S_{\crit}^{q}$.
For this reason precise
information about the
poles and the growth
of $\zeta_{\mu}^{q}$
on
line segments that
cross the critical strip $\Cal S_{\crit}^{q}$ is needed.
Such estimates are provided by the next result, i\.e\. Theorem 5.3 below.
This result
says that there is a sequence of horizontal line segments 
crossing the critical strip $\Cal S_{\crit}^{q}$
without
hitting the poles of $\zeta_{\mu}^{q}$
and such that
$\zeta_{\mu}^{q}$ is {\it uniformly}
bounded on these line segments.

\bigskip

\proclaim{Theorem 5.3. Growth estimates of $\zeta_{\mu}^{q}$
inside the critical strip 
$\alpha(q)\le\real(s)\le\beta(q)$
and the sequence $(t_{q,n})_{n}$}
Fix $q\in\Bbb R$.
Then there is an increasing sequence $(t_{q,n})_{n}$
of positive real numbers with $t_{q,n}\to\infty$
satisfying the following:
for all real numbers $c$, there
is a constant
$k_{c}$ such that
for all $\sigma\le c$ and all $n$,
we have
 $$
 \align
 |\zeta_{\mu}^{q}(\sigma\pm\ima t_{q,n})|
&\le
 k_{c}\,.
 \endalign
 $$
\endproclaim

\bigskip

\noindent
Theorem 5.3 is proven in Section 13.
The statement in
Theorem 5.3
is 
deep and relies on a number of very delicate estimates.
The reader is also referred to  [JoLaGo] for related results in
very general and abstract settings
and to
[Lap-vF2] for somewhat related results for $q=0$.

\bigskip

{\bf 5.4. Symbolic multifractal tubes of self-similar measures:
the
first explicit formula.}
Using the Mellin transform, we
first
obtain
an explicit formula for $V_{\mu,r}^{q,\sym}(K)$
expressing $V_{\mu,r}^{q,\sym}(K)$
as 
a complex contour integral
of $\zeta_{\mu}^{q}$.
Schematically this part can be represented as follows:

\bigskip

{\it
\centerline{The Mellin transform}

\centerline{$\Downarrow$}

\centerline{
$V_{\mu,r}^{q,\sym}(K)$
equals a complex contour integral of $\zeta_{\mu}^{q}$}
}

\bigskip

\noindent
More precisely, using
the Mellin transform technique,
we 
obtain
the
first explicit formula for $V_{\mu,r}^{q,\sym}(K)$, 
namely Theorem 5.4 below, 
expressing $V_{\mu,r}^{q,\sym}(K)$
as 
a complex contour integral
of the zeta-function $\zeta_{\mu}^{q}$.

\bigskip

\proclaim{Theorem 5.2. The first explicit formula for
$V_{\mu,r}^{q,\sym}(K)$} 
Fix $q\in\Bbb R$.
For
$c>\max\big(\,-dq,1-dq,\dots,d-dq,\beta(q)\,\big)$
and
$0<r<r_{\min}$, we have
 $$
 \align
 V_{\mu,r}^{q,\sym}(K)
&=
 \sum_{l}
  \kappa_{\mu}^{q,l}(K)
  \,
 \sigma_{q,l}r^{-l+dq}
  \,\,
 +
 \,\,
  \frac{1}{2\pi\ima}
 \int_{c-\ima \infty}^{c+\ima \infty}
Z_{\mu}^{q}(s)
 \,
 r^{-s}
 \,ds\,.
 \endalign
 $$
\endproclaim

\bigskip

\noindent
The proof of Theorem 5.4 is given in Section 14.

\bigskip

{\bf 5.5. Symbolic multifractal tubes of self-similar measures:
the
second explicit formula.}
Next,
using
the
first explicit formula for $V_{\mu,r}^{q,\sym}(K)$
expressing $V_{\mu,r}^{q,\sym}(K)$
as 
a complex contour integral
of $\zeta_{\mu}^{q}$
and the residue theorem,
we
obtain
an explicit formula for $V_{\mu,r}^{q,\sym}(K)$,
namely Theorem 5.5 below,
expressing $V_{\mu,r}^{q,\sym}(K)$ as 
a sum of residues
of $\zeta_{\mu}^{q}$.
Schematically this part can be represented as follows:

\bigskip

{\it
\centerline{
$V_{\mu,r}^{q,\sym}(K)$
equals a complex contour integral of $\zeta_{\mu}^{q}$}

\centerline{\&}

\centerline{
the residue theorem}

\centerline{$\Downarrow$}

\centerline{
$V_{\mu,r}^{q,\sym}(K)$
equals a sum of residues of $\zeta_{\mu}^{q}$}
}

\bigskip

\noindent
More precisely,  we now use 
the
first explicit formula for $V_{\mu,r}^{q,\sym}(K)$, i\.e\.  Theorem 5.4, 
and the residue theorem
together with the growth estimate provided by Theorem 5.3
to
obtain
the second explicit formula for $V_{\mu,r}^{q,\sym}(K)$,
namely Theorem 5.5 below,
expressing $V_{\mu,r}^{q,\sym}(K)$ as 
the sum of
a series involving the residues
of $\zeta_{\mu}^{q}$.

\bigskip

\proclaim{Theorem 5.5. The second explicit formula for
$V_{\mu,r}^{q,\sym}(K)$: a multifractal Steiner formula}
Fix $q\in\Bbb R$.
 Assume that
$\beta(q)
 \not\in
 \{
 -dq,1-dq,\ldots,d-dq\}$.
  Let $(t_{q,n})_{n}$
 be the sequence from Theorem 5.3.
For all $0<r< r_{\min}$, we have
 $$
 \align
  V_{\mu,r}^{q,\sym}(K)
&=
 \lim_{n}	
 \sum
 \Sb
 \omega
 \in P(\zeta_{\mu}^{q})\\
 {}\\
 |\imag(\omega)|\le t_{q,n}
 \endSb
 \res
 \Big(
 s\to
 Z_{\mu}^{q}(s)
 \,
 r^{-s};
 \omega
 \Big)\,. 
 \endalign
 $$
\endproclaim

\bigskip

\noindent
The proof of Theorem 5.5 is given in Section 15.
It is clear that
Theorem 5.5
may be viewed as 
a multifractal Steiner formula for 
the
symbolic multifractal Minkowski volume
of $V_{\mu,r}^{q,\sym}(K)$.
This viewpoint is perhaps even more clear in the following special
case.

\bigskip

\proclaim{Corollary  5.6. The second explicit formula for
$V_{\mu,r}^{q,\sym}(K)$: a multifractal Steiner formula}
Fix $q\in\Bbb R$.
 Assume that
$\beta(q)
 \not\in
 \{
 -dq,1-dq,\ldots,d-dq\}$
 and that all the poles of $\zeta_{\mu}^{q}$ are simple.
  Let $(t_{q,n})_{n}$
 be the sequence from Theorem 5.3.
For all $0<r< r_{\min}$, we have
 $$
 \align
  V_{\mu,r}^{q,\sym}(K)
&=
 \lim_{n}	
 \sum
 \Sb
 \omega
 \in P(\zeta_{\mu}^{q})\\
 {}\\
 |\imag(\omega)|\le t_{q,n}
 \endSb
 \Bigg(
 \sum_{l}
 \frac
 {\kappa_{\mu}^{q,l}(K)\,(\sigma_{q,l}-1)}
 {\omega-(l-dq)}
 \Bigg)
  \res
 (
 \zeta_{\mu}^{q}
 ;
 \omega
 )
 \,
 r^{-\omega}\,. 
 \endalign
 $$
\endproclaim 
\noindent{\it  Proof}\newline
\noindent
If $\omega$ is a simple pole of $\zeta_{\mu}^{q}$, then clearly
$\res
 (
 s\to
 Z_{\mu}^{q}(s)
 \,
 r^{-s};
 \omega
 )
 =
\res
 (
 s\to 
  (
 \sum_{l}
 \frac
 {\kappa_{\mu}^{q,l}(K)\,(\sigma_{q,l}-1)}
 {s-(l-dq)}
 )\allowmathbreak
 \,
 \zeta_{\mu}^{q}(s)
  \,
 r^{-s};
 \omega
 )
 =
  (
 \sum_{l}
 \frac
 {\kappa_{\mu}^{q,l}(K)\,(\sigma_{q,l}-1)}
 {\omega-(l-dq)}
 )
  \res
 (
 \zeta_{\mu}^{q}
 ;
 \omega
 )
 \,
 r^{-\omega}$. 
The desired result follows immediately from this observation and Theorem 5.5.
\hfill$\square$

\bigskip

\noindent
The 
formula in Corollary 5.6
has an even closer resemblance
to Steiner's formula.
Namely,
the
symbolic multifractal Minkowski volume
of $V_{\mu,r}^{q,\sym}(K)$
is written as a \lq\lq sum"
of powers of $r$.

\bigskip

{\bf 5.6. Symbolic multifractal tubes of self-similar measures:
the
third explicit formula.}
Finally, we
use
the second explicit formula for $V_{\mu,r}^{q,\sym}(K)$, i\.e\. Theorem 5.3,
expressing $V_{\mu,r}^{q,\sym}(K)$
as 
a sum of residues
of $\zeta_{\mu}^{q}$
together with 
 a very careful analysis of the residues of
$\zeta_{\mu}^{q}$
to 
prove Theorem 5.7
providing 
explicit formulas
for the limiting behaviour of $V_{\mu,r}^{q,\sym}(K)$
as $r\searrow0$.
Schematically this part can be represented as follows:

\bigskip

{\it
\centerline{
$V_{\mu,r}^{q,\sym}(K)$
equals a sum of residues of $\zeta_{\mu}^{q}$}

\centerline{\&}

\centerline{
a careful analysis of the residues of $\zeta_{\mu}^{q}$}

\centerline{$\Downarrow$}

\centerline{
an explicit formula for $\lim\limits_{r\searrow0}V_{\mu,r}^{q,\sym}(K)$}
}

\bigskip

\noindent
Recall, that
 the
second explicit formula for $V_{\mu,r}^{q,\sym}(K)$, i\.e\. Theorem 5.5, 
expressing $V_{\mu,r}^{q,\sym}(K)$
as 
a sum of residues
of $\zeta_{\mu}^{q}$ 
says that
if $q\in\Bbb R$
and
$\beta(q)
 \not\in
 \{
 -dq,1-dq,\ldots,d-dq
 \}$,
then 
 $$
 \align
 V_{\mu,r}^{q,\sym}(K)
 &=
   \lim_{n}	
 \sum
 \Sb
 \omega
 \in P(\zeta_{\mu}^{q})\\
 {}\\
 |\imag(\omega)|\le t_{q,n}
 \endSb
 \res
 \Big(
 s\to
 Z_{\mu}^{q}(s)
 \,
 r^{-s};
 \omega
 \Big)
 \endalign
 $$
for all $0<r<r_{\min}$,
and so
 $$
 \align
 \frac{1}{r^{-\beta(q)}}V_{\mu,r}^{q,\sym}(K)
&=
 \frac{1}{r^{-\beta(q)}}
   \lim_{n}	
 \sum
 \Sb
 \omega
 \in P(\zeta_{\mu}^{q})\\
 {}\\
 |\imag(\omega)|\le t_{q,n}
 \endSb
 \res
 \Big(
 s\to
 Z_{\mu}^{q}(s)
 \,
 r^{-s};
 \omega
 \Big)\,.
 \tag5.20
 \endalign
 $$
 for all $0<r<r_{\min}$.
The following heuristics  
suggests a plausible approach for analyzing
the asymptotic behaviour of
$ \frac{1}{r^{-\beta(q)}}V_{\mu,r}^{q,\sym}(K)$ as 
$r\searrow 0$.
Since all poles $\omega$ of $\zeta_{\mu}^{q}$ 
lie on or to the left of the critical line $\Cal L_{\crit}^{q}$ (see Section 5.3), 
i\.e\. since
$\real(\omega)\le\beta(q)$,
it is tempting to
split the sum in (5.20)
into following two parts, namely,
the sum of those 
$\omega\in P(\zeta_{\mu}^{q})$ for which $\real(\omega)<\beta(q)$
and the sum
of those 
$\omega\in P(\zeta_{\mu}^{q})$ for which $\real(\omega)=\beta(q)$,
i\.e\.
we attempt to write
 $$
 \align
 \frac{1}{r^{-\beta(q)}}V_{\mu,r}^{q,\sym}(K)
&=
  \frac{1}{r^{-\beta(q)}}
 \lim_{n}	
 \sum
 \Sb
 \omega
 \in P(\zeta_{\mu}^{q})\\
 {}\\
 |\imag(\omega)|\le t_{q,n}
 \endSb
 \res
 \Big(
 s\to
 Z_{\mu}^{q}(s)
 \,
 r^{-s};
 \omega
 \Big)\\
&=
 \pi_{q}^{\sym}(r)
 \,\,
 +
 \,\,
 E_{q}^{\sym}(r)
 \tag5.21
 \endalign
 $$
where
 $$
 \align
  \pi_{q}^{\sym}(r)
&= 
 \frac{1}{r^{-\beta(q)}}
  \lim_{n}	
 \sum
 \Sb
 \omega
 \in P(\zeta_{\mu}^{q})\\
 {}\\
 |\imag(\omega)|\le t_{q,n}\\
 {}\\
\real(\omega)=\beta(q)
 \endSb
 \res
 \Big(
 s\to
 Z_{\mu}^{q}(s)
 \,
 r^{-s};
 \omega
 \Big)\,,
 \tag5.22\\
&{}\\
 E_{q}^{\sym}(r)
&= 
 \frac{1}{r^{-\beta(q)}}
  \lim_{n}	
 \sum
 \Sb
 \omega
 \in P(\zeta_{\mu}^{q})\\
 {}\\
 |\imag(\omega)|\le t_{q,n}\\
 {}\\
\real(\omega)<\beta(q)
 \endSb
 \res
 \Big(
 s\to
 Z_{\mu}^{q}(s)
 \,
 r^{-s};
 \omega
 \Big)\,;
 \tag5.23
 \endalign
 $$
of course, at the moment we do not even know if the limits (5.22) and (5.23)
exist.
For  
 poles $\omega$ of $\zeta_{\mu}^{q}$
 with
 $\real(\omega)<\beta(q)$
 we have
 $\real(\beta(q)-\omega)>0$, 
 and it therefore seems
 plausible that
 if
$\omega$ is a pole of $\zeta_{\mu}^{q}$
 with
 $\real(\omega)<\beta(q)$, then
 $$
 \res
 \Big(
 s\to
 Z_{\mu}^{q}(s)
 \,
 r^{\beta(q)-s};
 \omega
 \Big)\to 0
 \,\,\,\,\text{ as $r\searrow0$,}
 $$
  suggesting that
 $$
 \align
 E_{q}^{\sym}(r)
&= 
 \frac{1}{r^{-\beta(q)}}
  \lim_{n}	
 \sum
 \Sb
 \omega
 \in P(\zeta_{\mu}^{q})\\
 {}\\
 |\imag(\omega)|\le t_{q,n}\\
 {}\\
 \real(\omega)<\beta(q)
 \endSb
 \res
 \Big(
 s\to
 Z_{\mu}^{q}(s)
 \,
 r^{-s};
 \omega
 \Big)\\
&= 
  \lim_{n}	
 \sum
 \Sb
 \omega
 \in P(\zeta_{\mu}^{q})\\
 {}\\
 |\imag(\omega)|\le t_{q,n}\\
 {}\\
 \real(\omega)<\beta(q)
 \endSb
 \res
 \Big(
 s\to
 Z_{\mu}^{q}(s)
 \,
 r^{\beta(q)-s};
 \omega
 \Big)\\
&\to
  \lim_{n}	
 \sum
 \Sb
 \omega
 \in P(\zeta_{\mu}^{q})\\
 {}\\
 |\imag(\omega)|\le t_{q,n}\\
 {}\\
\real(\omega)<\beta(q)
 \endSb
 0\\
&=
 0
 \,\,\,\, 
 \text{as $r\searrow0$.}
 \tag5.24
 \endalign
 $$
Finally, 
 combining (5.21) and (5.24)
indicates that
the asymptotic behaviour of 
 $$
 \align
 \frac{1}{r^{-\beta(q)}}V_{\mu,r}^{q,\sym}(K)
&=
  \pi_{q}^{\sym}(r)
 \,\,
 +
 \,\,
 E_{q}^{\sym}(r)
 \endalign
 $$
is determined  by
$ \pi_{q}^{\sym}(r)$.
In Section 16
we will prove that the heuristic argument
outlined above is
correct.
As suggested by the 
expressions
for
$  \pi_{q}^{\sym}(r)$
and
 $E_{q}^{\sym}(r)$,
 a rigorous analysis depends on 
 a very careful analysis of
 the structure of the poles of $\zeta_{\mu}^{q}$
 on and \lq\lq near"
 the 
 \lq\lq critical line"
  $\Cal L_{\crit}^{q}=\{s\in\Bbb C\,|\, \real(s)=\beta(q)\}$.
  Indeed, in order to
  prove that
  the limit $  \pi_{q}^{\sym}(r)$
 exists, a good understanding of the 
  poles $\omega$ of $\zeta_{\mu}^{q}$
  with $\real(\omega)=\beta(q)$ is needed
  and
 in order to
  prove that
  $E_{q}^{\sym}(r)$
tends to zero as $r\searrow 0$,
 a good understanding of the 
  poles $\omega$ of $\zeta_{\mu}^{q}$
  with $\real(\omega)$ close to but not equal to $\beta(q)$
   is needed.
  While a good understanding  of the 
  poles $\omega$ of $\zeta_{\mu}^{q}$
  with $\real(\omega)=\beta(q)$
  is easily obtained
  (and is, in fact, provided by Proposition 5.2),
  it is highly non-trivial to 
  obtain 
  a good understanding
    of the 
  poles $\omega$ of $\zeta_{\mu}^{q}$
  with $\real(\omega)$ close to but not equal to $\beta(q)$.
Indeed, a very substantial part of Section 13
is devoted to this problem.
  In particular, Proposition 13.2, Proposition 13.3 and Theorem 13.5
provide
detailed information about the structure 
of the 
poles and residues of $\zeta_{\mu}^{q}$
 near the critical line.
Finally, in Section 16 these results are used to prove 
Theorem 5.7
showing
that
  $E_{q}^{\sym}(r)$
tends to zero as $r\searrow 0$.
More precisely,  Theorem 5.7 below provides
a
complete description
of the
limiting behaviour of $V_{\mu,r}^{q,\sym}(K)$
as $r\searrow0$.
Namely,
if
the set
$\{\log r_{1}^{-1},\ldots,\log r_{N}^{-1}\}$
is not contained in a discrete additive subgroup of $\Bbb R$,
then 
$\frac{1}{r^{-\beta}}V_{\mu,r}^{q,\sym}(K)$
behaves asymptotically
as a multiplicatively periodic function
$\pi_{q}^{\sym}$
and we provide an explicit formula for
$\pi_{q}^{\sym}$,
and 
if
the set
$\{\log r_{1}^{-1},\ldots,\log r_{N}^{-1}\}$
is  contained in a discrete additive subgroup of $\Bbb R$,
then 
$\frac{1}{r^{-\beta}}V_{\mu,r}^{q,\sym}(K)$
converges to a constant
$c_{q}^{\sym}$
and we provide an explicit formula for
$c_{q}^{\sym}$.
We will now give the precise statement of
Theorem 5.7.
 In Theorem 5.7 we  write
 $\fraction(x)$ for the fractional part of a real number $x$.

\bigskip

\proclaim{Theorem 5.7}
Fix $q\in\Bbb R$.
Assume that
$\beta(q)
 \not\in
 \{
 -dq,1-dq,\ldots,d-dq\}$.

 \roster
\item"(1)"
The non-arithmetic case:
If 
the set
$\{\log r_{1}^{-1},\ldots,\log r_{N}^{-1}\}$
is not contained in a discrete additive subgroup of $\Bbb R$,
then
 $$
 \frac{1}{r^{-\beta(q)}}
 \,
 V_{\mu,r}^{q,\sym}(K)
 =
 c_{q}^{\sym}
 +
 \varepsilon_{q}(r)
 $$
for all $0<r<r_{\min}$ 
where
$c_{q}^{\sym}\in\Bbb R$ is the constant given by
 $$
 c_{q}^{\sym}
 =
 -
 \frac
 {1}
 {\sum_{i}p_{i}^{q}r_{i}^{\beta(q)}\log r_{i}}
\,
\sum_{l}\frac{\kappa_{\mu}^{q,l}(K)\,(\sigma_{q,l}-1)}
{\beta(q)-(l-dq)}
 $$ 
and
$\varepsilon_{q}^{\sym}(r)
 \to0$
 as $r\searrow0$.
 In addition, $K$ is
 $(q,\beta(q))$ symbolic multifractal Minkowski
 measurable with respect to $\mu$ with
  $$
  M_{\mu}^{q,\beta(q),\sym}(K)
  =
   -
 \frac
 {1}
 {\sum_{i}p_{i}^{q}r_{i}^{\beta(q)}\log r_{i}}
\,
\sum_{l}\frac{\kappa_{\mu}^{q,l}(K)\,(\sigma_{q,l}-1)}
{\beta(q)-(l-dq)}\,.
  $$
\item"(2)"
The arithmetic case:
If 
the set
$\{\log r_{1}^{-1},\ldots,\log r_{N}^{-1}\}$
is contained in a discrete additive subgroup of $\Bbb R$
and
$\langle\log r_{1}^{-1},\ldots,\log r_{N}^{-1}\rangle=u\Bbb Z$
with $u>0$,
then 
 $$
 \frac{1}{r^{-\beta(q)}}
 \,
 V_{\mu,r}^{q,\sym}(K)
 =
 \pi_{q}^{\sym}(r)
 +
 \varepsilon_{q}^{\sym}(r)
 $$
for all $0<r<r_{\min}$  
where
$\pi_{q}^{\sym}:(0,\infty)\to\Bbb R$ is the multiplicatively periodic function
with period equal to $e^{u}$
(i\.e\.
$\pi_{q}^{\sym}(e^{u}r)
 =
 \pi_{q}^{\sym}(r)$
for all 
$r\in(0,\infty)$)
given by
 $$
 \align
 \qquad\qquad
 \pi_{q}^{\sym}(r)
&=
 -
 \frac
 {1}
 {\sum_{i}p_{i}^{q}r_{i}^{\beta(q)}\log r_{i}}
 \,\\
&\qquad\qquad
   \qquad
 \times 
 \,
 u
 \,
 \sum_{l=0,1,\ldots,d}
 \,
 \Bigg(
 \,
 \frac
 {\kappa_{\mu}^{q,l}(K)\,(\sigma_{q,l}-1)}
 {1-e^{-u(\beta(q)-(l-dq))}}
 \\
&{}\\ 
&\qquad\qquad
   \qquad\qquad
 \times
 \cases
 \frac{e^{-u(\beta(q)-(l-dq))    }+1}{2}
&\quad
 \text
 {for
 $r\in e^{\Bbb Z u}$;
 }\\
 e^{-u(\beta(q)-(l-dq))\fraction(-\frac{\log r}{u})      }
&\quad
 \text
 {for 
 $r\not\in e^{\Bbb Z u}$
 }
 \endcases
 \Bigg)\,.
 \endalign
 $$
and
$\varepsilon_{q}(r)
 \to0$
 as $r\searrow0$.
 In addition, $K$ is
 $(q,\beta(q))$ average symbolic multifractal Minkowski
 measurable with respect to $\mu$ with
  $$
  M_{\mu,\ave}^{q,\beta(q),\sym}(K)
  =
   -
 \frac
 {1}
 {\sum_{i}p_{i}^{q}r_{i}^{\beta(q)}\log r_{i}}
\,
\sum_{l}\frac{\kappa_{\mu}^{q,l}(K)\,(\sigma_{q,l}-1)}
{\beta(q)-(l-dq)}\,.
  $$
\endroster
\endproclaim

\bigskip

\noindent
Theorem 5.7
is proved in Section 16.
Theorem 5.7
is clearly a \lq\lq symbolic" version of Theorem 3.3.
However, the key difference between Theorem 5.7
and Theorem 3.3
is that Theorem 5.7 (in contrast to Theorem 3.3)
has explicit 
formulas for the 
constant $c_{q}^{\sym}$ and the function $\pi_{q}^{\sym}$.

\newpage


${}$

\bigskip
\bigskip
\bigskip
\bigskip
\bigskip
\bigskip
\bigskip
\bigskip

\centerline{\bigletter Part 2:}

\bigskip

\centerline{\bigletter Proofs of the Results from Section 3}


\bigskip
\bigskip
\bigskip
\bigskip

\heading{6. Proving that
 $$
 \lambda_{q,m}(r)
 \le
 \sum
  \Sb
  |\bold i|=|\bold j|=m\\
  \bold i\not=\bold j
  \endSb
 Q_{\bold i,\bold j}^{q}(r)
 $$
}\endheading

\bigskip

The main purpose of this section is to prove Proposition 6.3.
However,
we begin by introducing some notation.
Let
$\Sigma
 =
 \{1,\ldots,N\}$ and write
 $$
 \aligned
 \Sigma^{m}
&=\{1,\dots,n\}^{m}\,,\\
 \Sigma^{*}
&=\bigcup_{m}\Sigma_{m}\,,\\
 \Sigma^{\Bbb N}
&=
 \{1,\dots,n\}^{\Bbb N}\,,
 \endaligned
 \tag6.1
 $$
i\.e\. $\Sigma^{m}$ is the family of all
strings
$\bold i=i_{1}\ldots i_{m}$
of length $m$ 
with $i_{j}\in\{1,\ldots,N\}$,
$\Sigma^{*}$ is the family of all finite strings
$\bold i=i_{1}\ldots i_{m}$
with $i_{j}\in\{1,\ldots,N\}$
and
$\Sigma^{*}$ is the family of all infinite strings
$\bold i=i_{1}i_{2}\ldots$
with $i_{j}\in\{1,\ldots,N\}$.
For $\bold i\in\Sigma^{m}$, we write $|\bold i|=m$
for the length of $\bold i$.
and for a positive integer $n$ with $n\le m$,
we 
write
$\bold i|n=i_{1}\ldots i_{n}$
for the
truncation
of $\bold i$ to the $n$-th place.
Also, for 
$\bold i=i_{1}\dots i_{m},
\bold j=j_{1}\dots j_{n}\in\Sigma^{*}$,
let
$\bold i\bold j=
i_{1}\dots i_{m}j_{1}\dots j_{n}$
denote the concatenation of $\bold i$ and $\bold j$.
Next, if $\bold i=i_{1}\dots i_{m}\in\Sigma^{*}$, we will write
$S_{\bold i}=S_{i_{1}}\circ\cdots\circ S_{i_{m}}$,
$r_{\bold i}=r_{i_{1}}\cdots r_{i_{m}}$ and
$p_{\bold i}=p_{i_{1}}\cdots p_{i_{m}}$.
Finally, write
 $$
 \align
 r_{\min}
&=
 \min_{i=1,\ldots,N}r_{i}\,,\\
r_{\max}
&=
 \max_{i=1,\ldots,N}r_{i}\,.
 \endalign
 $$

We now introduce the two 
key quantities in this (and the subsequent)
sections, namely,
$Q_{\bold i,\bold j}^{q}(r)$
and
$\lambda_{q,m}(r)$.
For
$E\subseteq\Bbb R^{d}$ and $r>0$,
recall that
$B(E,r)$
denotes the (open) $r$-neighbourhood of $E$,
i\.e\.
$B(E,r)
 =
 \{x\in\Bbb R^{d}\,|\,\dist(x,E)< r\}$.
For $q\in\Bbb R$ and
$\bold i,\bold j\in\Sigma^{*}$ and $r>0$, write
 $$
 Q_{\bold i,\bold j}^{q}(r)
 =
 \frac{1}{r^{d}}
 \,
 \int
 \limits_{
 B(S_{\bold i}K,r)\cap B(S_{\bold j}K,r)
 }
 \mu(B(x,r))^{q}
 \,
 d\Cal L^{d}(x)\,.
 \tag6.2
 $$
Next, for
$q\in\Bbb R$ and a positive integer $m$
and $r>0$, write
 $$
 \align
 \lambda_{q,m}(r)
&=
 \Cal I_{\mu,r}^{q}
 \big(\,B(K,r)\,\big)
 \,
 -
 \,
 \sum_{|\bold i|=m}
 p_{\bold i}^{q}
  \,
 \bold 1_{(0,r_{\bold i}]}(r)
\,
 \Cal I_{\mu,r_{\bold i}^{-1}r}^{q}
 \big(\,B(K,r_{\bold i}^{-1}r)\,\big)\\
&=
V_{\mu,r}^{q}(K)
 -
 \sum_{|\bold i|=m}p_{\bold i}^{q}
 \,
 \bold 1_{(0,r_{\bold i}]}(r)
\,
V_{\mu,r_{\bold i}^{-1}r}^{q}(K)\,,
\tag6.3
\endalign
$$                                             
If $m=1$, then we will
write
$\lambda_{}(r)$ for $\lambda_{q,1}(r)$, i\.e\. we will write
$$
\align
{}
\,\,
 \lambda_{q}(r)
&=
 \Cal I_{\mu,r}^{q}
 \big(\,B(K,r)\,\big)
 \,
 -
 \,
 \sum_{|i|=1}
 p_{\bold i}^{q}
  \,
 \bold 1_{(0,r_{i}]}(r)
\,
 \Cal I_{\mu,r_{i}^{-1}r}^{q}
 \big(\,B(K,r_{i}^{-1}r)\,\big)\\
&=
V_{\mu,r}^{q}(K)
 -
 \sum_{|i|=1}p_{i}^{q}
 \,
 \bold 1_{(0,r_{i}]}(r)
\,
V_{\mu,r_{i}^{-1}r}^{q}(K)\,,
\tag6.4
\endalign
$$

The main result in the section
is Proposition 6.3
providing an upper bound for
the difference
$\lambda_{q,m}(r)
 =
V_{\mu,r}^{q}(K)
 -
 \sum_{|\bold i|=m}p_{\bold i}^{q}
 \,
 \bold 1_{(0,r_{\bold i}]}(r)
\,
V_{\mu,r_{\bold i}^{-1}r}^{q}(K)$            
in terms of 
$Q_{\bold i,\bold j}^{q}(r)$; 
namely, in Proposition 6.3
we will prove that if $r>0$ is sufficiently small, than
 $$
 |\lambda_{q,m}(r)|
 \le
 \sum
  \Sb
  |\bold i|=|\bold j|=m\\
  \bold i\not=\bold j
  \endSb
 Q_{\bold i,\bold j}^{q}(r)\,.
 \tag6.5
 $$
We now turn towards the proof of Proposition 6.3.
We begin with two lemmas.

\bigskip

\proclaim{Lemma 6.1}
Fix $q\in\Bbb R$. Let $m\in\Bbb N$.
\roster
\item"(1)" For $r>0$,
we have
 $$
 \Cal I_{\mu,r}^{q}
 \big(\,B(K,r)\,\big)
 \le
 \sum_{|\bold i|=m}
 \Cal I_{\mu,r}^{q}
 \big(\,B(S_{\bold i}K,r)\,\big)\,.
 $$
\item"(2)" For $r>0$,
we have
 $$ 
 -
 \sum
 \Sb
 |\bold i|=|\bold j|=m\\
 \bold i\not=\bold j
 \endSb
 Q_{\bold i,\bold j}^{q}(r)
 \,
 +
 \,
 \sum_{|\bold i|=m}
 \Cal I_{\mu,r}^{q}
 \big(\,B(S_{\bold i}K,r)\,\big)
 \le
 \Cal I_{\mu,r}^{q}
 \big(\,B(K,r)\,\big)\,.
 $$
\endroster
\endproclaim 
\noindent{\it  Proof}\newline
\noindent 
(1)
Fix $r>0$.
Since 
$K=\cup_{|\bold i|=m}S_{\bold i}K$ we obtain
 $$
 \align
 \Cal I_{\mu,r}^{q}
 \big(\,B(K,r)\,\big)
&=
 \frac{1}{r^{d}}
 \,
 \int
 \limits_{
 B(K,r)
 }
 \mu(B(x,r))^{q}
 \,
 d\Cal L^{d}(x)\\
&= 
 \frac{1}{r^{d}}
 \,
 \int
 \limits_{
 B(\cup_{|\bold i|=m}S_{\bold i}K,r)
 }
 \mu(B(x,r))
 \,
 d\Cal L^{d}(x)\\
&\le
 \frac{1}{r^{d}}
 \,
 \int
 \limits_{
 \cup_{|\bold i|=m}B(S_{\bold i}K,r)
 }
 \mu(B(x,r))^{q}
 \,
 d\Cal L^{d}(x)\\
&\le
 \sum_{|\bold i|=m}
 \frac{1}{r^{d}}
 \,
 \int
 \limits_{
 B(S_{\bold i}K,r)
 }
 \mu(B(x,r))^{q}
 \,
 d\Cal L^{d}(x)\\
&\le
 \sum_{|\bold i|=m}
 \Cal I_{\mu,r}^{q}
 \big(\,B(S_{\bold i}K,r)\,\big)\,.
 \endalign
 $$

\noindent (2)
Fix $r>0$.
For $\bold i\in\Sigma^{*}$ with 
$|\bold i|=m$, write
 $$
 G_{\bold i}
 =
 B(S_{\bold i}K,r)
 \,
 \setminus
 \,
 \bigcup
 \Sb
 |\bold j|=|\bold i|\\
 \bold j\not=\bold i
 \endSb
 B(S_{\bold j}K,r)\,.
 $$
Also, 
for $\bold i,\bold j\in\Sigma^{*}$ with 
$|\bold i|=|\bold j|=m$, write
 $$
 H_{\bold i,\bold j}
 =
 B(S_{\bold i}K,r)
 \cap
 B(S_{\bold j}K,r)\,.
 $$

Since clearly $(G_{\bold i})_{|\bold i|=m}$
is a family of pairwise disjoint sets with
$\cup_{|\bold i|=m}G_{\bold i}\subseteq B(K,r)$, 
we conclude
that
 $$
 \align
 \Cal I_{\mu,r}^{q}
 \big(\,B(K,r)\,\big)
&=
 \frac{1}{r^{d}}
 \,
 \int
 \limits_{
 B(K,r)
 }
 \mu(B(x,r))^{q}
 \,
 d\Cal L^{d}(x)\\
&\ge
 \frac{1}{r^{d}}
 \,
 \int
 \limits_{
 \cup_{|\bold i|=m}G_{\bold i}
 }
 \mu(B(x,r))^{q}
 \,
 d\Cal L^{d}(x)\\
&=
 \sum_{|\bold i|=m}
 \frac{1}{r^{d}}
 \,
 \int
 \limits_{
 G_{\bold i}
 }
 \mu(B(x,r))^{q}
 \,
 d\Cal L^{d}(x)\,.
 \tag6.6
 \endalign
 $$

Next, note that
 $$
 B(S_{\bold i}K,r)
 \subseteq
 G_{\bold i}
 \,
 \cup
 \,
 \bigcup
 \Sb
 |\bold j|=|\bold i|\\
 \bold j\not=\bold i
 \endSb
 H_{\bold i,\bold j}\,.
 $$
It follows from this that
 $$
 \align
&\frac{1}{r^{d}}
 \,
 \int
 \limits_{
 G_{\bold i}
 }
 \mu(B(x,r))^{q}
 \,
 d\Cal L^{d}(x)
 \,
 +
 \,
 \sum
 \Sb
 |\bold j|=|\bold i|\\
 \bold j\not=\bold i
 \endSb
 \frac{1}{r^{d}}
 \,
 \int
 \limits_{
 H_{\bold i,\bold j}
 }
 \mu(B(x,r))^{q}
 \,
 d\Cal L^{d}(x)\\
&\qquad\qquad
 \qquad\qquad
 \qquad\qquad
 \qquad\qquad
 \ge
 \frac{1}{r^{d}}
 \,
 \int
 \limits_{
 G_{\bold i}
 \cup
 \cup_{|\bold j|=|\bold i|\,,\,\bold j\not=\bold i}
 H_{\bold i,\bold j}
 }
 \mu(B(x,r))^{q}
 \,
 d\Cal L^{d}(x)\\
&\qquad\qquad
 \qquad\qquad
 \qquad\qquad
 \qquad\qquad
 \ge
 \frac{1}{r^{d}}
 \,
 \int
 \limits_{
 B(S_{\bold i}K,r) 
 }
 \mu(B(x,r))^{q}
 \,
 d\Cal L^{d}(x)\,,
 \endalign
 $$
whence
 $$
 \align
&\frac{1}{r^{d}}
 \,
 \int
 \limits_{
 G_{\bold i}
 }
 \mu(B(x,r))^{q}
 \,
 d\Cal L^{d}(x)\\
&\qquad\qquad
 \ge
 \frac{1}{r^{d}}
 \,
 \int
 \limits_{
 B(S_{\bold i}K,r) 
 }
 \mu(B(x,r))^{q}
 \,
 d\Cal L^{d}(x)
 \,
 -
 \,
 \sum
 \Sb
 |\bold j|=|\bold i|\\
 \bold j\not=\bold i
 \endSb
 \frac{1}{r^{d}}
 \,
 \int
 \limits_{
 H_{\bold i,\bold j}
 }
 \mu(B(x,r))^{q}
 \,
 d\Cal L^{d}(x)\,.
 \tag6.7
 \endalign
 $$

Finally, combining (6.6) and (6.7) gives
 $$
 \align
 \Cal I_{\mu,r}^{q}
 \big(\,B(K,r)\,\big)
&\ge 
 \sum_{|\bold i|=m}
 \frac{1}{r^{d}}
 \,
 \int
 \limits_{
 G_{\bold i}
 }
 \mu(B(x,r))^{q}
 \,
 d\Cal L^{d}(x)\\
&\ge
 \sum_{|\bold i|=m}
 \left(
 \frac{1}{r^{d}}
 \,
 \int
 \limits_{
 B(S_{\bold i}K,r) 
 }
 \mu(B(x,r))^{q}
 \,
 d\Cal L^{d}(x)
 \,
 -
 \,
 \sum
 \Sb
 |\bold j|=|\bold i|\\
 \bold j\not=\bold i
 \endSb
 \frac{1}{r^{d}}
 \,
 \int
 \limits_{
 H_{\bold i,\bold j}
 }
 \mu(B(x,r))^{q}
 \,
 d\Cal L^{d}(x)
 \right)\\
&=
 \sum_{|\bold i|=m}
 \frac{1}{r^{d}}
 \,
 \int
 \limits_{
 B(S_{\bold i}K,r) 
 }
 \mu(B(x,r))^{q}
 \,
 d\Cal L^{d}(x)
 \,
 -
 \,
 \sum
 \Sb
 |\bold j|=|\bold i|=m\\
 \bold j\not=\bold i
 \endSb
 \frac{1}{r^{d}}
 \,
 \int
 \limits_{
 H_{\bold i,\bold j}
 }
 \mu(B(x,r))^{q}
 \,
 d\Cal L^{d}(x)\\
&=
 \sum_{|\bold i|=m}
 \Cal I_{\mu,r}^{q}
 \big(\,B(S_{\bold i}K,r)\,\big)
 \,
 -
 \,
 \sum
 \Sb
 |\bold i|=|\bold j|=m\\
 \bold i\not=\bold j
 \endSb
 Q_{\bold i,\bold j}^{q}(r)\,.
 \endalign
 $$
This completes the proof. 
\hfill$\square$

\bigskip

\proclaim{Lemma 6.2}
Fix $q\in\Bbb R$.
Let $\bold i\in\Sigma^{*}$.
\roster
\item"(1)" For $r>0$,
we have
 $$
 \Cal I_{\mu,r}^{q}
 \big(\,B(S_{\bold i}K,r)\,\big)
 \le
 p_{\bold i}^{q}
 \Cal I_{\mu,r_{\bold i}^{-1}r}^{q}
 \big(\,B(K,r_{\bold i}^{-1}r)\,\big)
 \,
 +
 \,
 \sum
 \Sb
 |\bold j|=|\bold i|\\
 \bold j\not=\bold i
 \endSb
 Q_{\bold i,\bold j}^{q}(r)\,.
 $$
\item"(2)" For $0\le q$ and $r>0$,
we have
 $$ 
 p_{\bold i}^{q}
 \Cal I_{\mu,r_{\bold i}^{-1}r}^{q}
 \big(\,B(K,r_{\bold i}^{-1}r)\,\big)
 \le
 \Cal I_{\mu,r}^{q}
 \big(\,B(S_{\bold i}K,r)\,\big)\,.
 $$
\endroster
\endproclaim 
\noindent{\it  Proof}\newline
\noindent 
(1)
Fix $r>0$.
Write
 $$
 G
 =
 B(S_{\bold i}K,r)
 \setminus
 \bigcup
 \Sb
 |\bold j|=|\bold i|\\
 \bold j\not=\bold i
 \endSb
 B(S_{\bold j}K,r)\,.
 $$
Also,
for $\bold j\in\Sigma^{*}$ with
 $|\bold j|=|\bold i|$, write
 $$
 H_{\bold j}
 =
 B(S_{\bold i}K,r)
 \cap
 B(S_{\bold j}K,r).
 $$

Note that
 $$
 B(S_{\bold i}K,r)
 \subseteq
 G
 \,
 \cup
 \,
 \bigcup
 \Sb
 |\bold j|=|\bold i|\\
 \bold j\not=\bold i
 \endSb
 H_{\bold j}\,.
 $$
It follows from this 
and the equation
$\mu(B(x,r))
=
\sum_{|\bold j|=|\bold i|}p_{\bold j}\mu(S_{\bold j}^{-1}B(x,r))$,
that
 $$
 \align
 \Cal I_{\mu,r}^{q}
&\big(\,B(S_{\bold i}K,r)\,\big)\\
&=
 \frac{1}{r^{d}}
 \,
 \int
 \limits_{
 B(S_{\bold i}K,r)
 }
 \mu(B(x,r))^{q}
 \,
 d\Cal L^{d}(x)\\
&\le
 \frac{1}{r^{d}}
 \,
 \int
 \limits_{
 G
 \cup
 \cup_{
 |\bold j|=|\bold i|,
 \bold j\not=\bold i
 }
 H_{\bold j}
 }
 \mu(B(x,r))^{q}
 \,
 d\Cal L^{d}(x)\\
&\le 
 \frac{1}{r^{d}}
 \,
 \int
 \limits_{
 G
 }
 \mu(B(x,r))^{q}
 \,
 d\Cal L^{d}(x)
 \,
 +
 \,
 \sum
 \Sb
 |\bold j|=|\bold i|\\
 \bold j\not=\bold i
 \endSb
 \frac{1}{r^{d}}
 \,
 \int
 \limits_{
 H_{\bold j}
 }
 \mu(B(x,r))^{q}
 \,
 d\Cal L^{d}(x)\\
&= 
 \frac{1}{r^{d}}
 \,
 \int
 \limits_{
 G
 }
 \left(
 \sum_{|\bold j|=|\bold i|}p_{\bold j}\mu(S_{\bold j}^{-1}B(x,r))
 \right)^{q}
 \,
 d\Cal L^{d}(x)
 \,
 +
 \,
 \sum
 \Sb
 |\bold j|=|\bold i|\\
 \bold j\not=\bold i
 \endSb
 \frac{1}{r^{d}}
 \,
 \int
 \limits_{
 H_{\bold j}
 }
 \mu(B(x,r))^{q}
 \,
 d\Cal L^{d}(x)\,.
 \tag6.8  
 \endalign
 $$
Now
observe that
$S_{\bold j}^{-1}B(x,r)=\varnothing$ for all $x\in G$ and all
$\bold j\in\Sigma^{*}$
with
$|\bold j|=|\bold i|$ and $\bold j\not=\bold i$.
It follows from this 
that
$\sum_{|\bold j|=|\bold i|}p_{\bold j}\mu(S_{\bold j}^{-1}B(x,r))
=
p_{\bold i}\mu(S_{\bold i}^{-1}B(x,r),r)$,
and (6.8) therefore simplifies to
 $$
 \align
 \Cal I_{\mu,r}^{q}
 \big(\,B(S_{\bold i}K,r)\,\big)
&\le
 p_{\bold i}^{q}
 \,
 \frac{1}{r^{d}}
 \,
 \int
 \limits_{
 G
 }
 \mu(S_{\bold i}^{-1}B(x,r))^{q}
 \,
 d\Cal L^{d}(x)
 \,
 +
 \,
 \sum
 \Sb
 |\bold j|=|\bold i|\\
 \bold j\not=\bold i
 \endSb
 \frac{1}{r^{d}}
 \,
 \int
 \limits_{
 H_{\bold j}
 }
 \mu(B(x,r))^{q}
 \,
 d\Cal L^{d}(x)\\
&=
 p_{\bold i}^{q}
 \,
 \frac{1}{r^{d}}
 \,
 \int
 \limits_{
 G
 }
 \mu(S_{\bold i}^{-1}B(x,r))^{q}
 \,
 d\Cal L^{d}(x)
 \,
 +
 \,
 \sum
 \Sb
 |\bold j|=|\bold i|\\
 \bold j\not=\bold i
 \endSb
 Q_{\bold i,\bold j}(r)\\
&=
 p_{\bold i}^{q}
 \,
 \frac{1}{r^{d}}
 \,
 \int
 \limits_{
 G
 }
 \mu(B(S_{\bold i}^{-1}x,r_{\bold i}^{-1}r))^{q}
 \,
 d\Cal L^{d}(x)
 \,
 +
 \,
 \sum
 \Sb
 |\bold j|=|\bold i|\\
 \bold j\not=\bold i
 \endSb
 Q_{\bold i,\bold j}(r)\\ 
&=
 p_{\bold i}^{q}
 \,
 \frac{1}{r^{d}}
 \,
 r_{\bold i}^{d}
 \int
 \limits_{
 S_{\bold i}^{-1}G
 }
 \mu(B(x,r_{\bold i}^{-1}r))^{q}
 \,
 d\Cal L^{d}(x)
 \,
 +
 \,
 \sum
 \Sb
 |\bold j|=|\bold i|\\
 \bold j\not=\bold i
 \endSb
 Q_{\bold i,\bold j}(r)\,.
 \tag6.9
 \endalign
 $$
Finally, using the fact that
$S_{\bold i}^{-1}G\subseteq
S_{\bold i}^{-1}B(S_{\bold i}K,r)\subseteq
B(K,r_{\bold i}^{-1}r)$, we conclude from (6.9) that
 $$
 \align
 \Cal I_{\mu,r}^{q}
 \big(\,B(S_{\bold i}K,r)\,\big)
&\le
 p_{\bold i}^{q}
 \,
 \frac{1}{(r_{\bold i}^{-1}r)^{d}}
 \,
 \int
 \limits_{
 B(K,r_{\bold i}^{-1}r)
 }
 \mu(B(x,r_{\bold i}^{-1}r))^{q}
 \,
 d\Cal L^{d}(x)
 \,
 +
 \,
 \sum
 \Sb
 |\bold j|=|\bold i|\\
 \bold j\not=\bold i
 \endSb
 Q_{\bold i,\bold j}(r)\\
&=
 p_{\bold i}^{q}
 \Cal I_{\mu,r_{\bold i}^{-1}r}^{q}
 \big(\,B(K,r_{\bold i}^{-1}r)\,\big)
 \,
 +
 \,
 \sum
 \Sb
 |\bold j|=|\bold i|\\
 \bold j\not=\bold i
 \endSb
 Q_{\bold i,\bold j}^{q}(r)\,.
 \endalign
 $$

\noindent (2)
Fix $r>0$.
Write
 $$
 G
 =
 B(S_{\bold i}K,r)
 \setminus
 \bigcup
 \Sb
 |\bold j|=|\bold i|\\
 \bold j\not=\bold i
 \endSb
 B(S_{\bold j}K,r)\,.
 $$
Also, write
 $$
 H
 =
 B(S_{\bold i}K,r)
 \,
 \cap
 \,
 \bigcup
 \Sb
 |\bold j|=|\bold i|\\
 \bold j\not=\bold i
 \endSb
 B(S_{\bold j}K,r)\,.
 $$

Note that $B(S_{\bold i}K,r)=G\cup H$
and
$G\cap H=\varnothing$.
It follows from this 
and the identity
$\mu(B(x,r))
=
\sum_{|\bold j|=|\bold i|}p_{\bold j}\mu(S_{\bold j}^{-1}B(x,r))$,
that
 $$
 \align
 \Cal I_{\mu,r}^{q}
 \big(\,B(S_{\bold i}K,r)\,\big)
&=
 \frac{1}{r^{d}}
 \,
 \int
 \limits_{
 B(S_{\bold i}K,r)
 }
 \mu(B(x,r))^{q}
 \,
 d\Cal L^{d}(x)\\
&=
 \frac{1}{r^{d}}
 \,
 \int
 \limits_{
 G
 }
 \mu(B(x,r))^{q}
 \,
 d\Cal L^{d}(x)
 \,
 +
 \,
 \frac{1}{r^{d}}
 \,
 \int
 \limits_{
 H
 }
 \mu(B(x,r))^{q}
 \,
 d\Cal L^{d}(x)\\ 
&=
 \frac{1}{r^{d}}
 \,
 \int
 \limits_{
 G
 }
 \left(
 \sum_{|\bold j|=|\bold i|}p_{\bold j}\mu(S_{\bold j}^{-1}B(x,r))
 \right)^{q}
 \,
 d\Cal L^{d}(x)
 \,
 +
 \,
 \frac{1}{r^{d}}
 \,
 \int
 \limits_{
 H
 }
 \mu(B(x,r))^{q}
 \,
 d\Cal L^{d}(x)\,.\\
&{}
 \tag6.10
  \endalign
 $$
Now
observe that
$S_{\bold j}^{-1}B(x,r)=\varnothing$ for all $x\in G$ and all
$\bold j\in\Sigma^{*}$
with
$|\bold j|=|\bold i|$ and $\bold j\not=\bold i$.
It follows from this 
that
$\sum_{|\bold j|=|\bold i|}p_{\bold j}\mu(S_{\bold j}^{-1}B(x,r))
=
p_{\bold i}\mu(S_{\bold i}^{-1}B(x,r))$,
and (6.10) therefore simplifies to 
 $$
 \align
 \Cal I_{\mu,r}^{q}
 \big(\,B(S_{\bold i}K,r)\,\big)
&=
 p_{\bold i}^{q}
 \frac{1}{r^{d}}
 \,
 \int
 \limits_{
 G
 }
 \mu(S_{\bold i}^{-1}B(x,r))^{q}
 \,
 d\Cal L^{d}(x)
 \,
 +
 \,
 \frac{1}{r^{d}}
 \,
 \int
 \limits_{
 H
 }
 \mu(B(x,r))^{q}
 \,
 d\Cal L^{d}(x)\,.\\
&{} 
 \tag6.11
 \endalign
 $$ 
Once more using the fact that
that $B(S_{\bold i}K,r)=G\cup H$
and
$G\cap H=\varnothing$, we conclude from (6.11) that
 $$
 \align
 \Cal I_{\mu,r}^{q}
 \big(\,B(S_{\bold i}K,r)\,\big)
&=
 p_{\bold i}^{q}
 \frac{1}{r^{d}}
 \,
 \int
 \limits_{
 G\cup H
 }
 \mu(S_{\bold i}^{-1}B(x,r))^{q}
 \,
 d\Cal L^{d}(x)\\
&\qquad\qquad 
 \,
 -
 \,
 p_{\bold i}^{q}
 \frac{1}{r^{d}}
 \,
 \int
 \limits_{
 H
 }
 \mu(S_{\bold i}^{-1}B(x,r))^{q}
 \,
 d\Cal L^{d}(x)
 \,
 +
 \,
 \frac{1}{r^{d}}
 \,
 \int
 \limits_{
 H
 }
 \mu(B(x,r))^{q}
 \,
 d\Cal L^{d}(x)\\
&=
 p_{\bold i}^{q}
 \frac{1}{r^{d}}
 \,
 \int
 \limits_{
 B(S_{\bold i}K,r)
 }
 \mu(S_{\bold i}^{-1}B(x,r))^{q}
 \,
 d\Cal L^{d}(x)\\
&\qquad\qquad 
 \,
 +
 \,
 \frac{1}{r^{d}}
 \,
 \int
 \limits_{
 H
 }
 \big(
 \mu(B(x,r))^{q}
 -
 p_{\bold i}^{q}\mu(S_{\bold i}^{-1}B(x,r))^{q}
 \big)
 \,
 d\Cal L^{d}(x)\,. 
 \tag6.12
 \endalign
 $$ 
Finally, note that
since $0\le q$, we conclude that
$\mu(B(x,r))^{q}
 =
 (
 \sum_{|\bold j|=|\bold i|}p_{\bold j}\mu(S_{\bold j}^{-1}B(x,r))
 )^{q}
 \allowmathbreak
 \ge
 p_{\bold i}^{q}\mu(S_{\bold i}^{-1}B(x,r))^{q}$
for all $x$ and all $r>0$. 
It therefore follows from (6.12) that
 $$
 \align
 \Cal I_{\mu,r}^{q}
 \big(\,B(S_{\bold i}K,r)\,\big)
&\ge
 p_{\bold i}^{q}
 \frac{1}{r^{d}}
 \,
 \int
 \limits_{
 B(S_{\bold i}K,r)
 }
 \mu(S_{\bold i}^{-1}B(x,r))^{q}
 \,
 d\Cal L^{d}(x)\\
&=
 p_{\bold i}^{q}
 \frac{1}{r^{d}}
 \,
 \int
 \limits_{
 B(S_{\bold i}K,r)
 }
 \mu(B(S_{\bold i}^{-1}x,r_{\bold i}^{-1}r))^{q}
 \,
 d\Cal L^{d}(x)\\ 
&=
 p_{\bold i}^{q}
 \frac{1}{r^{d}}
 \,
 r_{\bold i}^{d}
 \,
 \int
 \limits_{
 S_{\bold i}^{-1}B(S_{\bold i}K,r)
 }
 \mu(B(x,r_{\bold i}^{-1}r))^{q}
 \,
 d\Cal L^{d}(x)\\  
&=
 p_{\bold i}^{q}
 \frac{1}{(r_{\bold i}^{-1}r)^{d}}
 \,
 \int
 \limits_{
 B(K,r_{\bold i}^{-1}r)
 }
 \mu(B(x,r_{\bold i}^{-1}r))^{q}
 \,
 d\Cal L^{d}(x)\\   
&=
 \Cal I_{\mu,r_{\bold i}^{-1}r}^{q}
 \big(\,B(K,r_{\bold i}^{-1}r)\,\big)\,.
 \endalign
 $$ 
This completes the proof.
\hfill$\square$

\bigskip

\proclaim{Proposotion 6.3}
Fix
$q\in\Bbb R$.
Let $m\in\Bbb N$.
\roster
\item"(1)"
If $0\le q$ and the OSC is satisfied, then
for $r>0$, we have
 $$
 \left|
 \Cal I_{\mu,r}^{q}
 \big(\,B(K,r)\,\big)
 \,
 -
 \,
 \sum_{|\bold i|=m}
 p_{\bold i}^{q}
 \Cal I_{\mu,r_{\bold i}^{-1}r}^{q}
 \big(\,B(K,r_{\bold i}^{-1}r)\,\big)
 \right|
 \,\,
 \le
 \,\,
 \sum
 \Sb
 |\bold i|=|\bold j|=m\\
 \bold i\not=\bold j
 \endSb
 Q_{\bold i,\bold j}^{q}(r)\,.
$$
In particular, for 
$r>0$ with $r<r_{\min}^{m}$, we have
 $$
 |\lambda_{q,m}(r)|
 \,\,
 \le
 \,\,
 \sum
 \Sb
 |\bold i|=|\bold j|=m\\
 \bold i\not=\bold j
 \endSb
 Q_{\bold i,\bold j}^{q}(r)\,.
$$
\item"(2)"
If the SSC is satisfied, then
for $r>0$ with 
$r<\frac{1}{2}\min_{|\bold i|=|\bold j|=m,\bold i\not=\bold j}\dist(S_{\bold i}K,S_{\bold j}K)$,
 we have
 $$
 \Cal I_{\mu,r}^{q}
 \big(\,B(K,r)\,\big)
 \,
 -
 \,
 \sum_{|\bold i|=m}
 p_{\bold i}^{q}
 \Cal I_{\mu,r_{\bold i}^{-1}r}^{q}
 \big(\,B(K,r_{\bold i}^{-1}r)\,\big)
 \,\,
 =
 \,\,
 0\,.
 {}
 \qquad\qquad
 \,\,\,\,\,\,
$$
In particular, for 
$r>0$ with 
$r<
\min\big(\,
r_{\min}^{m}
\,,\,
\frac{1}{2}\min_{|\bold i|=|\bold j|=m,\bold i\not=\bold j}\dist(S_{\bold i}K,S_{\bold j}K)
\,\big)$, we have
 $$
 \lambda_{q,m}(r)
 \,\,
 =
 \,\,
 0\,.
 {}\qquad\qquad
 \,\,\,\,\,\,
$$
\endroster
\endproclaim 
\noindent{\it  Proof}\newline
\noindent 
(1)
This
follows immediately from Lemma 6.1 and Lemma 6.2.

\noindent
(2)
Since $r<\frac{1}{2}\min_{|\bold i|=|\bold j|=m,\bold i\not=\bold j}\dist(S_{\bold i}K,S_{\bold j}K)$
and
$K=\cup_{|\bold i|=m}S_{\bold i}K$,
we clearly have
 $$
 \align
 \Cal I_{\mu,r}^{q}
 \big(\,B(K,r)\,\big)
&=
 \frac{1}{r^{d}}
 \,
 \int
 \limits_{
 B(K,r)
 }
 \mu(B(x,r))^{q}
 \,
 d\Cal L^{d}(x)\\
&= 
 \frac{1}{r^{d}}
 \,
 \int
 \limits_{
 B(\cup_{|\bold i|=m}S_{\bold i}K,r)
 }
 \mu(B(x,r))
 \,
 d\Cal L^{d}(x)\\
&=
 \frac{1}{r^{d}}
 \,
 \int
 \limits_{
 \cup_{|\bold i|=m}B(S_{\bold i}K,r)
 }
 \mu(B(x,r))^{q}
 \,
 d\Cal L^{d}(x)\\
&=
 \sum_{|\bold i|=m}
 \frac{1}{r^{d}}
 \,
 \int
 \limits_{
 B(S_{\bold i}K,r)
 }
 \mu(B(x,r))^{q}
 \,
 d\Cal L^{d}(x)\\
&= 
\sum_{|\bold i|=m}
 \frac{1}{r^{d}}
 \,
 \int
 \limits_{
 B(S_{\bold i}K,r)
 }
  \left(
 \sum_{|\bold j|=|\bold i|}p_{\bold j}\mu(S_{\bold j}^{-1}B(x,r))
 \right)^{q}
 \,
 d\Cal L^{d}(x)\,.
 \tag6.13
 \endalign
 $$
Next,
observe that
since
$r<\frac{1}{2}\min_{|\bold i|=|\bold j|=m,\bold i\not=\bold j}\dist(S_{\bold i}K,S_{\bold j}K)$,
we conclude that
$S_{\bold j}^{-1}B(x,r)=\varnothing$ for all 
$x\in B(S_{\bold i}K,r)$ for all
$\bold i,\bold j\in\Sigma^{*}$
with
$|\bold j|=|\bold i|$ and $\bold j\not=\bold i$.
It follows from this 
that
$\sum_{|\bold j|=|\bold i|}p_{\bold j}\mu(S_{\bold j}^{-1}B(x,r))
=
p_{\bold i}\mu(S_{\bold i}^{-1}B(x,r))$,
and (6.13) therefore simplifies to
 $$
 \align
 \Cal I_{\mu,r}^{q}
 \big(\,B(K,r)\,\big)
&=
 \sum_{|\bold i|=m}
 p_{\bold i}^{q}
 \,
 \frac{1}{r^{d}}
 \,
 \int
 \limits_{
 B(S_{\bold i}K,r)
 }
 \mu(S_{\bold i}^{-1}B(x,r))^{q}
 \,
 d\Cal L^{d}(x)\\
&=
 \sum_{|\bold i|=m}
 p_{\bold i}^{q}
 \,
 \frac{1}{r^{d}}
 \,
 \int
 \limits_{
 B(S_{\bold i}K,r)
 }
 \mu(B(S_{\bold i}^{-1}x,r_{\bold i}^{-1}r))^{q}
 \,
 d\Cal L^{d}(x)\\ 
&=
 \sum_{|\bold i|=m}
 p_{\bold i}^{q}
 \,
 \frac{1}{r^{d}}
 \,
 r_{\bold i}^{d}
 \int
 \limits_{
 S_{\bold i}^{-1}B(S_{\bold i}K,r)
 }
 \mu(B(x,r_{\bold i}^{-1}r))^{q}
 \,
 d\Cal L^{d}(x)\\
&=
 \sum_{|\bold i|=m}
 p_{\bold i}^{q}
 \,
 \frac{1}{(r_{\bold i}^{-1}r)^{d}}
 \int
 \limits_{
 B(K,r_{\bold i}^{-1}r)
 }
 \mu(B(x,r_{\bold i}^{-1}r))^{q}
 \,
 d\Cal L^{d}(x)\\
&=
 \sum_{|\bold i|=m}
 p_{\bold i}^{q}
 \Cal I_{\mu,r_{\bold i}^{-1}r}^{q}
 \big(\,B(K,r_{\bold i}^{-1}r)\,\big)\,. 
 \endalign
 $$
This completes the proof.
\hfill$\square$

\bigskip
\bigskip


\heading{7. Proving that
 $$
 Q_{\bold i,\bold j}^{q}(r)
 \le
 \text{\rm constant}
 \,
 Z_{m}^{q}(r)$$
}\endheading

\bigskip

The main purpose of this section is to prove
Proposition 7.8.
However,
we begin by introduction some notation.
For $\bold i,\bold h\in\Sigma^{*}$, we write
 $$
 \bold i\prec\bold h
 $$
if and only if $\bold i$ is  a substring of $\bold h$, i\.e\.
if and only if there are strings 
$\bold s,\bold t\in\Sigma^{*}$
such that
 $$
 \bold h
 =
 \bold s\bold i\bold t\,.
 $$

Next,
if $(S_{1},\dots,S_{n})$ satisfies the OSC, then it follows from
a result by Schief [Schi]
that there exists an open, bounded and non-empty 
subset $U$ of $\Bbb R^{d}$ with
 $$
 \aligned
&\cup_{i}S_{i}U\subseteq U\,,\\
&S_{i}U\cap S_{j}U=\varnothing 
 \quad\text{for all $i,j$ with $i\not=j$,}\\
&U\cap K\not=\varnothing\,.
 \endaligned
 \tag7.1
 $$ 
In addition, it is easily seen that
$S_{\bold i}K\subseteq\overline{S_{\bold i}U}$ for all 
$\bold i\in\Sigma^{*}$,
and that
$S_{\bold i}K\cap S_{\bold j}U=\varnothing$
for all $\bold i,\bold j\in\Sigma^{*}$
with $|\bold i|=|\bold j|$ and 
$\bold i\not=\bold j$,
cf\. [Hu]. 
Also, since $U$ is open and bounded there exist $\rho_{1},\rho_{2}>0$
such that $U$ contains a ball of radius $\rho_{1}$, and $U$ is contained in a 
ball of radius $\rho_{2}$.
Since $U\cap K\not=\varnothing$, we can choose $\bold l\in\Sigma^{*}$ such that
 $$
 S_{\bold l}K\subseteq U\,,
 \tag7.2
 $$
and the compactness of $S_{\bold l}K$ now implies that
 $$
 d_{0}=\dist(S_{\bold l}K,\Bbb R^{d}\setminus U)>0\,.
 \tag7.3
 $$
For brevity write $D_{0}=\diam K$.
Also,
for a positive integer $m\in\Bbb N$, let $M_{m}\in\Bbb N$
 be chosen such that
  $$
  \frac{1}{r_{\max}^{M_{m}-1}}
  \ge
  2
  \,
  \frac{D_{0}}{d_{0}}
  \,
 \frac{1}{r_{\min}^{m-1}}
 \tag7.4
 $$
(recall that
$r_{\min}=\min_{i=1,\ldots,N}r_{i}$
and
$r_{\max}=\max_{i=1,\ldots,N}r_{i}$).
Now put
 $$
 \align
 r_{m}
&=
 r_{\min}^{M_{m}+m+|\bold l|}\,,
 \tag7.5\\
 a_{m}
&=
 \frac{1}{D_{0}}\,\frac{r_{\min}}{r_{\max}^{M_{m}+m}}\,,
 \tag7.6\\
 b_{m}
&=
 \frac{1}{D_{0}}\,\frac{1}{r_{\min}^{M_{m}+m}}\,,
 \tag7.7
 \endalign
 $$ 
and define 
$Z_{m}^{q}:(0,\infty)\to\Bbb R$ by
 $$
 Z_{m}^{q}(r)
 =
 \sum
 \Sb
  \bold h\in\Sigma^{*}\\
  |\bold h|\ge|\bold l|\\
  a_{m}r\le r_{\bold h}\le b_{m}r\\
  \bold l\not\prec\bold h
 \endSb
 p_{\bold h}^{q}\,.
 \tag7.8
 $$ 
If $m=1$,
then we will write
 $Z^{q}(r)$ for $Z_{1}^{q}(r)$, i\.e\.
 we will write
  $$
  Z^{q}(r)
 =
 \sum
 \Sb
  \bold h\in\Sigma^{*}\\
  |\bold h|\ge|\bold l|\\
  a_{1}r\le r_{\bold h}\le b_{1}r\\
  \bold l\not\prec\bold h
 \endSb
 p_{\bold h}^{q}\,.
 \tag7.9
 $$ 
The main purpose of this section to prove
Proposition 7.8
saying
if $m$ is a positive integer and 
$\bold i,\bold j\in\Sigma^{*}$ with
$|\bold i|=|\bold j|=m$ 
and $\bold i\not=\bold j$,
then there is a constant $c_{\bold i,\bold j,m}>0$ such that
for $r>0$
with
$r<\frac{1}{2}r_{m}$, we have
 $$
 Q_{\bold i,\bold j}^{q}(r)
 \le
 \,
 \cases
 c_{\bold i,\bold j,m}
 \,
 Z_{m}^{q}(\tfrac{1}{2}r)
&\quad
 \text{for $q<0$;}\\ 
&{}\\ 
 c_{\bold i,\bold j,m}
 \,
 Z_{m}^{q}(2r)
&\quad
 \text{for $0\le q$.}
 \endcases
 $$

We now turn towards the proof of Proposition 7.8.
However, we first make the following definition. Namely,
for a string $\bold i=i_{1}\ldots i_{m}\in\Sigma^{*}$,
let
 $$
 \hat{\bold i}=i_{1}\ldots i_{m-1}
 $$
denote the \lq\lq parent" of $\bold i$,
and 
for $r>0$ write
 $$
 \Sigma_{r}^{*}
 =
 \Big\{\bold i\in\Sigma^{*}
 \,\Big|\,
 r_{\bold i}\diam K<r\le r_{\hat{\bold i}}\diam K
 \Big\}\,.
 $$

\bigskip

\proclaim{Lemma 7.1}
Assume that the OSC is satisfied.

Let $\bold l$ and $d_{0}$ be as in (7.2) and (7.3).
Assume that
$u,w\in\Sigma$
and
$\bold u,\bold w,\bold h\in\Sigma^{*}$
satisfy the following conditions:
\roster
\item"(i)" $u\not=w$:
\item"(ii)" 
$\dist(\, S_{u\bold u\bold h}K \,,\, S_{w\bold w}K  \,)\le r_{u\bold u\bold h}d_{0}$.
\endroster
Then
$\bold l\not\prec\bold h$.
\endproclaim
\noindent{\it  Proof}\newline
\noindent
Let $U$ be the open set in (7.1).
Assume, in order to obtain a contradiction,
 that
$\bold l$ is
a substring of $\bold h$, 
i\.e\.
there exist strings $\bold s,\bold t\in\Sigma^{*}$ 
such that
$\bold h=\bold s\bold l\bold t$. 
Hence 
 $$
 \align
 \dist
 \big(
 \,S_{u\bold u\bold h}K
 \,,\,
 \Bbb R^{d}\setminus S_{u\bold u\bold s}U\,
 \big)
&= 
 \dist
 \big(
 \,S_{u\bold u\bold s\bold l\bold t}K
 \,,\,
 \Bbb R^{d}\setminus S_{u\bold u\bold s}U\,
 \big)\\
&\ge
 \dist
 \big(
 \,S_{u\bold u\bold s\bold l}K
 \,,\,
 \Bbb R^{d}\setminus S_{u\bold u\bold s}U\,
 \big)\\
&=
 r_{u\bold u\bold s}
 \dist
 \big(
 \,S_{\bold l}K
 \,,\,
 \Bbb R^{d}\setminus U\,
 \big)\\  
&\ge
 \dist
 \big(
 \,S_{\bold u\bold l}K\,,\,
 \Bbb R^{d}\setminus S_{\bold u}U\,
 \big)\\
&=
 r_{u\bold u\bold s}d_{0}\\
&>
 r_{u\bold u\bold h}d_{0}\,.
 \tag7.10
 \endalign
 $$
Also, $u\not =w$ (by (i)), we conclude that
$S_{u\bold u\bold s}U\cap S_{w\bold w}K=\varnothing$, i\.e\.
$S_{w\bold w}K
\subseteq
\Bbb R^{d}\setminus S_{u\bold u\bold s}U$, whence (using (ii))
 $$
 \align
 \dist
 \big(
 \,S_{u\bold u\bold h}K
 \,,\,
 \Bbb R^{d}\setminus S_{u\bold u\bold s}U\,
 \big) 
&\le
 \dist
 \big(
 \,S_{u\bold u\bold h}K
 \,,\,
S_{w\bold w}K\,
 \big)\\ 
 &\le
 r_{u\bold u\bold h}d_{0}
 \tag7.11
 \endalign
 $$
Inequalities (7.10) and (7.11) give the desired contradiction.
\hfill$\square$

\bigskip

\proclaim{Lemma 7.2}
Assume that the OSC is satisfied.

Let $\bold l$ be as in (7.2).
Let $m\in\Bbb N$.
Let $\bold i,\bold j\in\Sigma^{*}$ with
$|\bold i|=|\bold j|=m$ and
$\bold i\not=\bold j$.
Assume that
$r>0$, $x\in\Bbb R^{d}$ and 
$\bold k\in\Sigma^{*}$ satisfy the following conditions:
\roster
\item"(i)" 
$ 0<r<r_{m}$;
\item"(ii)"
$x\in B(S_{\bold i}K,r)\cap B(S_{\bold j}K,r)$;
\item"(iii)"
$\bold k\in\Sigma_{r}^{*}$
and
$x\in B(S_{\bold k}K,r)$.
\endroster
Then there are strings $\bold u,\bold h,\bold v\in\Sigma^{*}$ 
with
$\bold k=\bold u\bold h\bold v$
such that
 $$
 \align
 |\bold u|
&=m\,,\\
 |\bold v|
&=M_{m}\,,\\
 |\bold h|
&\ge|\bold l|\,,\\
 \bold l
&\not\prec
 \bold h\,.
 \endalign
 $$
\endproclaim
\noindent{\it  Proof}\newline
\noindent
Recall, that $D_{0}=\diam K$.
Since
$\bold k\in\Sigma_{r}^{*}$, we conclude that
$r_{\min}^{|\bold k|}D_{0}
\le
r_{\bold k}D_{0}
\le
r
\le
r_{m}
=
r_{\min}^{M_{m}+m+|\bold l|}$,
whence
$|\bold k|
\ge
M_{m}+m+|\bold l|$. It follows from this that
there
are $\bold u,\bold h,\bold v\in\Sigma^{*}$ 
with
$\bold k=\bold u\bold h\bold v$
such that
$|\bold u|=m$,
$|\bold v|=M_{m}$
and
$|\bold h|\ge|\bold l|$.
We must now show that
 $$
 \bold l
 \not\prec
 \bold h\,.
 $$
Note that $|\bold i|=m=|\bold u|\le|\bold k|$
and
$|\bold j|=m=|\bold u|\le|\bold k|$.
Hence, since
$\bold i\not=\bold j$, we can find
$\bold w\in\{\bold i,\bold j\}$ such that
 $$
 \bold w\not=\bold k|m=\bold u\,.
 $$
Write
$\bold u=u_{1}\ldots u_{m}$
and
$\bold w=w_{1}\ldots w_{m}$.
Since
$\bold w\not=\bold u$, there is $s\in\{1,\ldots,m\}$ such that
 $$
 \align
 w_{i}
&=
 u_{i}
 \,\,\,\,\text{for $i=1,\ldots,s-1$,}\\
 w_{s}
&\not=
 u_{s}\,.
 \endalign
 $$
Next, we prove the following two claims.

\bigskip

\noindent{\it Claim 1:  We have
$\dist
 (
 \,
 S_{u_{s}\ldots u_{m}\bold h}K
 \,,\,
 S_{w_{s}\ldots w_{m}}K
 \,
 )
\le
 2\frac{1}{r_{\min}^{m-1}}r$.
}
 
\noindent
{\it Proof of Claim 1.}
Since
$w_{1}=u_{1},\ldots,w_{s-1}=u_{s-1}$,
we deduce that
 $$
 \align
 \dist
 \big(
 \,
 S_{u_{s}\ldots u_{m}\bold h}K
 \,,\,
 S_{w_{s}\ldots w_{m}}K
 \,
 \big)
&=
 \frac{1}{r_{u_{1}\ldots u_{s-1}}}
 \dist
 \big(
 \,
 S_{u_{1}\ldots u_{m}\bold h}K
 \,,\,
 S_{w_{1}\ldots w_{m}}K
 \,
 \big)\\
&\le
 \frac{1}{r_{\min}^{s-1}}
 \dist
 \big(
 \,
 S_{u_{1}\ldots u_{m}\bold h}K
 \,,\,
 S_{w_{1}\ldots w_{m}}K
 \,
 \big)\\ 
&\le
 \frac{1}{r_{\min}^{m-1}}
 \dist
 \big(
 \,
 S_{u_{1}\ldots u_{m}\bold h}K
 \,,\,
 S_{w_{1}\ldots w_{m}}K
 \,
 \big)\\  
&=
 \frac{1}{r_{\min}^{m-1}}
 \dist
 \big(
 \,
 S_{\bold u\bold h}K
 \,,\,
 S_{\bold w}K
 \,
 \big)\,. 
 \tag7.12
 \endalign
 $$
Next, since
$S_{\bold u\bold h\bold v}K\subseteq S_{\bold u\bold h}K$, it follows 
from (7.12) that
 $$
 \align
 \dist
 \big(
 \,
 S_{u_{s}\ldots u_{m}\bold h}K
 \,,\,
 S_{w_{s}\ldots w_{m}}K
 \,
 \big)
&\le
 \frac{1}{r_{\min}^{m-1}}
 \dist
 \big(
 \,
 S_{\bold u\bold h\bold v}K
 \,,\,
 S_{\bold w}K
 \,
 \big)\\
&=
 \frac{1}{r_{\min}^{m-1}}
 \dist
 \big(
 \,
 S_{\bold k}K
 \,,\,
 S_{\bold w}K
 \,
 \big)\\ 
&\le
 \frac{1}{r_{\min}^{m-1}}
 \Big(
 \dist
 \big(
 \,
 S_{\bold k}K
 \,,\,
 x
 \,
 \big)
 +
 \dist
 \big(
 \,
 x
 \,,\,
 S_{\bold w}K
 \,
 \big)
 \Big)\,.\\
&{} 
 \tag7.13
 \endalign
 $$
However,
since
$x\in B(S_{\bold i}K,r)\cap B(S_{\bold j}K,r)
\subseteq
B(S_{\bold w}K,r)$,
we have
$\dist(x,S_{\bold w}K)\le r$.
Similarly, since
$x\in B(S_{\bold k}K,r)$,
we have
$\dist(S_{\bold k}K,x)\le r$.
It therefore follows from (7.13) that
 $$
 \align
 \dist
 \big(
 \,
 S_{u_{s}\ldots u_{m}\bold h}K
 \,,\,
 S_{w_{s}\ldots w_{m}}K
 \,
 \big)
&\le
 2\frac{1}{r_{\min}^{m-1}}r\,.
 \endalign
 $$
This completes the proof of Claim 1.

\bigskip

\noindent{\it Claim 2:  We have
$d_{0}r_{u_{s}\ldots u_{m}\bold h}
\ge
 \frac{d_{0}}{D_{0}}
 \frac{1}{r_{max}^{M_{m}-1}}r$.
}
 
\noindent
{\it Proof of Claim2.}
We have
(using the fact that $\bold k\in\Sigma_{r}^{*}$)
 $$
 \align
 d_{0}r_{u_{s}\ldots u_{m}\bold h}
&\ge
 d_{0}
 r_{\bold u\bold h}\\
&\ge
 d_{0}
 \frac{1}{r_{max}^{M_{m}}}r_{\bold u\bold h\bold v}\\ 
&\ge
 \frac{d_{0}}{D_{0}}
 \frac{1}{r_{max}^{M_{m}-1}}r\,.
 \endalign
 $$
This completes the proof of Claim 2.

\bigskip

It now follows from Claim 1, Claim 2 and (7.4) that
 $$
 \dist
 \big(
 \,
 S_{u_{s}\ldots u_{m}\bold h}K
 \,,\,
 S_{w_{s}\ldots w_{m}}K
 \,
 \big) 
 \le
 r_{u_{s}\ldots u_{m}\bold h}d_{0}\,,
 $$
and we therefore deduce from Lemma 7.1 that
$\bold l
 \not\prec
 \bold h$.
This completes the proof.
\hfill$\square$

\bigskip

\proclaim{Lemma 7.3}
Assume that the OSC is satisfied.
For $\bold i\in\Sigma^{*}$, we have
 $$
 \mu(S_{\bold i}K)=p_{\bold i}\,.
 $$
\endproclaim
\noindent{\it  Proof}\newline
\noindent
This lemma is proved in [Graf].
\hfill$\square$

\bigskip

\noindent
Lemma 7.4 below is a slight modification of a result due to 
Hutchinson [Hu] and the proof is therefore omitted.
Moreover, Lemma 7.5
is a standard result and the proof of Lemma 7.5 is therefore 
also omitted.

\bigskip

\proclaim{Lemma 7.4}
Let $r,k,k_{1},k_{2}>0$, and let $(V_{i})_{i}$ be a family of pairwise disjoint
open 
subsets 
of $\Bbb R^{d}$.
Assume that each set $V_{i}$ contains a ball of radius 
$k_{1}r$ and is contained in a ball of radius $k_{2}r$.
Then
 $$
 \big|\,
 \{i\,|\, \overline {V_{i}}\cap B(x,kr)\not=\varnothing\}
 \,\big|
 \le
 \left(
 \frac{k+2k_{2}}{k_{1}}
 \right)^{d}
 $$
for all $x\in\Bbb R^{d}$. 
\endproclaim

\bigskip

\proclaim{Lemma 7.5}
Fix $q\in\Bbb R$.
Let $c>0$ and $(a_{i})_{i\in I}$ be a family of real numbers with $|I|\le c$.
Then
$\left(
 \sum_{i\in I}a_{i}
 \right)^{q}
 \le
 \max(1,c^{q-1})
 \sum_{i\in I}a_{i}^{q}$.
\endproclaim

\bigskip

\noindent
Before stating and proving the next proposition we need the following definition.
Namely,
for $r>0$,
we will say that a subset $F$ of $\Bbb R^{d}$ is $r$-separated
if
 $$
 B(x,r)\cap B(y,r)=\varnothing
 $$
for all $x,y\in F$.

\bigskip

\proclaim{Proposition 7.6}
Fix $q\in\Bbb R$
and
assume that the OSC is satisfied.

Let
$m\in\Bbb N$.
Let $\bold i,\bold j\in\Sigma^{*}$ with
$|\bold i|=|\bold j|=m$
and
$\bold i\not=\bold j$.

There exists a constant $k_{\bold i,\bold j,m}>0$ such that
if $r>0$ with $r<r_{m}$ and
$F\subseteq B(S_{\bold i}K,r)\cap B(S_{\bold j}K,r)$
is a $(1+\frac{\rho_{2}}{D_{0}})r$-separated
set, then we have
 $$
 \sum_{x\in F}\mu(B(x,r))^{q}
 \le
 \,
 k_{\bold i,\bold j,m}
 \,
 Z_{m}^{q}(r)
 $$ 
(recall, that  
$\rho_{2}$
is a positive real number such that
the non-empty and open set $U$ from
(7.1) contains a ball of radius equal to $\rho_{2}$, 
and that
$Z_{m}^{q}(r)$ is defined in (7.8)).
\endproclaim
\noindent{\it  Proof}\newline
\noindent 
Let $U$, $\bold l$ and $d_{0}$ be as in (7.1), (7.2) and (7.3), respectively.
Fix $0<r<r_{m}$.
For each $x\in F$ we may choose $\bold k(x)\in\Sigma_{r}^{*}$ such that
$x\in S_{\bold k(x)}K$. 
We clearly have
 $$
 S_{\bold k(x)}K
 \subseteq
 B(x,r)\cap K
 \subseteq
 \bigcup
   \Sb
   \bold k\in\Sigma_{r}^{*}\\
   \dist(x,S_{\bold k}K)\le r
   \endSb 
 S_{\bold k}K\,,
 $$
for all $x\in F$, whence
 $$
 \align
  \mu(B(x,r))^{q}
&\le
 \cases
 \mu(S_{\bold k(x)}K)^{q}
 &\qquad
 \text{for $q\le 0$;}\\
 &{}\\ 
 {\dsize
 \mu
 \left(
 \bigcup
   \Sb
   \bold k\in\Sigma_{r}^{*}\\
   \dist(x,S_{\bold k}K)\le r
   \endSb 
 S_{\bold k}K
 \right)^{q}
 }
 &\qquad
 \text{for $0\le q$;}
 \endcases\\
&{}\\ 
&\le
 \cases
 \mu(S_{\bold k(x)}K)^{q}
 &\qquad
 \text{for $q\le 0$;}\\
 &{}\\ 
 {\dsize
 \left(
 \sum
   \Sb
   \bold k\in\Sigma_{r}^{*}\\
   \dist(x,S_{\bold k}K)\le r
   \endSb 
 \mu(S_{\bold k}K)
 \right)^{q}
 }
 &\qquad
 \text{for $0\le q$,}
 \endcases\\
 \endalign
 $$
for all $x\in F$.
This implies that
 $$
 \align
  \sum_{x\in F}\mu(B(x,r))^{q}
&\le
 \cases
 \sum_{x\in F}
 \mu(S_{\bold k(x)}K)^{q}
 &\qquad
 \text{for $q\le 0$;}\\
 &{}\\ 
 {\dsize
 \sum_{x\in F}
 \left(
 \sum
   \Sb
   \bold k\in\Sigma_{r}^{*}\\
   \dist(x,S_{\bold k}K)\le r
   \endSb 
 \mu(S_{\bold k}K)
 \right)^{q}
 }
 &\qquad
 \text{for $0\le q$.}
 \endcases
 \tag7.14
 \endalign
 $$

Next, we prove the following two claims.

\bigskip

\noindent
{\it Claim 1.
There is a constant $C_{1}>0$ such that
 $$
 \left(
 \sum
   \Sb
   \bold k\in\Sigma_{r}^{*}\\
   \dist(x,S_{\bold k}K)\le r
   \endSb 
 \mu(S_{\bold k}K)
 \right)^{q}
 \le
 C_{1}
 \sum
   \Sb
   \bold k\in\Sigma_{r}^{*}\\
   \dist(x,S_{\bold k}K)\le r
   \endSb 
 \mu(S_{\bold k}K)^{q}
 $$
for all $r>0$ and all $x\in\Bbb R^{d}$. 
}

\noindent{\it Proof of Claim 1.}
Recall that $U$ contains a ball of radius $\rho_{1}$ and is contained 
in a ball of radius $\rho_{2}$.

For $\bold k\in\Sigma_{r}^{*}$, we therefore conclude that
$S_{\bold k}U$ 
contains a ball of radius
$r_{\bold k}\rho_{1}$
and that
$r_{\bold k}\rho_{1}
\ge
r_{\hat\bold k}r_{\min}\rho_{1}
\ge
\frac{r_{\min}\rho_{1}}{D_{0}}r$
(because $\bold k\in\Sigma_{r}^{*}$).
We deduce from this that
$S_{\bold k}U$ 
contains a ball of radius
$\frac{r_{\min}\rho_{1}}{D_{0}}r$.

For $\bold k\in\Sigma_{r}^{*}$, we also conclude that
$S_{\bold k}U$ 
is contained
in a ball of radius
$r_{\bold k}\rho_{2}$
and that
$r_{\bold k}\rho_{2}
\le
\frac{\rho_{2}}{D_{0}}r$
(because $\bold k\in\Sigma_{r}^{*}$).
We deduce from this that
$S_{\bold k}U$ 
is contained in
a ball of radius
$\frac{\rho_{2}}{D_{0}}r$.

Next,
since $(S_{\bold k}U)_{\bold k\in\Sigma_{r}^{*}}$ is a pairwise disjoint 
family of sets with
$S_{\bold k}K\subseteq\overline{S_{\bold k}U}$, 
Lemma 7.4 therefore
 implies that
 $$
 \align
 \Bigg|
 \,
 \Big\{
 \bold k\in\Sigma_{r}^{*}
  \,\Big|\,
   \dist(x,S_{\bold k}K)\le r
 \Big\}
 \,
 \Bigg|
&\le
 \Bigg|
 \,
 \Big\{
 \bold k\in\Sigma_{r}^{*}
  \,\Big|\,
  S_{\bold k}K\cap B(x,r)
  \not=
  \varnothing
 \Big\}
 \,
 \Bigg|\\
&\le
 \Bigg|
 \,
 \Big\{
 \bold k\in\Sigma_{r}^{*}
  \,\Big|\,
  \overline{S_{\bold k}U}\cap B(x,r)
  \not=
  \varnothing
 \Big\}
 \,
 \Bigg|\\
&\le
 C_{0}
 \tag7.15
 \endalign
 $$
for all $x$,
where
$C_{0}=((1+2\frac{\rho_{2}}{D_{0}})/\frac{r_{\min}\rho_{1}}{D_{0}})^{d}$.

Finally, we deduce from (7.15) and 
Lemma 7.5 that 
 $$
 \left(
 \sum
   \Sb
   \bold k\in\Sigma_{r}^{*}\\
   \dist(x,S_{\bold k}K)\le r
   \endSb 
 \mu(S_{\bold k}K)
 \right)^{q}
 \le
 C_{1}
 \sum
   \Sb
   \bold k\in\Sigma_{r}^{*}\\
   \dist(x,S_{\bold k}K)\le r
   \endSb 
 \mu(S_{\bold k}K)^{q}\,,
 $$
where 
$C_{1}=\max(1,C_{0}^{q-1})$.
This completes the proof of Claim 1.

\bigskip

\noindent{\it Claim 2.
There is a constant $C_{2}>0$ such that
$$
 \sum_{x\in F}
 \sum
 \Sb
   \bold k\in\Sigma_{r}^{*}\\
   \dist(x,S_{\bold k}K)\le r
   \endSb 
 \mu(S_{\bold k}K)^{q}
 \le
 C_{2}
 \,
 \sum_{|\bold u|=m}
 \,
 \sum_{|\bold v|=M_{m}}
 \,
 \sum
 \Sb
  \bold h\in\Sigma^{*}\\
    \bold u\bold h\bold v\in\Sigma_{r}^{*}\\
  |\bold h|\ge|\bold l|\\
  \bold l\not\prec\bold h
  \endSb 
 \mu(S_{\bold u\bold h\bold v}K)^{q}
 $$
for all $r>0$.}

\noindent{\it Proof of Claim 2.}
Again, recall that $U$ contains a ball of radius $\rho_{1}$ and is contained 
in a ball of radius $\rho_{2}$.

Fx $y\in\Bbb R^{d}$.

We first  prove  that
if $y\in\Bbb R^{d}$, then
 $$
 \bigcup
  \Sb
  \bold k\in\Sigma_{r}^{*}\\
  \dist(y,S_{\bold k}K)\le r
  \endSb
 S_{\bold k}U
 \subseteq
 B(y,(1+\tfrac{\rho_{2}}{D_{0}})r)\,.
 \tag7.16
 $$
Indeed, 
this follows from the fact that
$S_{\bold k}U$ 
is contained
in a ball of radius
$r_{\bold k}\rho_{2}$
and 
$r_{\bold k}\rho_{2}
\le
\frac{\rho_{2}}{D_{0}}r$
for all  $\bold k\in\Sigma_{r}^{*}$.
This proves (7.16).

\medskip

Next,
for $\bold k\in\Sigma_{r}^{*}$, we conclude that
$S_{\bold k}U$ 
contains a ball of radius
$r_{\bold k}\rho_{1}$
and that
$r_{\bold k}\rho_{1}
\ge
r_{\hat\bold k}r_{\min}\rho_{1}
\ge
\frac{r_{\min}\rho_{1}}{D_{0}}r$
(because $\bold k\in\Sigma_{r}^{*}$).
We deduce from this that:
 $$
 \bigcup
  \Sb
  \bold k\in\Sigma_{r}^{*}\\
  \dist(y,S_{\bold k}K)\le r
  \endSb
 S_{\bold k}U
 \quad
 \text{
 contains a ball of radius $\frac{r_{\min}\rho_{1}}{D_{0}}r$.
 } 
 \tag7.17
$$

Also, we deduce from (7.16)  that:
 $$
 \bigcup
  \Sb
  \bold k\in\Sigma_{r}^{*}\\
  \dist(y,S_{\bold k}K)\le r
  \endSb
 S_{\bold k}U
 \quad
 \text{
 is contained in a ball of radius $\left(1+\frac{\rho_{2}}{D_{0}}\right)r$.
 } 
 \tag7.18
$$

Next, we prove that:
 $$
 \left(
  \bigcup
  \Sb
  \bold k\in\Sigma_{r}^{*}\\
  \dist(y,S_{\bold k}K)\le r
  \endSb
 S_{\bold k}U
\right)_{y\in F}
\quad
\text{
is a pairwise disjoint family of open sets.
}
\tag7.19
$$
Indeed, this follows from (7.16) and the
fact that
$F$ is $(1+\frac{\rho_{2}}{D_{0}})r$-separated. This proves (7.19).

Finally, we deduce from (7.17), (7.18), (7.19) and Lemma 7.4 that
if $x\in F$, then
 $$
 \align
 &
 \left|
 \,
 \left\{
 y\in F
 \,
 \left|
 \,
 \left(
  \bigcup
  \Sb
  \bold k\in\Sigma_{r}^{*}\\
  \dist(y,S_{\bold k}K)\le r
  \endSb
 S_{\bold k}K
 \right)
 \cap
 \left(
  \bigcup
  \Sb
  \bold k\in\Sigma_{r}^{*}\\
  \dist(x,S_{\bold k}K)\le r
  \endSb
 S_{\bold k}K
 \right)
 \not=\varnothing
 \right.
 \right\}
 \,
 \right|\\
&{} \\
&\qquad\qquad
  \qquad\qquad
 \le
  \left|
 \,
 \left\{
 y\in F
 \,
 \left|
 \,
 \left(
  \bigcup
  \Sb
  \bold k\in\Sigma_{r}^{*}\\
  \dist(y,S_{\bold k}K)\le r
  \endSb
 S_{\bold k}K
 \right)
 \cap
 B(x,(1+\tfrac{\rho_{2}}{D_{0}})r)
 \not=\varnothing
 \right.
 \right\}
 \,
 \right|\\
&{} \\
&\qquad\qquad
  \qquad\qquad
 \le
  \left|
 \,
 \left\{
 y\in F
 \,
 \left|
 \,
 \left(
  \bigcup
  \Sb
  \bold k\in\Sigma_{r}^{*}\\
  \dist(y,S_{\bold k}K)\le r
  \endSb
 \overline{S_{\bold k}U}
 \right)
 \cap
 B(x,(1+\tfrac{\rho_{2}}{D_{0}})r)
 \not=\varnothing
 \right.
 \right\}
 \,
 \right|\\
&{} \\
&\qquad\qquad
  \qquad\qquad
 \le
  \left|
 \,
 \left\{
 y\in F
 \,
 \left|
 \,
 \left(
 \overline{
  \bigcup
  \Sb
  \bold k\in\Sigma_{r}^{*}\\
  \dist(y,S_{\bold k}K)\le r
  \endSb
 S_{\bold k}U
 }
 \right)
 \cap
 B(x,(1+\tfrac{\rho_{2}}{D_{0}})r)
 \not=\varnothing
 \right.
 \right\}
 \,
 \right|\\
&{} \\
&\qquad\qquad
  \qquad\qquad
  \le
  C_{2}
  \tag7.20
 \endalign
 $$
where
$C_{2}
=
(
((1+\frac{\rho_{2}}{D_{0}})+2(1+\frac{\rho_{2}}{D_{0}}))
/
(\frac{r_{\min}\rho_{1}}{D_{0}})
)^{d}$.

We deduce from (7.20) that each term in the sum
 $$
 \sum_{x\in F}
 \sum
  \Sb
  \bold k\in\Sigma_{r}^{*}\\
  \dist(x,S_{\bold k}K)\le r
  \endSb
 \mu(S_{\bold k}K)^{q}
 $$ 
is repeated atmost
$C_{2}$ times.
This observation and 
Lemma 7.2
(which is applicable since $0<r<r_{m}$)
 now gives
 $$
 \sum_{x\in F}
 \sum
 \Sb
   \bold k\in\Sigma_{r}^{*}\\
   \dist(x,S_{\bold k}K)\le r
   \endSb 
 \mu(S_{\bold k}K)^{q}
 \le
 C_{2}
 \,
 \sum_{|\bold u|=m}
 \,
 \sum_{|\bold v|=M_{m}}
 \,
 \sum
 \Sb
  \bold h\in\Sigma^{*}\\
  \bold u\bold h\bold v\in\Sigma_{r}^{*}\\
  |\bold h|\ge|\bold l|\\
  \bold l\not\prec\bold h
  \endSb 
 \mu(S_{\bold u\bold h\bold v}K)^{q}\,.
 $$
This completes the proof of Claim 2.

\bigskip

Combining (7.14), Claim 1 and Claim 2 now yields
 $$
 \align
  \sum_{x\in F}\mu(B(x,r))^{q}
&\le
 \cases
 {\dsize
 \sum_{x\in F}
 \mu(S_{\bold k(x)}K)^{q}
 }
 &\qquad
 \text{for $q\le 0$;}\\
 &{}\\ 
 {\dsize
 \sum_{x\in F}
 \left(
 \sum
   \Sb
   \bold k\in\Sigma_{r}^{*}\\
   \dist(x,S_{\bold k}K)\le r
   \endSb 
 \mu(S_{\bold k}K)
 \right)^{q}
 }
 &\qquad
 \text{for $0\le q$;}
 \endcases\\
&{}\\ 
&\le
 \cases
 {\dsize
 \sum_{x\in F}
 \sum
   \Sb
   \bold k\in\Sigma_{r}^{*}\\
   \dist(x,S_{\bold k}K)\le r
   \endSb 
 \mu(S_{\bold k}K)^{q}
 }
 &\qquad
 \text{for $q\le 0$;}\\
 &{}\\ 
 {\dsize
 C_{1}
 \sum_{x\in F}
 \sum
   \Sb
   \bold k\in\Sigma_{r}^{*}\\
   \dist(x,S_{\bold k}K)\le r
   \endSb 
 \mu(S_{\bold k}K)^{q}
 }
 &\qquad
 \text{for $0\le q$;}
 \endcases\\
&{}\\
&\le
 C_{1}C_{2}
 \,
 \sum_{|\bold u|=m}
 \,
 \sum_{|\bold v|=M_{m}}
 \,
 \sum
 \Sb
  \bold h\in\Sigma^{*}\\
  \bold u\bold h\bold v\in\Sigma_{r}^{*}\\
  |\bold h|\ge|\bold l|\\
  \bold l\not\prec\bold h
  \endSb 
 \mu(S_{\bold u\bold h\bold v}K)^{q}\\
&{}\\
&=
 C_{1}C_{2}
 \,
 \sum_{|\bold u|=m}
 \,
 \sum_{|\bold v|=M_{m}}
 \,
 \sum
 \Sb
  \bold h\in\Sigma^{*}\\
  \bold u\bold h\bold v\in\Sigma_{r}^{*}\\
  |\bold h|\ge|\bold l|\\
  \bold l\not\prec\bold h
  \endSb 
 p_{\bold u}^{q}p_{\bold h}^{q}p_{\bold v}^{q}\,,
 \tag7.21
 \endalign
 $$
where the last equality is due to Lemma 7.3. 
However, if
$\bold u,\bold h,\bold v\in\Sigma^{*}$
with
$|\bold u|=m$, $|\bold v|=M_{m}$
and
$\bold u\bold h\bold v\in\Sigma_{r}^{*}$, then
$r_{\bold h}
 =
 \frac{r_{\bold u\bold h\bold v}}{r_{\bold u}r_{\bold v}}
 \le
 \frac{r}{D_{0}r_{\min}^{m+M_{m}}}
 =
 b_{m}r$
and
$r_{\bold h}
 =
 \frac{r_{\bold u\bold h\bold v}}{r_{\bold u}r_{\bold v}}
 \ge
 \frac{rr_{\min}}{D_{0}r_{\max}^{m+M_{m}}}
 =
 a_{m}r$.
We deduce from this and (7.21) that
 $$
 \align
  \sum_{x\in F}\mu(B(x,r))^{q}
&\le
 C_{1}C_{2}
 \,
 \sum_{|\bold u|=m}
 \,
 \sum_{|\bold v|=M_{m}}
 \,
 \sum
 \Sb
  \bold h\in\Sigma^{*}\\
  a_{m}\le r_{\bold h}\le b_{m}r\\
  |\bold h|\ge|\bold l|\\
  \bold l\not\prec\bold h
  \endSb 
 p_{\bold u}^{q}p_{\bold h}^{q}p_{\bold v}^{q}\\
&{}\\ 
&\le
 C_{1}
 \,
 C_{2}
 \,
 N^{m+M_{m}}
 \,
 \big(\max_{i}p_{i}^{q}\big)^{m+M_{m}}
 \,
 \sum
 \Sb
  \bold h\in\Sigma^{*}\\
  a_{m}\le r_{\bold h}\le b_{m}r\\
  |\bold h|\ge|\bold l|\\
  \bold l\not\prec\bold h
  \endSb 
 p_{\bold h}^{q}\\
&=
 k_{\bold i,\bold j,m}\,Z_{m}^{q}(r)
 \endalign
 $$
where
$k_{\bold i,\bold j,m}
=
C_{1}
 \,
 C_{2}
 \,
 N^{m+M_{m}}
 \,
 (\max_{i}p_{i}^{q})^{m+M_{m}}$.
\hfill$\square$

 \bigskip

\noindent
We now turn towards the proof of the main result in this section, namely, Proposition 7.8
However, in order to 
deduce Proposition 7.8 from 
Proposition 7.6 we need the following simple covering lemma.

\bigskip

\proclaim{Lemma 7.7}
Let $d\in\Bbb N$ and $u\ge 1$.
Then there exists a positive integer $\chi\in\Bbb N$
satisfying the following:
if $E\subseteq\Bbb R^{d}$ is a bounded set and $s>0$,
then there are sets
$F_{1},\ldots,F_{\chi}\subseteq E$ such that
\roster
\item"(1)"
The set $F_{i}$ is $us$-separated for each $i$.
\item"(2)"
 We have
$E\subseteq\cup_{i=1}^{\chi}\cup_{y\in F_{i}}B(y,s)$.
\endroster
\endproclaim
\noindent{\it  Proof}\newline
\noindent
This is easily proved and the proof is therefore omitted.
\hfill$\square$

\bigskip

\noindent
We can now prove Proposition 7.8.

 \bigskip

\proclaim{Proposition 7.8}
Fix $q\in\Bbb R$
and
assume that the OSC is satisfied.

Let
$m\in\Bbb N$.
Let $\bold i,\bold j\in\Sigma^{*}$ with
$|\bold i|=|\bold j|=m$
and
$\bold i\not=\bold j$.

There exists a constant $c_{\bold i,\bold j,m}>0$ such that
if $r>0$ with $r<\frac{1}{2}r_{m}$, then we have
 $$
 Q_{\bold i,\bold j}^{q}(r)
 \le
 \,
 \cases
 c_{\bold i,\bold j,m}
 \,
 Z_{m}^{q}(\tfrac{1}{2}r)
&\quad
 \text{for $q<0$;}\\ 
&{}\\ 
 c_{\bold i,\bold j,m}
 \,
 Z_{m}^{q}(2r)
&\quad
 \text{for $0\le q$.}
 \endcases
 $$ 
(recall, that  $Z_{m}^{q}(r)$ is defined in (7.8)).
\endproclaim
\noindent{\it  Proof}\newline
\noindent 
It follows from Lemma 7.7
that there is a positive integer $\chi$
such that for all $r>0$
we can find
sets
$F_{r,1},\ldots,F_{r,\chi}\subseteq
B(S_{\bold i}K,r)\cap B(S_{\bold j}K,r)$
satisfying:
 $$ 
\text{
the set
$F_{r,i}$
is
$(1+\tfrac{\rho_{2}}{D_{0}})2r$-separated
for each $i$,}
\tag7.22
$$
and
 $$
 B(S_{\bold i}K,r)\cap B(S_{\bold j}K,r)
 \subseteq
 \bigcup_{i=1}^{\chi}
 \,
 \bigcup_{y\in F_{r,i}}
 \,
 B(y,\tfrac{1}{2}r)\,.
 \tag7.23
 $$

Fix $0<r<\frac{1}{2}r_{m}$.
It follows from (7.23) that
 $$
 \align
 Q_{\bold i,\bold j}^{q}(r)
&=
\frac{1}{r^{d}}
\int\limits_{B(S_{\bold i}K,r)\cap B(S_{\bold j}K,r)}
\mu(B(x,r))^{q}
\,
d\Cal L^{d}(x)\\
&\le
 \frac{1}{r^{d}}
 \,
 \sum_{i=1}^{\chi}
 \,
 \sum_{y\in F_{r,i}}
 \,
 \int\limits_{B(y,\frac{1}{2}r)}
\mu(B(x,r))^{q}
\,
d\Cal L^{d}(x)\,.
\tag7.24
\endalign
$$
Now, note that
for $x\in B(y,\tfrac{1}{2}r)$, we clearly have
$B(y,\frac{1}{2}r)
\subseteq
B(x,r)
\subseteq
B(y,2r)$, whence
 $$
 \mu(B(x,r))^{q}
 \le
 \cases
 \mu(B(y,\frac{1}{2}r))^{q}
&\quad
 \text{for $q<0$;}\\
 \mu(B(y,2r))^{q}
&\quad
 \text{for $0\le q$.}  
\endcases
\tag7.25
$$
Next,
writing
$\Omega_{d}=\Cal L^{d}(B(0,1))$ for the volume of the unit ball in $\Bbb R^{d}$
and combining (7.24) and (7.25) gives
 $$
 \align
 Q_{\bold i,\bold j}^{q}(r)
&\le
 \cases
 {\dsize
 \frac{1}{r^{d}}
 \,
 \sum_{i=1}^{\chi}
 \,
 \sum_{y\in F_{r,i}}
 \,
 \int\limits_{B(y,\frac{1}{2}r)}
\mu(B(y,\tfrac{1}{2}r))^{q}
\,
d\Cal L^{d}(x)
}
&\quad
\text{for $q<0$;}
 \\
 {\dsize
 \frac{1}{r^{d}}
 \,
 \sum_{i=1}^{\chi}
 \,
 \sum_{y\in F_{r,i}}
 \,
 \int\limits_{B(y,\frac{1}{2}r)}
\mu(B(y,2r))^{q}
\,
d\Cal L^{d}(x)
}
&\quad
\text{for $0\le q$;}
\endcases\\
&{}\\
&\le
 \cases
 {\dsize
 \frac{1}{r^{d}}
 \,
 \sum_{i=1}^{\chi}
 \,
 \sum_{y\in F_{r,i}}
 \,
\mu(B(y,\tfrac{1}{2}r))^{q}
\,
\Cal L^{d}(B(y,\tfrac{1}{2}r))
}
&\quad
\text{for $q<0$;}
 \\
 {\dsize
 \frac{1}{r^{d}}
 \,
 \sum_{i=1}^{\chi}
 \,
 \sum_{y\in F_{r,i}}
 \,
\mu(B(y,2r))^{q}
\,
\Cal L^{d}(B(y,\tfrac{1}{2}r))
}
&\quad
\text{for $0\le q$;}
\endcases\\
&{}\\
&\le
 \cases
 {\dsize
 \frac{\Omega_{d}}{2^{d}}
 \,
 \sum_{i=1}^{\chi}
 \,
 \sum_{y\in F_{r,i}}
 \,
\mu(B(y,\tfrac{1}{2}r))^{q}}
&\quad
\text{for $q<0$;}
\\
 {\dsize
\frac{\Omega_{d}}{2^{d}}
 \,
 \sum_{i=1}^{\chi}
 \,
 \sum_{y\in F_{r,i}}
 \,
\mu(B(y,2r))^{q}
}
&\quad
\text{for $0\le q$.}
\endcases
\tag7.26
\endalign
$$
However, since both $\frac{1}{2}r<r_{m}$ and $2r<r_{m}$ and, in addition,
the set
$F_{r,i}$ is 
$(1+\frac{\rho_{2}}{D_{0}})2r$-separated (cf\. (7.22)), 
and therefore, in particular, 
$(1+\frac{\rho_{2}}{D_{0}})\frac{1}{2}r$-separated,
we conclude from Proposition 7.6 that
 $$
 \sum_{y\in F_{r,i}}
 \,
\mu(B(y,\tfrac{1}{2}r))^{q}
\le
k_{\bold i,\bold j,m}
\,
Z_{m}^{q}(\tfrac{1}{2}r)
\tag7.27
$$
and
$$
 \sum_{y\in F_{r,i}}
 \,
\mu(B(y,2r))^{q}
\le
k_{\bold i,\bold j,m}
\,
Z_{m}^{q}(2r)
\tag7.28
$$
for all $i$
where $k_{\bold i,\bold j,m}$ is the constant in Proposition 7.6.
Finally, we deduce from (7.26), (7.27) and (7.28) that
$$
 \align
 Q_{\bold i,\bold j}^{q}(r)
&\le
 \cases
 {\dsize
 \frac{\Omega_{d}}{2^{d}}
 \,
 \sum_{i=1}^{\chi}
 \,
k_{\bold i,\bold j,m}
\,
Z_{m}^{q}(\tfrac{1}{2}r)}
&\quad
\text{for $q<0$;}
\\
 {\dsize
\frac{\Omega_{d}}{2^{d}}
 \,
 \sum_{i=1}^{\chi}
 \,
k_{\bold i,\bold j,m}
\,
Z_{m}^{q}(2r)}
&\quad
\text{for $0\le q$.}
\endcases\\
&{}\\
&=
 \cases
 {\dsize
{}\,
c_{\bold i,\bold j,m}
\,
Z_{m}^{q}(\tfrac{1}{2}r)}
&\quad
\text{for $q<0$;}
\\
 {\dsize
\,
c_{\bold i,\bold j,m}
\,
Z_{m}^{q}(2r)}
&\quad
\text{for $0\le q$,}
\endcases
\endalign
$$
where
$c_{\bold i,\bold j,m}
=
\frac{\Omega_{d}}{2^{d}}
 \,
 \chi
 \,
k_{\bold i,\bold j,m}$.
This completes the proof.
\hfill$\square$

\bigskip
\bigskip


\heading{8. Proving that
 $$
 Z_{m}^{q}(r)
 \le
 \text{\rm constant}
 \,
 r^{-\gamma(q)}
 $$
}\endheading

\bigskip

The main purpose of this section is to prove Proposition 8.2.  
Let $\bold l$ be as in (7.2) and fix $q\in\Bbb R$.
Observe that that function
$\Xi^{q}:
 s\to
 \sum_{
   |\bold i|=|\bold l|\,,\,
   \bold i\not=\bold l
  }
 p_{\bold i}^{q}r_{\bold i}^{s}$
is continuous and strictly decreasing with
$\lim_{s\to-\infty}\Xi^{q}(s)=\infty$ and 
$\lim_{s\to\infty}\Xi^{q}(s)=0$. Hence, there exists a unique
$\gamma(q)\in\Bbb R$ such that
 $$
 \sum
 \Sb
   |\bold i|=|\bold l|\\
   \bold i\not=\bold l
   \endSb 
 p_{\bold i}^{q}r_{\bold i}^{\gamma(q)}
 =
 1\,.
 \tag8.1
 $$
Also, note that 
since
$\sum_{|\bold i|=|\bold l|,\bold i\not=\bold l}
 p_{\bold i}^{q}r_{\bold i}^{\gamma(q)}
 =
 1
 =
 \sum_{|\bold i|=|\bold l|}
 p_{\bold i}^{q}r_{\bold i}^{\beta(q)}$
it follows that 
 $$
 \gamma(q)<\beta(q)\,.
 \tag8.2
 $$
The main purpose of this section is to prove Proposition 8.2 saying that
for each positive integer $m$,
there is a constant $c_{m}$ such that
for all $r>0$, we have
 $$
Z_{m}^{q}(r)
  \,
  \le
  \,
  c_{m}r^{-\gamma(q)}
 $$
(recall, that $Z_{m}^{q}(r)$ is defined in (7.8)).
We begin with a small lemma.

\bigskip

\proclaim{Lemma 8.1}
Fix $q\in\Bbb R$.
Let $m\in\Bbb N$.
For $c>0$, we have
 $$
 \sup_{c\le r}
 \,
 r^{\gamma(q)}\,Z_{m}^{q}(r)<\infty
 $$
(recall, that $Z_{m}^{q}(r)$ is defined in (7.8)).
\endproclaim
\noindent{\it  Proof}\newline
\noindent 
Let
 $$
 \Gamma_{m}
 =
 \Big\{
 \bold h\in\Sigma^{*}
 \,\Big|\,
 a_{m}c\le r_{\bold h}
 \Big\}\,.
 $$ 
Observe that if $\bold h\in\Gamma_{m}$, then
$a_{m}c\le r_{\bold h}\le r_{\max}^{|\bold h|}$, whence
$|\bold h|\le \frac{\log a_{m}c}{\log r_{\max}}$,
and since there are only finitely many strings $\bold h\in\Sigma^{*}$ with
$|\bold h|\le \frac{\log a_{m}c}{\log r_{\max}}$, we deduce from this that
 $$
 |\Gamma_{m}|<\infty\,.
 \tag8.3
 $$

Next, note that if $\bold h\in\Sigma^{*}$ satisfies
$a_{m}r\le r_{\bold h}\le b_{m}r$, then
$r^{\gamma(q)}\le\frac{1}{a_{m}^{\gamma(q)}}r_{\bold h}^{\gamma(q)}$
if $\gamma(q)\ge 0$
and
$r^{\gamma(q)}\le\frac{1}{b_{m}^{\gamma(q)}}r_{\bold h}^{\gamma(q)}$
if $\gamma(q)\le 0$.
We conclude from this that  if $\bold h\in\Sigma^{*}$ satisfies
$a_{m}r\le r_{\bold h}\le b_{m}r$, then
 $$
 r^{\gamma(q)}
 \le
 \max\Big(
 \tfrac{1}{a_{m}^{\gamma(q)}}\,,\,
 \tfrac{1}{b_{m}^{\gamma(q)}}
 \Big)
 \,
 r_{\bold h}^{\gamma(q)}
 \tag8.4
$$
for all $q$.

Now combining (8.3) and (8.4) gives
 $$
 \align
 \sup_{c\le r}
  \,
 r^{\gamma(q)}\,Z_{m}^{q}(r)
&\le
  \sup_{c\le r}
  \,
  \sum
 \Sb
  \bold h\in\Sigma^{*}\\
  a_{m}r\le r_{\bold h}\le b_{m}r\\
  |\bold h|\ge|\bold l|\\
  \bold l\not\prec\bold h
  \endSb 
 r^{\gamma(q)} 
 p_{\bold h}^{q}\\
&\le
 \max\Big(
 \tfrac{1}{a_{m}^{\gamma(q)}}\,,\,
 \tfrac{1}{b_{m}^{\gamma(q)}}
 \Big)
 \,
  \sup_{c\le r}
  \,
  \sum
 \Sb
  \bold h\in\Sigma^{*}\\
  a_{m}r\le r_{\bold h}\le b_{m}r\\
  |\bold h|\ge|\bold l|\\
  \bold l\not\prec\bold h
  \endSb 
  r_{\bold h}^{\gamma(q)} 
 p_{\bold h}^{q}\\
&\le
 \max\Big(
 \tfrac{1}{a_{m}^{\gamma(q)}}\,,\,
 \tfrac{1}{b_{m}^{\gamma(q)}}
 \Big)
 \,
   \sum
 \Sb
  \bold h\in\Sigma^{*}\\
  a_{m}c\le r_{\bold h}
  \endSb 
  r_{\bold h}^{\gamma(q)} 
 p_{\bold h}^{q}\\ 
&\le
 \max\Big(
 \tfrac{1}{a_{m}^{\gamma(q)}}\,,\,
 \tfrac{1}{b_{m}^{\gamma(q)}}
 \Big)
 \,
   \sum_{\bold h\in\Gamma_{m}}
  r_{\bold h}^{\gamma(q)} 
 p_{\bold h}^{q}\\ 
&<
 \infty\,,
 \endalign
 $$
where the sum 
$ \sum_{\bold h\in\Gamma_{m}}
  r_{\bold h}^{\gamma(q)} 
 p_{\bold h}^{q}$ is finite since the set $\Gamma_{m}$ is finite by (8.3).
\hfill$\square$

\bigskip

\proclaim{Proposition 8.2}
Fix $q\in\Bbb R$.
Let $\bold l$ be as in (7.2).
Let $m\in\Bbb N$.
Then there exists a constant $c_{m}>0$ such that
for $r>0$, we have
 $$
 Z_{m}^{q}(r)
  \,
  \le
  \,
  c_{m}r^{-\gamma(q)}
 $$
(recall, that $Z_{m}^{q}(r)$ is defined in (7.8)).
\endproclaim
\noindent{\it  Proof}\newline
\noindent 
Choose $\delta_{m}>0$ such that
$b_{m}\frac{\delta_{m}}{r_{\min}^{|\bold l|}}\le r_{\min}^{|\bold l|}$.
Next,
define $W_{m}^{q}:(0,\infty)\to[0,\infty)$ by
$W_{m}^{q}(r)
=
 r^{\gamma(q)}Z_{m}^{q}(r)$.

Observe that 
for all $r>0$, we have 
 $$
 \align
 Z_{m}^{q}(r)
&=
 \sum
 \Sb
  \bold h\in\Sigma^{*}\\
  a_{m}r\le r_{\bold h}\le b_{m}r\\
  |\bold h|\ge|\bold l|\\
  \bold l\not\prec\bold h
  \endSb 
 p_{\bold h}^{q}\\
&=
\sum
  \Sb
  |\bold i|=|\bold l|\\
  \bold i\not=\bold l
  \endSb
 \, 
 \sum
 \Sb
  \bold j\in\Sigma^{*}\\
  a_{m}r\le r_{\bold i\bold j}\le b_{m}r\\
  \bold l\not\prec\bold i\bold j
  \endSb 
 p_{\bold i\bold j}^{q}\\
&\le
 \sum
  \Sb
  |\bold i|=|\bold l|\\
  \bold i\not=\bold l
  \endSb
 p_{\bold i} 
 \, 
 \sum
 \Sb
  \bold j\in\Sigma^{*}\\
  a_{m}\frac{r}{r_{\bold i}}\le r_{\bold j}\le b_{m}\frac{r}{r_{\bold i}}\\
  \bold l\not\prec\bold j
  \endSb 
 p_{\bold j}^{q}\,. 
 \tag8.5
 \endalign
 $$  
However,
if
$0<r<\delta_{m}$
and
$\bold i,\bold j\in\Sigma^{*}$ 
satisfy
$|\bold i|=|\bold l|$
and
$r_{\bold j}\le b_{m}\frac{r}{r_{\bold i}}$, then
$r_{\min}^{|\bold j|}
\le
r_{\bold j}
\le
b_{m}\frac{r}{r_{\bold i}}
\le
b_{m}\frac{\delta_{m}}{r_{\bold i}}
\le
b_{m}\frac{\delta_{m}}{r_{\min}^{|\bold i|}}
=
b_{m}\frac{\delta_{m}}{r_{\min}^{|\bold l|}}
\le
r_{\min}^{|\bold l|}$, whence
$|\bold j|\ge|\bold l|$.
We conclude from this and (8.5) that
if $0<r<\delta_{m}$, then
 $$
 \align
 Z_{m}^{q}(r)
&\le
 \sum
  \Sb
  |\bold i|=|\bold l|\\
  \bold i\not=\bold l
  \endSb
 p_{\bold i} 
 \, 
 \sum
 \Sb
  \bold j\in\Sigma^{*}\\
  |\bold j|\ge|\bold l|\\
  a_{m}\frac{r}{r_{\bold i}}\le r_{\bold j}\le b_{m}\frac{r}{r_{\bold i}}\\
  \bold l\not\prec\bold j
  \endSb 
 p_{\bold j}^{q}\\
&=
 \sum
  \Sb
  |\bold i|=|\bold l|\\
  \bold i\not=\bold l
  \endSb
 p_{\bold i} 
 Z_{m}^{q}(\tfrac{r}{r_{\bold i}})\,. 
 \endalign
 $$  
Hence, for $0<r<\delta_{m}$ we obtain
 $$
 \align
 W_{m}^{q}(r)
&=
 r^{\gamma(q)}Z_{m}^{q}(r)\\ 
&\le
 \sum
   \Sb
   |\bold i|=|\bold l|\\
   \bold i\not=\bold l
   \endSb
  p_{\bold i}^{q} 
  r_{\bold i}^{\gamma(q)}
  (\tfrac{r}{r_{\bold i}})^{\gamma(q)}
  Z_{m}^{q}(\tfrac{r}{r_{\bold i}})\\
&=
 \sum
   \Sb
   |\bold i|=|\bold l|\\
   \bold i\not=\bold l
   \endSb
  p_{\bold i}^{q} 
  r_{\bold i}^{\gamma(q)}
  W_{m}^{q}(\tfrac{r}{r_{\bold i}})\,.
  \tag8.6
  \endalign
 $$

Let 
$\Delta=r_{\max}^{|\bold l|}$.
It follows from (8.6) 
and definition (8.1) of $\gamma(q)$
that, if $0<a<\delta_{m}$, then
 $$
 \align
 \sup_{a\Delta\le r<\delta_{m}}W_{m}^{q}(r)
&\le
 \sup_{a\Delta\le r<\delta_{m}}
  \sum
   \Sb
   |\bold i|=|\bold l|\\
   \bold i\not=\bold l
   \endSb
  p_{\bold i}^{q} 
  r_{\bold i}^{\gamma(q)}
  W_{m}^{q}(\tfrac{r}{r_{\bold i}})\\
&\le
  \sum
   \Sb
   |\bold i|=|\bold l|\\
   \bold i\not=\bold l
   \endSb
  p_{\bold i}^{q} 
  r_{\bold i}^{\gamma(q)}
 \sup_{a\Delta\le r<\delta_{m}}
 W_{m}^{q}(\tfrac{r}{r_{\bold i}})\\
&\le
  \sum
   \Sb
   |\bold i|=|\bold l|\\
   \bold i\not=\bold l
   \endSb
  p_{\bold i}^{q} 
  r_{\bold i}^{\gamma(q)}
 \sup_{a\le s}
 W_{m}^{q}(s)\\
&=
 \sup_{a\le s}W_{m}^{q}(s)\,,
 \tag8.7
 \endalign
 $$
We deduce from (8.7) that if $0<a<\delta_{m}$, then
 $$
 \align
 \sup_{a\Delta\le r}W_{m}^{q}(r)
&\le
 \max
 \Bigg(
 \,
 \sup_{a\Delta\le r<\delta_{m}}W_{m}^{q}(r)
 \,,\,
 \sup_{\delta_{m}< r}W_{m}^{q}(r)
 \,
 \Bigg)\\
&\le
 \max
 \Bigg(
 \,
 \sup_{a\le r}W_{m}^{q}(r)
 \,,\,
 \sup_{\delta_{m}< r}W_{m}^{q}(r)
 \,
 \Bigg)\\
&=
  \sup_{a\le r}W_{m}^{q}(r)\,.
  \qquad\qquad
  \qquad\qquad
  \qquad
  \text{[since $a<\delta_{m}$]}
 \tag8.8
 \endalign
 $$

Next, choose a positive integer $k_{m}$ such that
$\Delta^{k_{m}}<\delta_{m}$. 
Repeated applications of (8.8) now gives
  $$
  \align
  \sup_{0<r}W_{m}^{q}(r)
 &=
  \sup_{k\in\Bbb N}
  \,\,\,\,
  \sup_{\Delta^{k+k_{m}}\le r}
  \,\,\,\,
  W_{m}^{q}(r)\\
&\le
  \sup_{k\in\Bbb N}
  \,
  \sup_{\Delta^{k+k_{m}-1}\le r}
  \,
  W_{m}^{q}(r)\\
&\,\,\,
  \vdots\\
&\le
  \sup_{k\in\Bbb N}
  \,\,\,\,\,\,
  \sup_{\Delta^{k_{m}}\le r}
  \,\,\,\,\,\,\,
  W_{m}^{q}(r)\\
&\le
  \sup_{\Delta^{k_{m}}\le r}
  \,
  W_{m}^{q}(r)\\
&=
  \sup_{\Delta^{k_{m}}\le r}
  \,
 r^{\gamma(q)}\,Z_{m}^{q}(r)
  \tag8.9
  \endalign
  $$
Finally, combining (8.9) and Lemma 8.1 shows that
$r^{\gamma(q)}Z_{m}^{q}(r)
=
W_{m}^{q}(r)
\le
\sup_{s>0}W_{m}^{q}(s)
\le
\sup_{\Delta^{k_{m}}\le s}
  \,
 s^{\gamma(q)}\,Z_{m}^{q}(s)<\infty$.
This completes the proof of Proposition 8.2.
\hfill$\square$

\bigskip
\bigskip


\heading{9. Proof of Theorem 3.3
}\endheading

\bigskip

The purpose of this section is to prove Theorem 3.3.
The proof are based on renewal theory and, in particular, 
a recent reneval theorem by 
Levitin \& Vassiliev [LeVa].
Below we state Levitin \& Vassiliev's renewal theorem.

\bigskip

\proclaim{Theorem 9.1.
Levitin \& Vassiliev's renewal theorem [LeVa]}
Let
$t_{1},\dots, t_{N}>0$
and
$p_{1},\ldots,p_{N}>0$ with $\sum_{i}p_{i}=1$.
Define the probability measure $P$ by
 $$
 P
 =
 \sum_{i}p_{i}\delta_{t_{i}}\,.
 $$
 
Let $\lambda,\Lambda:\Bbb R\to\Bbb R$
be real valued functions 
satisfying the following conditions:
\roster 
\item"(i)"
The function $\lambda$ is piecewise continuous;
\item"(ii)"
There are constants 
$c,k>0$ such that
 $$
 |\lambda(t)|
 \le 
 ce^{-k|t|}
 $$
for all $t\in\Bbb R$;
\item"(iii)"
We have
 $$
 \Lambda(t)
 \to
 0
 \,\,
 \text{as $t\to-\infty$;}
 $$
\item"(iv)"
We have
 $$
 \Lambda(t)
 =
 \int\Lambda(t-s)\,dP(s)
 +
 \lambda(t)
 $$
for all $t\in\Bbb R$.
\endroster

The the following holds:
\roster
\item"(1)"
The
non-arithmetic case:
If $\{t_{1},\ldots,t_{N}\}$
is not contained in a discrete additive subgroup of $\Bbb R$,
then
 $$
 \Lambda(t)
 =
 c
 +
 \varepsilon(t)
 $$
for all $t\in\Bbb R$ where
 $$
 \align
 c
&=
 \frac{1}{\int s\,dP(s)}
 \int\lambda(s)\,ds\,,\\
 \varepsilon(t)
& 
 \to
 0
 \,\,
 \text{as $t\to\infty$.}
\endalign
$$
In addition,
 $$
  \frac{1}{T}\int_{0}^{T}\Lambda(t)\,dt
 \to
  c
 =
 \frac{1}{\int s\,dP(s)}
 \int\lambda(s)\,ds
 \,\,\text{as $T\to\infty$.}
 \tag9.1
 $$

\item"(2)"
The
arithmetic case:
If $\{t_{1},\ldots,t_{N}\}$
is contained in a discrete additive subgroup of $\Bbb R$
and
$\langle t_{1},\ldots,t_{N}\rangle=u\Bbb Z$
with $u>0$,
then
 $$
 \Lambda(t)
 =
 \pi(t)
 +
 \varepsilon(t)
 $$
for all $t\in\Bbb R$ where
 $$
 \align
 \pi(t)
&=
 \frac{1}{\int s\,dP(s)}
 u
 \sum_{n\in\Bbb Z}
 \lambda(t+nu)\,,\\
 \varepsilon(t)
& 
 \to
 0
 \,\,
 \text{as $t\to\infty$.}
\endalign
$$
In addition, the function $\pi$ is $u$-periodic
and
 $$
  \frac{1}{T}\int_{0}^{T}\Lambda(t)\,dt
 \to
  c
 =
 \frac{1}{\int s\,dP(s)}
 \int\lambda(s)\,ds
 \,\,\text{as $T\to\infty$.}
 \tag9.2
 $$
 \endroster
\endproclaim
\noindent{\it  Proof}\newline
\noindent
All statements, except (9.1) and (9.2), follow 
[LeVa].
Below
we prove (9.1) and (9.2).
Indeed,
(9.1) follows immediately
and (9.2) is proved as follows.
Namely,
since $\pi$ is $u$-periodic we conclude that
 $$
 \align
 \frac{1}{T}\int_{0}^{T}\Lambda(t)\,dt
&=
\frac{1}{T}\int_{0}^{T}\pi(t)\,dt
+
\frac{1}{T}\int_{0}^{T}\varepsilon(t)\,dt\\
&\to
 \frac{1}{u}\int_{0}^{u}\pi(t)\,dt\\
&=
  \frac{1}{\int t\,dP(t)}
   \int_{0}^{u}
 \sum_{n\in\Bbb Z}
 \lambda(t+nu)
 \,dt\,.
\tag9.3
 \endalign
 $$ 
Now observe that since
$|\lambda(t)|
 \le 
 ce^{-k|t|}$
for all $t\in\Bbb R$
and
$\int ce^{-k|t|}\,dt<\infty$, it follows from 
two applications of Lebesgue's Dominated Convergence Theorem 
and the fact that $\pi$ is $u$-periodic
that
 $$
 \align
  \int_{0}^{u}
 \sum_{n\in\Bbb Z}
 \lambda(t+nu)
 \,dt
&= 
  \sum_{n\in\Bbb Z}
  \int_{0}^{u}
 \lambda(t+nu)
 \,dt\\
&=
  \sum_{n\in\Bbb Z}
  \int_{nu}^{(n+1)u}
 \lambda(t)
 \,dt\\
&=
  \sum_{n\in\Bbb Z}
  \int
 \bold 1_{[nu,(n+1)u)}(t)\,\lambda(t)
 \,dt\\
&=
  \int
  \sum_{n\in\Bbb Z}
  \bold 1_{[nu,(n+1)u)}(t)\,\lambda(t)
 \,dt\\ 
&=
  \int
 \lambda(t)
 \,dt\,.
 \tag9.4
 \endalign
 $$
Finally, combining (9.3) and (9.4)
shows that
 $$ 
 \align
 \frac{1}{T}\int_{0}^{T}\Lambda(t)\,dt
&\to
  \frac{1}{\int t\,dP(t)}
  \int
 \lambda(t)
 \,dt\,.
 \endalign
 $$
This completes the proof.
\hfill$\square$

\bigskip

In order to prove Theorem 3.3,
we will apply
Levitin \& Vassiliev's renewal theorem
to the probability measure $P=P_{q}$
and the
functions
$\lambda=\lambda_{q}^{0}$
and
$\Lambda=\Lambda_{q}^{0}$
defined as follows.
Namely,
first
recall that $\lambda_{q}:(0,\infty)\to\Bbb R$ is defined by
 $$
 \align
 \lambda_{q}(r)
&=
V_{\mu,r}^{q}(K)
 -
 \sum_{i}p_{i}^{q}
 \,
 \bold 1_{(0,r_{i}]}(r)
\,
V_{\mu,r_{i}^{-1}r}^{q}(K)
\endalign
$$
and
define $\Lambda_{q}:(0,\infty)\to\Bbb R$ by
 $$
 \Lambda_{q}(r)
 =
 V_{\mu,r}^{q}(K)\,.
 \qquad\qquad
 \qquad\qquad
 \qquad\qquad
 \,\,\,\,{}
 $$
Also, define
$\lambda_{q}^{0}:\Bbb R\to\Bbb R$ by
 $$
 \lambda_{q}^{0}(t)
 =
 \bold 1_{[0,\infty)}(t)
\,
e^{-t\beta(q)}
 \lambda_{q}(e^{-t})
 \qquad\qquad
 \qquad\,\,\,
 {}
 $$
and define $\Lambda_{q}^{0}:\Bbb R\to\Bbb R$ by
 $$
 \Lambda_{q}^{0}(t)
= 
 \bold 1_{[0,\infty)}(t)
\,
e^{-t\beta(q)}
 \Lambda_{q}(e^{-t})\,.
 \qquad\qquad
 \qquad
 {}
 $$
Finally,  define the probability measure $P_{q}$ by
 $$
 P_{q}
 =
 \sum_{i}p_{i}^{q}r_{i}^{\beta(q)}\delta_{\log r_{i}^{-1}}\,.
 $$
Next, we show 
(in Propositions 9.2--9.5)
that
probability measure $P=P_{q}$
and the
functions
$\lambda=\lambda_{q}^{0}$
and
$\Lambda=\Lambda_{q}^{0}$
satisfy conditions (i)--(iv)
in 
Levitin \& Vassiliev's renewal theorem. 
 
\bigskip

\proclaim{Proposition 9.2}
Fix $q\in\Bbb R$
and
assume 
that one of the following conditions is satisfied:
\roster
\item"(i)"
The OSC is satisfied and $0\le q$;
\item"(ii)"
The SSC is satisfied.
\endroster
Then
the function $\lambda_{q}^{0}$ is piecewise continuous.
\endproclaim
\noindent{\it  Proof}\newline
\noindent
Define $f:(0,\infty)\to\Bbb R$ by
 $$
 f(r)
 =
 \int\limits_{B(K,r)}
 \mu(B(x,r))^{q}
 \,d\Cal L^{d}(x)\,.
 $$
It clearly suffices to show that $f$ is continuous.
We therefore fix $r_{0}>0$
and prove that $f$ is continuous at $r_{0}$.
For $0<h<1$, define
$\varphi,\varphi_{h}^{-},\varphi_{h}^{+}:\Bbb R^{d}\to\Bbb R$ by
 $$
 \align
 \varphi(x)
&=
 \bold 1_{B(K,r_{0})}(x)
 \,
 \mu(B(x,r_{0}))^{q}\,,\\
\varphi_{h}^{-}(x)
&=
 \bold 1_{B(K,r_{0}-h)}(x)
 \,
 \mu(B(x,r_{0}-h))^{q}\,,\\
 \varphi_{h}^{+}(x)
&=
 \bold 1_{B(K,r_{0}+h)}(x)
 \,
 \mu(B(x,r_{0}+h))^{q}\,.
 \endalign
 $$

Since
$\cup_{0<h<1}B(x,r_{0}-h)=B(x,r_{0})$,
we conclude that
$\mu(B(x,r_{0}-h))\to\mu(B(x,r_{0}))$ as $h\searrow0$, whence
$\varphi_{h}^{-}(x)\to\varphi(x)$
for all $x$ as $h\searrow0$.
We also have
$|\varphi_{h}^{-}(x)|\le\bold 1_{B(K,r_{0})}(x)$
for all $x$ where
$\int\bold 1_{B(K,r_{0})}(x)\,d\Cal L^{d}(x)=\Cal L^{d}(B(K,r_{0}))<\infty$,
and we therefore conclude from 
Lebesgue's Dominated Convergence Theorem that
 $$
 f(r_{0}-h)
 =
 \int\varphi_{h}^{-}(x)\,d\Cal L^{d}(x)
 \to
 \int\varphi(x)\,d\Cal L^{d}(x)
 =
 f(r_{0})
 \tag9.5
 $$
as $h\searrow0$.

Similarly,
since
$\cap_{0<h<1}B(x,r_{0}+h)=\overline{B(x,r_{0})}$,
we conclude that
$\mu(B(x,r_{0}+h))\to\mu(\overline{B(x,r_{0})})$ as $h\searrow0$, whence
$\varphi_{h}^{+}(x)
\to
\bold 1_{\overline{B(K,r_{0})}}(x)
 \,
 \mu(\overline{B(x,r_{0})})^{q}
$
for all $x$ as $h\searrow0$.
We also have
$|\varphi_{h}^{+}(x)|\le\bold 1_{B(K,r_{0}+1)}(x)$
for all $x$ where
$\int\bold 1_{B(K,r_{0}+1)}(x)\,d\Cal L^{d}(x)=\Cal L^{d}(B(K,r_{0}+1))<\infty$,
and we therefore conclude from the Dominated Convergence Theorem that
 $$
 f(r_{0}+h)
 =
 \int\varphi_{h}^{-}(x)\,d\Cal L^{d}(x)
 \to
 \int\bold 1_{\overline{B(K,r)}}(x)
 \,
 \mu(\overline{B(x,r)})^{q}
 \,d\Cal L^{d}(x)
 \tag9.6
 $$
as $h\searrow0$. 
Since clearly
$\Cal L^{d}(\overline{B(K,r_{0})})
=
\Cal L^{d}(B(K,r_{0}))$, we deduce from (9.6) that
 $$
 f(r_{0}+h)
  \to
 \int\bold 1_{B(K,r)}(x)
 \,
 \mu(\overline{B(x,r)})^{q}
 \,d\Cal L^{d}(x)
 \tag9.7
 $$
as $h\searrow0$. 
Next, it is proved in 
[Mat]
that
either $K$ lies in a $l$-dimensional affine subspace 
of $\Bbb R^{d}$ with $l<d$ 
or
$\mu(K\cap\Gamma)=0$ for every $l$ dimensional $C^{1}$-submanifold
$\Gamma\subseteq\Bbb R^{d}$ with $0<l<d$. 
This implies that
$\mu(\overline{B(x,r_{0})}\setminus B(x,r_{0}))=0$ for all $x$, and we therefore deduce from 
(9.7) that
 $$
 f(r_{0}+h)
  \to
 \int\bold 1_{B(K,r_{0})}(x)
 \,
 \mu(B(x,r_{0}))^{q}
 \,d\Cal L^{d}(x)
 =
 f(r_{0})
 \tag9.8
 $$
as $h\searrow0$.

Finally, it follows from (9.5) and (9.8) that $f$ is continuous at $r_{0}$. 
\hfill$\square$

\bigskip

\proclaim{Proposition 9.3}
Fix $q\in\Bbb R$
and
assume 
that one of the following conditions is satisfied:
\roster
\item"(i)"
The OSC is satisfied and $0\le q$;
\item"(ii)"
The SSC is satisfied.
\endroster
Then
there is a constant $c>0$ such that
$|\lambda_{q}^{0}(t)|
\le
c
e^{-(\beta(q)-\gamma(q))|t|}$
for all $t\in\Bbb R$.
\endproclaim
\noindent{\it  Proof}\newline
\noindent
For a positive integer $m$ and $\bold i,\bold j\in\Sigma^{*}$, let
$r_{m}$ be the number defined in (7.5).
Also,
let
$c_{\bold i,\bold j,m}$ be the constant in Proposition 7.8
and
let 
$c_{m}$ be the constant in Proposition 8.2.

Choose $t_{0}>0$ such that
$e^{-t}<\min(r_{\min},\frac{1}{2}r_{1})$ for $t\ge t_{0}$,
and
observe that for $t\ge t_{0}$, we have
 $$
 \align
 |\lambda_{q}^{0}(t)|
&=
 \bold 1_{[0,\infty)}(t)
\,
e^{-t\beta(q)}
\,
|\lambda_{q}(e^{-t})|\\ 
&=
 e^{-t\beta(q)}
 \Bigg|
V_{\mu,e^{-t}}^{q}(K)
 -
 \sum_{i}p_{i}^{q}
 \,
 \bold 1_{(0,r_{i}]}(e^{-t})
\,
V_{\mu,r_{i}^{-1}e^{-t}}^{q}(K)
 \Bigg|\\
&=
 e^{-t\beta(q)}
 \Bigg|
V_{\mu,e^{-t}}^{q}(K)
 -
 \sum_{i}p_{i}^{q}
V_{\mu,r_{i}^{-1}e^{-t}}^{q}(K)
 \Bigg|
 \qquad\qquad
 \qquad\quad
 \text{[since $e^{-t}<r_{\min}\le r_{i}$]}
 \\
&\le
  e^{-t\beta(q)}
 \sum
 \Sb
  |i|=|j|=1\\
  i\not=j
 \endSb
 Q_{i,j}^{q}(e^{-t})
 \qquad\qquad
 \qquad\qquad
 \qquad\qquad
 \quad\,\,\,\,\,
 \text{[by Proposition 6.3]}\\
&{}\\ 
&\le
 \cases
 {\dsize
  e^{-t\beta(q)}
 \sum
 \Sb
  |i|=|j|=1\\
  i\not=j
 \endSb
 c_{i,j,1}
 Z_{1}^{q}(\tfrac{1}{2}e^{-t})
 }
&\quad\qquad
  \,
 \text{for $q<0$;}
 \\ 
  {\dsize
  e^{-t\beta(q)}
 \sum
 \Sb
  |i|=|j|=1\\
  i\not=j
 \endSb
 c_{i,j,1}
 Z_{1}^{q}(2e^{-t})
 }
&\quad\qquad
  \,
 \text{for $0\le q$}
\endcases 
 \qquad \text{[by Proposition 7.8]}\\
&{}\\ 
&\le
 \cases
 {\dsize
  e^{-t\beta(q)}
 \sum
 \Sb
  |i|=|j|=1\\
  i\not=j
 \endSb
 c_{i,j,1}
 \,
 c_{1}
 \,
 (\tfrac{1}{2}e^{-t})^{-\gamma(q)}
 }
&\quad
 \text{for $q<0$;}
 \\ 
  {\dsize
  e^{-t\beta(q)}
 \sum
 \Sb
  |i|=|j|=1\\
  i\not=j
 \endSb
 c_{i,j,1}
 \,
 c_{1}
 \,
 (2e^{-t})^{-\gamma(q)}
 }
&\quad
 \text{for $0\le q$}
\endcases 
 \qquad
 \text{[by Proposition 8.2]}\\
&{}\\ 
&=
 ce^{-(\beta(q)-\gamma(q))t} 
 \tag9.9
 \endalign 
$$
where
$c=
c_{1}
\,
\max((\frac{1}{2})^{-\gamma(q)},2^{-\gamma(q)})
\,
\sum_{|i|=|j|=1,i\not=j}
 c_{i,j,1}$.

Next, since $\lambda_{q}^{0}$ is piecewise continuous 
(by Proposition 9.2), 
we conclude that $\lambda_{q}^{0}$ is bounded on
the compact interval
$[0,t_{0}]$, i\.e\. there is a constant $M_{0}$ such that 
$|\lambda_{q}^{0}(t)|\le M_{0}$ for all $t\in[0,t_{0}]$.
It follows from this and (9.9) that
 $$
 |\lambda_{q}^{0}(t)|
 \le
 \max
 \Big(
 \tfrac{M_{0}}{e^{-(\beta(q)-\gamma(q))t_{0}}}
 \,,\,
 c
 \Big)
 \,e^{-(\beta(q)-\gamma(q))t}
 \tag9.10
 $$ 
for all $t\ge 0$.

The statement now follows from (9.10) and the fact that
$\lambda_{q}^{0}(t)=0$ for all $t<0$.
 \hfill$\square$

\bigskip

\proclaim{Proposition 9.4}
Fix $q\in\Bbb R$.
Then
$\Lambda_{q}^{0}(t)
 \to
 0$
as $t\to-\infty$.
\endproclaim
\noindent{\it  Proof}\newline
\noindent
Indeed, this follows trivially from the fact that
 $\Lambda_{q}^{0}(t)=0$ for all $t<0$.
 \hfill$\square$

\bigskip

\proclaim{Proposition 9.5}
Fix $q\in\Bbb R$.
Then
$\Lambda_{q}^{0}(t)
=
 \int\Lambda_{q}^{0}(t-s)\,dP_{q}(s)
\,
+
\,
\lambda_{q}^{0}(t)$
for all $t\in\Bbb R$.
\endproclaim
\noindent{\it  Proof}\newline
\noindent
It follows immediately from the definitions of
$\lambda_{q}^{0}$, $\Lambda_{q}^{0}$ and $P_{q}$ that
 $$
 \align
  \Lambda_{q}^{0}(t)
&=
 \bold 1_{[0,\infty)}(t)
\,
e^{-t\beta(q)}
 \Lambda_{q}(e^{-t})\\    
&=
  \bold 1_{[0,\infty)}(t)
\,
e^{-t\beta(q)} 
 \Bigg(
 \sum_{i}p_{i}^{q}
 \,
 \bold 1_{(0,r_{i}]}(e^{-t})
\,
V_{\mu,r_{i}^{-1}e^{-t}}^{q}(K) 
 \,
 +
 \,
 \lambda_{q}(e^{-t})
 \Bigg)\\
&=
 \sum_{i}p_{i}^{q}
 \,
 e^{-t\beta(q)} 
 \,
 \bold 1_{(0,r_{i}]}(e^{-t})
 \,
  \bold 1_{[0,\infty)}(t)
\,
V_{\mu,r_{i}^{-1}e^{-t}}^{q}(K) 
\,
+
\,
\lambda_{q}^{0}(t)\\
&=
  \sum_{i}
 p_{i}^{q}r_{i}^{\beta(q)}
 \,
 \bold 1_{[0,\infty)}(t-\log r_{i}^{-1})
 \,
  \bold 1_{[0,\infty)}(t)
\,
e^{-\beta(q)(t-\log r_{i}^{-1})}
V_{\mu,e^{-(t-\log r_{i}^{-1})}}^{q}(K)
+
\lambda_{q}^{0}(t)\\
&=
  \sum_{i}
 p_{i}^{q}r_{i}^{\beta(q)}
 \,
 \bold 1_{[0,\infty)}(t-\log r_{i}^{-1})
\,
e^{-\beta(q)(t-\log r_{i}^{-1})}
V_{\mu,e^{-(t-\log r_{i}^{-1})}}^{q}(K)
\,
+
\,
\lambda_{q}^{0}(t)\\
&=
  \sum_{i}
 p_{i}^{q}r_{i}^{\beta(q)}
 \Lambda_{q}^{0}(t-\log r_{i}^{-1})
\,
+
\,
\lambda_{q}^{0}(t)\\
&=
 \int\Lambda_{q}^{0}(t-s)\,dP_{q}(s)
\,
+
\,
\lambda_{q}^{0}(t)
\endalign
$$
for all $t\in\Bbb R$.
\hfill$\square$

\bigskip

We can now prove Theorem 3.3.

\bigskip

\noindent{\it  Proof of Theorem 3.3}\newline
\noindent
It follows from Propositions 9.2--9.5 that
Theorem 9.1 can be applied 
to
the probability measure $P=P_{q}$
and the
functions
$\lambda=\lambda_{q}^{0}$
and
$\Lambda=\Lambda_{q}^{0}$.
We divide the 
proof into two cases.

\bigskip

\noindent
{\it Case 1: 
If $\{\log r_{1}^{-1},\dots,\log r_{N}^{-1}\}$, is 
not contained in a discrete additive subgroup of $\Bbb R$.}
If $\{\log r_{1}^{-1},\dots,\log r_{N}^{-1}\}$, is 
not contained in a discrete additive subgroup of $\Bbb R$, then
Theorem 9.1 implies that
 $$
 \Lambda_{q}^{0}(t)
 =
 c_{q}+\varepsilon_{q}^{0}(t)
 $$
where
$c_{q}\in\Bbb R$ is the constant given by
$$
 \align
 c_{q} 
&=
 \frac{1}{\int s\,dP_{q}(s)}
 \int \lambda_{q}^{0}(s)\,ds\\
&=
 \frac{1}{-\sum_{i}p_{i}^{q}r_{i}^{\beta(q)}\log r_{i}	}
 \,
 \int_{0}^{\infty}
 e^{-s\beta(q)}
 \lambda_{q}(e^{-s})
 \,ds\\ 
 &=
 \frac{1}{-\sum_{i}p_{i}^{q}r_{i}^{\beta(q)}\log r_{i}	}
 \,
 \int_{0}^{1}
 r^{\beta(q)}
 \lambda_{q}(r) 
 \,
 \frac{dr}{r}
 \endalign
 $$
and
 $$ 
 \varepsilon_{q}^{0}(t)
\to
 0
 \,\,\text{as $t\to\infty$.}  
 \qquad\qquad
 \qquad\qquad
 \qquad\qquad
 \quad\,\,\,\,
  $$
In particular, we have
 $$
  \align
  r^{\beta(q)}\,V_{\mu,r}^{q}(K)
 &=
  \Lambda_{q}^{0}(\log\tfrac{1}{r})
  =
  c_{q}
  +
  \varepsilon_{q}(r)
  \tag9.11
  \endalign
  $$ 
where
$  \varepsilon_{q}(r)
=
\varepsilon_{q}^{0}(\log\tfrac{1}{r})
\to 0$ as $r\searrow 0$.

Finally,
it follows  from 
(9.11) that
 $$
 r^{\beta(q)}\,V_{\mu,r}^{q}(K)
 \to
 c_{q}
 \,\,
 \text{as $r\searrow0$.}
 $$

 \bigskip

 \noindent
{\it Case 2:
If $\{\log r_{1}^{-1},\dots,\log r_{N}^{-1}\}$ is 
contained in a discrete additive subgroup of $\Bbb R$.}
If $\{\log r_{1}^{-1},\dots,\log r_{N}^{-1}\}$ is 
contained in a discrete additive subgroup of $\Bbb R$
and
$\langle t_{1},\ldots,t_{N}\rangle=u\Bbb Z$
with $u>0$, 
then
Theorem 9.1 implies that
 $$
 \Lambda_{q}^{0}(t)
 =
 \pi_{q}^{0}(r)+\varepsilon_{q}^{0}(t)
 $$
where
$\pi_{q}^{0}:\Bbb R\to\Bbb R$
is the function  given by 
$$
 \align
 \pi_{q}^{0}(t) 
&=
 \frac{1}{\int s\,dP_{q}(s)}
 \,
 u
 \sum_{n\in\Bbb Z}\lambda_{q}^{0}(t+nu)\\
&=
 \frac{1}{-\sum_{i}p_{i}^{q}r_{i}^{\beta(q)}\log r_{i}	} \,
 u
 \sum_{n\in\Bbb Z}\lambda_{q}^{0}(t+nu)
 \endalign
 $$
and
 $$ 
  \varepsilon_{q}^{0}(t)
 \to
 0
 \,\,\text{as $t\to\infty$.}  
 \qquad\qquad
 \qquad
 $$
Moreover, we have
 $$
 \pi_{q}^{0}(t+u)
 =
 \pi_{q}^{0}(t)
 $$
for all $t\in\Bbb R$, i\.e\.
$\pi_{q}^{0}$ is additively periodic
with period equal to $u$.
In particular, we have
 $$
  \align
  r^{\beta(q)}\,V_{\mu,r}^{q}(K)
&=
  \Lambda_{q}^{0}(\log\tfrac{1}{r})
  =
  \pi_{q}(r)
  +
  \varepsilon_{q}(r)
  \endalign
  $$ 
where
$\pi_{q}:\Bbb R\to\Bbb R$
is the function  given by 
 $$
 \align
 \pi_{q}(r) 
&=
 \pi_{q}^{0}(\log\tfrac{1}{r})\\ 
&=
 \frac{1}{-\sum_{i}p_{i}^{q}r_{i}^{\beta(q)}\log r_{i}	}
  \,
 u
 \sum_{n\in\Bbb Z}\lambda_{q}^{0}(\log\tfrac{1}{r}+nu)\\
&=
 \frac{1}{-\sum_{i}p_{i}^{q}r_{i}^{\beta(q)}\log r_{i}	}
  \,
 u
 \sum_{n\in\Bbb Z}
 \bold 1_{[0,\infty)}(\log\tfrac{1}{r}+nu)
 \,
 e^{-\beta(q)(\log\frac{1}{r}+nu)}
 \,
 \lambda_{q}(e^{-(\log\frac{1}{r}+nu)})\\
&=
 \frac{1}{-\sum_{i}p_{i}^{q}r_{i}^{\beta(q)}\log r_{i}	}
  \,
 u
 \sum
  \Sb
  n\in\Bbb Z\\
  {}\\
  re^{nu}\le 1
  \endSb
 (re^{nu})^{\beta(q)}
 \,
 \lambda_{q}(re^{nu})\\
  \endalign
 $$
and
$  \varepsilon_{q}(r)
=
\varepsilon_{q}^{0}(\log\tfrac{1}{r})
\to 0$ as $r\searrow 0$.
Moreover, 
since 
$\pi_{q}^{0}$ is additively periodic
with period equal to $u$, we have
 $$
 \pi_{q}(e^{u}r)
 =
 \pi_{q}^{0}(\log\tfrac{1}{e^{u}r})
 =
\pi_{q}^{0}(\log\tfrac{1}{r}-u)
= 
 \pi_{q}^{0}(\log\tfrac{1}{r})
 =
 \pi_{q}(r)
 $$
for all $r>0$, i\.e\.
$\pi_{q}$ is multiplicatively
periodic
with period equal to $e^{u}$.

 Finally
 it
 follows from Theorem 9.1 that
 $$
 \frac{1}{T}\int_{0}^{T}\Lambda_{q}^{0}(t)\,dt
 \to
 c_{q}
 \,\,
 \text{as $T\to\infty$.}
 $$ 
However,
since
 $$
 \align 
\frac{1}{T}\int_{0}^{T}\Lambda_{q}^{0}(t)\,dt
&=
 \frac{1}{T}\int_{0}^{T}e^{-t\beta(q)}\,V_{\mu,e^{-t}}^{q}(K)\,dt\\
&=
 \frac{1}{-\log e^{-T}}
 \int_{e^{-T}}^{1}s^{\beta(q)}\,V_{\mu,s}^{q}(K)\,\frac{ds}{s}\,,
 \endalign
 $$
we now conclude that
 $$
  \frac{1}{-\log r}
 \int_{r}^{1}s^{\beta(q)}\,V_{\mu,s}^{q}(K)\,\frac{ds}{s}
 \to
 c_{q}
  \,\,
 \text{as $r\searrow 0$.}
 $$
This completes the proof of Theorem  3.3 in Case 2.
\hfill$\square$


\newpage

${}$

\bigskip
\bigskip
\bigskip
\bigskip
\bigskip
\bigskip
\bigskip
\bigskip

\centerline{\bigletter Part 3:}

\bigskip

\centerline{\bigletter Proofs of the Results from Section 4}


\bigskip
\bigskip
\bigskip
\bigskip

\heading{10. Analysis of
$\Cal H_{\mu}^{q,\beta(q)}$
}\endheading

The purpose of this section is 
twofold.
Firstly, we
prove Theorem 10.3.
Secondly, we apply
Theorem 10.3 to
obtain an explicit formula for the multifractal Hausdorff measure 
$ \Cal H_{\mu}^{q,\beta(q)}(S_{\bold i}K)$ of $S_{\bold i}K$, cf\. Proposition 10.5.(6);
this formula plays an important part in Sections 11--12
 when identifying the weak limit of the
(suitably
normalized)
tube measure
$\Cal I_{\mu,r}^{q}$.
We now turn towards the proof of Theorem 10.3.
We begin with two auxiliary lemmas.

\bigskip

\proclaim{Lemma 10.1}
Let $\mu$ be a probability measure on $\Bbb R^{d}$.
Fix $q,t\in\Bbb R$ and $E\subseteq\Bbb R^{d}$.
Let  $\varepsilon>0$.
Then there is a constant $c>0$ and a sequence $(r_{n})_{n}$ of 
positive real numbers with
$r_{n}\to0$ satisfying the following:
for all $F\subseteq E$ and $n$, 
the inequality below is satisfied, namely,
 $$
 \left.
 \matrix
 \text{for $q<0$, we have}
 \qquad
 \overline{\Cal H}_{\mu,2r_{n}}^{q,t}(F)\\
 \text{\,for $0\le q$, we have}
 \qquad
 \overline{\Cal H}_{\mu,\frac{1}{2}r_{n}}^{q,t}(F)\\
 \endmatrix
 \right\}
 \le
 c
 \,
 r_{n}^{t-(\underline\dim_{\Min,\mu}^{q}(E)+\varepsilon)}\,.
  $$
\endproclaim 
\noindent{\it  Proof}\newline
Since
$\liminf_{r\searrow0}\frac{\log V_{\mu,r}^{q}(E)}{-\log r}
 <
 \underline\dim_{\Min,\mu}^{q}(E)+\varepsilon$,
we can find 
a sequence $(r_{n})_{n}$ of real numbers with
$r_{n}\to0$
such that
$\frac{\log V_{\mu,r_{n}}^{q}(E)}{-\log r_{n}}
 <
 \underline\dim_{\Min,\mu}^{q}(E)+\varepsilon$
for all $n$, whence
 $$
 V_{\mu,r_{n}}^{q}(E)
 \le
 r_{n}^{-(\underline\dim_{\Min,\mu}^{q}(E)+\varepsilon)}
 \tag10.1
 $$
for all $n$.

Next, fix $F\subseteq E$
and let $\delta>0$.

It is clear that we can choose a countable centered covering
$(\,B(x_{n,i},r_{n})\,)_{i\in I_{n}}$
of $F$ such that there are 
$2^{d}$ subsets
$I_{n,1},\ldots, I_{n,2^{d}}$
of $I$ satisfying:
 $$
 \gather
 \bigcup_{k}I_{n,k}=I_{n}\,,\\
 I_{n,k}\cap I_{n,l}=\varnothing
 \,\,\,\,
 \text{for $k\not=l$,}\\
 \text{
 the sets
 $(\,B(x_{n,i},r_{n})\,)_{i\in I_{n,k}}$
 are pairwise disjoint.
 }
 \endgather 
 $$

Writing $\Omega_{d}$ for the Lebesgue measure of the unit ball in $\Bbb R^{d}$,
we have
 $$  
 \aligned
 \overline{\Cal H}_{\mu,2r_{n}}^{q,t}(F)
&\le
 \sum_{i\in I_{n}}
 \mu(B(x_{n,i},2r_{n}))^{q}(4r_{n})^{t}\\
&=
 \frac{4^{t}}{\Omega_{d}}
 \,
 \sum_{i\in I_{n}}
 \frac{1}{r_{n}^{d}}
 \,
 \int\limits_{B(x_{n,i},r_{n})}
 \mu(B(x_{n,i},2r_{n}))^{q}
 \,d\Cal L^{d}(x)
 \,
 r_{n}^{t}
 \quad
 \text{for $q<0$,}\\
&{}\\ 
 \overline{\Cal H}_{\mu,\frac{1}{2}r_{n}}^{q,t}(F)
&\le
 \sum_{i\in I_{n}}
 \mu(B(x_{n,i},\tfrac{1}{2}r_{n}))^{q}r_{n}^{t}\\
&=
 \frac{2^{d}}{\Omega_{d}}
 \,
 \sum_{i\in I_{n}}
 \frac{1}{r_{n}^{d}}
 \,
 \int\limits_{B(x_{n,i},r_{n})}
 \mu(B(x_{n,i},\tfrac{1}{2}r_{n}))^{q}
 \,d\Cal L^{d}(x)
 \,
 r_{n}^{t}
 \quad
 \text{for $0\le q$.}
 \endaligned
 \tag10.2
 $$
However,
if $x\in B(x_{n,i},r_{n})$, then clearly
$B(x_{n,i},\frac{1}{2})\subseteq B(x,r_{n})\subseteq B(x_{n,i},2r_{n})$,
whence 
 $$
 \aligned
 \mu(B(x_{n,i},2r_{n}))^{q}
&\le
 \mu(B(x,r_{n}))^{q}
 \quad
 \text{for $q<0$,}\\
\mu(B(x_{n,i},\tfrac{1}{2}r_{n}))^{q}
&\le
 \mu(B(x,r_{n}))^{q}
 \quad
 \text{for $0\le q$.}
 \endaligned 
 \tag10.3
 $$
It follows from (10.2) and (10.3) that
 $$  
 \aligned
 \overline{\Cal H}_{\mu,2r_{n}}^{q,t}(F)
&\le
 \frac{4^{t}}{\Omega_{d}}
 \,
 \sum_{i\in I_{n}}
 \frac{1}{r_{n}^{d}}
 \,
 \int\limits_{B(x_{n,i},r_{n})}
 \mu(B(x,r_{n}))^{q}
 \,d\Cal L^{d}(x)
 \,
 r_{n}^{t}\\ 
 &\le
 c_{0}
 \,
 \sum_{k=1}^{2^{d}}
 \,
 \sum_{i\in I_{n,k}}
 \frac{1}{r_{n}^{d}}
 \,
 \int\limits_{B(x_{n,i},r_{n})}
 \mu(B(x,r_{n}))^{q}
 \,d\Cal L^{d}(x)
 \,
 r_{n}^{t}
  \quad
 \text{for $q<0$,}\\
&{}\\ 
\overline{\Cal H}_{\mu,\frac{1}{2}r_{n}}^{q,t}(F)
&\le
 \frac{2^{d}}{\Omega_{d}}
 \,
 \sum_{i\in I_{n}}
 \frac{1}{r_{n}^{d}}
 \,
 \int\limits_{B(x_{n,i},r_{n})}
 \mu(B(x,r_{n}))^{q}
 \,d\Cal L^{d}(x)
 \,
 r_{n}^{t}\\ 
 &\le
 c_{0}
 \,
 \sum_{k=1}^{2^{d}}
 \,
 \sum_{i\in I_{n,k}}
 \frac{1}{r_{n}^{d}}
 \,
 \int\limits_{B(x_{n,i},r_{n})}
 \mu(B(x,r_{n}))^{q}
 \,d\Cal L^{d}(x)
 \,
 r_{n}^{t}
  \quad
 \text{for $0\le q$,}
 \endaligned
 \tag10.4
 $$
where
$c_{0}
=
\max(\frac{4^{t}}{\Omega_{d}}, \frac{2^{d}}{\Omega_{d}})$
Next, using the fact
that the sets
$(\,B(x_{n,i},r_{n})\,)_{i\in I_{n,k}}$
 are pairwise disjoint
 and $x_{n,i}\in F\subseteq E$, we conclude from (10.4) that:
 $$  
 \align
 \left.
 \matrix
 \text{for $q<0$, we have}
 \qquad
 \overline{\Cal H}_{\mu,2r_{n}}^{q,t}(F)\\
 \text{\,for $0\le q$, we have}
 \qquad
 \overline{\Cal H}_{\mu,\frac{1}{2}r_{n}}^{q,t}(F)\\
 \endmatrix
 \right\}
 &\le
c_{0}
 \,
 \sum_{k=1}^{2^{d}}
 \,
 \sum_{i\in I_{n,k}}
 \frac{1}{r_{n}^{d}}
 \,
 \int\limits_{B(x_{n,i},r_{n})}
 \mu(B(x,r_{n}))^{q}
 \,d\Cal L^{d}(x)
 \,
 r_{n}^{t}\\
&=
c_{0}
 \,
 \sum_{k=1}^{2^{d}}
 \,
 \frac{1}{r_{n}^{d}}
 \,
 \int\limits_{\cup_{i\in I_{n,k}}B(x_{n,i},r_{n})}
 \mu(B(x,r_{n}))^{q}
 \,d\Cal L^{d}(x)
 \,
 r_{n}^{t}\\
&\le
 c_{0}
 \,
 \sum_{k=1}^{2^{d}}
 \,
 \frac{1}{r_{n}^{d}}
 \,
 \int\limits_{B(E,r_{n})}
 \mu(B(x,r_{n}))^{q}
 \,d\Cal L^{d}(x)
 \,
 r_{n}^{t}\\
&=
c_{0}
 \,
 \sum_{k=1}^{2^{d}}
 \,
 V_{\mu,r_{n}}^{q}(E)
 \,
 r_{n}^{t}\,.
 \tag10.6
 \endalign
 $$
Finally, combining (10.1) and (10.5) shows that
 $$  
 \left.
 \matrix
 \text{for $q<0$, we have}
 \qquad
 \overline{\Cal H}_{\mu,2r_{n}}^{q,t}(F)\\
 \text{\,for $0\le q$, we have}
 \qquad
 \overline{\Cal H}_{\mu,\frac{1}{2}r_{n}}^{q,t}(F)\\
 \endmatrix
 \right\}
 \le
 c
 \,
 r_{n}^{t-(\underline\dim_{\Min,\mu}^{q}(E)+\varepsilon)}
$$
where $c=2^{d}c_{0}$.
\hfill$\square$

\bigskip

\proclaim{Lemma 10.2}
Let $\mu$ be a probability measure on $\Bbb R^{d}$.
Fix $q,t\in\Bbb R$ and $E\subseteq\Bbb R^{d}$.
Then
 $$
 \dim_{\Haus,\mu}^{q}(E)
 \le
 \underline\dim_{\Min,\mu}^{q}(E)
 \le
 \overline\dim_{\Min,\mu}^{q}(E)\,.
 $$
\endproclaim
\noindent{\it  Proof}\newline
\noindent
It clearly suffices to prove that 
$ \dim_{\Haus,\mu}^{q}(E)
 \le
 \underline\dim_{\Min,\mu}^{q}(E)$.
Let $\varepsilon>0$
and write
$t=\underline\dim_{\Min,\mu}^{q}(E)+\varepsilon$.
It follows 
from Lemma 10.1
that there
is a
constant $c>0$ and a sequence $(r_{n})_{n}$ of positive real numbers with
$r_{n}\to0$ satisfying the following:
for all $F\subseteq E$ and $n$, 
the inequality below is satisfied, namely,
 $$
 \align
 \left.
 \matrix
 \text{for $q<0$, we have}
 \qquad
 \overline{\Cal H}_{\mu,2r_{n}}^{q,t}(F)\\
 \text{\,for $0\le q$, we have}
 \qquad
 \overline{\Cal H}_{\mu,\frac{1}{2}r_{n}}^{q,t}(F)\\
 \endmatrix
 \right\}
&\le
 c
 \,
 r_{n}^{t-(\underline\dim_{\Min,\mu}^{q}(E)+\varepsilon)}\\
&=
 c\,r_{n}^{0}\\
&=
 c\,.
 \endalign
  $$

Hence for all $F\subseteq E$, we have
 $$
 \align
 \overline{\Cal H}_{\mu}^{q,t}(F)
&=
\cases
 \limsup_{n}
 \overline{\Cal H}_{\mu,2r_{n}}^{q,t}(F)
&\quad
 \text{for $q<0$;}\\
\limsup_{n}
 \overline{\Cal H}_{\mu,\frac{1}{2}r_{n}}^{q,t}(F)
&\quad
 \text{for $0\le q$}\\ 
 \endcases\\
&\le 
 c\,. 
 \endalign
 $$
Since, $F\subseteq E$ was arbitrary, this implies that
 $$
   \Cal H_{\mu}^{q,t}(E)
  =
  \sup_{F\subseteq E}
   \overline{\Cal H}_{\mu}^{q,t}(F) 
 \le
 c\,,
 $$
whence
$\dim_{\Haus,\mu}^{q}(E)
 \le
 t
 =
 \underline\dim_{\Min,\mu}^{q}(E)+\varepsilon$.
Finally, letting $\varepsilon\searrow0$ gives the desired result.
\hfill$\square$

\bigskip

\noindent
We can now prove Theorem 10.3.

\bigskip

\proclaim{Theorem 10.3}
Fix $q\in\Bbb R$
and
assume that the OSC is satisfied.

For
$\bold i,\bold j\in\Sigma^{*}$ with
$|\bold i|=|\bold j|$ and 
$\bold i\not=\bold j$, 
we have
 $$
 \align
 \dim_{\Haus,\mu}^{q}(S_{\bold i}K\cap S_{\bold j}K)
 \le
 \underline\dim_{\Min,\mu}^{q}(S_{\bold i}K\cap S_{\bold j}K)
&\le
 \overline\dim_{\Min,\mu}^{q}(S_{\bold i}K\cap S_{\bold j}K)\\
&\le
 \gamma(q)\\
&<
 \beta(q)\\
&=
 \dim_{\Haus,\mu}^{q}(K)\,.
 \endalign
 $$
\endproclaim
\noindent{\it  Proof}\newline
\noindent
For a positive integer $m$ and $\bold i,\bold j\in\Sigma^{*}$
with
$|\bold i|=|\bold j|=m$, let
$r_{m}$ be the number defined in (7.5).
Also, let
$c_{\bold i,\bold j,m}$ be the constant in Proposition 7.8
and 
let
$c_{m}$ be the constant in Proposition 8.2.

Now, fix $\bold i,\bold j\in\Sigma^{*}$
with
$|\bold i|=|\bold j|$
and write $m$ for the common value of
$|\bold i|$ and $|\bold j|$, i\.e\. write
$|\bold i|=|\bold j|=m$.
It clearly suffices to prove that
$\overline\dim_{\Min,\mu}^{q}(S_{\bold i}K\cap S_{\bold j}K)
 \le
 \gamma(q)$.
In order to prove 
$\overline\dim_{\Min,\mu}^{q}(S_{\bold i}K\cap S_{\bold j}K)
 \le
 \gamma(q)$, we note that
 for $r<r_{m}$, we have
 $$
 \align
 V_{\mu,r}^{q}
&(S_{\bold i}K\cap S_{\bold j}K)\\
&=
 \frac{1}{r^{d}}
 \int\limits_{
 B(\,S_{\bold i}K\cap S_{\bold j}K\,,r\,)
 } 
 \mu(B(x,r))^{q}
 \,
 d\Cal L^{d}(x)\\
&\le
 \frac{1}{r^{d}}
 \int\limits_{
 B(S_{\bold i}K,r)\cap B(S_{\bold j}K,r)
 } 
 \mu(B(x,r))^{q}
 \,
 d\Cal L^{d}(x) \\
&=
 Q_{\bold i,\bold j}^{q}(r)\\
&\le
 \cases
 {\dsize
  c_{\bold i,\bold j,m}
 \,
 Z_{m}^{q}(\tfrac{1}{2}r)
 }
&\quad\qquad
 \text{for $q<0$;}
 \\ 
  {\dsize
 c_{\bold i,\bold j,m}
 \,
 Z_{m}^{q}(2r)
 }
&\quad\qquad
 \text{for $0\le q$}
\endcases 
 \qquad\qquad
 \text{[by Proposition 7.8]}\\
&\le
 \cases
 {\dsize
 c_{\bold i,\bold j,m}
 \,
 c_{m}
 \,
 (\tfrac{1}{2}r)^{-\gamma(q)}
 }
&\quad
 \text{for $q<0$;}
 \\ 
  {\dsize
 c_{\bold i,\bold j,m}
 \,
 c_{m}
 \,
 (2r)^{-\gamma(q)}
 }
&\quad
 \text{for $0\le q$}
\endcases 
 \qquad\quad
 \,\,\,\,\,
 \text{[by Proposition 8.2]}\\
&\le
 c\,r^{-\gamma(q)} 
 \tag10.7
 \endalign 
$$
where
$c=
\max((\frac{1}{2})^{-\gamma(q)},2^{-\gamma(q)})
\,
c_{m}
\,
c_{\bold i,\bold j,m}$.
It follows immediately from inequality (10.7) that
$\overline\dim_{\Min,\mu}^{q}(S_{\bold i}K\cap S_{\bold j}K)
\le
 \gamma(q)$.
 \hfill$\square$

\bigskip

Next, we
 prove Proposition 10.5
 providing an explicit 
 formula for the multifractal Hausdorff measure 
$ \Cal H_{\mu}^{q,\beta(q)}(S_{\bold i}K)$ of $S_{\bold i}K$.
We begin with some definitions and an auxiliary lemma. 
Let $q\in\Bbb R$ and $E\subseteq\Bbb R^{d}$.
For a locally finite measure $\nu$ on $\Bbb R^{d}$ and a bi-measurable map
$T:\Bbb R^{d}\to\Bbb R^{d}$,
the lower and upper $q$-th order Jacobian of $T$ on $E$ with respect 
to $\nu$ are defined by
 $$
 \underline J_{\nu}^q(T,E)
 =
 \liminf_{r\searrow0}
 \inf_{x\in E}
 \left(
    \frac
	   {
	    \nu TB(x,r)
	   }
	   {
	    \nu B(x,r)
	   }
 \right)^q\,,
 $$
and
 $$ 
 \overline J_{\nu}^q(T,E)
 =
 \limsup_{r\searrow0}
 \sup_{x\in E}
 \left(
    \frac
	   {
	    \nu TB(x,r)
	   }
	   {
	    \nu B(x,r)
	   }
 \right)^q\,,
 $$	  	   
respectively.
If $\underline J_{\nu}^q(T,E)$ and $\overline J_{\nu}^q(T,E)$
coincide, we write 
$J_{\nu}^q(T,E)$ for the common value and call it the
$q$-th order Jacobian of $T$ on $E$ with respect 
to $\nu$.
The main importance of the Jacobians for our purpose is that
they determine the scaling behaviour  of
$\Cal H_{\nu}^{q,t}$ and $\Cal P_{\nu}^{q,t}$.
This is stated formally in the next lemma.

\bigskip		   
		   
\proclaim{Lemma 10.4} 
Let $\nu$ be a probability measure on $\Bbb R^d$
and let $T:\Bbb R^d\to\Bbb R^d$ be a similarity map, i.e.
there exists a constant $c\in(0,\infty)$ such that
$|Tx-Ty|=c|x-y|$ for all
$x,y\in\Bbb R^d$. Assume that
$T(\supp\nu)\subseteq\supp\nu$.
Let $q,t\in\Bbb R$ and $E\subseteq\supp\nu$. Then
 $$
 \underline J_{\nu}^q(T,E)c^t
             \Cal H_{\nu}^{q,t}(E)
		     \le
		     \Cal H_{\nu}^{q,t}(TE)
		     \le
		     \overline J_{\nu}^q(T,E)c^t
             \Cal H_{\nu}^{q,t}(E)\,.
             $$
\endproclaim
\noindent{\it  Proof}\newline
\noindent Follows easily from the definitions. See also [Ol1, Lemma 4.3].
\hfill$\square$

\bigskip

\proclaim{Proposition 10.5}
Fix $q\in\Bbb R$
and
assume 
that one of the following conditions is satisfied:
\roster
\item"(i)"
The OSC is satisfied and $0\le q$;
\item"(ii)"
The SSC is satisfied.
\endroster
For $\bold i\in\Sigma^{*}$
and
$\delta>0$ write
 $$
 \Delta_{\bold i,\delta}
 =
 \bigcup
   \Sb
   |\bold j|=|\bold i|\\
   \bold j\not=\bold i
   \endSb
 S_{\bold i}^{-1}B(S_{\bold j}K,\delta)\,.
 $$ 
\roster
\item"(1)" 
For $\bold i,\bold j\in\Sigma^{*}$ with
$|\bold i|=|\bold j|$ and 
$\bold i\not=\bold j$, we have
 $$
 \Cal H_{\mu}^{q,\beta(q)}(S_{\bold i}K\cap S_{\bold j}K)
 =
 0\,.
 $$
\item"(2)"
For $\bold i,\bold j\in\Sigma^{*}$ with
$|\bold i|=|\bold j|$ and 
$\bold i\not=\bold j$, we have
 $$
 \Cal H_{\mu}^{q,\beta(q)}(S_{\bold i}^{-1}S_{\bold j}K)
  =
 0\,.
 $$
\item"(3)"
For  $\bold i\in\Sigma^{*}$, we have
  $$
 \Cal H_{\mu}^{q,\beta(q)}(\Delta_{\bold i,\delta})
 \to 0
 \quad
 \text{as $\delta\searrow 0$.}
 $$
\item"(4)"
For $\bold i\in\Sigma^{*}$, we have
 $$
 \Cal H_{\mu}^{q,\beta(q)}(S_{\bold i}\Delta_{\bold i,\delta})
 \to 0
 \quad
 \text{as $\delta\searrow 0$.}
 $$
\item"(5)"
For $\bold i\in\Sigma^{*}$ and $\delta>0$, we have
 $$
 \underline J_{\mu}^{q}(S_{\bold i},K\setminus \Delta_{\bold i,\delta})
 =
 \overline J_{\mu}^{q}(S_{\bold i},K\setminus \Delta_{\bold i,\delta})
 =
 p_{\bold i}^{q}
 \,.
 $$
\item"(6)"
For $q\ge 0$ and $\bold i\in\Sigma^{*}$, we have
 $$
 \Cal H_{\mu}^{q,\beta(q)}(S_{\bold i}K)
 =
 r_{\bold i}^{\beta(q)}
 p_{\bold i}^{q}
 \,
 \Cal H_{\mu}^{q,\beta(q)}(K)\,.
 $$
\endroster
\endproclaim
\noindent{\it  Proof}\newline
\noindent
(1) It follows from (1) that
 $\dim_{\Haus,\mu}^{q}(S_{\bold i}K\cap S_{\bold j}K)
 <
 \beta(q)$,
 whence
$\Cal H_{\mu}^{q,\beta(q)}(S_{\bold i}K\cap S_{\bold j}K)
 =
 0$.

\noindent
(2)
We divide the proof into two cases.

\bigskip

\noindent
{\it Case 1: The OSC is satisfied and $0\le q$.}
It follows from Lemma 10.4 that
 $$
 \align
 \Cal H_{\mu}^{q,\beta(q)}(S_{\bold i}^{-1}S_{\bold j}K)
&=
 \Cal H_{\mu}^{q,\beta(q)}(S_{\bold i}^{-1}(S_{\bold i}K\cap S_{\bold j}K))\\
&\le
 r_{i}^{-\beta(q)}\overline J_{\mu}^{q}(S_{\bold i}^{-1},K)
 \,
 \Cal H_{\mu}^{q,t}(S_{\bold i}K\cap S_{\bold j}K)\,.
 \endalign
 $$
Since $\Cal H_{\mu}^{q,\beta(q)}(S_{\bold i}K\cap S_{\bold j}K)=0$
(by Theorem 10.3), it
therefore 
suffices to show that
$\overline J_{\mu}^{q}(S_{\bold i}^{-1},K)<\infty$. 
Below we prove that
 $\overline J_{\mu}^{q}(S_{\bold i}^{-1},K)<\infty$.
Indeed, 
(using the fact that $q\ge 0$) we have
 $$
 \align
 \overline J_{\mu}^{q}(S_{\bold i}^{-1},K)
&=
 \limsup_{r\searrow 0}\,
 \sup_{x\in K}
 \left(
 \frac{\mu(S_{\bold i}^{-1}B(x,r))}{\mu B(x,r)}
 \right)^{q}\\
&=
 \limsup_{r\searrow 0}\,
 \sup_{x\in K}
 \left(
 \frac
 {p_{\bold i}\mu(S_{\bold i}^{-1}B(x,r))}
 {\sum_{|\bold j|=|\bold i|}p_{\bold j}\mu(S_{\bold j}^{-1}B(x,r))}
 \right)^{q}
 p_{\bold i}^{-q}\\
&\le
 p_{\bold i}^{-q}\\
&<
 \infty\,.
 \endalign
 $$

\bigskip

\noindent
{\it Case 2: The SSC is satisfied.}
Indeed, if the SSC is satisfied, then
$S_{\bold i}^{-1}S_{\bold j}K$
for all
$\bold i,\bold j\in\Sigma^{*}$ with
$|\bold i|=|\bold j|$ and 
$\bold i\not=\bold j$,
and the result therefore follows immediately.

\noindent (3)
Since 
$(\,\Delta_{\bold i,\delta}\,)_{\delta>0}$
is a decreasing family of sets
(i\.e\. if $0<\delta\le \rho$, then
$\Delta_{\bold i,\delta}
\subseteq
\Delta_{\bold i,\rho}$)
with
$\cap_{\rho>0}\Delta_{\bold i,\rho}
=
\cup_{|\bold j|=|\bold i|,\bold j\not=\bold i}
S_{\bold i}^{-1}S_{\bold j}K$
(because 
$\cup_{|\bold j|=|\bold i|,\bold j\not=\bold i}
S_{\bold i}^{-1}S_{\bold j}K$
is closed),
it follows from (2) that
$$
 \align
 \Cal H_{\mu}^{q,\beta(q)}(\Delta_{\bold i,\delta})
&\to
 \Cal H_{\mu}^{q,\beta(q)}(\cap_{\rho>0}\Delta_{\bold i,\rho})\\
&=
 \Cal H_{\mu}^{q,\beta(q)}
 (\cup_{|\bold j|=|\bold i|,\bold j\not=\bold i}S_{\bold i}^{-1}S_{\bold j}K)\\
&=
 0\,.
 \endalign
 $$

\noindent (4)
Since 
$(\,S_{\bold i}\Delta_{\bold i,\delta}\,)_{\delta>0}$
is a decreasing family of sets
(i\.e\. if $0<\delta\le \rho$, then
$S_{\bold i}\Delta_{\bold i,\delta}
\subseteq
S_{\bold i}\Delta_{\bold i,\rho}$)
with
$\cap_{\rho>0}S_{\bold i}\Delta_{\bold i,\rho}
=
\cap_{\rho>0}
\big(
\cup_{|\bold j|=|\bold i|,\bold j\not=\bold i}
\big(S_{\bold i}K\cap B(S_{\bold j}K,\rho)\big)
\big)
=
\cup_{|\bold j|=|\bold i|,\bold j\not=\bold i}
\big(S_{\bold i}K\cap S_{\bold j}K\big)$
(because 
$\cup_{|\bold j|=|\bold i|,\bold j\not=\bold i}
\big(S_{\bold i}K\cap S_{\bold j}K\big)$
is closed),
it follows from (1) that
$$
 \align
 \Cal H_{\mu}^{q,\beta(q)}(S_{\bold i}\Delta_{\bold i,\delta})
&\to
 \Cal H_{\mu}^{q,\beta(q)}(\cap_{\rho>0}S_{\bold i}\Delta_{\bold i,\rho})\\
&=
 \Cal H_{\mu}^{q,\beta(q)}
 \big(
 \cup_{|\bold j|=|\bold i|,\bold j\not=\bold i}
\big(S_{\bold i}K\cap S_{\bold j}K\big)
\big)\\
&=
 0\,.
 \endalign
 $$

\noindent (5)
Let $x\in K\setminus \Delta_{\bold i,\delta}$. Then
$S_{\bold i}x\not\in
\cup_{|\bold j|=|\bold i|,\bold j\not=\bold i}B(S_{\bold j}K,\delta)$, whence
$S_{\bold i}B(x,r)\cap K
\subseteq S_{\bold i}K\setminus 
\cup_{|\bold j|=|\bold i|,\bold j\not=\bold i}S_{\bold j}K$
for $0<r<\delta$, and so
 $$
 \frac
 {\mu S_{\bold i}B(x,r)}{\mu B(x,r)}
 =
 \frac
 {\sum_{|\bold j|=|\bold i|}
 \mu(S_{\bold j}^{-1}S_{\bold i}B(x,r))}
 {\mu B(x,r)}
 =
 \frac
 {p_{\bold i}\mu B(x,r)}{\mu B(x,r)}
 =
 p_{\bold i}
 $$
for all $0<r<\delta$.

\noindent
(6) We have
 $$
 \align
 \Cal H_{\mu}^{q,\beta(q)}(S_{\bold i}K)
&=
 \lim_{\delta\searrow 0}
 \Cal H_{\mu}^{q,\beta(q)}(S_{\bold i}(K\setminus \Delta_{\bold i,\delta}))
 \qquad\qquad
 \qquad\qquad
 \quad
 \text{[by (4)]}
\\&=
 \lim_{\delta\searrow 0}
 r_{\bold i}^{\beta(q)}
 J_{\mu}^{q}(S_{\bold i},K\setminus \Delta_{\bold i,\delta})
 \,
 \Cal H_{\mu}^{q,\beta(q)}(K\setminus \Delta_{\bold i,\delta})
 \quad\,\,\,\,
 \text{[by Lemma 10.4 and (5)]}\\
&=
 \lim_{\delta\searrow 0}
 r_{\bold i}^{\beta(q)}
 p_{\bold i}^{q}
 \,
 \Cal H_{\mu}^{q,\beta(q)}(K\setminus \Delta_{\bold i,\delta})
 \qquad\qquad
 \qquad\quad\,\,\,\,
 \text{[by (5)]}\\
&=
 r_{\bold i}^{\beta(q)}
 p_{\bold i}^{q}
 \,
 \Cal H_{\mu}^{q,\beta(q)}(K)\,.
 \qquad\qquad
 \qquad\qquad
 \qquad\quad\,\,\,\,
 \text{[by (3)]} 
 \endalign
 $$
This completes the proof.  
\hfill$\square$

\bigskip
\bigskip


\heading{11. Proof of Theorem 4.1}\endheading

\bigskip

The purpose of this section is to prove Theorem 4.1.
We begin with some technical lemmas.
For $\bold i\in\Sigma^{*}$,
we define $h_{\bold i}^{q}:(0,\infty)\to\Bbb R$ by
 $$
 h_{\bold i}^{q}(r)
 =
 \frac{V_{\mu,r_{\bold i}^{-1}r}^{q}(K)} {(r_{\bold i}^{-1}r)^{-\beta(q)}}
 \Bigg/
 \frac{V_{\mu,r}^{q}(K)} {r^{-\beta(q)}}\,.
 $$

\bigskip

\proclaim{Lemma 11.1}
Fix $q\in\Bbb R$
and
assume 
that one of the following conditions is satisfied:
\roster
\item"(i)"
The OSC is satisfied and $0\le q$;
\item"(ii)"
The SSC is satisfied.
\endroster
Let $\bold i\in\Sigma^{*}$.
Then
 $$
 h_{\bold i}^{q}(r)
 \to
 1\quad\text{as $r\searrow 0$.}
 $$
\endproclaim.
\noindent{\it  Proof}\newline
\noindent
This statement follows immediately from Theorem 3.3 in the non-arithmetic case.
In the arithmetic case there exists for each $i=1,\dots,N$ a positive 
integer $n_{i}$ such that
$\log\frac{1}{r_{i}}=u n_{i}$, and hence, if we write  
$\bold i=i_{1}\dots i_{m}$ and
$m_{i}=|\{k\mid i_{k}=i\}|$ for $i=1,\dots,m$, then  
Theorem 3.3 implies that
$h_{\bold i}^{q}(r)
=
\frac
{\pi_{q}(r_{\bold i}^{-1}r)+\varepsilon_{q}(r_{\bold i}^{-1}r)}
{\pi_{q}(r)+\varepsilon_{q}(r)}
=
\frac
{\pi_{q}((\prod_{i}r^{-m_{i}})r)+\varepsilon_{q}(r_{\bold i}^{-1}r)}
{\pi_{q}(r)+\varepsilon_{q}(r)}
=
\frac
{\pi_{q}(e^{u\sum_{i}n_{i}m_{i}}r)+\varepsilon_{q}(r_{\bold i}^{-1}r)}
{\pi_{q}(r)+\varepsilon_{q}(r)}
=
\frac
{\pi_{q}(r)+\varepsilon_{q}(r_{\bold i}^{-1}r)}
{\pi_{q}(r)+\varepsilon_{q}(r)}
\to 1$ as $r\searrow 0$.
\hfill$\square$

\bigskip

\proclaim{Lemma 11.2}
Fix $q\in\Bbb R$
and
assume 
that one of the following conditions is satisfied:
\roster
\item"(i)"
The OSC is satisfied and $0\le q$;
\item"(ii)"
The SSC is satisfied.
\endroster
Then we have:
\roster
\item"(1)"
For $\bold i\in\Sigma^{*}$, we have
 $$
 \frac
 {1}
 {\Cal I_{\mu,r}^{q}(\Bbb R^{d})}
 \Cal I_{\mu,r}^{q}\big(\,B(S_{\bold i}K,r)\,\big)
 \to
 p_{\bold i}^{q}r_{\bold i}^{\beta(q)}\,.
 $$
\item"(2)"
For $\bold i\in\Sigma^{*}$, we have
 $$
 \frac
 {1}
 {\Cal H_{\mu}^{q,\beta(q)}(K)}
 \Cal H_{\mu}^{q,\beta(q)}\big(\,B(S_{\bold i}K,r)\,\big)
 \to
 p_{\bold i}^{q}r_{\bold i}^{\beta(q)}\,.
 $$ 
\endroster
\endproclaim
\noindent{\it  Proof}\newline
\noindent
(1)
Since
$V_{\mu,r_{\bold i}^{-1}r}^{q}(K)
 =
 \Cal I_{\mu,r_{\bold i}^{-1}r}^{q}
 \big(\,B(K,r_{\bold i}^{-1}r)\,\big)$,
it follows from Lemma 6.2 that
 $$
 \align
 \Bigg|
 \Cal I_{\mu,r}^{q}\big(\,B(S_{\bold i}K,r)\,\big)
 -
 p_{\bold i}^{q}
 V_{\mu,r_{\bold i}^{-1}r}^{q}(K)
 \Bigg|
&=
 \Bigg|
 \Cal I_{\mu,r}^{q}\big(\,B(S_{\bold i}K,r)\,\big)
 -
 p_{\bold i}^{q}
 \Cal I_{\mu,r_{\bold i}^{-1}r}^{q}
 \big(\,B(K,r_{\bold i}^{-1}r)\,\big)
 \Bigg|\\
&\le
 \sum
   \Sb
   |\bold j|=|\bold i|\\
   \bold j\not=\bold i
   \endSb
 Q_{\bold i,\bold j}^{q}(r)\,.
 \endalign
 $$
Next,
since
$V_{\mu,r_{\bold i}^{-1}r}^{q}(K)
 =
 \Cal I_{\mu,r}^{q}
 \big(\,\Bbb R^{d}\,\big)$,
this implies that
 $$
 \align
 \Bigg|
 \frac
 {1}
 {\Cal I_{\mu,r}^{q}(\Bbb R^{d})}
 \Cal I_{\mu,r}^{q}\big(\,B(S_{\bold i}K,r)\,\big)
 -
 p_{\bold i}^{q}
 \frac
 {
 V_{\mu,r_{\bold i}^{-1}r}^{q}(K)
 }
 {
 V_{\mu,r}^{q}(K)
 }
 \Bigg|
&\le
 \sum
   \Sb
   |\bold j|=|\bold i|\\
   \bold j\not=\bold i
   \endSb
 \frac
 {
 Q_{\bold i,\bold j}^{q}(r)
 }
 {
 V_{\mu,r}^{q}(K)
 }\,.
 \tag11.1
 \endalign
 $$
However,
$p_{\bold i}^{q}
 \frac
 {
 V_{\mu,r_{\bold i}^{-1}r}^{q}(K)
 }
 {
 V_{\mu,r}^{q}(K)
 }
 =
 p_{\bold i}^{q}r_{\bold i}^{\beta(q)}h_{\bold i}^{q}(r)$,
and it therefore follows from (11.1) that
 $$
 \align
 \Bigg|
 \frac
 {1}
 {\Cal I_{\mu,r}^{q}(\Bbb R^{d})}
 \Cal I_{\mu,r}^{q}\big(\,B(S_{\bold i}K,r)\,\big)
 -
 p_{\bold i}^{q}r_{\bold i}^{\beta(q)}h_{\bold i}^{q}(r)
 \Bigg|
&\le
 \sum
   \Sb
   |\bold j|=|\bold i|\\
   \bold j\not=\bold i
   \endSb
 \frac
 {
 Q_{\bold i,\bold j}^{q}(r)
 }
 {
 V_{\mu,r}^{q}(K)
 }\,.
 \tag11.2
 \endalign
 $$

Write $m=|\bold i|$.
Let
$r_{m}$ be the number defined in (7.5).
Also,
for $\bold j\in\Sigma^{*}$
with $|\bold i|=|\bold j|=m$, 
let
$c_{\bold i,\bold j,m}$ be the constant in Proposition 7.8
and
let 
$c_{m}$ be the constant in Proposition 8.2.
Observe that for $0<r<\frac{1}{2}r_{m}$, we have
 $$
 \align
 Q_{\bold i,\bold j}^{q}(r)
&\le
 \cases
 {\dsize
 c_{\bold i,\bold j,m}
 Z_{m}^{q}(\tfrac{1}{2}r)
 }
&\quad\qquad
  \,
 \text{for $q<0$;}
 \\ 
  {\dsize
 c_{\bold i,\bold j,m}
 Z_{m}^{q}(2r)
 }
&\quad\qquad
  \,
 \text{for $0\le q$}
\endcases 
 \qquad\,
  \text{[by Proposition 7.8]}\\
&{}\\ 
&\le
 \cases
 {\dsize
 c_{\bold i,\bold j,m}
 \,
 c_{m}
 \,
 (\tfrac{1}{2}r)^{-\gamma(q)}
 }
&\quad
 \text{for $q<0$;}
 \\ 
  {\dsize
 c_{\bold i,\bold j,m}
 \,
 c_{m}
 \,
 (2r)^{-\gamma(q)}
 }
&\quad
 \text{for $0\le q$}
\endcases 
 \qquad
 \text{[by Proposition 8.2]}\\
&{}\\ 
&=
 c
 \,
  c_{\bold i,\bold j,m}
 \,
 c_{m}
 \,
 r^{-\gamma(q)} \,.
 \endalign 
$$
where
$c=
\max((\frac{1}{2})^{-\gamma(q)},2^{-\gamma(q)})$.
Hence, (11.2)
implies that
 $$
 \align
&\Bigg|
 \frac
 {1}
 {\Cal I_{\mu,r}^{q}(\Bbb R^{d})}
 \Cal I_{\mu,r}^{q}\big(\,B(S_{\bold i}K,r)\,\big)
 -
 p_{\bold i}^{q}r_{\bold i}^{\beta(q)}
 \Bigg|\\
&\qquad\qquad
 \le
 \Bigg|
 \frac
 {1}
 {\Cal I_{\mu,r}^{q}(\Bbb R^{d})}
 \Cal I_{\mu,r}^{q}\big(\,B(S_{\bold i}K,r)\,\big)
 -
 p_{\bold i}^{q}r_{\bold i}^{\beta(q)}h_{\bold i}^{q}(r)
 \Bigg|
 +
 |
 p_{\bold i}^{q}r_{\bold i}^{\beta(q)}h_{\bold i}^{q}(r)
 -
 p_{\bold i}^{q}r_{\bold i}^{\beta(q)}
 |\\
&\qquad\qquad
 \le
 \sum
   \Sb
   |\bold j|=|\bold i|\\
   \bold j\not=\bold i
   \endSb
 \frac
 {
 Q_{\bold i,\bold j}^{q}(r)
 }
 {
 V_{\mu,r}^{q}(K)
 } 
 +
 p_{\bold i}^{q}r_{\bold i}^{\beta(q)}
 |
 h_{\bold i}^{q}(r)
 -
 1
 |\\  
&\qquad\qquad
 \le
 \sum
   \Sb
   |\bold j|=|\bold i|\\
   \bold j\not=\bold i
   \endSb
 \frac
 {
 c
 \,
 c_{\bold i,\bold j,m}
 \,
 c_{m}
 \,
 (2r)^{-\gamma(q)}
 }
 {
 V_{\mu,r}^{q}(K)
 } 
 +
 p_{\bold i}^{q}r_{\bold i}^{\beta(q)}
 |
 h_{\bold i}^{q}(r)
 -
 1
 |\\    
&\qquad\qquad
 =
 \sum
   \Sb
   |\bold j|=|\bold i|\\
   \bold j\not=\bold i
   \endSb
   c
   \,
   c_{\bold i,\bold j,m}
   \,
   c_{m}
   \,
   r^{\beta(q)-\gamma(q)}
 \Bigg(
 \frac
 {
 1
 }
 {
 r^{-\beta(q)}
 } 
 V_{\mu,r}^{q}(K)
 \Bigg)^{-1}
 +
 p_{\bold i}^{q}r_{\bold i}^{\beta(q)}
 |
 h_{\bold i}^{q}(r)
 -
 1
 |\\    
 \tag11.3
 \endalign
 $$
for all $0<r<\frac{1}{2}r_{m}$.
Since
$r^{\beta(q)-\gamma(q)}\to 0$ as $r\searrow 0$ (because
$\beta(q)-\gamma(q)>0$ by (8.2)) and
$ \frac
 {
 1
 }
 {
 r^{-\beta(q)}
 } 
 V_{\mu,r}^{q}(K)$
is bounded away from $0$ for $r$ small enough by Theorem 3.3,
it follows 
from Lemma 11.1 and (11.3)
that
 $$
 \frac
 {1}
 {\Cal I_{\mu,r}^{q}(\Bbb R^{d})}
 \Cal I_{\mu,r}^{q}\big(\,B(S_{\bold i}K,r)\,\big)
 \to
 p_{\bold i}^{q}r_{\bold i}^{\beta(q)}\,.
 $$

\noindent 
(2)
For brevity write
$H
=
 \frac{1}{\Cal H_{\mu}^{q,\beta(q)}(K)}\Cal H_{\mu}^{q,\beta(q)}\floor K$.
Since $\big(\,B(S_{\bold i}K,r)\,\big)_{r>0}$ 
is a decreasing family of sets 
(i\.e\.
if $0<s\le t$, then 
$B(S_{\bold i}K,s)\subseteq B(S_{\bold i}K,t)$)
with
$\cap_{r>0}B(S_{\bold i}K,r)=S_{\bold i}K$
(because $S_{\bold i}K$ is closed)
and
$H(\,B(S_{\bold i}K,r)\,)\le 1<\infty$ for all $r>0$, 
we conclude that
 $$
 \align
 H\big(\,B(S_{\bold i}K,r)\,\big)
&\to
 H\Bigg(\,\bigcap_{r>0}B(S_{\bold i}K,r)\,\Bigg)\\
&=
 H(S_{\bold i}K)\,.
 \tag11.4
 \endalign 
 $$
Next, it follows from Proposition 10.5 that
$H(S_{\bold i}K)=
r_{\bold i}^{\beta(q)}
p_{\bold i}^{q}$.
Combining this and (11.4) shows that
 $$
 H\big(\,B(S_{\bold i}K,r)\,\big)
 \to 
 p_{\bold i}^{q}r_{\bold i}^{\beta(q)}\,.
 \tag11.5
 $$
This completes the proof.
\hfill$\square$

\bigskip

\proclaim{Lemma 11.3}
Fix $q\in\Bbb R$
and
assume 
that one of the following conditions is satisfied:
\roster
\item"(i)"
The OSC is satisfied and $0\le q$;
\item"(ii)"
The SSC is satisfied.
\endroster
For $m\in\Bbb R$ and $r>0$, let $(E_{\bold i}(r))_{|\bold i|=m}$
be a family of Borel sets such that
$E_{\bold i}(r)\subseteq B(S_{\bold i}K,r)$
for all $\bold i$
and
$\cup_{|\bold i|=m}E_{\bold i}(r)
=
\cup_{|\bold i|=m}B(S_{\bold i}K,r)$.
Then we have:
\roster
\item"(1)"
For $\bold i\in\Sigma^{*}$,
we have
 $$
 \frac
 {1}
 {\Cal I_{\mu,r}^{q}(\Bbb R^{d})}
 \Cal I_{\mu,r}^{q}\big(\,E_{\bold i}(r)\,\big) 
 \to
 p_{\bold i}^{q}r_{\bold i}^{\beta(q)}\,.
 $$
\item"(2)"
For $\bold i\in\Sigma^{*}$,
we have
 $$
 \frac
 {1}
 {\Cal H_{\mu}^{q,\beta(q)}(K)}
 \Cal H_{\mu}^{q,\beta(q)}\big(\,E_{\bold i}(r)\cap K\,\big) 
 \to
 p_{\bold i}^{q}r_{\bold i}^{\beta(q)}\,.
 $$
\endroster
\endproclaim
\noindent{\it  Proof}\newline
(1)
Recall, that we write
$\Cal V_{\mu,r}^{q}
=\frac
 {1}
 {\Cal I_{\mu,r}^{q}(\Bbb R^{d})}
 \Cal I_{\mu,r}^{q}$
for $r>0$, and note that Lemma 11.2 implies that
 $$
 \Cal V_{\mu,r}^{q}\big(\,B(S_{\bold i}(K,r)\,\big) 
 \to
 p_{\bold i}^{q}r_{\bold i}^{\beta(q)}
 \tag11.6
 $$
for all $\bold i$. 
 
Since clearly
$1
\le
\sum_{|\bold i|=m}\Cal V_{\mu,r}^{q}(E_{\bold i}(r))$
(because
$1
=
\Cal V_{\mu,r}^{q}(B(K,r))
\le
\Cal V_{\mu,r}^{q}(\cup_{|\bold i|=m}B(S_{\bold i}K,r))
=
\Cal V_{\mu,r}^{q}(\cup_{|\bold i|=m}E_{\bold i}(r))
\le
\sum_{|\bold i|=m}\Cal V_{\mu,r}^{q}(E_{\bold i}(r))$),
we conclude from (11.6) that
 $$
 \align
 \sum_{|\bold i|=m}
 \Big(
 \Cal V_{\mu,r}^{q}\big(\,B(S_{\bold i}(K,r)\,\big)
 -
 \Cal V_{\mu,r}^{q}\big(\,E_{\bold i}(r)\,\big) 
 \Big)
&=
 \sum_{|\bold i|=m}
 \Cal V_{\mu,r}^{q}\big(\,B(S_{\bold i}(K,r)\,\big)
 -
 \sum_{|\bold i|=m}
 \Cal V_{\mu,r}^{q}\big(\,E_{\bold i}(r)\,\big) \\
&\le
 \sum_{|\bold i|=m}
 \Cal V_{\mu,r}^{q}\big(\,B(S_{\bold i}(K,r)\,\big)
 -
 1\\
&\to
 \sum_{|\bold i|=m}
 p_{\bold i}^{q}r_{\bold i}^{\beta(q)}
 -
 1\\
&=
 0\,.
 \tag11.7
 \endalign
 $$
Next, since
$\Cal V_{\mu,r}^{q}(\,B(S_{\bold i}(K,r)\,)
 -
 \Cal V_{\mu,r}^{q}(\,E_{\bold i}(r)\,) 
 \ge
 0$
for all $\bold i$
(because
$E_{\bold i}(r)
\subseteq
B(S_{\bold i}K,r)$), it follows from (11.7) that
 $$
 \Cal V_{\mu,r}^{q}\big(\,B(S_{\bold i}(K,r)\,\big)
 -
 \Cal V_{\mu,r}^{q}\big(\,E_{\bold i}(r)\,\big) 
 \to
 0
 \tag11.8
 $$
for all $\bold i$.

Finally,
we deduce from (11.6) and (11.8) that
 $$
 \align
 \Cal V_{\mu,r}^{q}\big(\,E_{\bold i}(r)\,\big) 
&=
 \Cal V_{\mu,r}^{q}\big(\,B(S_{\bold i}(K,r)\,\big) 
 -
 \Big(
 \Cal V_{\mu,r}^{q}\big(\,B(S_{\bold i}(K,r)\,\big)
 -
 \Cal V_{\mu,r}^{q}\big(\,E_{\bold i}(r)\,\big) 
 \Big)\\
&\to
 p_{\bold i}^{q}r_{\bold i}^{\beta(q)}
 -
 0\\
&=
 p_{\bold i}^{q}r_{\bold i}^{\beta(q)}
 \endalign
 $$
for all $\bold i$.

\noindent
(2) 
The proof of (2) is similar to the proof of (1) and is therefore 
omitted.
\hfill$\square$

 \bigskip

We can now prove Theorem 4.1.

\bigskip

\noindent{\it  Proof of Theorem 4.1}\newline
\noindent
Fix $q\ge 0$.
Recall, that we write
$\Cal V_{\mu,r}^{q}
=
\frac
 {1}
 {\Cal I_{\mu,r}^{q}(\Bbb R^{d})}
 \Cal I_{\mu,r}^{q}$.
Also, for brevity write
$H
=
\frac
 {1}
 {\Cal H_{\mu}^{q,\beta(q)}(K)}
 \Cal H_{\mu}^{q,\beta(q)}\floor K$.

We must now prove that
 $$
 {}
 \qquad\quad\,\,
 \Cal V_{\mu,r}^{q} 
 \to
 H
 \,\,\,\,
 \text{weakly.}
 $$
We therefore
let $f:\Bbb R^{d}\to\Bbb R$ be a continuous
function with compact support.
We must now prove that
 $$
 \int f\,d\Cal V_{\mu,r}^{q}
 \to
 \int f\,dH\,.
 \tag11.10
 $$

Below we prove (11.10).

Let $\varepsilon>0$.
Since $f$ is uniformly continuous
(because $f$ is continuous with compact support)  there is a real number 
$\delta_{0}>0$ such that if $x,y\in\Bbb R^{d}$ satisfy 
$|x-y|\le\delta_{0}$, then
 $$
 |f(x)-f(y)|
 \le
 \tfrac{\varepsilon}{8}\,.
 \tag11.11
 $$

Next, since
$\max_{|\bold i|=m}\diam S_{\bold i}K
\le
r_{\max}^{m}\diam K
\to 0$ as $m\to\infty$,
there is a positive integer $m_{0}$ such that
$\diam S_{\bold i}K\le\frac{\delta_{0}}{4}$
for all $\bold i$ with $|\bold i|=m_{0}$.

Also, for each positive number $r>0$, we
may clearly choose a family 
$\big(\,E_{\bold i}(r)\,\big)_{|\bold i|=m_{0}}$
of pairwise disjoint Borel sets $E_{\bold i}(r)$
such that
$E_{\bold i}(r)\subseteq B(S_{\bold i}K,r)$
for all $\bold i$
and
$\cup_{|\bold i|=m_{0}}E_{\bold i}(r)
=
\cup_{|\bold i|=m_{0}}B(S_{\bold i}K,r)$.

Finally, fix $\bold i$ with $|\bold i|=m_{0}$.
It follows from Lemma 11.3
that
$H(E_{\bold i}(r))
\to
p_{\bold i}^{q}r_{\bold i}^{\beta(q)}$,
and we can therefore
choose a positive
real number $s_{\bold i}>0$ such that
for $0<r<s_{\bold i}$ we have
 $$
 \aligned
 \Big|
 H(\,E_{\bold i}(r)\,)
 -
 p_{\bold i}^{q}r_{\bold i}^{\beta(q)}
 \Big|
&\le
 \frac{\varepsilon}{4\|f\|_{\infty}N^{m_{0}}}\,,\\
 H(\,E_{\bold i}(r)\,)
&\le
 2p_{\bold i}^{q}r_{\bold i}^{\beta(q)}\,.
 \endaligned
 \tag11.12
 $$
Similarly, 
it also follows from Lemma 11.3
that
$\Cal V_{\mu,r}^{q}(E_{\bold i}(r))
\to
p_{\bold i}^{q}r_{\bold i}^{\beta(q)}$,
and we can therefore
choose a positive
real number $t_{\bold i}>0$ such that
for $0<r<t_{\bold i}$ we have
 $$
 \aligned
 \Big|
 \Cal V_{\mu,r}^{q}(\,E_{\bold i}(r)\,)
 -
 p_{\bold i}^{q}r_{\bold i}^{\beta(q)}
 \Big|
&\le
 \frac{\varepsilon}{4\|f\|_{\infty}N^{m_{0}}}\,,\\
 \Cal V_{\mu,r}^{q}(\,E_{\bold i}(r)\,)
&\le
 2p_{\bold i}^{q}r_{\bold i}^{\beta(q)}\,.
 \endaligned
 \tag11.13
 $$

Let
 $$
 r_{0}
 =
 \min
 \Big(
 \min_{|\bold i|=m_{0}}s_{\bold i}\,,\,
 \min_{|\bold i|=m_{0}}t_{\bold i}\,,\,
 \tfrac{\delta_{0}}{4}
 \Big)\,.\
 $$

We will now prove that if $0<r<r_{0}$, then
 $$ 
 \Bigg|
 \int f\,d\Cal V_{\mu,r}^{q}
 -
 \int f\,dH
 \Bigg|
 \le
 \varepsilon\,.
 $$

Fix $0<r<r_{0}$.
Since the family
$\big(\,E_{\bold i}(r)\,\big)_{|\bold i|=m_{0}}$
consists
of pairwise disjoint Borel sets $E_{\bold i}(r)$
such that
$E_{\bold i}(r)\subseteq B(S_{\bold i}K,r)$
for all $\bold i$
and
$\cup_{|\bold i|=m_{0}}E_{\bold i}(r)
=
\cup_{|\bold i|=m_{0}}B(S_{\bold i}K,r)$, we conclude that
 $$
 \align
 \Bigg|
 \int f\,d\Cal V_{\mu,r}^{q}
 -
 \int f\,dH
 \Bigg|
&=
 \left|
 \sum_{|\bold i|=m_{0}}\int_{E_{\bold i}(r)} f\,d\Cal V_{\mu,r}^{q}
 -
 \sum_{|\bold i|=m_{0}}\int_{E_{\bold i}(r)} f\,dH
 \right|\\
&\le
 \sum_{|\bold i|=m_{0}}
 \left|
 \int_{E_{\bold i}(r)} f\,d\Cal V_{\mu,r}^{q}
 -
 \int_{E_{\bold i}(r)} f\,dH
 \right|\,.
 \tag11.14
 \endalign
 $$
Next, for each $\bold i\in\Sigma^{*}$ with $|\bold i|=m_{0}$,
fix $x_{\bold i}\in E_{\bold i}(r)$.
We now have, using (11.12), (11.13) and (11.14) that
 $$
 \align
 \Bigg|
 \int f\,d\Cal V_{\mu,r}^{q}
 -
 \int f\,dH
 \Bigg|
&\le
 \sum_{|\bold i|=m_{0}}
 \left|
 \int_{E_{\bold i}(r)} f\,d\Cal V_{\mu,r}^{q}
 -
 \int_{E_{\bold i}(r)} f\,dH
 \right|\\
&\le
 \sum_{|\bold i|=m_{0}}
 \Bigg(
 \left|
 \int_{E_{\bold i}(r)} f\,d\Cal V_{\mu,r}^{q}
 -
 \int_{E_{\bold i}(r)} f(x_{\bold i})\,d\Cal V_{\mu,r}^{q}
 \right|\\
&\qquad\qquad
 \qquad\qquad
 +
 \left|
 \int_{E_{\bold i}(r)} f(x_{\bold i})\,d\Cal V_{\mu,r}^{q}
 -
 \int_{E_{\bold i}(r)} f(x_{\bold i})\,dH
 \right|\\
&\qquad\qquad
 \qquad\qquad
 +
 \left|
 \int_{E_{\bold i}(r)} f(x_{\bold i})\,dH
 -
 \int_{E_{\bold i}(r)} f\,dH
 \right|
 \Bigg)\\
&\le
 \sum_{|\bold i|=m_{0}}
 \Bigg(
 \int_{E_{\bold i}(r)} |f-f(x_{\bold i})|\,d\Cal V_{\mu,r}^{q}\\
&\qquad\qquad
 \qquad\qquad
 +
 |f(x_{\bold i})|
 \left|
 \int_{E_{\bold i}(r)} \,d\Cal V_{\mu,r}^{q}
 -
 \int_{E_{\bold i}(r)} \,dH
 \right|\\
&\qquad\qquad
 \qquad\qquad
 +
 \int_{E_{\bold i}(r)} |f(x_{\bold i})-f|\,dH
 \Bigg)\\ 
&\le
 \sum_{|\bold i|=m_{0}}
 \Bigg(
 \sup_{x,y\in E_{\bold i}(r)}
 |f(x)-f(y)|
 \,
 \Cal V_{\mu,r}^{q}(\,{E_{\bold i}(r)}\,)\\
&\qquad\qquad
 \qquad\qquad
 +
 \|f\|_{\infty}
 \Big|\Cal V_{\mu,r}^{q}(\,{E_{\bold i}(r)}\,)-H(\,{E_{\bold i}(r)}\,)\Big|\\
&\qquad\qquad
 \qquad\qquad
 +
 \sup_{x,y\in E_{\bold i}(r)}
 |f(x)-f(y)|
 \,
 H(\,{E_{\bold i}(r)}\,)
 \Bigg)\\ 
&\le
 \sum_{|\bold i|=m_{0}}
 \Bigg(
 2
 \sup_{x,y\in E_{\bold i}(r)}
 |f(x)-f(y)|
 \,
 2p_{\bold i}^{q}r_{\bold i}^{\beta(q)}\\
&\qquad\qquad
 \qquad\qquad
 +
 \|f\|_{\infty}
 \Big|\Cal V_{\mu,r}^{q}(\,{E_{\bold i}(r)}\,)-H(\,{E_{\bold i}(r)}\,)\Big|
 \Bigg)\\  
&\le
 4
 \sum_{|\bold i|=m_{0}}
 \sup_{x,y\in B(S_{\bold i}K,r)}
 |f(x)-f(y)|
 \,
 p_{\bold i}^{q}r_{\bold i}^{\beta(q)}\\
&\qquad\qquad
 +
 \|f\|_{\infty}
 \sum_{|\bold i|=m_{0}}
 \Bigg(
 \Big|\Cal V_{\mu,r}^{q}(\,{E_{\bold i}(r)}\,)-p_{\bold i}^{q}r_{\bold i}^{\beta(q)}\Big|
 +
 \Big|p_{\bold i}^{q}r_{\bold i}^{\beta(q)}-H(\,{E_{\bold i}(r)}\,)\Big|
 \Bigg)\\  
&\le
 4
 \sum_{|\bold i|=m_{0}}
 \sup_{x,y\in B(S_{\bold i}K,r)}
 |f(x)-f(y)|
 \,
 p_{\bold i}^{q}r_{\bold i}^{\beta(q)}\\
&\qquad\qquad
   \qquad\qquad
 +
 \|f\|_{\infty}
 \sum_{|\bold i|=m_{0}}
 \Bigg(
 \frac{\varepsilon}{4\|f\|_{\infty}N^{m_{0}}}
 +
 \frac{\varepsilon}{4\|f\|_{\infty}N^{m_{0}}}
 \Bigg)\\  
\allowdisplaybreak
&\le
 4
 \sum_{|\bold i|=m_{0}}
 \sup_{x,y\in B(S_{\bold i}K,r)}
 |f(x)-f(y)|
 \,
 p_{\bold i}^{q}r_{\bold i}^{\beta(q)}\\
&\qquad\qquad
   \qquad\qquad
 +
 \frac{\varepsilon}{2N^{m_{0}}}
 \sum_{|\bold i|=m_{0}}
 1\\ 
\allowdisplaybreak
&\le
 4
 \sum_{|\bold i|=m_{0}}
 \sup_{x,y\in B(S_{\bold i}K,r)}
 |f(x)-f(y)|
 \,
 p_{\bold i}^{q}r_{\bold i}^{\beta(q)}
 \,\,
 +
 \,\,
 \frac{\varepsilon}{2}\,.
 \tag11.15
 \endalign
 $$
However, for each $\bold i\in\Sigma^{*}$ with
$|\bold i|=m_{0}$ and
and for all $x,y\in B(S_{\bold i}K,r)$, we have
$|x-y|
\le
2(\diam S_{\bold i}K+r)
\le
2(\frac{\delta_{0}}{4}+r_{0})
\le
2(\frac{\delta_{0}}{2}+\frac{\delta_{0}}{2})
=
\delta_{0}$,
and it therefore follows from (11.11) that
$|f(x)-f(y)|\le\frac{\varepsilon}{8}$.
This and (11.15) imply that
 $$
 \align
 \Bigg|
 \int f\,d\Cal V_{\mu,r}^{q}
 -
 \int f\,dH
 \Bigg|
&\le
 4\,\frac{\varepsilon}{8}
 \sum_{|\bold i|=m_{0}}
 p_{\bold i}^{q}r_{\bold i}^{\beta(q)}
 +
 \frac{\varepsilon}{2}\\
&=
 \frac{\varepsilon}{2}
 +
 \frac{\varepsilon}{2}\\
&=
 \varepsilon\,.
 \endalign
 $$
This completes the proof.
\hfill$\square$

\bigskip
\bigskip


\heading{12. Proof of Theorem 4.2}\endheading

\bigskip

The purpose of this section is to prove Theorem 4.2. We begin we a small lemma.

\bigskip

\proclaim{Lemma 12.1}
Let
$\varphi,\Phi:(0,1)\to(0,\infty)$
be measurable functions such that
$\int_{r}^{1}\Phi(s)\,\frac{ds}{s}<\infty$ for all $r$
and
$\int_{r}^{1}\varphi(s)\Phi(s)\,ds<\infty$
for all $r$.
Let $c,C\ge 0$.
Assume that
 $$
 \align
 \frac{1}{-\log r}
 \int_{r}^{1}\Phi(s)\,\frac{ds}{s}
&\to 
 C
 \,\,\,\,
 \text{as $r\searrow 0$,}\\
 \varphi(r)
&\to 
 c
 \,\,\,\,\,\,
 \text{as $r\searrow 0$.}
 \endalign
 $$
Then the following holds.
 \roster
 \item"(1)"
We have 
 $$
 \frac{1}{-\log r}
 \int_{r}^{1}\varphi(s)\Phi(s)\,\frac{ds}{s}
\to 
 c\,C
 \,\,\,\,
 \text{as $r\searrow 0$.}
 $$
 \item"(2)"
We have 
 $$
 {}\quad\,\,\,\,\,
 \frac{1}{-\log r}
 \int_{r}^{1}\varphi(s)\,\frac{ds}{s}
\to 
 c
 \,\,\,\,\,\,\,\,\,
 \text{as $r\searrow 0$.}
 $$
\endroster 
\endproclaim 
\noindent{\it  Proof}\newline
\noindent
(1)
This follows by a standard argument and the proof is therefore omitted.

\noindent
(2) This follows from (1) by putting $\Phi=1$ and noticing
that
$\frac{1}{-\log r}
 \int_{r}^{1}\,\frac{ds}{s}  
 =
 1$ for all $r>0$.
 This completes the proof.
\hfill$\square$

\bigskip

We can now prove Theorem 4.2.

\bigskip

\noindent{\it  Proof of Theorem 4.2}\newline
\noindent
Recall, that
the $(q,\beta(q))$ multifractal Minkowski content 
$M_{\mu}^{q,\beta(q)}(K)$
and 
the $(q,\beta(q))$ average  multifractal Minkowski content 
$M_{\mu,\ave}^{q,\beta(q)}(K)$
are defined by
 $$
 \align
 M_{\mu}^{q,\beta(q)}(K)
&=
 \lim_{r\searrow 0}
 \frac{1}{r^{-\beta(q)}}
 V_{\mu,r}^{q}(K)\,,\\
 M_{\mu,\ave}^{q,\beta(q)}(K)
&=
 \lim_{r\searrow 0}
 \frac{1}{-\log r}
 \int_{r}^{1}
 \frac{1}{s^{-\beta(q)}}
 V_{\mu,s}^{q}(K)
 \,
 \frac{ds}{s}\,,
 \endalign
 $$
provided the limits exist.
Also, recall (from Theorem 3.2 or Theorem 3.3) that 
the multifractal Minkowski dimensions are given by
$\underline\dim_{\Min,\mu}^{q}(K)
 =
 \overline\dim_{\Min,\mu}^{q}(K)
 =
 \beta(q)$.
In particular, this implies that the 
the tube
measures
$\underline{\Cal S}_{\mu,r}^{q}$
and
$\overline{\Cal S}_{\mu,r}^{q}$
conincide,
and that
the average tube measures
$\underline{\Cal S}_{\mu,r,\ave}^{q}$
an
$\overline{\Cal S}_{\mu,r,\ave}^{q}$
coincide.
Below we write
$\Cal S_{\mu,r}^{q}$ for the 
common value of 
$\underline{\Cal S}_{\mu,r}^{q}$
and
$\overline{\Cal S}_{\mu,r}^{q}$,
and we write
$\Cal S_{\mu,r,\ave}^{q}$ for the 
common value of 
$\underline{\Cal S}_{\mu,r,\ave}^{q}$
an
$\overline{\Cal S}_{\mu,r,\ave}^{q}$, i\.e\. we write
 $$
 \align
 \Cal S_{\mu,r}^{q}
&=
 \frac{1}{r^{-\beta(q)}}
 \Cal I_{\mu,r}^{q}\,,\\
 \Cal S_{\mu,r,\ave}^{q,\beta(q)}
&=
 \frac{1}{-\log r}
 \int_{r}^{1}
 \frac{1}{s^{-\beta(q)}}
 \Cal I_{\mu,s}^{q}
 \,
 \frac{ds}{s}\,.
 \endalign
 $$
In particular, we note that
 $$
 \Cal S_{\mu,r,\ave}^{q}
 =
 \frac{1}{-\log r}
 \int_{r}^{1}
  \Cal S_{\mu,s}^{q}
 \,
 \frac{ds}{s}\,.
 $$
Finally, notice that if
$f:\Bbb R^{d}\to\Bbb R$ is  a continuous
function with compact support, then
a standard argument shows that
 $$
 \int
 f\,d\Cal S_{\mu,r,\ave}^{q}
 =
 \frac{1}{-\log r}
 \int_{r}^{1}
 \Bigg(
  \int f\,d\Cal S_{\mu,s}^{q}
 \Bigg)
 \,
 \frac{ds}{s}\,.
 \tag12.1
 $$
Recall, that we write
$\Cal V_{\mu,r}^{q}
=
\frac
 {1}
 {\Cal I_{\mu,r}^{q}(\Bbb R^{d})}
 \Cal I_{\mu,r}^{q}$.
 Also, for brevity write
$H
=
\frac
 {1}
 {\Cal H_{\mu}^{q,\beta(q)}(K)}
 \Cal H_{\mu}^{q,\beta(q)}\floor K$.
We can now prove the statements in Theorem 4.2.

\noindent
(1)
It follows from Theorem 3.2
that 
$\underline\dim_{\Min,\mu}^{q}(K)
 =
 \overline\dim_{\Min,\mu}^{q}(K)
 =
 \beta(q)$
and this clearly implies the desired statement.

\noindent
(2)
Let
$f:\Bbb R^{d}\to\Bbb R$ is  a continuous
function with compact support.
Since clearly
$\Cal S_{\mu,r}^{q}
=
\frac{1}{r^{-\beta(q)}}
\,
V_{\mu,r}^{q}(K)
\,
\Cal V_{\mu,r}^{q}$,
it now follows from Theorem 3.3 and Theorem 4.1 that
 $$
 \align
 \int f\,d\Cal S_{\mu,r}^{q}
&=
 \frac{1}{r^{-\beta(q)}}
 \,
 V_{\mu,r}^{q}(K)
 \,
 \int f\,d\Cal V_{\mu,r}^{q}\\
&\to
 M_{\mu}^{q,\beta(q)}(K)
 \int f\,dH\,.
 \tag12.2
 \endalign
 $$
We also deduce from 
Lemma 12.1.(2)
(applied to the function
$\varphi:(0,\infty)\to(0,\infty)$ defined by
$\varphi(r)= \int f\,d\Cal S_{\mu,r}^{q}$),
(12.1) and (12.2)
that
 $$
 \align
 \int f\,d\Cal S_{\mu,r,\ave}^{q}
&=
 \frac{1}{-\log r}
 \int_{r}^{1}
 \Bigg(
  \int f\,d\Cal S_{\mu,s}^{q}
 \Bigg)
 \,
 \frac{ds}{s} \\
&\to
 M_{\mu}^{q,\beta(q)}(K)
 \int f\,dH\\
&=
\,
 M_{\mu,\ave}^{q,\beta(q)}(K)
 \int f\,dH\,. 
\endalign 
 $$
This completes the proof.

\noindent
(3)
Let
$f:\Bbb R^{d}\to\Bbb R$ is  a continuous
function with compact support.
function.
As above, since
$\Cal S_{\mu,r}^{q}
=
\frac{1}{r^{-\beta(q)}}
\,
V_{\mu,r}^{q}(K)
\,
\Cal V_{\mu,r}^{q}$,
it now follows from 
Lemma 12.1.(1)
(applied to
the functions $\varphi,\Phi:(0,1)\to(0,\infty)$ defined by
$\varphi(r)
=
\int f\,d\Cal V_{\mu,r}^{q}$
and
$\Phi(r)
=
\frac{1}{r^{-\beta(q)}}\,V_{\mu,r}^{q}(K)$),
Theorem 3.3, Theorem 4.1
and (12.1) that
 $$
 \align
 \int
 f\,d\Cal S_{\mu,r,\ave}^{q}
&=
\frac{1}{-\log r}
 \int_{r}^{1}
 \Bigg(
  \int f\,d\Cal S_{\mu,s}^{q}
 \Bigg)
 \,
 \frac{ds}{s} \\
&=
 \frac{1}{-\log r}
 \int_{r}^{1}
 \Bigg(
\frac{1}{r^{-\beta(q)}}
\,
V_{\mu,r}^{q}(K)
 \int f\,d\Cal V_{\mu,s}^{q}
 \Bigg)
 \,
 \frac{ds}{s}\\
&\to
 M_{\mu,\ave}^{q,\beta(q)}(K)
 \int f\,dH\,.
 \endalign
 $$
This completes the proof. 
\hfill$\square$


\newpage

${}$

\bigskip
\bigskip
\bigskip
\bigskip
\bigskip
\bigskip
\bigskip
\bigskip

\centerline{\bigletter Part 4:}

\bigskip

\centerline{\bigletter Proofs of the Results from Section 5}


\bigskip
\bigskip
\bigskip
\bigskip

\heading{13. 
Analysis of the poles of $\zeta_{\mu}^{q}$}
\endheading

In this section we establish various technical 
growth
estimates
related to the
poles and
residues of 
the zeta function $\zeta_{\mu}^{q}$.
These estimates will play important parts in the proofs of Theorem 5.4, Theorem 5.5
and Theorem 5.7.
We begin by defining the number $\alpha(q)$;
this number plays an important part 
in describing
the location of the poles of $\zeta_{\mu}^{q}$.
Fix $q\in\Bbb R$
and define
$\alpha(q)$ by
 $$
 \alpha(q)
 =
 \inf
 \left\{
 t\in\Bbb R
 \,
 \left|
 \,
 \sum
 \Sb
  i\\
  r_{i}=r_{\min}
 \endSb 
 p_{i}^{q}r_{\min}^{t}
 \le
  1
 +
 \sum
 \Sb
  i\\
  r_{i}>r_{\min}
 \endSb 
 p_{i}^{q}r_{i}^{t}
 \right.
 \right\}\,.
 \tag13.1
 $$
Also, recall that $\beta(q)$ is defined by
 $$
 \sum_{i}p_{i}^{q}r_{i}^{\beta(q)}
=1\,.
\tag13.2
$$  
  Using the numbers $\alpha(q)$ and $\beta(q)$ we can now 
  describe the location of the poles of $\zeta_{\mu}^{q}$. 
  We first prove the statements in Proposition 5.2;
  however, for the benefit of the reader, before proving Proposition 5.2 we
  repeat the statements of the proposition.
Recall, that
if
$f$
is 
a meromorphic function, 
then $Z(f)$ denotes the set of zeros of $f$
and that
$P(f)$ denotes the set of poles of $f$.

\bigskip

\proclaim{Proposition 13.1 (i\.e\. statements (1)--(4.2) in Proposition 5.2).
The poles of $\zeta_{\mu}^{q}$}
Fix $q\in\Bbb R$.
\roster
\item"(1)"
We have
$-\infty<\alpha(q)\le\beta(q)<\infty$.
\item"(2)"
We have
 $$
 P(\zeta_{\mu}^{q})
 =
 Z
 \Big(
 s\to 1-\sum_{i}p_{i}^{q}r_{i}^{s}
 \Big)\,.
 \qquad\qquad
 \quad\,\,
  $$
\item"(3)"
We have
   $$
 P(\zeta_{\mu}^{q})
\subseteq
 \Big\{
 s\in\Bbb C
 \,\Big|\,
 \alpha(q)
 \le
 \real(s)
 \le
 \beta(q)
 \Big\}\,.
 \,\,\,\,\,\,\,\,\,
 $$
\item"(4.1)"
The poles $\omega$ with
$\real(\omega)=\beta(q)$
in the non-arithmetic case:
If 
the set
$\{\log r_{1}^{-1},\ldots,\log r_{N}^{-1}\}$
is not contained in a discrete additive subgroup of $\Bbb R$,
then
 $$
  P(\zeta_{\mu}^{q})
\cap \Big\{
 s\in\Bbb C
 \,\Big|\,
 \real(s)
 =
 \beta(q)
 \Big\}
 =
 \big\{
 \beta(q)
 \big\}\,.
 $$
\item"(4.2)"
The poles $\omega$ with
$\real(\omega)=\beta(q)$
in the arithmetic case:
If 
the set
$\{\log r_{1}^{-1},\ldots,\log r_{N}^{-1}\}$
is contained in a discrete additive subgroup of $\Bbb R$
and
$\langle\log r_{1}^{-1},\ldots,\log r_{N}^{-1}\rangle=u\Bbb Z$
with $u>0$,
then
 $$
 {}\qquad\,\,\,\,\,\,\,\,\,
  P(\zeta_{\mu}^{q})
\cap 
\Big\{
 s\in\Bbb C
 \,\Big|\,
 \real(s)
 =
 \beta(q)
 \Big\}
 =
 \,\,\,\,\beta(q)\,+\tfrac{2\pi}{u}\ima\Bbb Z\,,
 $$
and for each $i$, there is a unique integer $k_{i}$ such that
$\log r_{i}^{-1}=k_{i}u$ 
and, in addition,
 $$
 \align
 {}\qquad\qquad
   \qquad\qquad
   \qquad
 P(\zeta_{\mu}^{q})
&=
 \Big(\,\beta(q)+\tfrac{2\pi}{u}\ima\Bbb Z\,\Big)\\
&\qquad\qquad
 \,\,
 \cup
 \,\,
 \bigcup
 \Sb
 w\in Z(z\to 1-\sum_{i}p_{i}^{q}z^{k_{i}})\\
 {}\\
 w\not=e^{-u\beta(q)}
 \endSb
 \Big(
 \,
 -\tfrac{\log |w|}{u}
 -
 \tfrac{\Arg w}{u}
 \ima
 +
 \tfrac{2\pi}{u}\ima \Bbb Z
 \,
 \Big)
 \endalign
 $$
(where $\Arg z$ denotes the principal argument of $z\in\Bbb C$).

 \endroster
\endproclaim 
\noindent{\it  Proof}\newline
\noindent
(1)
It clearly suffices to prove that
$\alpha(q)>-\infty$ and $\alpha(q)\le\beta(q)$.
We first prove  that  $\alpha(q)>-\infty$. 
Indeed,
if we define 
$\Phi_{q}:\Bbb R\to\Bbb R$ by
$ 
\Phi_{q}(t)
=
 1
 +
 \sum_{i\,,\,r_{i}>r_{\min}}
 p_{i}^{q}r_{i}^{t}
-
 \sum_{i\,,\,r_{i}=r_{\min}}
 p_{i}^{q}r_{\min}^{t}$,
 then
 $\lim_{t\to-\infty}\Phi_{q}(t)=-\infty$
 and
  $\lim_{t\to\infty}\Phi_{q}(t)=1$,
  whence
  $\alpha(q)=\inf\{t\in\Bbb R\,|\,\Phi_{q}(t)\ge 0\}\in\Bbb R$.
  Next, we show that
  $\alpha(q)\le \beta(q)$.
  Indeed, it follows from the definition of $\beta(q)$ that
  $\sum_{i\,,\,r_{i}=r_{\min}}p_{i}^{q}r_{\min}^{\beta(q)}
  =
  1
  -
  \sum_{i\,,\,r_{i}>r_{\min}}p_{i}^{q}r_{\min}^{\beta(q)}
  \le
  1
  +
  \sum_{i\,,\,r_{i}>r_{\min}}p_{i}^{q}r_{\min}^{\beta(q)}$, and so
  $\beta(q)\ge\alpha(q)$.

\noindent
(2)
Since
$\zeta_{\mu}^{q}(s)
=
\frac{\sum_{i}p_{i}^{q}r_{i}^{s}}{1-\sum_{i}p_{i}^{q}r_{i}^{s}}$
for
$s\in\Bbb C\setminus
Z(s\to 1-\sum_{i}p_{i}^{q}r_{i}^{s})$,
 it follows immediately that
$P(\zeta_{\mu}^{q})
 =
 Z(
 s\to 1-\sum_{i}p_{i}^{q}r_{i}^{s})$.
 
 \noindent
 (3)
 Let $s=\sigma+\ima t\in\Bbb C$ with $\sigma,t\in\Bbb R$ be such that
$\sum_{i}p_{i}^{q}r_{i}^{s}
=1$.
We must now prove that
$\alpha(q)
 \le
 \sigma
 \le
 \beta(q)$.
 
We first prove that $\sigma\le\beta(q)$.
Indeed,
since
$1
=
|\sum_{i}p_{i}^{q}r_{i}^{s}|
\le
\sum_{i}p_{i}^{q}|r_{i}^{s}|
=
\sum_{i}p_{i}^{q}r_{i}^{\sigma}$,
we conclude immediately from the 
definition of $\beta(q)$ that
$\sigma\le \beta(q)$.

Next, we prove that $\alpha(q)\le\sigma$.
To prove this inequality  we note that
$\sum_{r_{i}=r_{\min}}p_{i}^{q}r_{i}^{s}
+
\sum_{r_{i}>r_{\min}}p_{i}^{q}r_{i}^{s}
=
\sum_{i}p_{i}^{q}r_{i}^{s}
=
1$,
from which we see that
$\sum_{r_{i}=r_{\min}}p_{i}^{q}r_{\min}^{\sigma}
=
\sum_{r_{i}=r_{\min}}p_{i}^{q}\,|r_{\min}^{s}|
=
|\sum_{r_{i}=r_{\min}}p_{i}^{q}r_{i}^{s}|
=
|1-\sum_{r_{i}>r_{\min}}p_{i}^{q}r_{i}^{s}|
\le
1+\sum_{r_{i}>r_{\min}}p_{i}^{q}\,|r_{i}^{s}|
=
1+\sum_{r_{i}>r_{\min}}p_{i}^{q}r_{i}^{\sigma}$.
It follows immediately from this and the definition of $\alpha(q)$ that
$\alpha(q)\le\sigma$.

\noindent
(4.1)
We must prove that
$P(\zeta_{\mu}^{q})
\cap 
\{s\in\Bbb C\,|\,
 \real(s)
 =
 \beta(q)\}
 =
 \{\beta(q)\}$.
We first note that
$\beta(q)
\in
Z(
 s\to 1-\sum_{i}p_{i}^{q}r_{i}^{s})
\cap 
\{
 s\in\Bbb C
 \,|\,
 \real(s)
 =
 \beta(q)
 \}
=
  P(\zeta_{\mu}^{q})
\cap 
\{
 s\in\Bbb C
 \,|\,
 \real(s)
 =
 \beta(q)
 \}$.
 Next,
 we prove that
 if 
 $\omega
 \in
  P(\zeta_{\mu}^{q})
\cap 
\{
 s\in\Bbb C
 \,|\,
 \real(s)
 =
 \beta(q)
 \}$, then
 $\omega=\beta(q)$.
 We therefore fix
  $\omega
 \in
  P(\zeta_{\mu}^{q})
\cap 
\{
 s\in\Bbb C
 \,|\,
 \real(s)
 =
 \beta(q)
 \}$.
 It follows that there is $t\in\Bbb R$
 such that
 $\omega=\beta(q)+\ima t$.
We must now show that $t=0$.
Observe that since
 $\beta(q)+\ima t
 =
 \omega
 \in
 P(\zeta_{\mu}^{q})=
 Z(
 s\to 1-\sum_{i}p_{i}^{q}r_{i}^{s})$,
 we have
 $1
 =\sum_{i}p_{i}^{q}r_{i}^{\beta(q)+\ima t}
 =\sum_{i}p_{i}^{q}r_{i}^{\beta(q)}e^{\ima t\log r_{i}}$.
 We therefore deduce that
  $1
  =|1|
 =|\sum_{i}p_{i}^{q}r_{i}^{\beta(q)}e^{\ima t\log r_{i}}|
 \le \sum_{i}|p_{i}^{q}r_{i}^{\beta(q)}e^{\ima t\log r_{i}}|
=\sum_{i}p_{i}^{q}r_{i}^{\beta(q)}
=1$, and so
 $$
 \Bigg|\sum_{i}p_{i}^{q}r_{i}^{\beta(q)}e^{\ima t\log r_{i}}\Bigg|
 =
  \sum_{i}|p_{i}^{q}r_{i}^{\beta(q)}e^{\ima t\log r_{i}}|\,.
  \tag13.3
 $$
 We 
 conclude from (13.3) that
 the ($2$-dimensional planar) vectors 
 given by
 $p_{1}^{q}r_{1}^{\beta(q)}e^{\ima t\log r_{1}},
 \allowmathbreak
 \ldots,
 \allowmathbreak
 p_{N}^{q}r_{N}^{\beta(q)}e^{\ima t\log r_{N}}$
 must be positive multiples
 of a common ($2$-dimensional planar) unit vector,
 i\.e\. there is a 
 ($2$-dimensional planar)
 unit vector $e^{\ima\theta}$ with $\theta\in[-\pi,\pi)$
 and positive numbers $\lambda_{1},\ldots,\lambda_{N}\ge 0$ such that
 $p_{i}^{q}r_{i}^{\beta(q)}e^{\ima t\log r_{i}}=\lambda_{i}e^{\ima\theta}$
 for all $i$.
Since $p_{i}^{q}r_{i}^{\beta(q)}>0$, this implies that
$p_{i}^{q}r_{i}^{\beta(q)}=\lambda_{i}$,
and consequently
 $e^{\ima t\log r_{i}}=e^{\ima \theta}$.
 It follows from this 
 that
 $t\log r_{i}-\theta\in2\pi\Bbb Z$ and we can therefore
  find integers $m_{i}$ such that
  $$
  t\log r_{1}=\theta+2\pi m_{1}\,,\ldots,\,
  t\log r_{N}=\theta+2\pi m_{N} \,.
  \tag13.4
  $$
  Next,
 since 
$1
 =\sum_{i}p_{i}^{q}r_{i}^{\beta(q)}e^{\ima t\log r_{i}}$, we
 conclude from (13.4) that
 $1
 =
 \sum_{i}p_{i}^{q}r_{i}^{\beta(q)}e^{\ima t\log r_{i}}
 =
 \sum_{i}p_{i}^{q}r_{i}^{\beta(q)}e^{\ima(\theta+2\pi m_{i})}
 =
 (\sum_{i}p_{i}^{q}r_{i}^{\beta(q)})e^{\ima\theta}
 =
 e^{\ima\theta}$, whence $\theta=0$ (because $\theta\in[-\pi,\pi)$).
This and another application of (13.4) shows that
  $$
  t\log r_{1}=2\pi m_{1}\,,\ldots,\,
  t\log r_{N}=2\pi m_{N} \,.
  \tag13.5
  $$
Finally, if $t\not=0$, then it follows from (13.5) that
 $\log r_{i}^{-1}=-\frac{2\pi}{t} m_{i}$ for all $i$, and
 the set 
 $\{\log r_{1}^{-1},\ldots,\log r_{N}^{-1}\}$ is therefore contained in the 
 discrete additive subgroup $-\frac{2\pi}{t} \Bbb Z$.
However, this contracts the fact that 
 $\{\log r_{1}^{-1},\ldots,\log r_{N}^{-1}\}$ is not contained
 in any
 discrete additive subgroup of $\Bbb R$.
 Consequently, we conclude that $t=0$.

\noindent
(4.2)
Since $\log r_{i}^{-1}=k_{i}u$, we deduce that if $s\in \Bbb C$, then
$1-\sum_{i}p_{i}^{q}r_{i}^{s}
=
1-\sum_{i}p_{i}^{q}(e^{-us})^{k_{i}}$, whence
 $$
 \align
  P(\zeta_{\mu}^{q})
&=
 Z
 \Big(
 s\to 1-\sum_{i}p_{i}^{q}r_{i}^{s}
 \Big)\\
&=
 Z
 \Big(
 s\to 1-\sum_{i}p_{i}^{q}(e^{-us})^{k_{i}}
 \Big)\\ 
&=
\bigcup
 \Sb
 w\in Z(z\to 1-\sum_{i}p_{i}^{q}z^{k_{i}})\\
 \endSb
\Big\{
s\in\Bbb C\,\Big|\,w=e^{-us}
\Big\}\\
&=
 \Big(\,\beta(q)+\tfrac{2\pi}{u}\ima\Bbb Z\,\Big)
 \,\,
 \cup
 \,\,
 \bigcup
 \Sb
 w\in Z(z\to 1-\sum_{i}p_{i}^{q}z^{k_{i}})\\
 {}\\
 w\not=e^{-u\beta(q)}
 \endSb
 \Big(
 \,
 -\tfrac{\log |w|}{u}
 -
 \tfrac{\Arg w}{u}
 \ima
 +
 \tfrac{2\pi}{u}\ima \Bbb Z
 \,
 \Big)\,.
\endalign
$$
This completes the proof.
\hfill$\square$

\bigskip

\noindent
The next two propositions, i\.e\.
Proposition 13.2 and Proposition 13.3,
contain detailed information about the poles of $\zeta_{\mu}^{q}$
near the \lq\lq critical line" $\real(s)=\beta(q)$.

\bigskip

\proclaim{Proposition 13.2. The poles of $\zeta_{\mu}^{q}$ near the 
\lq\lq critical line" $\real(s)=\beta(q)$} 
 Fix $q\in\Bbb R$.
Then there is a number $b(q)\in\Bbb R$
with the following properties:

\roster
\item"(1)"
We have
$b(q)<\beta(q)$.

\item
"(2)"
If:
 $$
 \text{
 $\omega$ is a pole of $\zeta_{\mu}^{q}$
 with
 }
 \,\,
 \omega
 \in
 \Big\{
  s\in\Bbb C\,\Big|\,b(q)\le \real(s)\le \beta(q)
  \Big\}\,,
  $$

\noindent
then:   

\smallskip

${}$\,${}$\,${}$
 (i)
 $\omega$ is a simple pole of $\zeta_{\mu}^{q}$;
 
 \smallskip
 
${}$
(ii)
$
\res(\zeta_{\mu}^{q};\omega)
=
\frac{1}{-\sum_{i}p_{i}^{q}r_{i}^{\omega}\log r_{i}}$;

(iii) 
$|\res(\zeta_{\mu}^{q};\omega)|
 \le
 -\tfrac{1}{\log r_{\max}}$.

\endroster  
\endproclaim
\noindent{\it  Proof}\newline
\noindent
Choose $i_{\star}$ such that 
  $$
 r_{i_{\star}}
 =
 \cases
 r_{\max}
&\quad
 \text{if $\beta(q)\le 0$;}\\
 r_{\min}
&\quad
 \text{if $0<\beta(q)$,}\\ 
 \endcases
 $$
and define $b_{\star}(q)\in\Bbb R$ by
 $$
 \sum_{i\not=i_{\star}}
 p_{i}^{q}r_{i}^{b_{\star}(q)}
 =
 1\,.
 $$
Observe that
$\sum_{i\not=i_{\star}}
 p_{i}^{q}r_{i}^{b_{\star}(q)}
 =
 1
 =
 \sum_{i}
 p_{i}^{q}r_{i}^{\beta(q)}
 >
 \sum_{i\not=i_{\star}}
 p_{i}^{q}r_{i}^{\beta(q)}$, whence
$b_{\star}(q)
 <
 \beta(q)$. 
Next, since
$b_{\star}(q)<\beta(q)$ 
we can choose a real number $b(q)$
with
 $$
  b_{\star}(q)
<
b(q)
 <
 \beta(q)
 $$  
such that
$b(q)<\beta(q)\le 0$ if $\beta(q)\le 0$
and 
$0\le b(q)<\beta(q)$  if $0<\beta(q)$.

We must now prove that if 
$\omega$ is a pole of $\zeta_{\mu}^{q}$
 with
$
 \omega
 \in
 \{
  s\in\Bbb C\,|\,b(q)\le \real(s)\le \beta(q)
 \}$,
then    
$\omega$ is a simple pole of $\zeta_{\mu}^{q}$
with
$
\res(\zeta_{\mu}^{q};\omega)
=
\frac{1}{-\sum_{i}p_{i}^{q}r_{i}^{\omega}\log r_{i}}$
and
$|\res(\zeta_{\mu}^{q};\omega)|
 \le
 -\tfrac{1}{\log r_{\max}}$.  
We therefore fix
a pole 
$\omega$  of $\zeta_{\mu}^{q}$
 with
$
 \omega
 \in
 \{
  s\in\Bbb C\,|\,b(q)\le \real(s)\le \beta(q)
 \}$.
Define $f,g:\Bbb C\to\Bbb C$  
by
$g(s)=\sum_{i}p_{i}^{q}r_{i}^{s}$
and
$f(s)=1-\sum_{i}p_{i}^{q}r_{i}^{s}$,
and note that
$\zeta_{\mu}^{q}
 =
 \frac{g}{f}$.
Also, note that since $\omega$ is pole of  $\zeta_{\mu}^{q}$, we conclude that
$f(\omega)=0$ and, consequently, $g(\omega)=1-f(\omega)=1$.
We now prove the following three claims.

\medskip

\noindent
{\it Claim 1.
$\real(p_{i}^{q}r_{i}^{\omega})\ge 0$ for all $i$.}

\noindent
{\it Proof of Claim 1.}
Indeed, otherwise there an index $j$ such that
$\real(p_{j}^{q}r_{j}^{\omega})< 0$, whence
(using the fact that
$g(\omega)=1$)
 $$
 \align
 1
&=
g(\omega)\\
&=
\real(g(\omega)\\
&= 
 \real\Bigg(\sum_{i}p_{i}^{q}r_{i}^{\omega}\Bigg)\\
&=
 \sum_{i}\real(p_{i}^{q}r_{i}^{\omega})\\
&<
 \sum_{i\not=j}\real(p_{i}^{q}r_{i}^{\omega})\\
&\le
 \sum_{i\not=j}p_{i}^{q}r_{i}^{\sigma}\,.
 \tag13.6
\endalign
$$
Next, note that 
if
$\beta(q)\le 0$, then
$r_{i_{\star}}=r_{\max}\ge r_{j}$
and
$\sigma=\real(\omega)\le\beta(q)\le 0$,
whence
$r_{i_{\star}}^{\sigma}\le r_{j}^{\sigma}$,
and so
$\sum_{i\not=j}p_{i}^{q}r_{i}^{\sigma}
\le
\sum_{i\not=i_{\star}}p_{i}^{q}r_{i}^{\sigma}$.
On the other hand,
if $0<\beta(q)$, then
$r_{i_{\star}}=r_{\min}\le r_{j}$
and
$0\le b(q)\le\real(\omega)=\sigma$,
whence
$r_{i_{\star}}^{\sigma}\le r_{j}^{\sigma}$,
and so
$\sum_{i\not=j}p_{i}^{q}r_{i}^{\sigma}
\le
\sum_{i\not=i_{\star}}p_{i}^{q}r_{i}^{\sigma}$.
Consequently,
we always have
$\sum_{i\not=j}p_{i}^{q}r_{i}^{\sigma}
\le
\sum_{i\not=i_{\star}}p_{i}^{q}r_{i}^{\sigma}$.
We conclude from this and (13.6) that
 $$
 \align
 1
 &<
 \sum_{i\not=j}p_{i}^{q}r_{i}^{\sigma}\\
&\le
 \sum_{i\not=i_{\star}}p_{i}^{q}r_{i}^{\sigma}\,.
 \tag13.7
\endalign
$$
However, since
$\sigma=\real(\omega)\ge b(q)>b_{\star}(q)$, we deduce that
$r_{i}^{\sigma}<r_{i}^{b_{\star}(q)}$.
The definition of $b_{\star}(q)$ and (13.7)
therefore imply that
 $$
 \align
 1
&<
 \sum_{i\not=i_{\star}}p_{i}^{q}r_{i}^{\sigma}\\
&\le 
 \sum_{i\not=i_{\star}}p_{i}^{q}r_{i}^{b_{\star}(q)}\\
&=
1\,.
\tag13.8
\endalign
$$
The desired contradiction follows immediately from (13.8). This completes the proof of Claim 1.

\medskip

\noindent{\it
Claim 2.
We have $g(\omega)=1$. In particular, $g(\omega)\not=0$.}

\noindent{\it
Proof of Claim 2.}
Indeed, this has already been observed above. This completes the proof of Claim 2.

 \medskip

\noindent{\it
Claim 3.
We have $\real(f'(\omega))\ge -\log r_{\max}$. In particular, 
$f(\omega)=0$
and
$f'(\omega)\not=0$.}

\noindent{\it
Proof of Claim 3.}
We have
$f'(\omega)
=
-\sum_{i}p_{i}^{q}r_{i}^{\omega}\log r_{i}$ and so
 $$
 \align
 \real(f'(\omega))
&=
 -\sum_{i}\real(p_{i}^{q}r_{i}^{\omega})\,\log r_{i}\,.
 \tag13.9
\endalign
 $$
Now observe that
it follows from Claim 1 that $\real(p_{i}^{q}r_{i}^{\omega})\ge 0$
for all $i$, and (13.9) 
therefore implies that
 $$
 \align
 \real(f'(\omega))
&=
 -\sum_{i}\real(p_{i}^{q}r_{i}^{\omega})\,\log r_{i}\\
&\ge
 -\sum_{i}\real(p_{i}^{q}r_{i}^{\omega})\,\log r_{\max}\\
&= 
-
 \real
 \Bigg(
 \sum_{i}p_{i}^{q}r_{i}^{\omega}
 \Bigg)
 \,
 \log r_{\max}\\
&= 
-
 \real
 (g(\omega))
 \,
 \log r_{\max}\\
&= 
-
 \real
 (1)
 \,
 \log r_{\max}\\
&= 
-
 \log r_{\max}\,.
\endalign
 $$
This completes the proof of Claim 3.

\medskip

Using the fact that $\zeta_{\mu}^{q}=\frac{g}{f}$,
we conclude from Claim 2 and Claim 3
that
$\omega$ is a simple pole of $\zeta_{\mu}^{q}$
with
$\res(\zeta_{\mu}^{q};\omega)=\frac{g(\omega)}{f'(\omega)}$, 
whence
$
 |\res(\zeta_{\mu}^{q};\omega)|
 =
 \frac{|g(\omega)|}{|f'(\omega)|}
 \le
 \frac{1}{|\real(f'(\omega))|}
 \le
 -
 \frac{1}{\log r_{\max}}$.
This completes the proof of Proposition 13.2.
\hfill$\square$

\bigskip

\noindent
Next, we state and prove Proposition 13.3.
The construction of the contour $\Gamma$ in Proposition 13.3
is motivated by  an argument in
[DeKo\"OzRa\"Ur].

\bigskip

\proclaim{Proposition 13.3.
The poles of $\zeta_{\mu}^{q}$ near the 
\lq\lq critical line"
$\real(s)=\beta(q)$: 
construction of $\Gamma$}
Fix $q\in\Bbb R$.
Assume that
$
 \beta(q)
 \not\in  
  \{
  -dq,1-dq,\ldots,d-dq
  \}$.
Let $b(q)$ be as in Proposition 13.2.  
Then there are two real numbers
$b_{0}(q)$ and $\beta_{0}(q)$,
and two real valued 
sequences
$(\,u_{n}(q)\,)_{n\in\Bbb Z}$
and
$(\,v_{n}(q)\,)_{n\in\Bbb Z}$
satisfying the following
conditions:
\roster
\item"(1)"
We have
$b(q)<b_{0}(q)<\beta_{0}(q)<\beta(q)$.
\item"(2)"
We have
$\{
  -dq,1-dq,\ldots,d-dq
  \}
  \cap
  [b_{0}(q),\beta(q)]
  =
  \varnothing$.

\item"(3)"
We have 
$0\not\in[b_{0}(q),\beta_{0}(q)]$.

\item"(4)"
We have
$u_{-n}(q)=-v_{n}(q)$ and $v_{-n}(q)=-u_{n}(q)$ for all $n$.
In addition,
$$
\ldots 
<
u_{-1}(q)
<
v_{-1}(q)
<
u_{0}(q)
<
0
<
v_{0}(q)
<
u_{1}(q)
<
v_{1}(q)
<
\dots
$$
and
$$
\gather
\,
\lim_{n\to\infty}\,\,u_{n}(q)
=
\,\,
\lim_{n\to\infty}\,\,v_{n}(q)
=
\,\,\,
\infty\,,\,\\
\lim_{n\to-\infty}u_{n}(q)
=
\lim_{n\to-\infty}v_{n}(q)
=
-\infty\,.
\endgather
$$

\item"(5)" We have
$u_{n+1}(q)-u_{n}(q)\ge-\frac{\pi}{\log r_{\min}}$
for all $n$.

\endroster

\noindent
Let
 $$
 \align
&\text{
 $\Pi_{n}^{+}$
 be the 
 directed
 horizontal
 line segment from}\\
&\qquad\qquad
  \qquad\qquad
  \qquad\qquad
  \qquad\qquad
 \text{  
 $b_{0}(q)+\ima v_{n}(q)$
 \quad
 to
 $\beta_{0}(q)+\ima v_{n}(q)$,
 }\\
&\text{
 $\Gamma_{n}^{+}$
 be the
 directed
 \lq\lq left" vertical
 line segment from}\\
&\qquad\qquad
  \qquad\qquad
  \qquad\qquad
  \qquad\qquad
 \text{  
 $b_{0}(q)+\ima u_{n}(q)$
 \quad
 to
 $b_{0}(q)+\ima v_{n}(q)$,
 }\\
&\text{
 $\Pi_{n}^{-}$
 be the
 directed
 horizontal
 line segment from}\\
&\qquad\qquad
  \qquad\qquad
  \qquad\qquad
  \qquad\qquad
 \text{  
 $\beta_{0}(q)+\ima u_{n}(q)$
 \quad
 to
 $b_{0}(q)+\ima u_{n}(q)$,
 }\\
&\text{
 $\Gamma_{n}^{-}$
 be the
 directed
  \lq\lq right" vertical
 line segment from}\\
&\qquad\qquad
  \qquad\qquad
  \qquad\qquad
  \qquad\qquad
 \text{  
 $\beta_{0}(q)+\ima v_{n-1}(q)$
 to
 $\beta_{0}(q)+\ima u_{n}(q)$,
 }\\
 \endalign
 $$
and
 $$
 \text{
 $\Gamma_{n}$ be the concatenation of
 $\Gamma_{n}^{-}$,
 $\Pi_{n}^{-}$,
 $\Gamma_{n}^{+}$
 and
 $\Pi_{n}^{+}$,
 }
 \qquad\qquad
 \qquad\qquad
 \qquad\,\,
 $$
   $$
 \text{
 $\Gamma$\,\,\,\,\, be the 
 concatenation of
  $\ldots,
  \Gamma_{-1},
  \Gamma_{0},
  \Gamma_{1},
 \ldots$.}
 \qquad\qquad
 \qquad\qquad
 \qquad\quad\,
 $$ 
\medskip 
\noindent
Then the following holds:
\roster
\item"(6)"
The path $\Gamma$ does not intersect the sets of poles of $\zeta_{\mu}^{q}$.
\item"(7)"
We have
$\sup_{s\in\Gamma}|\zeta_{\mu}^{q}(s)|<\infty$.
\item"(8)"
If $M$ is a bounded subset of $\Bbb R$, then
$\sup_{n}\sup_{s\in (M+u_{n}(q)\ima)\cup(M+v_{n}(q)\ima)}
|\zeta_{\mu}^{q}(s)|<\infty$.
\endroster

\noindent
Below we sketch the contour $\Gamma$.

\goodbreak

\newpage

\midinsert


 \vspace{130mm}
 \centerline{\hbox{\hskip -50mm\special{pdf=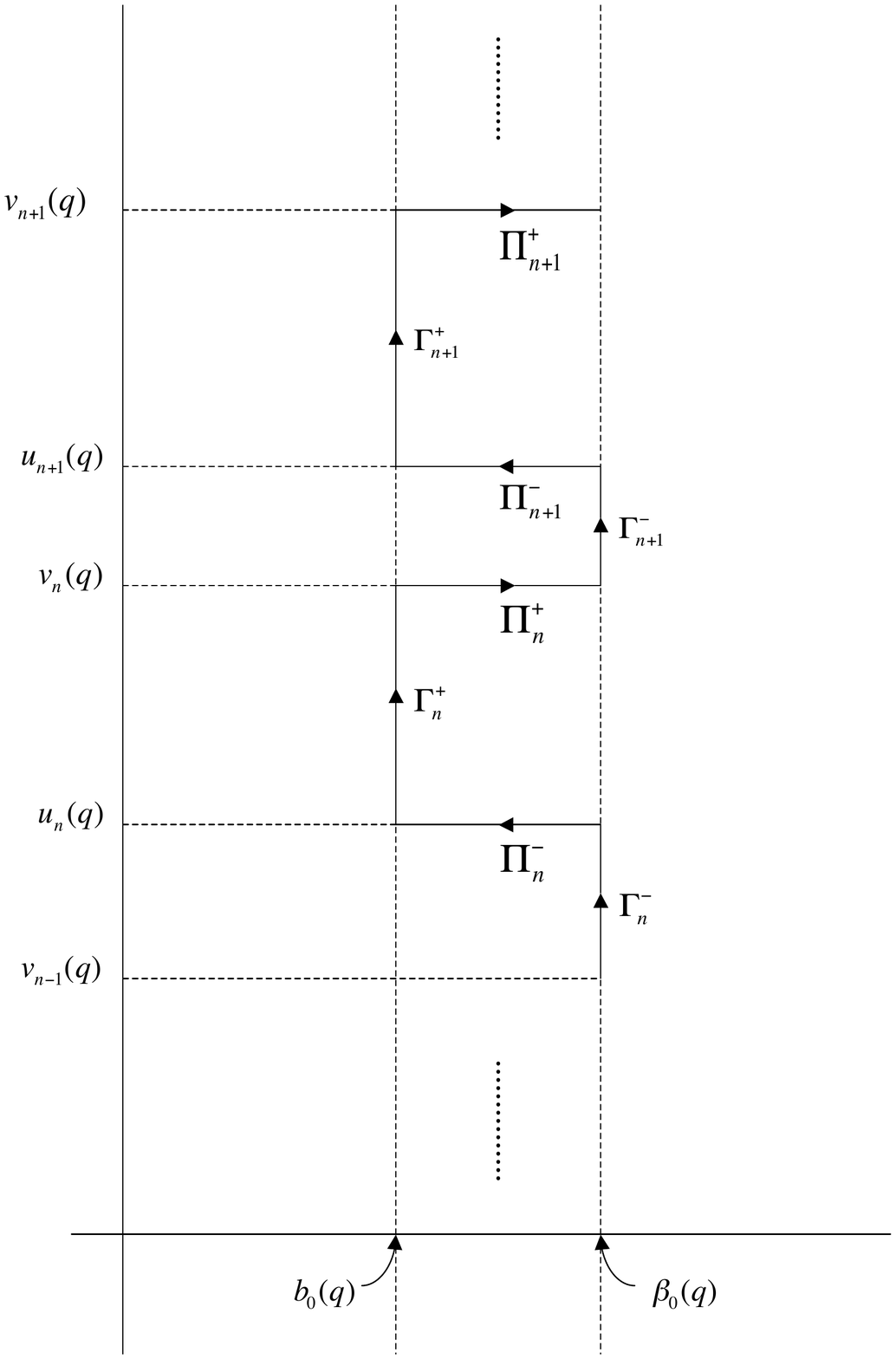
 scale=0.38}}}

\botcaption{\eightpoint \bf Fig\. 13.1}
\eightpoint 
The contour $\Gamma$.
\endcaption
\endinsert

\endproclaim

\noindent{\it  Proof}\newline
\noindent
Since 
$
 \beta(q)
 \not\in  
  \{
  -dq,1-dq,\ldots,d-dq
  \}$,
 we can choose $b_{0}(q)$ such that
$b(q)<b_{0}(q)<\beta(q)$,
$\{
  -dq,1-dq,\ldots,d-dq
  \}
  \cap
  [b_{0}(q),\beta(q)]
  =
  \varnothing$
and  
$0\not\in[b_{0}(q),\beta(q))$.
As $b_{0}(q)<\beta(q)$, we conclude that
$\sum_{i}p_{i}^{q}r_{i}^{b_{0}(q)}>1$. This implies that there is a positive number $d_{0}>0$ such that
 $$
 \sum_{i}p_{i}^{q}r_{i}^{b_{0}(q)}
 =
 1+d_{0}\,.
 \tag13.10
 $$
Next, note that we can choose $\theta_{0}>0$ such that
the following four conditions are satisfied:
$3\theta_{0}\frac{\log r_{\min}}{\log r_{\max}}
 \le
 (
 \frac{d_{0}}{1+d_{0}}
 )^{\frac{1}{2}}$,
$ \frac{\theta_{0}^{2}\min_{i}(p_{i}^{q}r_{i}^{\beta(q)})}{8}
 <
 d_{0}$,
$\frac{4\theta_{0}}{-\log r_{\max}}
 <
 \frac{2\pi-2\theta_{0}}{-\log r_{\min}}$
 and
$\theta_{0}
 \le
 \frac{\pi}{2}$.
Now, let
$\delta_{0}
 =
 \frac{\theta_{0}^{2}\min_{i}(p_{i}^{q}r_{i}^{\beta(q)})}{8}$
 and define $\beta_{0}(q)$ by
  $$
  \sum_{i}p_{i}^{q}r_{i}^{\beta_{0}(q)}
  =
  1+\delta_{0}\,.
  \tag13.11
  $$
The statements in the proposition
follow from the seven claims below.

\bigskip

\noindent{\it
Claim 1.
We have
$b(q)<b_{0}(q)<\beta_{0}(q)<\beta(q)$.}

\noindent
{\it Proof of Claim 1.}
It follows from the definition of $b_{0}(q)$ 
that
$b(q)<b_{0}(q)$,
and
since
$0
<
\delta_{0}
 =
 \frac{\theta_{0}^{2}\min_{i}(p_{i}^{q}r_{i}^{\beta(q)})}{8}
<
d_{0}$, we conclude that
$b_{0}(q)<\beta_{0}(q)<\beta(q)$.
This concludes the proof of Claim 1.

\bigskip

Next, for each $i$, we define
$\Theta_{i}:\Bbb R\to\Bbb R$ by
 $$
 \Theta_{i}(t)
 =
 \Arg(r_{i}^{t})
 =
 \Arg(e^{\ima t\log r_{i}})
 $$
(here we write $\Arg(z)$ for the principal argument of 
a complex number $z$, i\.e\. $\Arg(z)$
is the unique argument of $z$
belonging to the interval $[-\pi,\pi)$). 
Write
 $$
 G_{i}
 =
 \Big\{t\in\Bbb R\,\Big|\,
 |\Theta_{i}(t)|<\theta_{0}
 \Big\}
 $$
and note that
 $$
 G_{i}
 =
 \bigcup_{n\in\Bbb Z}
 B\Big(-n\pi\tfrac{1}{\log r_{i}},\rho_{i}\Big)
 $$
where $\rho_{i}=-\theta_{0}\frac{1}{\log r_{i}}$. 
Also write
 $$
 \align
 G
&=
 \Big\{t\in\Bbb R\,\Big|\,
 \max_{i}|\Theta_{i}(t)|<\theta_{0}
 \Big\}
 \endalign 
 $$ 
and note that
 $$
 \align
 G
&=
 \bigcap_{i}G_{i}\,.
 \endalign 
 $$ 
We now prove the following claim describing the 
structure of the set $G$.
Below we use the following notation, namely, if $M$ is a subset of $\Bbb R$, then
we will write
$-M=\{-x\,|\,x\in M\}$.

\bigskip

\noindent{\it
Claim 2.
\roster
\item"(1)" The set $G$ is open and $0\in G$.
\item"(2)" The set $G$ is unbounded.
\item"(3)" The set $\Bbb R\setminus G$ is unbounded.
\item"(4)" $G=-G$.
\endroster

\noindent
Proof of Claim 2.}
(1) This is clear since $G_{i}$ is open for all $i$ with $0\in G_{i}$.

\noindent
(2) In order to prove this statement we use Dirichlet's theorem on simultaneous
Diophantine approximation.
For the benefit of the reader we will now state this result.

\medskip

\noindent
{\it
Dirichlet's theorem on simultaneous
Diophantine approximation [Te, Lemma 14.1].
Let $n$ be a positive integer and 
let $\xi_{1},\ldots,\xi_{n}\in\Bbb R$.
Fix $Q\in\Bbb N$ and $\tau>0$.
Then there is $t\in[\tau,\tau Q^{n}]$ and 
there are $k_{1},\ldots,k_{n}\in\Bbb Z$ such that
$
|t\xi_{i}-k_{i}|\le\frac{1}{Q}
 $
for all $i$. 
}
 
\medskip

\noindent
We can now prove that $G$ is unbounded.
Let $\tau>0$. We must now show that there is $t\in G$ with $t\ge \tau$.
We first choose $Q\in\Bbb N$ with $\frac{1}{Q}<\frac{\theta_{0}}{2\pi}$.
Next, it follows from
Dirichlet's theorem on simultaneous
Diophantine approximation
(applied to $\xi_{i}=\frac{1}{2\pi}\log r_{i}$ for $i=1,\ldots,N$)
that there is 
$t\in[\tau,\tau Q^{N}]$ and that
there are $k_{1},\ldots,k_{N}\in\Bbb Z$ such that
$|\frac{t}{2\pi}\log r_{i}-k_{i}|\le\frac{1}{Q}<\frac{\theta_{0}}{2\pi}$
for all $i$, i\.e\.
$|t\log r_{i}-2\pi k_{i}|\le\theta_{0}$
for all $i$. 
This clearly implies that
$|\Theta_{i}(t)|
=
|\Arg(e^{\ima t\log r_{i}})|
<
\theta_{0}$ for all $i$, and so
$t\in G_{i}$ for all $i$, whence
$t\in\cap_{i}G_{i}=G$.
Also, since 
$t\in[\tau,\tau Q^{N}]$,
we see that $t\ge \tau$.

\noindent 
(3)
Since the set
$\Bbb R\setminus G_{i}$ is unbounded for all $i$, it follows that
$\Bbb R\setminus G
=
\Bbb R\setminus (\cap_{i} G_{i})
=
\cup_{i}(\Bbb R\setminus G_{i})$
is unbounded.

\noindent
(4)
This follows from the fact that
$G_{i}=-G_{i}$ for all $i$.
This completes the proof of Claim 2.

\bigskip

 \noindent
It follows from Claim 2 that
$0\in G$, $G=-G$
and that
$G$ is an open unbounded set with 
an infinite number of connected components
and that all connected components of $G$ are bounded.
We therefore conclude that
there
are numbers
$\ldots
<
a_{-2}<b_{-2}
<
a_{-1}<b_{-1}
<
a_{0}<0<b_{0}
<
a_{1}<b_{1}
<
a_{2}<b_{2}
<
\ldots$
such that
 $$
 G
 =
 \bigcup_{n}\,\,(a_{n},b_{n})\,.
 $$
and
 $$
 a_{n}\to\infty\,\,\,\,
 \text{as $n\to\infty$}
 \quad
 \text{and}
 \quad
 a_{n}\to-\infty\,\,\,\,
 \text{as $n\to-\infty$.}
 $$ 
In addition, since $G=-G$, we see that
$a_{-n}=-b_{n}$ and $b_{-n}=-a_{n}$ for all $n$.
Now put
 $$
 \align
 u_{n}(q)
&=
 a_{n}-\Big(-2\theta_{0}\tfrac{1}{\log r_{\max}}\Big)\,,\\
  v_{n}(q)
&=
 b_{n}+\Big(-2\theta_{0}\tfrac{1}{\log r_{\max}}\Big)\,.
\endalign
$$

\bigskip

\noindent{
\it
Claim 3.
Let $n$ be an integer.
Then the following statements hold.
\roster
\item"(1)"
$
\ldots 
<
u_{-1}(q)
<
v_{-1}(q)
<
u_{0}(q)
<
0
<
v_{0}(q)
<
u_{1}(q)
<
v_{1}(q)
<
\dots
$
\item"(2)"
$u_{-n}(q)=-v_{n}(q)$ and $v_{-n}(q)=u_{n}(q)$.
\item"(3)"
$u_{n+1}-u_{n}\ge-\frac{\pi}{\log r_{\min}}$.
\endroster
}
\noindent
{\it Proof of Claim 3.}
First note that
 $$a_{n+1}-b_{n}
\ge
\min_{i}\Big(-2\pi\tfrac{1}{\log r_{i}}-2\rho_{i}\Big)
=
\min_{i}\Big(-2\pi\tfrac{1}{\log r_{i}}-\Big(-2\theta_{0}\tfrac{1}{\log r_{i}}\Big)\Big)
=
-2(\pi-\theta_{0})\tfrac{1}{\log r_{\min}}\,.
\tag13.12
$$
(1) It is clear that
$u_{n}(q)<a_{n}<b_{n}<v_{n}(q)$
and since
$-\frac{4\theta_{0}}{\log r_{\max}}
<
-\frac{2\pi-2\theta_{0}}{\log r_{\min}}$,
(13.12) implies that
$-\frac{4\theta_{0}}{\log r_{\max}}
<
-\frac{2\pi-2\theta_{0}}{\log r_{\min}}
\le
a_{n+1}-b_{n}$, whence
$v_{n}(q)
=
 b_{n}+(-2\theta_{0}\frac{1}{\log r_{\max}})
 <
 a_{n+1}-(-2\theta_{0}\frac{1}{\log r_{\max}})
 =
 u_{n+1}(q)$.

\noindent
(2)
Since $a_{-n}=-b_{n}$ and $b_{-n}=-a_{n}$,
this follows immediately from the definitions
of $u_{n}(q)$ and $v_{n}(q)$.

\noindent
(3)
Since $\theta_{0}\le\frac{\pi}{2}$, we conclude from (13.12) that
$u_{n+1}(q)-u_{n}(q)
=
a_{n+1}-a_{n}
\ge
a_{n+1}-b_{n}
\ge
-\frac{2\pi-2\theta_{0}}{\log r_{\min}}
\ge
-\frac{\pi}{\log r_{\min}}$.
This completes the proof of Claim 3.

\bigskip

\noindent{
\it
Claim 4.
Let $i=1,\ldots,N$ and $n$ be an integer.
Then the following statements hold.
\roster
\item"(1)"
$\,
-(\frac{d_{0}}{1+d_{0}})^{\frac{1}{2}}
\le
\,\,\,\,\,\,
\Theta_{i}\big(t\big)
\,\,\,\,\,\,
\le
(\frac{d_{0}}{1+d_{0}})^{\frac{1}{2}}$
for $t\in[u_{n}(q),v_{n}(q)]$.
\item"(2)"
$\qquad\quad\,\,
\theta_{0}
\le
\Theta_{i}\big(u_{n}(q)\big)
\le
(\frac{d_{0}}{1+d_{0}})^{\frac{1}{2}}$.
\item"(3)"
$\,
-(\frac{d_{0}}{1+d_{0}})^{\frac{1}{2}}
\le
\Theta_{i}\big(v_{n}(q)\big)
\le
-\theta_{0}$.
\endroster
}
\noindent
{\it Proof of Claim 4.}
We first observe that if
$t=s+\varepsilon$ with
$s\in[a_{n},b_{n}]$ and $\varepsilon\in\Bbb R$,
then the following holds.
Namely,
since
$s\in[a_{n},b_{n}]$, we conclude that
$|\Arg(e^{\ima s\log r_{i}})|=|\Theta_{i}(s)|<\theta_{0}$ for all $i$.
This implies that for each $i$, we can find an integer $m_{i}\in\Bbb Z$ such that
$-\theta_{0}
 \le
 s\log r_{i}-2\pi m_{i}
 \le
 \theta_{0}$,
 and since $t=s+\varepsilon$, we deduce from this that
 $$
 -\theta_{0}+\varepsilon\log r_{\min}
 \le
 t\log r_{i}-2\pi m_{i}
 \le
 \theta_{0}+\varepsilon\log r_{\min}
 \tag13.13
$$
for all $i$.

\noindent
(1)
Let $t\in[u_{n}(q),v_{n}(q)]$.
It is clear that
$t=s+\varepsilon$
with
$s\in[a_{n},b_{n}]$ and  $|\varepsilon|\le-2\theta_{0}\frac{1}{\log r_{\max}}$.
It follows from this and (13.13) that 
for each $i$,
we can find an integer $m_{i}\in\Bbb Z$ such that
$-\theta_{0}+\varepsilon\log r_{\min}
 \le
 t\log r_{i}-2\pi m_{i}
 \le
 \theta_{0}+\varepsilon\log r_{\min}$, whence
$| t\log r_{i}-2\pi m_{i}|
\le
\theta_{0}-|\varepsilon|\log r_{\min}
\le
\theta_{0}+2\theta_{0}\frac{\log r_{\min}}{\log r_{\max}}
\le
3\theta_{0}\frac{\log r_{\min}}{\log r_{\max}}
\le
(\frac{d_{0}}{1+d_{0}})^\frac{1}{2}$.
This  clearly implies that
$|\Theta_{i}(t)|\le (\frac{d_{0}}{1+d_{0}})^\frac{1}{2}$.

\noindent
(2)
We have $u_{n}(q)=a_{n}+\varepsilon$ where
$\varepsilon=2\theta_{0}\frac{1}{\log r_{\max}}$.
It follows from this and (13.13) that 
for each $i$,
we can find an integer $m_{i}\in\Bbb Z$ such that
$-\theta_{0}+\varepsilon\log r_{\min}
 \le
 u_{n}(q)\log r_{i}-2\pi m_{i}
 \le
 \theta_{0}+\varepsilon\log r_{\min}$.
Since
$-\theta_{0}+\varepsilon\log r_{\min}
=
-\theta_{0}+2\theta_{0}\frac{\log r_{\min}}{\log r_{\max}}
\ge
-\theta_{0}+2\theta_{0}
=
\theta_{0}$
and
$\theta_{0}-\varepsilon\log r_{\min}
=
\theta_{0}+2\theta_{0}\frac{\log r_{\min}}{\log r_{\max}}
\le
3\theta_{0}\frac{\log r_{\min}}{\log r_{\max}}
\le
(\frac{d_{0}}{1+d_{0}})^\frac{1}{2}$, we therefore deduce that
$\theta_{0}
 \le
 u_{n}(q)\log r_{i}-2\pi m_{i}
 \le
(\frac{d_{0}}{1+d_{0}})^\frac{1}{2}$.
This clearly implies that
$\theta_{0}\le\Theta_{i}\big(u_{n}(q)\big)\le (\frac{d_{0}}{1+d_{0}})^\frac{1}{2}$.

\noindent
(3)
We have $v_{n}(q)=b_{n}+\varepsilon$ where
$\varepsilon=-2\theta_{0}\frac{1}{\log r_{\max}}$.
It follows from this and (13.13) that 
for each $i$,
we can find an integer $m_{i}\in\Bbb Z$ such that
$-\theta_{0}+\varepsilon\log r_{\min}
 \le
 v_{n}(q)\log r_{i}-2\pi m_{i}
 \le
 \theta_{0}+\varepsilon\log r_{\min}$.
Since
$-\theta_{0}+\varepsilon\log r_{\min}
=
-\theta_{0}-2\theta_{0}\frac{\log r_{\min}}{\log r_{\max}}
\ge
-3\theta_{0}\frac{\log r_{\min}}{\log r_{\max}}
\ge
-(\frac{d_{0}}{1+d_{0}})^\frac{1}{2}$
and
$\theta_{0}+\varepsilon\log r_{\min}
=
\theta_{0}-2\theta_{0}\frac{\log r_{\min}}{\log r_{\max}}
\le
\theta_{0}-2\theta_{0}
=
-\theta_{0}$, we therefore deduce that
$-(\frac{d_{0}}{1+d_{0}})^\frac{1}{2}
 \le
  v_{n}(q)\log r_{i}-2\pi m_{i}
 \le
-\theta_{0}$.
This clearly implies that
$-(\frac{d_{0}}{1+d_{0}})^\frac{1}{2}
\le
\Theta_{i}\big(v_{n}(q)\big)
\le 
-\theta_{0}$.
This completes the proof of Claim 4.

\bigskip

\noindent
In the claims below we use the following notation, namely, we define
the function
$f:\Bbb C\to\Bbb C$ by
$f(s)=1-\sum_{i}p_{i}^{q}r_{i}^{s}$.

\bigskip

\noindent{
\it
Claim 5.
Let $n$ be an integer.
Then the following statements hold.
\roster
\item"(1)"
$\inf_{s\in \Gamma_{n}^{-}}\real(f(s))\ge\delta_{0}$.
\item"(2)"
$\sup_{s\in \Gamma_{n}^{-}}|\zeta_{\mu}^{q}(s)|\le 
 \frac{\sum_{i}p_{i}^{q}r_{i}^{\beta_{0}}}{\delta_{0}}$.\endroster
}
\noindent
{\it Proof of Claim 5.}
(1)
Let $s\in \Gamma_{n}^{-}$. Consequently, there is $t\in[v_{n-1}(q),u_{n}(q)]$ such that
$s=\beta_{0}(q)+\ima t$.
As $t\in[v_{n-1}(q),u_{n}(q)]$, we deduce that
$t\not\in G$, whence
$ \max_{i}|\Theta_{i}(t)|\ge\theta_{0}$
and we can therefore find $i_{0}$ such that
$|\Theta_{i_{0}}(t)|\ge\theta_{0}$.
Using the fact that
$\cos \theta\le 1-\frac{1}{4}\theta^{2}$ for all $\theta$ with $|\theta|\le \frac{\pi}{2}$, this
implies that
 $$
 \align
 \cos\Theta_{i_{0}}(t)
&\le
 \cos\theta_{0}\\
&\le
 1-\frac{1}{4}\theta_{0}^{2} \\
&=
 1-2\delta_{0}\,\frac{1}{\min_{i}(p_{i}^{q}r_{i}^{\beta(q)})}\\
&\le
 1-2\delta_{0}
 \,
 \frac{1}{p_{i_{0}}^{q}r_{i_{0}}^{\beta(q)}}  
 \endalign
 $$
and so
 $$
 p_{i_{0}}^{q}r_{i_{0}}^{\beta_{0}(q)}\cos\Theta_{i_{0}}(t)
\le
p_{i_{0}}^{q}r_{i_{0}}^{\beta_{0}(q)}
-
2\delta_{0}\frac{r_{i_{0}}^{\beta_{0}(q)}}{r_{i_{0}}^{\beta(q)}}  \,.
\tag13.14
$$
Next, since $\beta_{0}(q)<\beta(q)$, we conclude that
$r_{i_{0}}^{\beta(q)}<r_{i_{0}}^{\beta_{0}(q)}$,
whence
$-\frac{r_{i_{0}}^{\beta_{0}(q)}}{r_{i_{0}}^{\beta(q)}} \le -1$.
It therefore follows from (13.14) that
 $$
 \align
  p_{i_{0}}^{q}r_{i_{0}}^{\beta_{0}(q)}\cos\Theta_{i_{0}}(t)
&\le
p_{i_{0}}^{q}r_{i_{0}}^{\beta_{0}(q)}-2\delta_{0}\,.
\tag13.15
\endalign
$$
Inequality (13.15) now  implies that
 $$
 \align
 \real
 \Bigg(
 \sum_{i}p_{i}^{q}r_{i}^{s}
 \Bigg)
&=
  \real
 \Big(
 p_{i_{0}}^{q}r_{i_{0}}^{\beta_{0}(q)}e^{\ima t\log r_{i}}
 \Big)
+
  \real
 \Bigg(
 \sum_{i\not= i_{0}}p_{i}^{q}r_{i}^{\beta_{0}(q)}e^{\ima t\log r_{i}}
 \Bigg)\\
&\le
 p_{i_{0}}^{q}r_{i_{0}}^{\beta_{0}(q)}\cos\Theta_{i_{0}}(t)
+
 \sum_{i\not= i_{0}}p_{i}^{q}r_{i}^{\beta_{0}(q)}
 \\
&\le
p_{i_{0}}^{q}r_{i_{0}}^{\beta_{0}(q)}-2\delta_{0} 
+
 \sum_{i\not= i_{0}}p_{i}^{q}r_{i}^{\beta_{0}(q)}
 \\
 &=
 \sum_{i}p_{i}^{q}r_{i}^{\beta_{0}(q)}-2\delta_{0} 
 \\
&=
 1+\delta_{0}-2\delta_{0} \\
&=
1-\delta_{0}\,. 
\tag13.16
\endalign
$$
We see from (13.16) that
$\real(f(s))
=
\real(1-\sum_{i}p_{i}^{q}r_{i}^{s})
=
1-\real(\sum_{i}p_{i}^{q}r_{i}^{s})
\ge
1-(1-\delta_{0})
=
\delta_{0}$.

\noindent
(2)
Since $\zeta_{\mu}^{q}=\frac{1-f}{f}$, we conclude from (1) that
$|\zeta_{\mu}^{q}(s)|
\le
 \frac{|1-f(s)|}{|f(s)|}
 \le
 \frac{\sum_{i}p_{i}^{q}r_{i}^{\beta_{0}(q)}}{|\real(f(s))|}
 \le
  \frac{\sum_{i}p_{i}^{q}r_{i}^{\beta_{0}(q)}}{\delta_{0}}$
    for all $s\in\Gamma_{n}^{-}$.
This completes the proof of Claim 5.

\bigskip

\noindent{
\it
Claim 6.
Let $n$ be an integer.
Then the following statements hold.
\roster
\item"(1)"
$\sup_{s\in \Gamma_{n}^{+}}\real(f(s))\le-\frac{d_{0}}{2}$.
\item"(2)"
$\sup_{s\in \Gamma_{n}^{+}}|\zeta_{\mu}^{q}(s)|\le 
 2\frac{\sum_{i}p_{i}^{q}r_{i}^{b_{0}}}{\delta_{0}}$.\endroster
}
\noindent
{\it Proof of Claim 6.}
(1)
Let $s\in \Gamma_{n}^{+}$. Consequently, there is $t\in[u_{n}(q),v_{n}(q)]$ such that
$s=b_{0}(q)+\ima t$.
As $t\in[u_{n}(q),v_{n}(q)]$, we conclude from statement (1) in Claim 4
that
$-(\frac{d_{0}}{1+d_{0}})^{\frac{1}{2}}
\le
\Theta_{i}\big(t\big)
\le
(\frac{d_{0}}{1+d_{0}})^{\frac{1}{2}}$
whence
$\cos \Theta_{i}\big(t\big)
\ge
\cos((\frac{d_{0}}{1+d_{0}})^{\frac{1}{2}})$
for all $i$.
Using the fact that
$\cos \theta\ge 1-\frac{1}{2}\theta^{2}$ for all $\theta$, this
implies that
 $$
 \align
 \cos\Theta_{i}(t)
&\ge
 \cos \Bigg( \bigg(\frac{d_{0}}{1+d_{0}}\bigg)^{\frac{1}{2}}\Bigg)\\
&\ge
 1-\frac{1}{2}\,\frac{d_{0}}{1+d_{0}}
 \endalign
$$
for all $i$,
and so
 $$
 \align
 \real
 \Bigg(
 \sum_{i}p_{i}^{q}r_{i}^{s}
 \Bigg)
&=
  \real
 \Bigg(
 \sum_{i}p_{i}^{q}r_{i}^{b_{0}(q)}e^{\ima t\log r_{i}}
 \Bigg)\\
&=
 \sum_{i}p_{i}^{q}r_{i}^{b_{0}(q)}\cos\Theta_{i}(t)
 \\
&\ge
 \sum_{i}p_{i}^{q}r_{i}^{b_{0}(q)}
 \Bigg(1-\frac{1}{2}\,\frac{d_{0}}{1+d_{0}}\Bigg)
 \\
 &=
(1+d_{0})
 \Bigg(1-\frac{1}{2}\,\frac{d_{0}}{1+d_{0}}\Bigg)
 \\
&=
 1+\frac{d_{0}}{2}\,.
 \tag13.17
 \endalign
$$
We see from (13.17) that
$\real(f(s))
=
\real(1-\sum_{i}p_{i}^{q}r_{i}^{s})
=
1-\real(\sum_{i}p_{i}^{q}r_{i}^{s})
\le
1-(1+\frac{d_{0}}{2})
=
-\frac{d_{0}}{2}$.

\noindent
(2)
Since $\zeta_{\mu}^{q}=\frac{1-f}{f}$, we conclude from (1) that
$|\zeta_{\mu}^{q}(s)|
\le
 \frac{|1-f(s)|}{|f(s)|}
 \le
 \frac{\sum_{i}p_{i}^{q}r_{i}^{b_{0}(q)}}{|\real(f(s))|}
 \le
  2\frac{\sum_{i}p_{i}^{q}r_{i}^{b_{0}(q)}}{d_{0}}$
  for all $s\in\Gamma_{n}^{+}$.
This completes the proof of Claim 6.

\bigskip

\noindent{
\it
Claim 7.
Let $n$ be an integer.
Let $M$ be a subset of $\Bbb R$
and write
$c=\inf_{x\in M}\sum_{i}p_{i}^{q}r_{i}^{x}$
and
$C=\sup_{x\in M}\sum_{i}p_{i}^{q}r_{i}^{x}$.
Then the following statements hold.
\roster
\item"(1)"
$\inf_{s\in (M+u_{n}(q)\ima)\cup(M+v_{n}(q)\ima)}
|\imag(f(s))|\ge c\sin\theta_{0}$.
\item"(2)"
$\sup_{s\in (M+u_{n}(q)\ima)\cup(M+v_{n}(q)\ima)}
|\zeta_{\mu}^{q}(s)|\le 
 \frac{C}{c\sin\theta_{0}}$
 .\endroster
In particular, if $M=[b_{0}(q),\beta_{0}(q)]$, then
 $M+u_{n}(q)\ima=\Pi_{n}^{-}$ and
 $M+v_{n}(q)\ima=\Pi_{n}^{+}$,
 and it therefore follows from (1) and (2) that
}
\roster
\item"(3)"
$\inf_{s\in \Pi_{n}^{-}\cup\Pi_{n}^{+}}
|\imag(f(s))|\ge\sin\theta_{0}$.
\item"(4)"
$\sup_{s\in \Pi_{n}^{-}\cup\Pi_{n}^{+}}
|\zeta_{\mu}^{q}(s)|\le 
 \frac{\sum_{i}p_{i}^{q}r_{i}^{b_{0}(q)}}{\sin\theta_{0}}$
 .\endroster

\noindent
{\it Proof of Claim 7.}
(1)
Let $s\in  (M+u_{n}(q)\ima)\cup(M+v_{n}(q)\ima)$. 
Consequently, there is $\sigma\in M$ 
and $t\in\{u_{n}(q),v_{n}(q)\}$ such that
$s=\sigma+\ima t$,
and so
 $$
 \align
 \imag
 \Bigg(
 \sum_{i}p_{i}^{q}r_{i}^{s}
 \Bigg)
&=
  \imag
 \Bigg(
 \sum_{i}p_{i}^{q}r_{i}^{\sigma}e^{\ima t\log r_{i}}
 \Bigg)\\
&=
 \sum_{i}p_{i}^{q}r_{i}^{\sigma}\sin \Theta_{i}(t)\,.
 \tag13.18
\endalign
$$
We now make the following two observations.
Namely,
it follows from
statement (2) in Claim 4
that
$\theta_{0}
\le
\Theta_{i}\big(u_{n}(q)\big)
\le
(\frac{d_{0}}{1+d_{0}})^{\frac{1}{2}}$
whence
$\sin \Theta_{i}\big(u_{n}(q)\big)
\ge
\sin\theta_{0}$
for all $i$,
and it
follows 
from statement (3) in Claim 4
that
$-(\frac{d_{0}}{1+d_{0}})^{\frac{1}{2}}
\le
\Theta_{i}\big(v_{n}(q)\big)
\le
-\theta_{0}$
whence
$\sin \Theta_{i}\big(v_{n}(q)\big)
\le
\sin(-\theta_{0})=-\sin\theta_{0}$
for all $i$.
These observations and (13.18) imply that
 $$
 \imag
 \Bigg(
 \sum_{i}p_{i}^{q}r_{i}^{s}
 \Bigg)
=
\cases
\dsize
\ge
\,\,\,\,\,
 \sum_{i}p_{i}^{q}r_{i}^{\sigma}
\sin \theta_{0}
&\quad
\text{for $t=u_{n}(q)$;}\\
\dsize
\le
 -
 \sum_{i}p_{i}^{q}r_{i}^{\sigma}
\sin \theta_{0}
&\quad
\text{for $t=v_{n}(q)$.}
\endcases
\tag13.19
$$
Next, since $\sigma\in M$, it follows that
$ \sum_{i}p_{i}^{q}r_{i}^{\sigma}\ge c$, and we therefore conclude from (13.19) that
 $$
 \imag
 \Bigg(
 \sum_{i}p_{i}^{q}r_{i}^{s}
 \Bigg)
 =
\cases
\dsize
\ge
\,\,\,\,\,
c\sin \theta_{0}
&\quad
\text{for $t=u_{n}(q)$;}\\
\dsize
\le
 -
c\sin \theta_{0}
&\quad
\text{for $t=v_{n}(q)$.}
\endcases
$$
It follows from this that
$|\imag(f(s))|
=
|\imag(1-\sum_{i}p_{i}^{q}r_{i}^{s})|
=
|-\imag(\sum_{i}p_{i}^{q}r_{i}^{s})|
\ge
c\sin\theta_{0}$.

\noindent
(2)
Since $\zeta_{\mu}^{q}=\frac{1-f}{f}$, we deduce from (1) that
$|\zeta_{\mu}^{q}(s)|
\le
 \frac{|1-f(s)|}{|f(s)|}
 \le
  \frac{C}{c\sin\theta_{0}}$
  for all $s\in (M+u_{n}(q)\ima)\cup(M+v_{n}(q)\ima)$.
This completes the proof of Claim 7.

\bigskip

The proof  now follows from Claims 1,3,5--7.
\hfill$\square$

\bigskip

\noindent
The next two results, i\.e Lemma 13.4 and Theorem 13.5,
give growth estimates
of the
zeta function $\zeta_{\mu}^{q}$ outside and inside the critical  strip, respectively.
For $i$ with $r_{i}>r_{\min}$, write $s_{i}=\frac{r_{\min}}{r_{i}}$ and put
$s_{\max}=\max_{r_{i}>r_{\min}}s_{i}$.

\bigskip

\proclaim{Lemma 13.4. Growth estimates of $\zeta_{\mu}^{q}$
outside the critical strip 
$\alpha(q)\le\real(s)\le\beta(q)$}
Fix $q\in\Bbb R$.
Define $f:\Bbb C\to\Bbb C$ by
$$
f(s)
=
1-\sum_{i}p_{i}^{q}r_{i}^{s}\,.
$$
There are two real numbers $A(q)$ and $B(q)$
and two  constants $c>0$ and $k>0$
satisfying the following:
\roster
\item"(1)"
We have
$A(q)
\le
\alpha(q)
\le
\beta(q)
\le
B(q)$.

\item"(2)"
For all $\sigma\ge B(q)$ and $t\in\Bbb R$, we have
 $
 |
 \frac{f'(\sigma+\ima t)}{f(\sigma+\ima t)}
 |
\le
cr_{\max}^{\sigma}$.

\item"(3)"
For all $\sigma\le A(q)$ and $t\in\Bbb R$, we have
 $
 |
 \frac{f'(\sigma+\ima t)}{f(\sigma+\ima t)}
 -
 \log r_{\min}
 |
\le
cs_{\max}^{-\sigma}$.

\item"(4)"
For all $\sigma\ge B(q)$ and $t\in\Bbb R$, we have
$
\frac{1}{2}\le|f(\sigma+\ima t)|\le\frac{3}{2}$.

\item"(5)"
For all $\sigma\le A(q)$ and $t\in\Bbb R$, we have
$
|\zeta_{\mu}^{q}(\sigma+\ima t)|
\le
k$.
\endroster
\endproclaim

\noindent{\it  Proof}\newline
\noindent
(1)--(4)
These results follows from straight forward estimates
and  we are therefore omitted the proofs.

\noindent
(5)
It suffices to show that there are numbers $\sigma_{q}$ and $c_{q}$ such that
for
all $\sigma\le \sigma_{q}$ and all $t\in\Bbb R$, we have
$|\zeta_{\mu}^{q}(\sigma+\ima t)|\le c_{q}$.
First, 
observe that
if we 
for a real number $\sigma$,
define $f_{\sigma}:\Bbb R\to \Bbb C$ by
$f_{\sigma}(t)
=
\sum_{i}p_{i}^{q}r_{i}^{\ima t}(\frac{r_{i}}{r_{\min}})^{\sigma}$,
then
 $$
 \align
 \left|
 1-\sum_{i}p_{i}^{q}r_{i}^{\sigma+\ima t}
 \right|
&=
r_{\min}^{\sigma}
 \Bigg|f_{\sigma}(t)- \frac{1}{r_{\min}^{\sigma}}\Bigg|\,.
 \tag13.20
 \endalign 
 $$

Next, note that 
 $$
 \align
 \sup_{t\in\Bbb R}
 \Bigg|
 \,
  \Bigg|f_{\sigma}(t)- \frac{1}{r_{\min}^{\sigma}}\Bigg|
  -
  \sum_{r_{i}=r_{\min}}p_{i}^{q}
  \,
  \Bigg|
&=
 \sup_{t\in\Bbb R}
  \Bigg|
 \,
  \Bigg|f_{\sigma}(t)- \frac{1}{r_{\min}^{\sigma}}\Bigg|
  -
 \Bigg|\sum_{r_{i}=r_{\min}}p_{i}^{q}r_{\min}^{\ima t}\Bigg|  
 \,
  \Bigg|\\
&\le
 \sup_{t\in\Bbb R}
 \Bigg|
 \,\,
 \,
  f_{\sigma}(t)- \frac{1}{r_{\min}^{\sigma}}
  \,\,
  -
  \,\,\,
  \sum_{r_{i}=r_{\min}}p_{i}^{q}r_{\min}^{\ima t}
  \,
  \,\,
  \Bigg|
  \\\   
&\le
 \sup_{t\in\Bbb R}
 \left|
 \sum_{i}p_{i}^{q}r_{i}^{\ima t}(\tfrac{r_{i}}{r_{\min}})^{\sigma}
 -
  \sum_{r_{i}=r_{\min}}p_{i}^{q}r_{\min}^{\ima t}
  \,\,
 \right|  
 +
 \frac{1}{r_{\min}^{\sigma}}\,.\\
& 
 \tag13.21
 \endalign
 $$
Since, clearly
$\sum_{i}p_{i}^{q}r_{i}^{\ima t}(\tfrac{r_{i}}{r_{\min}})^{\sigma}
 -
  \sum_{r_{i}=r_{\min}}p_{i}^{q}r_{\min}^{\ima t}
 =
 \sum_{r_{i}>r_{\min}}p_{i}^{q}r_{i}^{\ima t}(\tfrac{r_{i}}{r_{\min}})^{\sigma}$, we now deduce form 
(13.21) that
 $$
 \align
 \! \! \! \! \! \! \! \! \! \! \!
 \! \! \! \! \! \! \! \! \! \! \!
 \! \! \! \! \! \! \! \! \! \! \!
  \! \! \! \!
 \sup_{t\in\Bbb R}
 \Bigg|
 \,
  \Bigg|f_{\sigma}(t)- \frac{1}{r_{\min}^{\sigma}}\Bigg|
  -
  \sum_{r_{i}=r_{\min}}p_{i}^{q}
  \,
  \Bigg|
&=
 \sup_{t\in\Bbb R}
 \left|
 \sum_{r_{i}>r_{\min}}p_{i}^{q}r_{i}^{\ima t}(\tfrac{r_{i}}{r_{\min}})^{\sigma}
 \right|
 +
 \frac{1}{r_{\min}^{\sigma}}\\ 
&\le
 \sup_{t\in\Bbb R}
 \sum_{r_{i}>r_{\min}}p_{i}^{q}\,|r_{i}^{\ima t}|\,(\tfrac{r_{i}}{r_{\min}})^{\sigma}
 +
 \frac{1}{r_{\min}^{\sigma}}\\  
&\le
 \sum_{r_{i}>r_{\min}}p_{i}^{q}(\tfrac{r_{i}}{r_{\min}})^{\sigma}
 +
 \frac{1}{r_{\min}^{\sigma}}\\   
&\to
 0
 \,\,\,\,
 \text{as $\sigma\to-\infty$.}
 \tag13.22
 \endalign
 $$
It follows from (13.22) that we can choose a  number $\sigma_{q}$ with $\sigma_{q}<0$ such that
for all $\sigma\le \sigma_{q}$, we have
$\sup_{t\in\Bbb R}
 \big|
 \,
  |f_{\sigma}(t)-\frac{1}{r_{\min}^{\sigma}}|-\sum_{r_{i}=r_{\min}}p_{i}^{q}
  \,
  \big|
  \le
  \frac{1}{2}\sum_{r_{i}=r_{\min}}p_{i}^{q}$.
 This clearly implies that for 
for all $\sigma\ge \sigma_{q}$ and all $t\in\Bbb R$, we have
 $$
 \Bigg|
 f_{\sigma}(t)-\frac{1}{r_{\min}^{\sigma}}
 \Bigg|
 \ge
 \frac{1}{2}\sum_{r_{i}=r_{\min}}p_{i}^{q}\,.
 \tag13.23
 $$

Finally, combining (13.20) and (13.23) shows that
for
all $\sigma\le \sigma_{q}$ and all $t\in\Bbb R$, we have
 $$
   \align
  |\zeta_{\mu}^{q}(\sigma+\ima t)|
 &=
  \frac
  {
  |
  \sum_{i}p_{i}^{q}r_{i}^{\sigma+\ima t}
  |
  }
  {
  |
 1-\sum_{i}p_{i}^{q}r_{i}^{\sigma+\ima t}
 |
  }\\
 &\le
  \frac
  { 
 \sum_{i}p_{i}^{q}r_{i}^{\sigma}
 }
  {
 r_{\min}^{\sigma}
 |f_{\sigma}(t)-\frac{1}{r_{\min}^{\sigma}}|
  }\\
  &\le
  \frac
  { 
 \sum_{i}p_{i}^{q}r_{\min}^{\sigma}
 }
  {
r_{\min}^{\sigma}
 \frac{1}{2}\sum_{r_{i}=r_{\min}}p_{i}^{q}\
  }\\ 
 &=
 2
  \frac
  { 
 \sum_{i}p_{i}^{q}
 }
  {
 \sum_{r_{i}=r_{\min}}p_{i}^{q}\
  }\,.
 \endalign
 $$
This completes the proof
\hfill$\square$

\bigskip

\proclaim{Theorem 13.5 (i\.e\. Theorem 5.3). Growth estimates of $\zeta_{\mu}^{q}$
inside the critical strip 
$\alpha(q)\le\real(s)\le\beta(q)$}
Fix $q\in\Bbb R$.
Then there is an increasing a sequence $(t_{q,n})_{n}$
of positive real numbers with $t_{q,n}\to\infty$
satisfying the following:
for all real numbers $c$, there
is a constant
$k_{c}$ such that
for all $\sigma\le c$ and all $n$,
we have
 $$
 \align
 |\zeta_{\mu}^{q}(\sigma\pm\ima t_{q,n})|
&\le
 k_{c}\,.
 \endalign
 $$

\endproclaim
\noindent{\it  Proof}\newline
\noindent
For $n\in\Bbb Z$,
let $u_{n}(q)$ and $v_{n}(q)$ be defined as is Proposition 13.3.
Now
define $t_{q,n}$ by
 $$
 t_{q,n}=u_{n}(q)
 \tag13.24
 $$ 
for $n\in\Bbb N$ and 
observe that
it follows from Proposition 13.3
that
 $$
 -t_{q,n}
=
-u_{n}(q)
=
v_{-n}(q)
\tag13.25
$$
for $n\in\Bbb N$.
It also follows from Proposition 13.3 that
$t_{q,n}\to\infty$.
Next, fix a real number $c$.
We must now prove that
$\sup_{
 \sigma\le c\,,\,n\in\Bbb N
 }
 |\zeta_{\mu}^{q}(\sigma\pm\ima t_{q,n})|
<
\infty$.
Letting $A(q)$ be as 
is Lemma 13.4,
it follows
from 
Lemma 13.4.(4)
and
Proposition 13.3.(8)
(applied to the 
set
$M=[A(q),c]$)
that
$ \sup_{
 \sigma\le A(q),t\in\Bbb R
 }
 |\zeta_{\mu}^{q}(\sigma\pm\ima t)|<\infty$
 and
$  \sup_{n}\sup_{s\in([A(q),c]+u_{n}(q)\ima )\cup([A(q),c]+v_{n}(q)\ima)}
 |\zeta_{\mu}^{q}(s)|
 <
 \infty$,
 respectively, whence
 (using (13.24) and (13.25))
 $$
 \align
 \sup
 \Sb
 \sigma\le c\\
 n\in\Bbb N
 \endSb
 |\zeta_{\mu}^{q}(\sigma\pm\ima t_{q,n})|
&=
\max
\left(
\,
 \sup
 \Sb
 \sigma\le A(q)\\
 n\in\Bbb N
 \endSb
 |\zeta_{\mu}^{q}(\sigma\pm\ima t_{q,n})|
 \,,\,
  \sup
 \Sb
 \sigma\in[A(q),c]\\
 n\in\Bbb N
 \endSb
 |\zeta_{\mu}^{q}(\sigma+\ima t_{q,n})|
 \right.
 \,,\,\\
&\qquad\qquad
  \qquad\qquad
  \qquad\qquad 
  \qquad\qquad
  \left.
  \sup
 \Sb
 \sigma\in[A(q),c]\\
 n\in\Bbb N
 \endSb
 |\zeta_{\mu}^{q}(\sigma-\ima t_{q,n})| 
\right)\\
&\le
\max
\left(
\,
 \sup
 \Sb
 \sigma\le A(q)\\
t\in\Bbb R
 \endSb
 |\zeta_{\mu}^{q}(\sigma\pm\ima t)|
 \,,\,
 \sup_{n\in\Bbb Z}
  \sup
  \Sb
  s\in [A(q),c]+u_{n}(q)\ima \\
  {}
  \endSb
 |\zeta_{\mu}^{q}(s)|
 \right.
 \,,\,\\
&\qquad\qquad
  \qquad\qquad
  \qquad\qquad 
  \qquad\qquad
\left.
 \sup_{n\in\Bbb Z}
  \sup
   \Sb
  s\in [A(q),c]+v_{n}(q)\ima \\
  {}\\
  {}
  \endSb
 |\zeta_{\mu}^{q}(s)|
\right)\\
&<
\infty\,.
\endalign
$$
This completes the proof.
\hfill$\square$

\bigskip

The final result in this section, i\.e\. Theorem 13.8, provides an
estimate for the
density of the poles of the zeta function $\zeta_{\mu}^{q}$.
We first list two variants of Jensen's formula from complex analysis
(see [Con, p\. 280]) that
will be needed in order to prove
Theorem 13.8.
We begin with a definition.
For a holomorphic function $F$ and $R>0$ write
 $$
M_{F}(R)
=
\sup
 \Sb
 s\in\Bbb C\\
 |s|=R
 \endSb
 |F(s)|\,.
 $$
We can now state the two results from complex analysis needed
to prove Theorem 13.8. 
 
\bigskip

\proclaim{Proposition 13.6 [Te, Corollary 11.2]}
Let $F$ 
be a holomorphic function with $F(0)=1$
and
let $R$ and $\rho$ be positive real numbers with $0<\rho<\frac{1}{2}R$. Then
 $$
 \sup
 \Sb
 s\in\Bbb C\setminus Z(F)\\
 |s|\le \rho
 \endSb
 \left|
 \frac{F'(s)}{F(s)}
 -
 \sum
 \Sb
 \omega\in Z(F)\\
 |\omega|\le\frac{1}{2}R
 \endSb
 \frac{1}{s-\omega}
 \right|
 \le
 \frac{4R}{(R-2\rho)^{2}}
 \log
 M_{F}(R)\,.
 $$
 \endproclaim

 \bigskip

\proclaim{Proposition 13.7 [Te, Section 11; Rud, p\. 309]}
Let $F$ 
be a holomorphic function with $F(0)=1$
and
let $R$ be a positive real number.
Then
 $$
 \Big|
 \,
 \big\{\omega\in Z(F)
 \,\big|\,
 |\omega|\le \tfrac{1}{2}R
 \big\}
 \,
 \Big|
 \le
 \frac{1}{\log 2}
  \log
M_{F}(R)\,.
 $$
\endproclaim

\bigskip

\proclaim{Theorem 13.8 (i\.e\. statement (5) in Proposition 15.2).
Density of poles of $\zeta_{\mu}^{q}$}
Fix $q\in\Bbb R$.
Write
$\gamma = -\frac{1}{\pi}\log r_{\min}$.
We have
 $$
 \Big|
 \,
 \big\{
 \omega\in P(\zeta_{\mu}^{q})
 \,\big|\,
 |\imag(\omega)|\le t
 \big\}
 \,
 \Big|
 \,
 =
 \,
\gamma t
 +
 \Cal O(\log t)\,.
 $$
\endproclaim
\noindent{\it  Proof}\newline
For brevity write
$N_{t}
=
 |
 \,
 \{
 \omega\in P(\zeta_{\mu}^{q})
 \,|\,
 |\imag(\omega)|\le t
 \}
 \,
 |$.
Define $f:\Bbb C\to\Bbb C$ by
$f(s)=1-\sum_{i}p_{i}^{q}r_{i}^{s}$ and note that
$\zeta_{\mu}^{q}=\frac{1-f}{f}$. In particular, 
it follows from this that
$P(\zeta_{\mu}^{q})
 =
 Z(f)$, whence
$N_{t}
 =
 |
 \,
 \{
 \omega\in Z(f)
 \,|\,
 |\imag(\omega)|\le t
 \}
 \,
 |$.

Let the constants $r_{\max}$ and $s_{\max}$
be as in Lemma 13.4,
and
recall that 
it follows Lemma 13.4  that
there is a constant $c>0$ such that
if $\sigma\ge B(q)$ and $t\in\Bbb R$, then
 $$
 \Bigg|
 \frac{f'(\sigma+\ima t)}{f(\sigma+\ima t)}
 \Bigg|
\le
cr_{\max}^{\sigma}\,,
$$
and if $\sigma\le A(q)$ and $t\in\Bbb R$, then 
$$
\Bigg|
 \frac{f'(\sigma+\ima t)}{f(\sigma+\ima t)}
 -
 \log r_{\min}
 \Bigg|
\le
cs_{\max}^{-\sigma}\,.
$$
Also, 
for $t>0$,
we
choose $\sigma_{t}^{+}$ such that
$cr_{\max}^{\sigma_{t}^{+}}=\frac{1}{t}$  
and choose 
we
$\sigma_{t}^{-}$ such that
$cs_{\max}^{-\sigma_{t}^{-}}=\frac{1}{t}$.
It is clear that $\sigma_{t}^{+}\to\infty$ as $t\to\infty$
and that
$\sigma_{t}^{-}\to-\infty$ as $t\to\infty$ as $t\to\infty$.
We can therefore find $t_{0}>0$ such that
$\sigma_{t}^{+}\ge B(q)$ for $t\ge t_{0}$
and
$\sigma_{t}^{-}\le A(q)$ for $t\ge t_{0}$.

Next,
we
fix a positive real number $t$ with $t\ge t_{0}$
such that
$f(\sigma+\ima t)\not=0$
for all $\sigma\in\Bbb R$,
and 
define paths
$\Sigma_{t}^{+},\Sigma_{t}^{-},\Lambda_{t},\Delta_{t}$ in $\Bbb C$
by
 $$
 \align
 \Sigma_{t}^{+}
&\quad
 \text{
 is the directed line-segment from
 \quad\,
 $\sigma_{t}^{+}+\ima t$
 to
 $\sigma_{t}^{-}+\ima t$
 }\,,\\
 \Sigma_{t}^{-}
&\quad
 \text{
 is the directed line-segment from
 \quad\,
 $\sigma_{t}^{-}-\ima t$
 to
 $\sigma_{t}^{+}-\ima t$
 }\,,\\
  \Lambda_{t}\,\,
&\quad
 \text{
 is the directed line-segment from
\quad\,
  $\sigma_{t}^{+}-\ima t$
 to
 $\sigma_{t}^{+}+\ima t$
 }\,,\\
 \Delta_{t}
&\quad
 \text{
 is the directed line-segment from
\quad\,
  $\sigma_{t}^{-}+\ima t$
 to
 $\sigma_{t}^{-}-\ima t$
 }\,,
 \endalign
 $$
and let $C_{t}$ denote the simple closed path
obtained by concatenating 
the line segments
$\Lambda_{t}$,
$\Sigma_{t}^{+}$,  $\Delta_{t}$  and $\Sigma_{t}^{-}$.
Since $C_{t}$ does not pass through any of the zeros of $f$
(because 
$t$ is chosen such that
$f(\sigma+\ima t)\not=0$
for all $\sigma\in\Bbb R$)
and
since
all zeros $\omega$ of $f$ satisfy
$\real(\omega)
\in
[\alpha(q),\beta(q)]
\subseteq
[A(q),B(q)]$ (see  Lemma 13.4),
it follows from the 
Argument Principle (see [Con, p\. 123])
 that
$2\pi\ima N_{t}
 =
 2\pi\ima
 |
 \,
 \{
 \omega\in Z(f)
 \,|\,
 |\imag(\omega)|\le t
 \}
 \,
 |
 =
 \int_{C_{t}}\frac{f'(s)}{f(s)}\,ds$, whence
 $$
 \align
 2\pi
 N_{t}
&=
\imag
\int_{C_{t}}\frac{f'(s)}{f(s)}\,ds\\
&=
\imag
 \int_{\Sigma_{t}^{+}}\frac{f'(s)}{f(s)}\,ds
 +
 \imag
  \int_{\Sigma_{t}^{-}}\frac{f'(s)}{f(s)}\,ds
  +
  \imag
  \int_{\Lambda_{t}}\frac{f'(s)}{f(s)}\,ds
  +
  \imag
   \int_{\Delta_{t}}\frac{f'(s)}{f(s)}\,ds\,.
   \tag13.26
 \endalign
$$
As
$\int_{\Delta_{t}}\, ds=-2 t\ima $, 
and so
$\int_{\Delta_{t}}(-\log r_{\min})\, ds=-2\pi\ima\gamma t$,
we conclude from (13.26) that
 $$
 \align
 2\pi
 |N_{t}
&-\gamma t|\\
&=
\Bigg|
 \imag
 \int_{\Sigma_{t}^{+}}\frac{f'(s)}{f(s)}\,ds
 +
 \imag
  \int_{\Sigma_{t}^{-}}\frac{f'(s)}{f(s)}\,ds
  +
  \imag
  \int_{\Lambda_{t}}\frac{f'(s)}{f(s)}\,ds
  +
  \imag
   \int_{\Delta_{t}}\frac{f'(s)}{f(s)}\,ds
   -2\pi\gamma t
   \Bigg|
   \\
&=
 \Bigg|
\imag
 \int_{\Sigma_{t}^{+}}\frac{f'(s)}{f(s)}\,ds
 +
 \imag
  \int_{\Sigma_{t}^{-}}\frac{f'(s)}{f(s)}\,ds\\
  &\qquad\qquad
   \qquad\qquad
   \qquad\qquad  
  +
  \imag
  \int_{\Lambda_{t}}\frac{f'(s)}{f(s)}\,ds
  +
  \imag
   \int_{\Delta_{t}}
   \Bigg(
   \frac{f'(s)}{f(s)}\,
   -\log r_{\min}
   \Bigg)\,ds
 \Bigg| \\
&\le
 \Bigg|\imag\int_{\Sigma_{t}^{+}}\frac{f'(s)}{f(s)}\,ds\Bigg|
 +
  \Bigg|\imag\int_{\Sigma_{t}^{-}}\frac{f'(s)}{f(s)}\,ds\Bigg|\\
 &\qquad\qquad
   \qquad\qquad
   \qquad\qquad 
  +
  \Bigg|\int_{\Lambda_{t}}\frac{f'(s)}{f(s)}\,ds\Bigg|
  +
  \Bigg|
   \int_{\Delta_{t}}
   \Bigg(
   \frac{f'(s)}{f(s)}\,
   -\log r_{\min}
   \Bigg)
   \,ds
   \Bigg|\,.\\
 &{}  
   \tag13.27
 \endalign
 $$

We now estimates the four integrals in (13.27).

\bigskip

\noindent
{\it Claim 1.
We have
 $
  |\int_{\Lambda_{t}}\frac{f'(s)}{f(s)}\,ds|
  \le
  2$.
 }

\noindent
{\it Proof of Claim 1.} 
Write
$\ell(\Lambda_{t})=2t$ for the length of $\Lambda_{t}$.
For $s=\sigma+\ima \tau\in\Lambda_{t}$ with $\sigma,\tau\in\Bbb R$, we have
$\sigma\ge\sigma_{t}^{+}$
whence
$|\frac{f'(s)}{f(s)}|\le cr_{\max}^{\sigma}\le  \frac{1}{t}$,
and consequently
$|
  \int_{\Lambda_{t}}\frac{f'(s)}{f(s)}\,ds
  |
  \le
  \ell(\Lambda_{t})
  \,
  \sup_{s\in\Lambda_{t}}
  |
  \frac{f'(s)}{f(s)}
  |
   \le
   2t\,\frac{1}{t}
   =
   2$.
This proves Claim 1.

\bigskip

\noindent
{\it Claim 2.
We have
$
 |
   \int_{\Delta_{t}}
   (
   \frac{f'(s)}{f(s)}\,
   -\log r_{\min}
   )
   \,ds
   |
  \le
  2$.
 }

\noindent
{\it Proof of Claim 2.} 
Write
$\ell(\Delta_{t})=2t$ for the length of $\Delta_{t}$.
For $s=\sigma+\ima \tau\in\Delta_{t}$ with $\sigma,\tau\in\Bbb R$, we have
$\sigma\le\sigma_{t}^{-}$
whence
$|\frac{f'(s)}{f(s)}-\log r_{\min}|\le cs_{\max}^{-\sigma}\le  \frac{1}{t}$,
and consequently
 $|
   \int_{\Delta_{t}}
   (
   \frac{f'(s)}{f(s)}\,
   -\log r_{\min}
   )
   \,ds
   |
  \le
  \ell(\Delta_{t})
   \,
  \sup_{s\in\Delta_{t}}
  |
   \frac{f'(s)}{f(s)}\,
   -\log r_{\min}
   |
   \le
  2t\,\frac{1}{t}
  =
  2$.
This proves Claim 2.

\bigskip

\noindent
{\it Claim 3.
There are positive  constants $c_{0}$ and $c_{1}$ such that
 $|
  \imag
   \int_{\Sigma^{\pm}_{t}}
   \frac{f'(s)}{f(s)}
   \,ds
   |
  \le
  c_{0}+c_{1}\log t$.
 }

\noindent
{\it Proof of Claim 3.}
Write 
$\rho_{t}
=
\max(\,|B(q)-\sigma_{t}^{-}|,|B(q)-\sigma_{t}^{+}|\,)$
and
$R_{t}=3\rho_{t}$.
Next,
define $F_{t}:\Bbb C\to\Bbb C$ by
$F_{t}(s)
 =
 \frac{f(s+B(q)\pm\ima t)}{f(B(q)\pm\ima t)}$
and note that
$\frac{F_{t}'(\sigma)}{F_{t}(\sigma)}
=
 \frac{f'(\sigma+B(q)\pm\ima t)}{f(\sigma+B(q)\pm\ima t)}$
 for $\sigma\in\Bbb R$.
 Hence
 $$
 \align
 \Bigg|
\imag
  & \int_{\Sigma^{\pm}_{t}}
  \frac{f'(s)}{f(s)}
   \,ds
   \Bigg|\\
&=
\Bigg|
\imag
 \int_{\sigma_{t}^{-}-B(q)}^{\sigma_{t}^{+}-B(q)}
 \frac{f'(\sigma+B(q)\pm\ima t)}{f(\sigma+B(q)\pm\ima t)}
 \,d\sigma
 \Bigg|\\
 &=
 \Bigg|
 \imag
 \int_{\sigma_{t}^{-}-B(q)}^{\sigma_{t}^{+}-B(q)}
 \frac{F_{t}'(\sigma)}{F_{t}(\sigma)}
  \,d\sigma
  \Bigg|\\ 
 &\le
 \left|
 \imag
 \int_{\sigma_{t}^{-}-B(q)}^{\sigma_{t}^{+}-B(q)}
 \left(
 \frac{F_{t}'(\sigma)}{F_{t}(\sigma)}
 -
 \sum
 \Sb
 \omega\in Z(F_{t})\\
 |\omega|\le\frac{1}{2}R_{t}
 \endSb
 \frac{1}{\sigma-\omega}
 \right)
  \,d\sigma   
 +
  \imag
 \int_{\sigma_{t}^{-}-B(q)}^{\sigma_{t}^{+}-B(q)}
 \sum
 \Sb
 \omega\in Z(F_{t})\\
 |\omega|\le\frac{1}{2}R_{t}
 \endSb
 \frac{1}{\sigma-\omega}
  \,d\sigma    
   \right|\\
&\le
 \int_{\sigma_{t}^{-}-B(q)}^{\sigma_{t}^{+}-B(q)}
 \left|
 \frac{F_{t}'(\sigma)}{F_{t}(\sigma)}
 -
 \sum
 \Sb
 \omega\in Z(F_{t})\\
 |\omega|\le\frac{1}{2}R_{t}
 \endSb
 \frac{1}{\sigma-\omega}
 \right|
  \,d\sigma   
 +
  \left|
  \imag
 \int_{\sigma_{t}^{-}-B(q)}^{\sigma_{t}^{+}-B(q)}
 \sum
 \Sb
 \omega\in Z(F_{t})\\
  |\omega|\le\frac{1}{2}R_{t}
 \endSb
 \frac{1}{\sigma-\omega}
  \,d\sigma    
   \right|\\
&\le
(\sigma_{t}^{+}-\sigma_{t}^{-})
\,
\sup_{\sigma\in[\sigma_{t}^{-}-B(q),\sigma_{t}^{+}-B(q)]}
 \left|
 \frac{F_{t}'(\sigma)}{F_{t}(\sigma)}
 -
 \sum
 \Sb
 \omega\in Z(F_{t})\\
 |\omega|\le\frac{1}{2}R_{t}
 \endSb
 \frac{1}{\sigma-\omega}
 \right|\\
&\qquad\qquad
 \qquad\qquad
  \qquad\qquad
  \qquad\qquad
  \qquad\qquad
 +
  \sum
 \Sb
 \omega\in Z(F_{t})\\
 |\omega|\le\frac{1}{2}R_{t}
 \endSb
  \left|
  \imag
  \int_{\sigma_{t}^{-}-B(q)}^{\sigma_{t}^{+}-B(q)}
 \frac{1}{\sigma-\omega}
  \,d\sigma    
   \right|\,.
  \endalign
    $$
 Using the fact that
 $[\sigma_{t}^{-}-B(q),\sigma_{t}^{+}-B(q)]
 \subseteq
 [-\rho_{t},\rho_{t}]$
and  
$\sigma_{t}^{+}-\sigma_{t}^{-}\le2\rho_{t}\le R_{t}$,
we therefore conclude that
$$
 \align
 \Bigg|
\imag
  & \int_{\Sigma^{\pm}_{t}}
  \frac{f'(s)}{f(s)}
   \,ds
   \Bigg|\\
&\le
R_{t}
\,
\sup_{\sigma\in [-\rho_{t},\rho_{t}]}
 \left|
 \frac{F_{t}'(\sigma)}{F_{t}(\sigma)}
 -
 \sum
 \Sb
 \omega\in Z(F_{t})\\
 |\omega|\le\frac{1}{2}R_{t}
 \endSb
 \frac{1}{\sigma-\omega}
 \right|
  +
  \sum
 \Sb
 \omega\in Z(F_{t})\\
 |\omega|\le\frac{1}{2}R_{t}
 \endSb
  \left|
  \imag
  \int_{\sigma_{t}^{-}-B(q)}^{\sigma_{t}^{+}-B(q)}
 \frac{1}{\sigma-\omega}
  \,d\sigma    
   \right|\,.\\
&{}   
   \tag13.28
  \endalign
    $$
We now note that it follows from Proposition 13.6 that
 $$
 \align
\sup_{\sigma\in [-\rho_{t},\rho_{t}]}
 \left|
 \frac{F_{t}'(\sigma)}{F_{t}(\sigma)}
 -
 \sum
 \Sb
 \omega\in Z(F_{t})\\
 |\omega|\le\frac{1}{2}R_{t}
 \endSb
 \frac{1}{\sigma-\omega}
 \right|
 & \le
  \frac{4R_{t}}{(R_{t}-2\rho_{t})^{2}}
 \log
 M_{F_{t}}(R_{t})\\
 &=
   \frac{9}{R_{t}}
 \log
 M_{F_{t}}(R_{t})\,.
 \tag13.29
 \endalign
$$
Combining (13.28) and (13.29) 
gives
 $$
 \Bigg|
 \imag
   \int_{\Sigma^{\pm}_{t}}
   \frac{f'(s)}{f(s)}
   \,ds
   \Bigg|
 \le
 9
 \log
 M_{F_{t}}(R_{t}) +
  \sum
 \Sb
 \omega\in Z(F_{t})\\
 |\omega|\le\frac{1}{2}R_{t}
 \endSb
  \left|
  \imag
  \int_{\sigma_{t}^{-}-B(q)}^{\sigma_{t}^{+}-B(q)}
 \frac{1}{\sigma-\omega}
  \,d\sigma    
   \right|\,.
   \tag13.30
    $$
However, a simple computation shows that
if $\omega\in\Bbb C$ with $\imag\omega\not=0$, then
$\imag
  \int_{\sigma_{t}^{-}-B(q)}^{\sigma_{t}^{+}-B(q)}
 \frac{1}{\sigma-\omega}
  \,d\sigma   
  =
  \int_{\frac{\sigma_{t}^{-}-B(q)-\real\omega}{\imag\omega}}
  ^{\frac{\sigma_{t}^{+}-B(q)-\real\omega}{\imag\omega}}
  \frac{1}{1+x^{2}}\,dx$,
  and so
  $|
  \imag
  \int_{\sigma_{t}^{-}-B(q)}^{\sigma_{t}^{+}-B(q)}
 \frac{1}{\sigma-\omega}
  \,d\sigma   
|
\le\pi$.
We therefore conclude from (13.30) that
 $$
 \Bigg|
 \imag
   \int_{\Sigma^{\pm}_{t}}
   \frac{f'(s)}{f(s)}
   \,ds
   \Bigg|
\le
 9
 \log
 M_{F_{t}}(R_{t})
  +
  \pi
  \,
   \Big|
 \,
 \big\{\omega\in Z(F_{t})
 \,\big|\,
 |\omega|\le \tfrac{1}{2}R_{t}
 \big\}
 \,
 \Big|\,.
 \tag13.31
    $$
Now    
an application of Proposition 13.7 shows that
$|
 \,
 \{\omega\in Z(F_{t})
 \,|\,
 |\omega|\le \tfrac{1}{2}R_{t}
 \}
 \,
 |
 \le
 \frac{1}{\log 2}\log M_{F_{t}}(R_{t})$,
 and
 it therefore follows from (13.31) that
 $$
 \Bigg|
 \imag
   \int_{\Sigma^{\pm}_{t}}
   \frac{f'(s)}{f(s)}
   \,ds
   \Bigg|
 \le
 (9+\tfrac{\pi}{\log 2})
 \log
 M_{F_{t}}(R_{t})\,.
 \tag13.32
  $$

Next, we estimate 
$M_{F_{t}}(R_{t})$.
Indeed, 
since
$|f(B(q)\pm\ima t)|\ge\frac{1}{2}$
(cf\. Lemma 13.4),
we conclude that
 $M_{F_{t}}(R_{t})
 =
 \sup_{|s|=R_{t}}
  \frac{|f(s+B(q)\pm\ima t)|}{|f(B(q)\pm\ima t)|}
  \le
2(1+
  \sup_{\sigma\in[-R_{t},R_{t}]}
  \sum_{i}p_{i}^{q}r_{i}^{\sigma+B(q)})
  =
2+2\sum_{i}p_{i}^{q}r_{i}^{B(q)}r_{i}^{-R_{t}}
\le
2+Cr_{\min}^{-R_{t}}$
  where
  $C=2 \sum_{i}p_{i}^{q}r_{i}^{B(q)}$.
Also,
 $cs_{\max}^{-\sigma_{t}^{-}}=\frac{1}{t}$
 and
 $cr_{\max}^{\sigma_{t}^{+}}=\frac{1}{t}$.
 This clearly implies that
 there are positive constants 
 $c_{0}^{-}$,  $c_{0}^{+}$, $c_{1}^{-}$ and $c_{1}^{+}$
 such that
 $-\sigma_{t}^{-}= c_{0}^{-}+c_{1}^{-}\log t$ and
 $\sigma_{t}^{+}= c_{0}^{+}+c_{1}^{+}\log t$.
 Consequently
 $R_{t}=3\rho_{t}
 =
3\max(\,|B(q)-\sigma_{t}^{-}|,|B(q)-\sigma_{t}^{+}|\,)
\le
3(\sigma_{t}^{+}-\sigma_{t}^{-})
=
3(c_{0}^{+}+c_{0}^{-})+3(c_{1}^{+}+c_{1}^{-})\log t$.
We therefore conclude that
 $M_{F_{t}}(R_{t})
\le
2+Cr_{\min}^{-R_{t}}
  \le
2+Cr_{\min}^{-(3(c_{0}^{+}+c_{0}^{-})+3(c_{1}^{+}+c_{1}^{-})\log t)}
  =
  2+
  C_{0}\,C_{1}^{\log t}$
  where
  $C_{0}
  =
  Cr_{\min}^{-3(c_{0}^{+}+c_{0}^{-}) }$ and 
  $C_{1}=r_{\min}^{-3(c_{1}^{+}+c_{1}^{-}) }$.
 It follows from this and (13.32)
 that 
 $$
 \align
 \Bigg|
 \imag
   \int_{\Sigma^{\pm}_{t}}
   \frac{f'(s)}{f(s)}
   \,ds
   \Bigg|
&\le
 (9+\tfrac{\pi}{\log 2})
 \log (2+C_{0}\,C_{1}^{\log t})\,.
   \endalign
    $$
 This completes the proof of Claim 3.

\bigskip

The proof of the theorem now follows from (13.27), Claim 1, Claim 2 and Claim 3. 
\hfill$\square$

\bigskip
\bigskip


\heading{14. Proof of Theorem 5.4}
\endheading

The purpose of this section is to prove Theorem 5.4.
Recall that for a real number $q$ and $l=0,1,\ldots,d$,  we write
 $$
 \sigma_{q,l}=\sum_{i=1}^{N}p_{i}^{q}r_{i}^{l-dq}\,,
 $$
and 
  $$
 C_{\mu,r}^{q,l,\sym}(K)
 =
   \sum
 \Sb
 \bold i\\
 r_{\bold i}<r<r_{\hat\bold i}
 \endSb
 p_{\bold i}^{q}r_{\bold i}^{l-dq}
 \,\,
 +
 \,\,
 \tfrac{1+\frac{1}{\sigma_{q,l}}}{2}
 \,\,
   \sum
 \Sb
 \bold i\\
 r=r_{\hat\bold i}
 \endSb
 p_{\bold i}^{q}r_{\bold i}^{l-dq}\,.
 $$ 
Also, recall
that 
we define the symbolic $q$ multifractal Minkowski volume of $\mu$ by
 $$
 \align
 V_{\mu,r}^{q,\sym}(K)
&=
 \frac{1}{r^{d}}
 \sum_{l}
 C_{\mu,r}^{q,l,\sym}(K)
\,
r^{(d-l)+dq}\\
&=
 \sum_{l}
 C_{\mu,r}^{q,l,\sym}(K)
\,
r^{-l+dq}\,.
\endalign
$$

The key technique
used for proving Theorem 5.4 is to apply the Mellin transform
to (a suitably rescaled version of) the function
$r\to C_{\mu,r}^{q,l,\sym}(K)$.
The Mellin transform is a general method
for expressing
functions (satisfying various growth conditions)
as 
complex
contour integrals.
The precise statement is given by the Mellin transform theorem below.

\bigskip

\proclaim{Theorem 14.1. The Mellin transform theorem [Pat]}
Let $a,b\in[-\infty,\infty]$ with $a< b$
and
let $f:(0,\infty)\to\Bbb R$ be a real valued function.
Assume that the following conditions are satisfied:
\roster
\item"(i)"
The function $f$ is piecewise continuous on all compact subintervals of 
$(0,\infty)$,
and at  all discontinuity points $x_{0}>0$ of $f$, we have
$f(x_{0})=\frac{\lim_{x\searrow x_{0}}f(x)+\lim_{x\nearrow x_{0}}f(x)}{2}$;
\item"(ii)"
If $s\in\Bbb C$ satisfies
$a<\real(s)<b$,
then
$\int_{0}^{\infty}
 |x^{s-1}f(x)|\,dx
 <
 \infty$.
\endroster
Then we have:
\roster
\item"(1)"
For 
$s\in\Bbb C$ with
$a<\real(s)<b$
the integral
 $$
 \int_{0}^{\infty}
 x^{s-1}f(x)\,dx
 $$
is well-defined. 
\endroster
It follows from (1) that  the
function
$\Mellin f:\{s\in\Bbb C\,|\,a<\real(s)<b\}\to\Bbb C$ given by
 $$
 (\Mellin f)(s)
 =
 \int_{0}^{\infty}
 x^{s-1}f(x)\,dx
 $$
is well-defined. 
The function $\Mellin f$ is called the Mellin transform of $f$. 
\roster
\item"(2)"
For $c\in\Bbb R$ with
$a<c<b$ 
and $x>0$
the integral
 $$
\int_{c-\ima\infty}^{c+\ima\infty}x^{-s}(\Mellin f)(s)\,ds
$$
is well-defined.
\item"(3)"
For $c\in\Bbb R$ with
$a<c<b$
and
$x>0$, we have
$$
f(x)
=
\frac{1}{2\pi\ima}
\int_{c-\ima\infty}^{c+\ima\infty}x^{-s}(\Mellin f)(s)\,ds\,.
$$
\endroster
\endproclaim

\bigskip

\noindent
In order to prove Theorem 5.4,
we apply the Mellin transform theorem to the function
$r\to C_{\mu,r} ^{q,l,\sym}(K)$.
However, 
before applying the 
Mellin transform
it is 
useful
to 
\lq\lq rescale" the function 
$r\to C_{\mu,r} ^{q,l,\sym}(K)$.
This is to ensure
that all points of discontinuity satisfy Condition (i) in the Mellin transform theorem.
In order to do this, we first
define $E:\Bbb R\to\Bbb R$ by
 $$
 E(t)
 =
 \cases
 0
&\quad
 \text{for $t<0$;}\\
 \frac{1}{2}
&\quad
 \text{for $t=0$;}\\
 1
&\quad
 \text{for $0<t$.}
 \endcases
 $$
For $q\in\Bbb R$ and $l=0,1,\ldots,d$, 
we now
define the
\lq\lq rescaled" version 
 $r\to B_{\mu}^{q,l}(r)$ 
 of the function
$r\to C_{\mu,r} ^{q,l,\sym}(K)$
by
 $$
 B_{\mu}^{q,l}(r)
 =
 \sum_{\bold i}p_{\bold i}^{q}
 \bigg(\frac{r_{\bold i}}{r}\bigg)^{l-dq}E(r_{\bold i}-r)\,.
 $$
Proposition 14.2 below 
shows
 that
$B_{\mu}^{q,l}(r)$ is, indeed, a
\lq\lq rescaled" version of $C_{\mu,r} ^{q,l,\sym}(K)$.

\bigskip

\proclaim{Proposition 14.2}
Fix $q\in\Bbb R$
and $l=0,1,\ldots,d$.
For $0<r<r_{\min}$, we have
 $$
 C_{\mu,r}^{q,l,\sym}(K)r^{-l+dq} 
=
 \sigma_{q,l}r^{-l+dq}
 +
 (\sigma_{q,l}-1)
 B_{\mu}^{q,l}(r)\,.
 $$
\endproclaim
\noindent{\it  Proof}\newline
\noindent 
This follows from a straight forward
but lengthy calculation which we 
are omitting.
\hfill$\square$

 \bigskip

\noindent
Next, we 
apply the Mellin Transform Theorem
to the
\lq\lq rescaled" function
$B_{\mu}^{q,l}$.

\bigskip

\proclaim{Proposition 14.3}
Fix $q\in\Bbb R$ and $l=0,1,\ldots,d$.
Write
$ \Bbb H^{q,l}
 =
 \{
 s\in\Bbb C
 \,|\,
 \real(s)>\max\big(\,l-dq,\beta(q)\,\big)
 \}$.
\roster
\item"(1)"
The function $B_{\mu}^{q,l}$ is piecewise continuous on all compact subintervals of 
$(0,\infty)$,
and at  all discontinuity points $r_{0}>0$ of $B_{\mu}^{q,l}$, we have
$B_{\mu}^{q,l}(r_{0})
=
\frac{\lim_{r\searrow r_{0}}B_{\mu}^{q,l}(r)+\lim_{r\nearrow r_{0}}B_{\mu}^{q,l}(r)}{2}$.
\item"(2)"
For $s\in\Bbb H^{q,l}$, we have
$\int_{0}^{\infty}
 |B_{\mu}^{q,l}(r)r^{s-1}|
 \,
 dr
 <
 \infty$.
\endroster 
It follows from (1) and (2) 
that the Mellin transform
$\Mellin B_{\mu}^{q,l}:\Bbb H^{q,l}\to\Bbb C$
of
$B_{\mu}^{q,l}$
given by
 $$
 (\Mellin B_{\mu}^{q,l})(s)
 =
 \int_{0}^{\infty}
 B_{\mu}^{q,l}(r)r^{s-1}
 \,
 dr
 $$
is well-defined.
\roster
\item"(3)"
For $s\in\Bbb H^{q,l}$,
we have
 $$
 (\Mellin B_{\mu}^{q,l})(s)
 =
 \frac{1}{s-(l-dq)}\zeta_{\mu}^{q,l}(s)\,.
 $$
\item"(4)"
For
$c>\max\big(\,l-dq,\beta(q)\,\big)$
and
$r>0$ the integral
 $$
  \lim_{t\to\infty}
 \int_{c-\ima t}^{c+\ima t}
 (\Mellin B_{\mu}^{q,l})(s)
 \,
 r^{-s}
 \,ds
 $$
is well-defined. 
\item"(5)"
For
$c>\max\big(\,l-dq,\beta(q)\,\big)$
and
$r>0$, we have
 $$
 B_{\mu}^{q,l}(r)
 =
 \lim_{t\to\infty}
 \frac{1}{2\pi\ima}
 \int_{c-\ima t}^{c+\ima t}
 (\Mellin B_{\mu}^{q,l})(s)
 \,
 r^{-s}
 \,ds\,.
 $$
In particular,
for
$c>\max\big(\,l-dq,\beta(q)\,\big)$
and
$r>0$, we have
 $$
 B_{\mu}^{q,l}(r)
 =
 \lim_{t\to\infty}
 \frac{1}{2\pi\ima}
 \int_{c-\ima t}^{c+\ima t}
 \frac{1}{s-(l-dq)}\zeta_{\mu}^{q}(s)
 \,
 r^{-s}
 \,ds\,.
 $$

\endroster 
\endproclaim
\noindent{\it  Proof}\newline
\noindent
\noindent
(1)
This follows immediately from the 
fact that $E$ 
is piecewise continuous on all compact subintervals of 
$(0,\infty)$
and that 
 all discontinuity points $r_{0}>0$ of $E$ satisfy
 $E(r_{0})
=
\frac{\lim_{r\searrow r_{0}}E(r)+\lim_{r\nearrow r_{0}}E(r)}{2}$.

\noindent
(2)
Write
$s=\sigma+\ima t$
with $\sigma,t\in\Bbb R$.
It follows from the definition of
$B_{\mu}^{q,l}(r)$ that
 $$
 \align
 \int_{0}^{\infty}
 |B_{\mu}^{q,l}(r)r^{s-1}|
 \,
 dr
&= 
 \int_{0}^{\infty}
 \Bigg|
  \sum_{\bold i}p_{\bold i}^{q}
 \bigg(\frac{r_{\bold i}}{r}\bigg)^{l-dq}E(r_{\bold i}-r)r^{s-1}
 \Bigg|
 \,
 dr\\
&\le
 \int_{0}^{\infty}
  \sum_{\bold i}p_{\bold i}^{q}
 \bigg(\frac{r_{\bold i}}{r}\bigg)^{l-dq}E(r_{\bold i}-r)\,|r^{s-1}|
 \,
 dr\\
&=
 \int_{0}^{\infty}
  \sum_{\bold i}p_{\bold i}^{q}
 \bigg(\frac{r_{\bold i}}{r}\bigg)^{l-dq}E(r_{\bold i}-r)\,r^{\sigma-1}
 \,
 dr\\
&=
  \sum_{\bold i}
  \int_{0}^{\infty}
 p_{\bold i}^{q}
 \bigg(\frac{r_{\bold i}}{r}\bigg)^{l-dq}E(r_{\bold i}-r)\,r^{\sigma-1}
 \,
 dr\\
&=
  \sum_{\bold i}
   p_{\bold i}^{q}
  \int_{0}^{r_{\bold i}}
 \bigg(\frac{r_{\bold i}}{r}\bigg)^{l-dq}\,r^{\sigma-1}
 \,
 dr\,.
 \tag14.1
 \endalign
 $$
Since $\sigma=\real(s)>l-dq$, we deduce 
that
$\int_{0}^{r_{\bold i}}
 (\frac{r_{\bold i}}{r})^{l-dq}\,r^{\sigma-1}\,dr
 =
 \frac{1}{\sigma-(l-dq)}
 r_{\bold i}^{\sigma}$, 
 and we therefore conclude from (14.1)
 that
  $$
 \align
\int_{0}^{\infty}
 |B_{\mu}^{q,l}(r)r^{s-1}|
 \,
 dr
&\le
 \frac{1}{\sigma-(l-dq)}
  \sum_{\bold i}
   p_{\bold i}^{q}r_{\bold i}^{\sigma}\\
&=
 \frac{1}{\sigma-(l-dq)}
 \sum_{k=1}^{\infty}
 \sum
 \Sb
  |\bold i|=k
 \endSb
 p_{\bold i}^{q}
 r_{\bold i}^{\sigma}\\
&=
 \frac{1}{\sigma-(l-dq)}
 \sum_{k=1}^{\infty}
 \Bigg(
 \sum_{i}
 p_{i}^{q}
 r_{i}^{\sigma}
 \Bigg)^{k}
 \,.
 \tag14.2
 \endalign
 $$
Finally, since $\sigma=\real(s)>\beta(q)$, we conclude that
$\sum_{i}
 p_{i}^{q}
 r_{i}^{\sigma}<1$, whence
$\sum_{k=1}^{\infty}
 (
 \sum_{i}
 p_{i}^{q}
 r_{i}^{\sigma}
 )^{k} 
 <
 \infty$,
and we therefore deduce from
the previous
inequality 
(14.2) that
$\int_{0}^{\infty}
 |B_{\mu}^{q,l}(r)r^{s-1}|
 \,
 dr
\le
 \frac{1}{\sigma-(l-dq)}
 \sum_{k=1}^{\infty}
 (
 \sum_{i}
 p_{i}^{q}
 r_{i}^{\sigma}
 )^{k}
 <
 \infty$.

\noindent
(3)
Since the series
$B_{\mu}^{q,l}(r)r^{s-1}
 =
 \sum_{\bold i}p_{\bold i}^{q}
 (\frac{r_{\bold i}}{r})^{l-dq}E(r_{\bold i}-r)r^{s-1}$
 only has finitely many non-zero terms we immediately conclude that
 $$
 \align
 (\Mellin B_{\mu}^{q,l})(s)
&=
 \int_{0}^{\infty}
 B_{\mu}^{q,l}(r)r^{s-1}
 \,dr\\
&=
 \int_{0}^{\infty}
  \sum_{\bold i}p_{\bold i}^{q}
 \bigg(\frac{r_{\bold i}}{r}\bigg)^{l-dq}E(r_{\bold i}-r)\,r^{s-1}
 \,dr\\
&=
  \sum_{\bold i}p_{\bold i}^{q}
  \int_{0}^{\infty}
   \bigg(\frac{r_{\bold i}}{r}\bigg)^{l-dq}E(r_{\bold i}-r)\,r^{s-1}
 \,dr\\
&=
   \sum_{\bold i}p_{\bold i}^{q}
   \int_{0}^{r_{\bold i}}
 \bigg(\frac{r_{\bold i}}{r}\bigg)^{l-dq}\,r^{s-1}
 \,dr\,. 
 \tag14.3
  \endalign
 $$
As $\real(s)>l-dq$, we deduce 
that
$\int_{0}^{r_{\bold i}}
 (\frac{r_{\bold i}}{r})^{l-dq}\,r^{s-1}\,dr
 =
 \frac{1}{s-(l-dq)}
 r_{\bold i}^{s}$, 
 and
it now follows  from
 (14.3) that
  $$
 \align
 (\Mellin B_{\mu}^{q,l})(s)
&=
 \frac{1}{s-(l-dq)}
 \sum_{\bold i}
 p_{\bold i}^{q}
 r_{\bold i}^{s}\\
&=
 \frac{1}{s-(l-dq)}
 \zeta_{\mu}^{q}(s)\,.
 \endalign
 $$
 
\noindent
(4) 
This follows immediately from the Mellin Transform Theorem.

\noindent
(5) 
This follows immediately from the Mellin Transform Theorem and (3).
\hfill$\square$

\bigskip

 \noindent
 We can now prove Theorem 5.4.

 \noindent{\it  Proof of Theorem 5.4}\newline
Fix
$c>\max\big(\,-dq,1-dq,\dots,d-dq,\beta(q)\,\big)$
and
$0<r<r_{\min}$.
Using Proposition 14.2 and Proposition 14.3, we have
 $$
 \align
 V_{\mu,r}^{q,\sym}(K)
&=
 \sum_{l}
  \kappa_{\mu}^{q,l}(K)
 \,
 C_{\mu,r}^{q,l,\sym}(K)
\,
r^{-l+dq}\\
&=
 \sum_{l}
  \kappa_{\mu}^{q,l}(K)
 \,
 \Bigg(
  \sigma_{q,l}r^{-l+dq}
 \,\,
 +
 \,\,
 (\sigma_{q,l}-1)
 \,
 \frac{1}{2\pi\ima}
 \int_{c-\ima \infty}^{c+\ima \infty}
 \frac{1}{s-(l-dq)}\zeta_{\mu}^{q}(s)
 \,
 r^{-s}
 \,ds
 \Bigg)\\
&=
 \sum_{l}
  \kappa_{\mu}^{q,l}(K)
  \,
 \sigma_{q,l}r^{-l+dq}
  \,\,
 +
 \,\,
  \frac{1}{2\pi\ima}
 \int_{c-\ima \infty}^{c+\ima \infty}
  \Bigg(
  \sum_{l}
 \frac{\kappa_{\mu}^{q,l}(K)\, (\sigma_{q,l}-1)  }{s-(l-dq)}
 \Bigg)
 \,
 \zeta_{\mu}^{q}(s)
 \,
 r^{-s}
 \,ds\\
&=
 \sum_{l}
  \kappa_{\mu}^{q,l}(K)
  \,
 \sigma_{q,l}r^{-l+dq}
  \,\,
 +
 \,\,
  \frac{1}{2\pi\ima}
 \int_{c-\ima \infty}^{c+\ima \infty}
Z_{\mu}^{q}(s)
 \,
 r^{-s}
 \,ds\,.
  \endalign
 $$
This completes the proof.
\hfill$\square$

\bigskip
\bigskip


\heading{15. 
Proof of Theorem 5.5}
\endheading

The purpose of this section is to prove Theorem 5.5.
We first use the estimates from Section 11together
with
the residue theorem to compute the 
complex
contour integral
$ \frac{1}{2\pi\ima}\int_{c-\ima \infty}^{c+\ima \infty}
Z_{\mu}^{q}(s)
 \,
 r^{-s}
 \,ds$
 appearing in Theorem 5.4.

\bigskip

\proclaim{Proposition 15.1}
Fix $q\in\Bbb R$
and
$c>\max\big(\,-dq,1-dq,\ldots,d-dq,\beta(q)\,\big)$.
Let $(t_{q,n})_{n}$ be the sequence from Theorem 13.5.
For all $0<r< 1$, we have
 $$
 \align
 \frac{1}{2\pi\ima}
&\int_{c-\ima \infty}^{c+\ima \infty}
Z_{\mu}^{q}(s)
 \,
 r^{-s}
 \,ds
  =
 \lim_{n}	
 \sum
 \Sb
 \omega
 \in P(
 s\to
 Z_{\mu}^{q}(s)
 \,
 r^{-s}
 )\\
 {}\\
 |\imag(\omega)|\le t_{q,n}
 \endSb
 \res
 \Big(
 s\to
 Z_{\mu}^{q}(s)
 \,
 r^{-s};
 \omega
 \Big)\,.
 \endalign
 $$
 \endproclaim
\noindent{\it  Proof}\newline
\noindent
Fix $l=0,1,\ldots,d$.
It clearly suffices to prove that
for all $0<r<r_{\min}$, we have
 $$
 \align
 \frac{1}{2\pi\ima}
&\int_{c-\ima \infty}^{c+\ima \infty}
 \frac{1}{s-(l-dq)}\zeta_{\mu}^{q}(s)
 \,
 r^{-s}
 \,ds\\
&\qquad
 =
 \lim_{n}	
 \sum
 \Sb
 \omega
 \in P(
 s\to
 \frac{1}{s-(l-dq)}\zeta_{\mu}^{q}(s)
 \,
 r^{-s}
 )\\
 {}\\
 |\imag(\omega)|\le t_{q,n}
 \endSb
 \res
 \Big(
 s\to
 \frac{1}{s-(l-dq)}\zeta_{\mu}^{q}(s)
 \,
 r^{-s};
 \omega
 \Big)\,.
 \endalign
 $$

Let $A(q)$ be the constant in Lemma 13.4. 
Note that it follows from 
Lemma 13.4 that there is a constant $k$ such that
if $\sigma\le A(q)$ and $t\in\Bbb R$, then
 $$
|\zeta_{\mu}^{q}(\sigma+\ima t)|\le k\,.
$$
Also, note that it follows from Theorem 13.5 that there
is a constant $k_{c}$ such that
if $\sigma\le c$ and $n\in\Bbb N$, then
 $$
|\zeta_{\mu}^{q}(\sigma+\ima t_{q,n})|\le k_{c}\,.
$$
Next,
for all positive integers $n$ and $m$ with $-m\le A(q)$, we define paths
$\Sigma_{n,m}^{+},\Sigma_{n,m}^{-},\Lambda_{n},\Delta_{n,m}$ in $\Bbb C$
by
 $$
 \align
 \Sigma_{n,m}^{+}
&\quad
 \text{
 is the directed line-segment from
 \quad\,
 $c+\ima t_{q,n}$
 to
 $-m+\ima t_{q,n}$
 }\,,\\
 \Sigma_{n,m}^{-}
&\quad
 \text{
 is the directed line-segment from
 $-m-\ima t_{q,n}$
 to
 \quad\,
 $c-\ima t_{q,n}$
 }\,,\\
  \Lambda_{n}\,\,
&\quad
 \text{
 is the directed line-segment from
 \quad\,
  $c-\ima t_{q,n}$
 to
  \quad\,
 $c-\ima t_{q,n}$
 }\,,\\
 \Delta_{n,m}
&\quad
 \text{
 is the directed line-segment from
  $-m+\ima t_{q,n}$
 to
 $-m-\ima t_{q,n}$
 }\,.
 \endalign
 $$
Below we sketh the paths
$U_{n,m}$,
$L_{n,m}$,
$\Gamma_{n}$ and $\Lambda_{n,m}$.

\goodbreak

\midinsert


 \vspace{130mm}
 \centerline{\hbox{\hskip -50mm\special{pdf=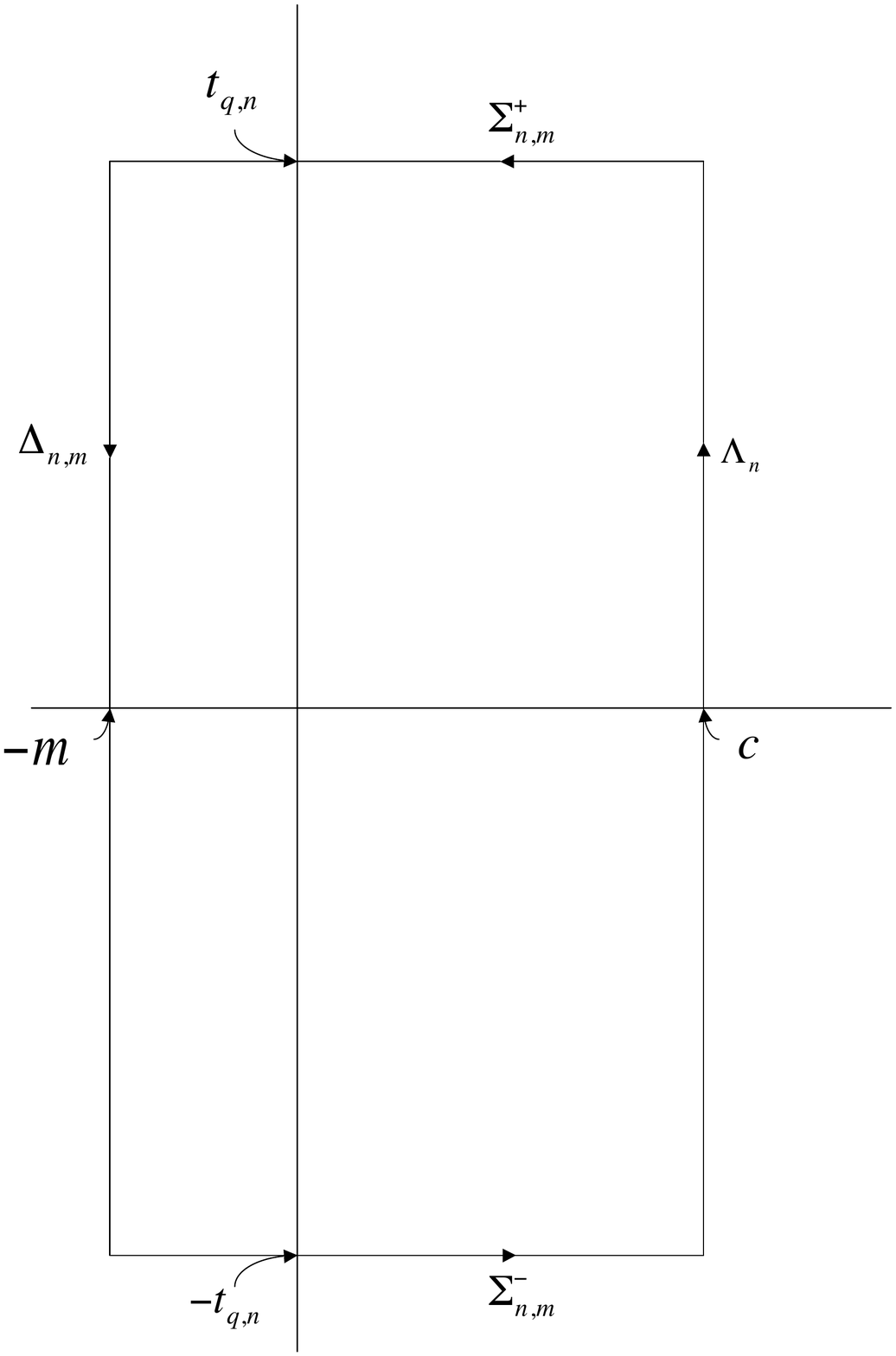
 scale=0.38}}}

\botcaption{\eightpoint \bf Fig\. 15.1}
\eightpoint 
The paths
$U_{n,m}$,
$L_{n,m}$,
$\Gamma_{n}$ and $\Lambda_{n,m}$.
\endcaption
\endinsert

\noindent
Since
the paths
$\Sigma_{n,m}^{+}$, 
$\Sigma_{n,m}^{-}$,
$\Lambda_{n}$ and
$\Delta_{n,m}$ 
enclose the region
$\{s\in\Bbb C\,|\,
-m
 <
 \real(s)
 <
 c\,,\, 
-t_{q,n}
 <
 \imag(s)
 <
 t_{q,n}\}$
and since all poles $\omega$ of
the function
$s\to
 \frac{1}{s-(l-dq)}\zeta_{\mu}^{q}(s)
 \,
 r^{-s}$
satisfy
$\real(\omega)
\in
[\alpha(q),\max(\,l-dq,\beta(q)\,)]
\subseteq
(-m,c)$,
it now follows from the 
residue theorem that
 $$
 \align
 2\pi\ima
&\sum
 \Sb
 \omega
 \in P(
 s\to
 \frac{1}{s-(l-dq)}\zeta_{\mu}^{q}(s)
 \,
 r^{-s}
 )\\
 -t_{q,n}
 <
 \imag(\omega)
 <
 t_{q,n} 
 \endSb
\res
 \Big(
 s\to
 \frac{1}{s-(l-dq)}\zeta_{\mu}^{q}(s)
 \,
 r^{-s};
 \omega
 \Big)\\
&\\ 
&\qquad\qquad
=
 2\pi\ima
 \sum
 \Sb
 \omega
 \in P(
 s\to
 \frac{1}{s-(l-dq)}\zeta_{\mu}^{q}(s)
 \,
 r^{-s}
 )\\
 -m
 <
 \real(\omega)
 <
 c\,\,\\
 -t_{q,n}
 <
 \imag(\omega)
 <
 t_{q,n} 
 \endSb
 \res
 \Big(
 s\to
 \frac{1}{s-(l-dq)}\zeta_{\mu}^{q}(s)
 \,
 r^{-s};
 \omega
 \Big)\\
&\\ 
&\qquad\qquad
= 
 \int_{c-\ima t_{q,n}}^{c+\ima t_{q,n}}
 \frac{1}{s-(l-dq)}\zeta_{\mu}^{q}(s)
 \,
 r^{-s}
 \,ds\\
&\qquad\qquad
\qquad
 +
 \,
 \int_{\Sigma_{n,m}^{+}}
 \frac{1}{s-(l-dq)}\zeta_{\mu}^{q}(s)
 \,
 r^{-s}
 \,ds\\
&\qquad\qquad
\qquad
 +
 \,
 \int_{\Delta_{n,m}}
 \frac{1}{s-(l-dq)}\zeta_{\mu}^{q}(s)
 \,
 r^{-s}
 \,ds\\
&\qquad\qquad
\qquad
 +
 \,
  \int_{\Sigma_{n,m}^{-}}
 \frac{1}{s-(l-dq)}\zeta_{\mu}^{q}(s)
 \,
 r^{-s}
 \,ds\,.
 \tag15.1
 \endalign
 $$

Using the estimates from Section 11, we will 
now provide estimates for 
the three integrals 
$ \int_{\Sigma_{n,m}^{+}}
 \frac{1}{s-(l-dq)}\zeta_{\mu}^{q}(s)
 \,
 r^{-s}
 \,ds$,
 $\int_{\Gamma_{n,m}^{+}}
 \frac{1}{s-(l-dq)}\zeta_{\mu}^{q}(s)
 \,
 r^{-s}
 \,ds$ and
 $\int_{\Sigma_{n,m}^{-}}
 \frac{1}{s-(l-dq)}\zeta_{\mu}^{q}(s)
 \,
 r^{-s}
 \,ds$
in (15.1)
 and show that they
 tend to $0$
by first fixing $n$ and letting $m\to\infty$,
and then letting $n\to\infty$.
Define $f_{n,m},f_{n},g_{n,m}:(0,1)\to\Bbb R$ by
 $$
 \align
 f_{n,m}(r)
&=
 \frac{k_{c}\,(r^{-c}-r^{m})}{-\log r}
 \,
 \frac{1}{t_{q,n}}\,,\\
 f_{n}(r)
&=
 \frac{k_{c}\,r^{-c}}{-\log r}
 \,
 \frac{1}{t_{q,n}}\,,\\
 g_{n,m}(r)
&=
 2k\,r^{m}
 \Big(
 \log\Big(\sqrt{t_{q,n}^{2}+(m+(l-dq))^{2}}+t_{q,n}\Big)-\log|m+(l-dq)|
 \Big)\,.
 \endalign
 $$
Below we estimate the integrals
$ \int_{\Sigma_{n,m}^{+}}
 \frac{1}{s-(l-dq)}\zeta_{\mu}^{q}(s)
 \,
 r^{-s}
 \,ds$,
 $\int_{\Gamma_{n,m}^{+}}
 \frac{1}{s-(l-dq)}\zeta_{\mu}^{q}(s)
 \,
 r^{-s}
 \,ds$ and
 $\int_{\Sigma_{n,m}^{-}}
 \frac{1}{s-(l-dq)}\zeta_{\mu}^{q}(s)
 \,
 r^{-s}
 \,ds$
using the functions $f_{n,m}$, $f_{n}$ and $g_{n,m}$.

\bigskip

\noindent{\it 
Claim 1.
For $0<r<1$, we have
 $$
 \align
 \Bigg|
 \int_{\Sigma_{n,m}^{\pm}}
 \frac{1}{s-(l-dq)}\zeta_{\mu}^{q}(s)
 \,
 r^{-s}
 \,ds
 \Bigg|
&\le
 f_{n,m}(r)\,.\\
 \endalign
 $$
}
\noindent{\it Proof of Claim 1.}
We have,
using Theorem 13.5,
 $$
 \align
 \Bigg|
 \int_{\Sigma_{n,m}^{\pm}}
&\frac{1}{s-(l-dq)}\zeta_{\mu}^{q}(s)
 \,
 r^{-s}
 \,ds
 \Bigg|\\
&\le
 \int_{-m}^{c}
 \frac{1}{|-\sigma-m+c\pm\ima t_{q,n}-(l-dq)|}
 \,
 |\zeta_{\mu}^{q}(-m-\sigma+c\pm\ima t_{q,n})|
 \,
 |r^{\sigma+m-c\mp\ima t_{q,n}}|
 \,d\sigma\\
&\le
 \int_{-m}^{c}
 \frac{1}{t_{q,n}}
 \,
 k_{c}
 \,
 r^{m+\sigma-c}
 \,d\sigma\\ 
&=
 f_{n,m}(r)\,.
 \endalign
 $$
This completes the proof of Claim 1.

\bigskip

\noindent{\it 
Claim 2.
For $0<r<1$, we have
 $$
 \Bigg|
 \int_{\Delta_{n,m}}
 \frac{1}{s}\zeta_{\mu}^{q}(s)
 \,
 r^{-s}
 \,ds
 \Bigg|
 \le
 g_{n,m}(r)\,.
 $$
}
\noindent{\it Proof of Claim 2.}
We have,
using Lemma 13.4,
 $$
 \align
 \Bigg|
 \int_{\Delta_{n,m}}
 \frac{1}{s-(l-dq)}\zeta_{\mu}^{q}(s)
 \,
 r^{-s}
 \,ds
 \Bigg|
&\le
 \int_{-m}^{c}
 \frac{1}{|-m-\ima t-(l-dq)|}
 \,
 |\zeta_{\mu}^{q}(-m-\ima t)|
 \,
 |r^{m+\ima t}|
 \,dt\\ 
&\le
 \int_{-m}^{c}
 \frac{1}{\sqrt{(m+(l-dq))^{2}+t^{2}}}
 \,
 k
 \,
 r^{m}
 \,dt\\ 
&=
 g_{n,m}(r)\,.
 \endalign
 $$
Ths completes the proof of Claim 2.

\bigskip

\noindent
Finally, we prove the following claim.

\bigskip

\noindent{\it 
Claim 3.
For $0<r<1$, we have
 $$
 2\pi\ima
 \sum
 \Sb
 \omega
 \in P(
 s\to
 \frac{1}{s-(l-dq)}\zeta_{\mu}^{q}(s)
 \,
 r^{-s}
 )\\
 -t_{q,n}
 <
 \imag(\omega)
 <
 t_{q,n} 
 \endSb
 \res
 \Big(
 s\to
 \frac{1}{s-(l-dq)}\zeta_{\mu}^{q}(s)
 \,
 r^{-s};
 \omega
 \Big)
 -
 \int_{c-\ima t_{q,n}}^{c+\ima t_{q,n}}
 \frac{1}{s-(l-dq)}\zeta_{\mu}^{q}(s)
 \,
 r^{-s}
 \,ds
 \to 
 0
 $$
}
\noindent{\it Proof of Claim 3.} 
Let $\varepsilon>0$.
Next, note that for each fixed $n\in\Bbb N$, we
have
$f_{n,m}(r)\to f_{n}(r)$
as $m\to\infty$
and
$g_{n,m}(r)\to 0$
as $m\to\infty$
(since $r<1$).
For each fixed $n\in\Bbb N$, we
can therefore choose a positive integer $M_{n}$ such that
if $m\ge M_{n}$, then
 $$
 \align
 f_{n,m}(r)
&\le
 2f_{n}(r)\,,\\
 g_{n,m}(r)
&\le
 \tfrac{\varepsilon}{3}\,.
 \endalign
 $$
Also, since
$f_{n}(r)\to 0$
as $n\to\infty$, we can choose a positive integer
$N_{0}$
such that if $n\ge N_{0}$, then
 $$
 f_{n}(r)\le\tfrac{\varepsilon}{6}\,.
 $$
Fix $n\ge N_{0}$.
Using (15.1), Claim 1 and Claim 2, we now conclude that
 $$
 \align
&\left|
 2\pi\ima
 \sum
 \Sb
 \omega
 \in P(
 s\to
 \frac{1}{s-(l-dq)}\zeta_{\mu}^{q}(s)
 \,
 r^{-s}
 )\\
 -t_{q,n}
 <
 \imag(\omega)
 <
 t_{q,n} 
 \endSb
 \res
 \Big(
 s\to
 \frac{1}{s-(l-dq)}\zeta_{\mu}^{q}(s)
 \,
 r^{-s};
 \omega
 \Big)
 -
 \int_{c-\ima t_{q,n}}^{c+\ima t_{q,n}}
 \frac{1}{s-(l-dq)}\zeta_{\mu}^{q}(s)
 \,
 r^{-s}
 \,ds
 \right|\\
&\qquad\qquad
 \qquad\qquad
 \qquad\qquad
 \qquad\qquad
 \qquad\qquad
 \qquad\qquad
 \le
 2f_{n,M_{n}}(r)
 +
 g_{n,M_{n}}(r)\\
&\qquad\qquad
 \qquad\qquad
 \qquad\qquad
 \qquad\qquad
 \qquad\qquad
 \qquad\qquad
 \le
 4f_{n}(r)
 +
 g_{n,M_{n}}(r)\\
&\qquad\qquad
 \qquad\qquad
 \qquad\qquad
 \qquad\qquad
 \qquad\qquad
 \qquad\qquad
 \le
 4\tfrac{\varepsilon}{6}
 +
 \tfrac{\varepsilon}{3}\\
&\qquad\qquad
 \qquad\qquad
 \qquad\qquad
 \qquad\qquad
 \qquad\qquad
 \qquad\qquad
 = 
 \varepsilon\,.
 \endalign
 $$
This proves Claim 3.

\bigskip

Finally, we deduce from Claim 3 that
for $0<r<1$, we have
 $$
 \align
 2\pi\ima
&\sum
 \Sb
 \omega
 \in P(
 s\to
 \frac{1}{s-(l-dq)}\zeta_{\mu}^{q}(s)
 \,
 r^{-s}
 )\\
 -t_{q,n}
 <
 \imag(\omega)
 <
 t_{q,n} 
 \endSb
 \res
 \Big(
 s\to
 \frac{1}{s-(l-dq)}\zeta_{\mu}^{q}(s)
 \,
 r^{-s};
 \omega
 \Big)\\
&=
 \left(
 \frac{1}{2\pi\ima}
 \sum
 \Sb
 \omega
 \in P(
 s\to
 \frac{1}{s-(l-dq)}\zeta_{\mu}^{q}(s)
 \,
 r^{-s}
 )\\
 -t_{q,n}
 <
 \imag(\omega)
 <
 t_{q,n} 
 \endSb
 \res
 \Big(
 s\to
 \frac{1}{s-(l-dq)}\zeta_{\mu}^{q}(s)
 \,
 r^{-s};
 \omega
 \Big)
 -
 \int_{c-\ima t_{q,n}}^{c+\ima t_{q,n}}
 \frac{1}{s-(l-dq)}\zeta_{\mu}^{q}(s)
 \,
 r^{-s}
 \,ds
 \right)\\
&\qquad\qquad
 \qquad\qquad
 \qquad\qquad
 \qquad\qquad
 +
 \int_{c-\ima t_{q,n}}^{c+\ima t_{q,n}}
 \frac{1}{s-(l-dq)}\zeta_{\mu}^{q}(s)
 \,
 r^{-s}
 \,ds\\
&\to
 0
 +
 \lim_{t\to\infty}\int_{c-\ima t}^{c+\ima t}
 \frac{1}{s-(l-dq)}\zeta_{\mu}^{q}(s)
 \,
 r^{-s}
 \,ds\\
&=
 \int_{c-\ima \infty}^{c+\ima \infty}
 \frac{1}{s-(l-dq)}\zeta_{\mu}^{q}(s)
 \,
 r^{-s}
 \,ds\,.
 \endalign
 $$
This completes the proof. 
\hfill$\square$

 \bigskip
 
 \noindent
 We can now prove the
  second explicit formula for $V_{\mu,r}^{q,\sym}(K)$, i\.e\. Theorem 5.5.

 \bigskip
 
\noindent{\it  Proof of Theorem 5.5}\newline
\noindent
It follows from Theorem 5.4 and Proposition 15.1
that
if $0<r< r_{\min}$, then
 $$
 \align
 V_{\mu,r}^{q,\sym}(K)
&=
 \sum_{l}
  \kappa_{\mu}^{q,l}(K)
  \,
 \sigma_{q,l}
 \,
 r^{-l+dq}
 \,\,
 +
 \,\,
 \frac{1}{2\pi\ima}
\int_{c-\ima \infty}^{c+\ima \infty}
Z_{\mu}^{q}(s)
 \,
 r^{-s}
 \,ds\\
&=
 \sum_{l}
  \kappa_{\mu}^{q,l}(K)
  \,
 \sigma_{q,l}
 \,
 r^{-l+dq}
 \,\,
 +
 \,\,
  \lim_{n}	
 \sum
 \Sb
 \omega
 \in P(
 s\to
 Z_{\mu}^{q}(s)
 \,
 r^{-s}
 )\\
 {}\\
 |\imag(\omega)|\le t_{q,n}
 \endSb
 \res
 \Big(
 s\to
 Z_{\mu}^{q}(s)
 \,
 r^{-s};
 \omega
 \Big)\,.\\
&{} 
 \tag15.2
 \endalign
 $$
 Also, if 
 $ \beta(q)\not=l-dq$, then
 clearly
 $P(
 s\to
 \frac{1}{s-(l-dq)}\zeta_{\mu}^{q}(s)
 \,
 r^{-s}
 )
 =
 P(\zeta_{\mu}^{q})\cup\{l-dq\}$
and
$P(\zeta_{\mu}^{q})\cap\{l-dq\}=\varnothing$.
It follows from this observation that
(15.2) can  be written as
 $$
 \align
 V_{\mu,r}^{q,\sym}(K)
&=
 \sum_{l}
  \kappa_{\mu}^{q,l}(K)
  \,
 \sigma_{q,l}
 \,
 r^{-l+dq}
 \,\,
 +
 \,\,
 \sum_{l}
 \res
 \Big(
 \,
 s 
 \to
  Z_{\mu}^{q}(s)
 \,
 r^{-s}
 ;
 l-dq
 \Big)\\
&\qquad\qquad
  \qquad\qquad
  \qquad\qquad
  \qquad
  +
  \,\, 
  \lim_{n}	
 \sum
 \Sb
 \omega
 \in P(
 \zeta_{\mu}^{q}
 )\\
 {}\\
 |\imag(\omega)|\le t_{q,n}
 \endSb
 \res
 \Big(
 s\to
 Z_{\mu}^{q}(s)
 \,
 r^{-s};
 \omega
 \Big)\,.\\
&{}
\tag15.3
 \endalign
 $$
Next,
if 
 $ \beta(q)\not\in\{-dq,1-dq,\ldots,d-dq\}$, then
 a simple calculation 
(using the fact that if $f$ and $g$ are meromorphic functions with
$f(\omega)\not=0$, $g(\omega)=0$ and $g'(\omega)\not=0$, then $\omega$
is a pole of $\frac{f}{g}$ and 
$\res(\frac{f}{g};\omega)=\frac{f(\omega)}{g'(\omega)}$)
shows that
 $$
 \align
  \res
 \Big(
 \,
 s 
 \to
  Z_{\mu}^{q}(s)
 \,
 r^{-s}
 ;
 l-dq
 \Big)
 &=
 \res
 \Bigg(
 \,
 s 
 \to
 \Bigg(
 \sum_{k}
  \frac{\kappa_{\mu}^{q,k}(K)\,(\sigma_{q,k}-1)}{s-(k-dq)}
  \Bigg)\zeta_{\mu}^{q}(s)
 \,
 r^{-s}
 ;
 l-dq
 \Bigg)\\
 &=
 \res
 \Big(
 \,
 s 
 \to
  \frac{\kappa_{\mu}^{q,l}(K)\,(\sigma_{q,l}-1)}{s-(l-dq)}
  \zeta_{\mu}^{q}(s)
 \,
 r^{-s}
 ;
 l-dq
 \Big)\\ 
 &=
  \kappa_{\mu}^{q,l}(K)\,(\sigma_{q,l}-1)
  \,
  \zeta_{\mu}^{q}(l-dq)
 \,
 r^{-(l-dq)}\\  
& =
-\kappa_{\mu}^{q,l}(K)
\,
 \sigma_{q,l}\,r^{-(l-dq)}\,.
 \tag15.4
 \endalign
 $$   
Combining (15.3) and (15.4) now yields 
  $$
 \align
 V_{\mu,r}^{q,\sym}(K)
&=
 \sum_{l}
  \kappa_{\mu}^{q,l}(K)
  \,
 \sigma_{q,l}
 \,
 r^{-l+dq}
 \,\,
 -
 \,\,
 \sum_{l}
  \kappa_{\mu}^{q,l}(K)
  \,
 \sigma_{q,l}
 \,
 r^{-l+dq}\\
&\qquad\qquad
  \qquad\qquad
  \qquad\qquad
  \qquad
  +
  \,\, 
  \lim_{n}	
 \sum
 \Sb
 \omega
 \in P(
 \zeta_{\mu}^{q}
 )\\
 {}\\
 |\imag(\omega)|\le t_{q,n}
 \endSb
 \res
 \Big(
 s\to
 Z_{\mu}^{q}(s)
 \,
 r^{-s};
 \omega
 \Big)\\
&=
 \lim_{n}	
 \sum
 \Sb
 \omega
 \in P(
 \zeta_{\mu}^{q}
 )\\
 {}\\
 |\imag(\omega)|\le t_{q,n}
 \endSb
 \res
 \Big(
 s\to
 Z_{\mu}^{q}(s)
 \,
 r^{-s};
 \omega
 \Big)\,.
 \endalign
 $$
This completes the proof.
\hfill$\square$

\bigskip
\bigskip


\heading{16. 
Proof of Theorem 5.7}
\endheading

The purpose of this section is to prove
Theorem 5.7.
Since the proof is somewhat long and involved we will now give a brief 
description of the main outline of the argument.
However, we first introduce some notation.
Fix $q\in\Bbb R$ and let $\Gamma$ denote the path defined in Proposition 13.3.
We now write
 $G^{q}$
 for the
 set of those $s\in\Bbb C$ such that
$s$ lies strictly to the right of $\Gamma$, i.e\.
 $$
 G^{q}
 =
 \left\{
 s\in\Bbb C
 \,\left|\,
 \sup
 \Sb
  z\in\Gamma\\
  \imag(z)=\imag(s)
  \endSb
  \real(z)
  <
  \real(s)
    \right.
  \right\}\,.
  \tag16.1
 $$
 Assuming that
$
 \beta(q)
 \not\in  
 \{
  -dq,1-dq,\ldots,d-dq
  \}$,
  the proof of Theorem 5.7 is divided into the following five parts:

\bigskip

\roster
\item"{\it Part 1.}"
{\it The behavior of  $Z_{\mu}^{q}(s) \,r^{-s}$ on $\Gamma$:}

\noindent
{\it Part 1.1:}
For all $0<r<1$,
the following integral exists, namely
 $$
  \int_{\Gamma}
 Z_{\mu}^{q}(s)
 \,
 r^{-s}
 \,
 ds\,.
 $$
{\it Part 1.2:} In addition
$$
\qquad\quad\,\,
 \frac{1}{r^{-\beta(q)}}
 \int_{\Gamma}
 Z_{\mu}^{q}(s)
 \,
 r^{-s}
 \,
 ds
\to
0
\,\,\,\,
\text{as $r\searrow 0$.}
 $$
Part 1 is proved in Lemma 16.1 and Theorem 16.2.

\bigskip

\item"{\it Part 2.}"
{\it The behavior of $Z_{\mu}^{q}(s)\,r^{-s}$
between $\Gamma$ and the critical line
$\real(s)=\beta(q)$:}

\noindent
{\it Part 2.1:}
For all $0<r<1$, the following limit exists, namely
 $$
  \lim_{n}
 \sum
 \Sb
 \omega\in P(\zeta_{\mu}^{q})\cap G^{q}\\
 |\imag(\omega)|\le t_{q,n}\\
 \real(\omega)<\beta(q)
 \endSb
 \res
 \Big(
 \,
 s 
 \to
  Z_{\mu}^{q}(s)
 \,
 r^{-s}
 ;
 \omega
 \Big)\,.
 $$
{\it Part 2.2:} In addition 
 $$
 \qquad\quad\,\,
 \frac{1}{r^{-\beta(q)}}
  \lim_{n}
 \sum
 \Sb
 \omega\in P(\zeta_{\mu}^{q})\cap G^{q}\\
 |\imag(\omega)|\le t_{q,n}\\
 \real(\omega)<\beta(q)
 \endSb
 \res
 \Big(
 \,
 s 
 \to
  Z_{\mu}^{q}(s)
 \,
 r^{-s}
 ;
 \omega
 \Big)
 \to
 0
 \,\,\,\,
 \text{as $r\searrow0$.}
 $$
Part 2 is proved in Theorem 16.3.

\bigskip

\item"{\it Part 3.}"
{\it The behavior of $Z_{\mu}^{q}(s) \,r^{-s}$
on the critical line
$\real(s)=\beta(q)$:}

\noindent
{\it Part 3.1:}
For all $r>0$, the following limit exists, namely
 $$
  \lim_{n}
 \sum
 \Sb
 \omega\in P(\zeta_{\mu}^{q})\cap G^{q}\\
 |\imag(\omega)|\le t_{q,n}\\
 \real(\omega)=\beta(q)
 \endSb
 \res
 \Big(
 \,
 s 
 \to
  Z_{\mu}^{q}(s)
 \,
 r^{-s}
 ;
 \omega
 \Big)\,.
 $$ 
%
{\it Part 3.2:}
In addition
 $$
 \frac{1}{r^{-\beta(q)}}
  \lim_{n}
 \sum
 \Sb
 \omega\in P(\zeta_{\mu}^{q})\cap G^{q}\\
 |\imag(\omega)|\le t_{q,n}\\
 \real(\omega)=\beta(q)
 \endSb
 \res
 \Big(
 \,
 s 
 \to
  Z_{\mu}^{q}(s)
 \,
 r^{-s}
 ;
 \omega
 \Big)
 \qquad\quad\,\,
 $$ 
is a multiplicatively periodic function of $r$. 
Part 3 is proved in Theorem 16.4.

 \newpage
 
\bigskip

\item"{\it Part 4.}"
{\it Computing
$ \frac{1}{2\pi\ima}
 \int_{c-\ima\infty}^{c+\ima\infty}
 Z_{\mu}^{q}(s)
 \,
 r^{-s}
\,
ds$
using 
Part 1.1, Part 2.1 and Part 3.1:}

\noindent
For all $0<r<1$, we have
 $$
 \align
 \frac{1}{2\pi\ima}
 \int_{c-\ima\infty}^{c+\ima\infty}
 Z_{\mu}^{q}(s)
 \,
 r^{-s}
\,
ds
&=
-
 \sum
 \Sb
 l=0,1,\ldots,d\\
\beta(q)<l-dq
 \endSb
 \kappa_{\mu}^{q,l}(K)
\,
 \sigma_{q,l}\,r^{-(l-dq)}\\
&\qquad
+
 \,
  \lim_{n}
 \sum
 \Sb
 \omega\in P(\zeta_{\mu}^{q})\cap G^{q}\\
 |\imag(\omega)|\le t_{q,n}\\
\real(\omega)<\beta(q)
 \endSb
 \res
 \Big(
 \,
 s 
 \to
  Z_{\mu}^{q}(s)
 \,
 r^{-s}
 ;
 \omega
 \Big)\\
&\qquad
+
\,
 \lim_{n}
 \sum
 \Sb
 \omega\in P(\zeta_{\mu}^{q})\cap G^{q}\\
 |\imag(\omega)|\le t_{q,n}\\
 \real(\omega)=\beta(q)
 \endSb
 \res
 \Big(
 \,
 s 
 \to
  Z_{\mu}^{q}(s)
 \,
 r^{-s}
 ;
 \omega
 \Big)\\ 
&\qquad
  + 
   \frac{1}{2\pi\ima}
\int_{\Gamma}
 \,\,
 Z_{\mu}^{q}(s)
 \,
 r^{-s}
 \,
 ds\,;
 \endalign 
 $$ 
observe that both of the two 
limits and the integral on the right hand side of the above equality
 are well-defined by Part 1.1, Part 2.1 and Part 3.1. 
Part 4 is proved in Theorem 16.5.

 \bigskip

\item"{\it Part 5.}"
{\it Proving Theorem 5.7
using 
Part 4, Part 1.2, Part 2.2 and Part 3.2:}

\noindent
Theorem 5.4 shows that
for all  $0<r<r_{\min}$, we have
$$
\align
\frac{1}{r^{-\beta(q)}}
  V_{\mu,r}^{q,\sym}(K)
&=
  \frac{1}{r^{-\beta(q)}}
 \sum_{l}
  \kappa_{\mu}^{q,l}(K)
  \,
 \sigma_{q,l}
 \,
 r^{-l+dq}
 \,\,
 \\
 &\qquad\qquad
 \qquad\qquad
  +
 \,\,
   \frac{1}{2\pi\ima}
   \frac{1}{r^{-\beta(q)}}
 \int_{c-\ima \infty}^{c+\ima \infty}
 Z_{\mu}^{q}(s)
 \,
 r^{-s}
 \,ds\,,
 \endalign
 $$
and Part 4, Part 1.2, Part 2.2 and Part 3.2
shows that
for all
$0<r<1$, we have
$$
\align
 &\frac{1}{r^{-\beta(q)}}
 \sum_{l}
  \kappa_{\mu}^{q,l}(K)
  \,
 \sigma_{q,l}
 \,
 r^{-l+dq}
 \,\,
 \\
 &\qquad\qquad
 \qquad\qquad
  +
 \,\,
   \frac{1}{2\pi\ima}
   \frac{1}{r^{-\beta(q)}}
 \int_{c-\ima \infty}^{c+\ima \infty}
 Z_{\mu}^{q}(s)
 \,
 r^{-s}
 \,ds
 =
 \pi_{q}^{\sym}(r)+\varepsilon_{q}^{\sym}(r)
 \endalign
 $$
where
$\pi_{q}^{\sym}$ is a multiplicatively periodic function
and
$\varepsilon_{q}^{\sym}(r)\to 0$ as $r\searrow 0$.
Consequently,
for all
$0<r<r_{\min}$, we have
$$
\frac{1}{r^{-\beta(q)}}
  V_{\mu,r}^{q,\sym}(K)
=
 \pi_{q}^{\sym}(r)+\varepsilon_{q}^{\sym}(r).
 $$
Part 5 is proved after the statement and proof of Theorem 16.5.

\endroster

\bigskip

After this brief outline, 
we now
state and prove the results in this section.

\bigskip

\proclaim{Lemma 16.1} 
\roster
\item"(1)"
Let $a_{1},\ldots,a_{n}$ be complex numbers with $\sum_{i}a_{i}=0$.
Let $x_{1},\ldots,x_{n}$ be real numbers and let $I$ be a compact interval with
$\{x_{1},\ldots,x_{n}\}\cap I=\varnothing$.
Then
 $$
 \sup
 \Sb
 z\in\Bbb C\\
 \real(z)\in I
 \endSb
 \Bigg|
 z^{2}
 \sum_{i}\frac{a_{i}}{z-x_{i}}
 \Bigg|
 <
 \infty\,.
 $$
 \item"(2)"
 Fix $q\in\Bbb R$.
 Assume that
 $
 \beta(q)
 \not\in  
  \{
  -dq,1-dq,\ldots,d-dq
  \}$.
 Let $b_{0}(q)$ be as in Proposition 13.3.
Then  
 $$
  \sup
 \Sb
 s\in\Bbb C\\
 \real(z)\in [b_{0}(q),\beta(q)]
 \endSb
 \Bigg|
 s^{2}
 \sum_{l=0}^{d}
 \frac
 {\kappa_{\mu}^{q,l}(K)\,(\sigma_{q,l}-1)}
 {s-(l-dq)}
 \Bigg|
 <
 \infty\,.
  $$
  \endroster
 \endproclaim
\noindent{\it  Proof}\newline
\noindent
(1) 
Below
we will use the following notation, namely, if $R$
is a 
polynomial, then we will write $\deg R$ for the degree of $R$.
Let $Q$ denote the polynomial defined by
$Q(z)=\prod_{i}(z-x_{i})$.
It is clear that there is a polynomial $P$ with $\deg P\le n-2$
such that
$\sum_{i}\frac{a_{i}}{z-x_{i}}
=
\frac{\sum_{i}a_{i}\prod_{j\not=i}(z-x_{j})}{\prod_{i}(z-x_{i})}
=
\frac{(\sum_{i}a_{i})z^{n-1}+P(z)}{Q(z)}$
for all $z\in\Bbb C\setminus \{x_{1},\ldots,x_{n}\}$,
and
since
$\sum_{i}a_{i}=0$, this shows that
 $z^{2}\sum_{i}\frac{a_{i}}{z-x_{i}}=\frac{z^{2}P(z)}{Q(z)}$
  for all
$z\in\Bbb C\setminus \{x_{1},\ldots,x_{n}\}$.
However,
since 
$\deg(z\to z^{2}P(z))\le n=\deg Q$
(because  $\deg P\le n-2$), we conclude that
$\limsup_{|z|\to\infty}|z^{2}\sum_{i}\frac{a_{i}}{z-x_{i}}|
=
\limsup_{|z|\to\infty}|\frac{z^{2}P(z)}{Q(z)}|<\infty$.
This clearly implies that 
there is a 
constant $A>0$ 
such that
$\sup_{|z|\ge A}
|z^{2}\sum_{i}\frac{a_{i}}{z-x_{i}}|<\infty$, whence
 $$
 \sup
 \Sb
 |z|\ge A\\
 \real(z)\in I
 \endSb
 \Bigg|
 z^{2}\sum_{i}\frac{a_{i}}{z-x_{i}}
\Bigg|
\le
\sup
 \Sb
 |z|\ge A\\
 \endSb
\Bigg|
z^{2}\sum_{i}\frac{a_{i}}{z-x_{i}}
\Bigg|
<
\infty\,.
\tag16.2
$$
Next,
let $C=\{z\in\Bbb C\,|\,|z|\le A\,,\,\real(z)\in I\}$. 
Since
$\{x_{1},\ldots,x_{n}\}\cap I=\varnothing$,
we conclude that the function $\Phi:C\to\Bbb C$ defined by
$\Phi(z)= z^{2}\sum_{i}\frac{a_{i}}{z-x_{i}}$
is well-defined and continuous on $C$.
Also, since 
 $C$ is compact (because $I$ is compact), 
the continuity of $\Phi$
implies that $\Phi$ is bounded on $C$,
i\.e\. 
$\sup
 \Sb
 |z|\le A\,,\,
 \real(z)\in I
 \endSb
 |
 z^{2}\sum_{i}\frac{a_{i}}{z-x_{i}}
|
=
 \sup
 \Sb
 |z|\le A\,,\,
 \real(z)\in I
 \endSb
 |\Phi(z)|
<
\infty$.
The desired result now follows from this and (16.2).

\noindent
(2)
This follows immediately
from (1) since
$
   \sum_{l}
 \kappa_{\mu}^{q,l}(K)\,(\sigma_{q,l}-1)
 =
 0$.
\hfill$\square$

\bigskip

\proclaim{Theorem 16.2. The behaviour of 
$Z_{\mu}^{q}(s)
 \,
 r^{-s}$ on $\Gamma$} 
Fix $q\in\Bbb R$.
Assume that
$
 \beta(q)
 \not\in  
 \{
  -dq,1-dq,\ldots,d-dq
  \}$.
Let $b_{0}(q)$, $\beta_{0}(q)$ and $\Gamma$ be as in Proposition 13.3.
\roster
\item"(1)"  
There is a constant
$c>0$ such that
for all $0<r<1$, we have
 $$
 \int_{\Gamma}
 \,\,
 |
 Z_{\mu}^{q}(s)
 \,
 r^{-s}
 |
 \,
 |ds|
 \le
 c
 \,
 r^{-\beta_{0}(q)}\,.
 $$
\item"(2)" 
For all $0<r<1$,
the function
$\Gamma\to\Bbb C
 \,:\,
 s
 \to
  Z_{\mu}^{q}(s)
 \,
 r^{-s}$
is integrable. 
In particular, 
for all $0<r<1$, the following integral
is well-defined, namely
 $$
 \int_{\Gamma}
 Z_{\mu}^{q}(s)
 \,
 r^{-s}
 \,
 ds\,.
 $$
\item"(3)"
We have
  $$
 \frac{1}{r^{-\beta(q)}}
 \int_{\Gamma}
 Z_{\mu}^{q}(s)
 \,
 r^{-s}
 \,
 ds
\to
0
\,\,\,\,
\text{as $r\searrow 0$.}
 $$
 \endroster
 \endproclaim
\noindent{\it  Proof}\newline
\noindent
(1) We first note that it follows from Lemma 16.1
and Proposition 13.3
that there is a constant $M$ such that
if $s\in\Bbb C$ with
$b_{0}(q)\le\real(s)\le \beta_{0}(q)$, then
 $$
 \Bigg|
 \sum_{l}
 \frac
 {\kappa_{\mu}^{q,l}(K)\,(\sigma_{q,l}-1)}
 {s-(l-dq)}
 \Bigg|
 \le
 M
 \frac{1}{|s|^{2}}\,,
 \tag16.3
 $$
and if $s\in\Gamma$, then
 $$
 |\zeta_{\mu}^{q}(s)|
 \le
 M\,.
 \tag16.4
 $$ 
For an integer $n$, let 
the numbers $u_{n}(q)$ and $v_{n}(q)$
and the paths
 $\Gamma_{n}^{-}$,
 $\Pi_{n}^{-}$,
 $\Gamma_{n}^{+}$
 and
 $\Pi_{n}^{+}$
 be defined as   in Proposition  13.3.
Inequalities now (16.3) and (16.4)  imply that
 $$
 \align
 \int_{\Gamma}
 |
  Z_{\mu}^{q}(s)
 \,
 r^{-s}
|
\,
|ds|
&= 
 \int_{\Gamma}
 \,\,
 \Bigg|
 \,
 \Bigg(
 \sum_{l}
 \frac
 {\kappa_{\mu}^{q,l}(K)\,(\sigma_{q,l}-1)}
 {s-(l-dq)}
 \Bigg)
 \,
 \zeta_{\mu}^{q}(s)
 r^{-s}
 \Bigg|
 \,
 |ds|\\
&\le
M^{2}
 \int_{\Gamma}
 \,\,
 \frac{1}{|s|^2}
 \,
 r^{-\real(s)}
 \,
 |ds|\\
&=
M^{2}
\,
\Bigg(
\sum_{n=-\infty}^{\infty}
\int_{\Pi_{n}^{-}}
 \,\,
 \frac{1}{|s|^2}
 \,
 r^{-\real(s)}
 \,
 |ds|\\
&\qquad\qquad
  \qquad
 +  
\sum_{n=-\infty}^{\infty}
\int_{\Pi_{n}^{+}}
 \,\,
 \frac{1}{|s|^2}
 \,
 r^{-\real(s)}
 \,
 |ds|\\
&\qquad\qquad
  \qquad
 + 
 \sum_{n=-\infty}^{\infty}
\,\,
\int_{\Gamma_{n}^{-}}
 \,\,
 \frac{1}{|s|^2}
 \,
 r^{-\real(s)}
 \,
 |ds|\\
&\qquad\qquad
  \qquad
 + 
\sum_{n=-\infty}^{\infty}
\,\,
\int_{\Gamma_{n}^{+}}
 \,\,
 \frac{1}{|s|^2}
 \,
 r^{-\real(s)}
 \,
 |ds|
 \Bigg)\,.
 \tag16.5
\endalign
$$

Below we analyse the 
sums
$\sum_{n=-\infty}^{\infty}
\int_{\Pi_{n}^{-}}
 \,\,
 \frac{1}{|s|^2}
 \,
 r^{-\real(s)}
 \,
 |ds|
 +  
\sum_{n=-\infty}^{\infty}
\int_{\Pi_{n}^{+}}
 \,\,
 \frac{1}{|s|^2}
 \,
 r^{-\real(s)}
 \,
 |ds|$
and
$
 \sum_{n=-\infty}^{\infty}
\,\,
\int_{\Gamma_{n}^{-}}
 \,\,
 \frac{1}{|s|^2}
 \,
 r^{-\real(s)}
 \,
 |ds|
 + 
\sum_{n=-\infty}^{\infty}
\,\,
\int_{\Gamma_{n}^{+}}
 \,\,
 \frac{1}{|s|^2}
 \,
 r^{-\real(s)}
 \,
 |ds|$ appearing on the left hand side of (16.5).

First we find an upper bound for the first of the two sums, namely, the following sum
$\sum_{n=-\infty}^{\infty}
\int_{\Pi_{n}^{-}}
 \,\,
 \frac{1}{|s|^2}
 \,
 r^{-\real(s)}
 \,
 |ds|
 +  
\sum_{n=-\infty}^{\infty}
\int_{\Pi_{n}^{+}}
 \,\,
 \frac{1}{|s|^2}
 \,
 r^{-\real(s)}
 \,
 |ds|$.
For brevity write $w=-\frac{\pi}{\log r_{\min}}$.
Since 
(see Proposition 13.3)
$\ldots 
<
u_{-1}(q)
<
v_{-1}(q)
<
u_{0}(q)
<
0
<
v_{0}(q)
<
u_{1}(q)
<
v_{1}(q)
<
\dots$
and
$u_{n+1}(q)-u_{n}(q)\ge w$ for all $n$,
 we conclude 
$v_{n}(q)\ge u_{n}(q)
\ge
(n-1)w+u_{1}(q)$ for all $n\ge 1$
and that
$|v_{-(n+1)}(q)|
\ge
|u_{-n}(q)|
\ge 
|-n|w+|u_{0}(q)|$
for all $n\ge 0$.
This clearly implies that
$\liminf_{n\to\pm\infty}\frac{|u_{n}(q)|}{|n|+1}\ge w$
and
$\liminf_{n\to\pm\infty}\frac{|v_{n}(q)|}{|n|+1}\ge w$.
Hence, if we write
$w_{0}=
\min(
\inf_{n}\frac{|u_{n}(q)|}{|n|+1},
\inf_{n}\frac{|v_{n}(q)|}{|n|+1}
)$, then
 $w_{0}>0$ and
 $$
 |u_{n}(q)|
 \ge
 w_{0}(|n|+1)\,,\,\,\,\,
 |v_{n}(q)|
 \ge
 w_{0}(|n|+1)
 \tag16.6
 $$
for all positive integers $n$.
Using (16.6), we now deduce that
if $0<r<1$, then
 $$
 \align
 \sum_{n=-\infty}^{\infty}
\int_{\Pi_{n}^{-}}
 \,\,
\frac{1}{|s|^2}
 \,
 r^{-\real(s)}
 \,
 |ds|
&+  
\sum_{n=-\infty}^{\infty}
\int_{\Pi_{n}^{+}}
 \,\,
 \frac{1}{|s|^2}
 \,
 r^{-\real(s)}
 \,
 |ds|\\
&=
 \sum_{n=-\infty}^{\infty}
  \int_{b_{0}(q)}^{\beta_{0}(q)}
 \frac{1}{\sigma^{2}+u_{n}(q)^{2}}
 r^{-\sigma}
 \,
 d\sigma\\
&\qquad
 +
\sum_{n=-\infty}^{\infty}
  \int_{b_{0}(q)}^{\beta_{0}(q)}
 \frac{1}{\sigma^{2}+v_{n}(q)^{2}}
 r^{-\sigma}
 \,
 d\sigma\\
&\le
\sum_{n=-\infty}^{\infty}
  \int_{b_{0}(q)}^{\beta_{0}(q)}
 \frac{1}{u_{n}(q)^{2}}
 r^{-\beta_{0}(q)}
 \,
 d\sigma\\
&\qquad
 +
\sum_{n=-\infty}^{\infty}
  \int_{b_{0}(q)}^{\beta_{0}(q)}
 \frac{1}{v_{n}(q)^{2}}
 r^{-\beta_{0}(q)}
 \,
 d\sigma\\
&\le
 2
 \frac{\beta_{0}(q)-b_{0}(q)}{w_{0}}
 \,
 \Bigg(
 \sum_{n=-\infty}^{\infty}
 \frac{1}{(|n|+1)^{2}}
  \Bigg)
  r^{-\beta_{0}(q)}\,.
  \tag16.7
 \endalign
 $$

Next, we find an upper bound for the second of the sums, 
namely,
$
 \sum_{n=-\infty}^{\infty}
\,\,
\int_{\Gamma_{n}^{-}}
 \,\,
 \frac{1}{|s|^2}
 \,
 r^{-\real(s)}
 \,
 |ds|
 + 
\sum_{n=-\infty}^{\infty}
\,\,
\int_{\Gamma_{n}^{+}}
 \,\,
 \frac{1}{|s|^2}
 \,
 r^{-\real(s)}
 \,
 |ds|$.
If $0<r<1$, then 
 $$
 \align
 \sum_{n=-\infty}^{\infty}
\,\,
\int_{\Gamma_{n}^{-}}
 \,\,
 \frac{1}{|s|^2}
 \,
 r^{-\real(s)}
 \,
 |ds|
&+ 
\sum_{n=-\infty}^{\infty}
\,\,
\int_{\Gamma_{n}^{+}}
 \,\,
 \frac{1}{|s|^2}
 \,
 r^{-\real(s)}
 \,
 |ds|
 \\
&
 =
\sum_{n=-\infty}^{\infty}
\,\,
\int_{v_{n-1}(q)}^{u_{n}(q)}
 \,\,
 \frac{1}{b_{0}(q)^{2}+t^{2}}
 \,
 r^{-b_{0}(q)}
 \,
 dt\\
&\qquad
 +
 \,
\sum_{n=-\infty}^{\infty}
\,\,
\int_{u_{n}(q)}^{v_{n}(q)}
 \,\,
 \frac{1}{\beta_{0}(q)^{2}+t^{2}}
 \,
 r^{-\beta_{0}(q)}
 \,
 dt\\
&
 \le
\sum_{n=-\infty}^{\infty}
\,\,
\int_{v_{n-1}(q)}^{v_{n}(q)}
 \,\,
 \frac{1}{\min(b_{0}(q)^{2},\beta_{0}(q)^{2})+t^{2}}
 \,
 dt
 \,\,
 \max\big(\, r^{-b_{0}(q)}\,,\, r^{-\beta_{0}(q)}\,\big)\\
&
 \le
\int_{-\infty}^{\infty}
 \frac{1}{\min(b_{0}(q)^{2},\beta_{0}(q)^{2})+t^{2}}
 \,
 dt
 \,\,
 r^{-\beta_{0}(q)}\\
&
 =
  \frac{\pi}{\min(|b_{0}(q)|,|\beta_{0}(q)|)}\,
 r^{-\beta_{0}(q)}\,.
 \tag16.8
 \endalign
 $$

Finally,
combining (16.5), (16.7) and (16.8) now gives
 $$
 \align
  \int_{\Gamma}
 \,\,
|
 Z_{\mu}^{q}(s)
 \,
 r^{-s}
| \,
 |ds|
&\le
 c
 \,
 r^{-\beta_{0}(q)}
 \endalign
 $$
where
$c
=
M^{2}
\big(
2
 \frac{\beta_{0}(q)-b_{0}(q)}{w_{0}}
 \sum_{n=-\infty}^{\infty}
 \frac{1}{(|n|+1)^{2}}
+
\frac{\pi}{\min(|b_{0}(q)|,|\beta_{0}(q)|)}
\big)$.

\noindent
(2)
This statement follows immediately from (1).

\noindent
(3)
Since $\beta_{0}(q)<\beta(q)$, it follows immediately from
(1) that
$| \frac{1}{r^{-\beta(q)}} 
\int_{\Gamma}
 Z_{\mu}^{q}(s)
 \,
 r^{-s}\,
 ds|
\le
 \frac{1}{r^{-\beta(q)}} 
\int_{\Gamma}
 |Z_{\mu}^{q}(s)
 \,
 r^{-s}
| \,
 |ds|
\le
 c
 r^{\beta(q)-\beta_{0}(q)}\to 0$ as $r\searrow0$.
\hfill$\square$

\bigskip

\proclaim{Theorem 16.3.
 The behavior of 
$Z_{\mu}^{q}(s)
 \,
 r^{-s}$ between
  $\Gamma$
  and the critical
  line
  $\real(s)=\beta(q)$}
Fix $q\in\Bbb R$.
Assume that
$\beta(q)
 \not\in  
  \{
 -dq,1-dq,\ldots,d-dq
   \}$.
   Let $(t_{q,n})_{n}$ be the 
   sequence from Theorem 13.5.
 \roster
 \item"(1)"  
 For $\omega\in P(\zeta_{\mu}^{q})\cap G^{q}$, define 
 $f_{\omega}:(0,1)\to\Bbb  C$ by
  $$
  f_{\omega}(r)
  =
   \frac{1}{r^{-\beta(q)}}\,
   \res
 \Big(
 \,
 s 
 \to
 Z_{\mu}^{q}(s)
 \,
 r^{-s}
 ;
 \omega
 \Big)\,.
 $$
Then
 $$
 \sum
 \Sb
 \omega\in P(\zeta_{\mu}^{q})\cap G^{q}\\
 \real(\omega)<\beta(q)
 \endSb 
 \|f_{\omega}\|_{\infty}
 <
 \infty\,.
 $$
 \item"(2)"
 For all $0<r<1$, the following limit exists, namely
  $$
 \frac{1}{r^{-\beta(q)}}
 \,\,\,\,
 \lim_{n}
 \sum
 \Sb
 \omega\in P(\zeta_{\mu}^{q})\cap G^{q}\\
 |\imag(\omega)|\le t_{q,n}\\
 \real(\omega)<\beta(q)
 \endSb
 \res
 \Big(
 \,
 s 
 \to
 Z_{\mu}^{q}(s)
 \,
 r^{-s}
 ;
 \omega
 \Big)\,.
  $$
\item"(3)"
 We have
 $$
 \qquad\qquad
 \qquad\,\,\,\,
 \frac{1}{r^{-\beta(q)}}
 \,\,\,\,
 \lim_{n}
 \sum
 \Sb
 \omega\in P(\zeta_{\mu}^{q})\cap G^{q}\\
 |\imag(\omega)|\le t_{q,n}\\
 \real(\omega)<\beta(q)
 \endSb
 \res
 \Big(
 \,
 s 
 \to
 Z_{\mu}^{q}(s)
 \,
 r^{-s}
 ;
 \omega
 \Big)
 \,
 \to
 \,
 0
 \,\,\,\,
 \text{as $r\searrow0$.}
 $$
 \endroster
 \endproclaim
\noindent{\it  Proof}\newline
\noindent
We first note that if
$\omega\in P(\zeta_{\mu}^{q})\cap G^{q}$, then
$\real(\omega)\in[b_{0}(q),\beta(q)]$, and
since
$\{-dq,1-dq,\ldots,d-dq\}\cap[b_{0}(q),\beta(q)]=\varnothing$,
we therefore deduce that
for all $l=0,1,\ldots,d$, we have
$\omega\not=l-dq$.
This implies that 
$$
\align
 \res(
 s 
 \to
 Z_{\mu}^{q}(s)
 \,
 r^{-s}
 ;
 \omega)
&=
 \res(
 s 
 \to
 \Bigg(\sum_{l}\frac{\kappa_{\mu}^{q,l}(K)\,(\sigma_{q,l}-1)}{s-(l-dq)}\Bigg)
 \,
 \zeta_{\mu}^{q}(s)
 \,
 r^{-s}
 ;
 \omega)\\
&=
 \Bigg(\sum_{l}\frac{\kappa_{\mu}^{q,l}(K)\,(\sigma_{q,l}-1)}{\omega-(l-dq)}\Bigg)
  \,
 r^{-\omega}
 \,
 \res(
 \zeta_{\mu}^{q}
 ;
 \omega)\,.
 \tag16.9
 \endalign
 $$

\noindent
(1)
Lemma 16.1 and Proposition 13.2 show that
 that there is a constant $M>0$
such that if $s\in\Bbb C$
with
$\real(s)\in[b_{0}(q),\beta(q)]$, then
 $$
 \Bigg|\sum_{l}\frac{\kappa_{\mu}^{q,l}(K)\,(\sigma_{q,l}-1)}{s-(l-dq)}\Bigg|
\le
M
\frac{1}{|s|^{2}}\,,
\tag16.10
$$
and
if $\omega\in P(\zeta_{\mu}^{q})\cap G^{q}$, then
 $$
 |\res(
 \zeta_{\mu}^{q}
 ;
 \omega)| 
 \le 
M\,.
\tag16.11
 $$
Combining (16.9) and the
inequalities (16.10) and (16.11) shows
that
 $$
 \align
 \sum
 \Sb
 \omega\in P(\zeta_{\mu}^{q})\cap G^{q}\\
 \real(\omega)<\beta(q)
 \endSb 
 \|f_{\omega}\|_{\infty}
 &=
 \lim_{n}
  \sum
 \Sb
 \omega\in P(\zeta_{\mu}^{q})\cap G^{q}\\
  |\imag(\omega)|\le t_{q,n}\\
 \real(\omega)<\beta(q)
 \endSb 
 \|f_{\omega}\|_{\infty}\\
&=
 \lim_{n}
 \sum
 \Sb
 \omega\in P(\zeta_{\mu}^{q})\cap G^{q}\\
 |\imag(\omega)|\le t_{q,n}\\
 \real(\omega)<\beta(q)
 \endSb
 \,\,\,\,
 \sup_{0<r<1}
 \Bigg|
 \frac{1}{r^{-\beta(q)}}
 \res
 \Big(
 \,
 s 
 \to
 Z_{\mu}^{q}(s)
 \,
 r^{-s}
 ;
 \omega
 \Big)
 \Bigg|\\
&\le
 \lim_{n}
 \sum
 \Sb
 \omega\in P(\zeta_{\mu}^{q})\cap G^{q}\\
 |\imag(\omega)|\le t_{q,n}\\
 \real(\omega)<\beta(q)
 \endSb
 \,\,\,\,
 \sup_{0<r<1}
 \,\,\,\,
 \frac{1}{r^{-\beta(q)}}
  \,
 \Bigg|\sum_{l}\frac{\kappa_{\mu}^{q,l}(K)\,(\sigma_{q,l}-1)}{\omega-(l-dq)}\Bigg|
  \,
 |r^{-\omega}|
 \,
  |\res(
 \zeta_{\mu}^{q}
 ;
 \omega)| 
\\
& \le
M^{2}
\,
 \lim_{n}
 \sum
 \Sb
 \omega\in P(\zeta_{\mu}^{q})\cap G^{q}\\
 |\imag(\omega)|\le t_{q,n}\\
 \real(\omega)<\beta(q)
 \endSb
  \,\,\,\,
 \sup_{0<r<1}
 \,\,\,\,
 \frac{1}{|\omega|^{2}}
 \,
 r^{\beta(q)-\real(\omega)}\\
&=
 M^{2}
 \,
 \lim_{n}
 \sum
 \Sb
 \omega\in P(\zeta_{\mu}^{q})\cap G^{q}\\
 |\imag(\omega)|\le t_{q,n}\\
 \real(\omega)<\beta(q)
 \endSb
  \,\,\,\,
 \frac{1}{|\omega|^{2}}\\
 \allowdisplaybreak
  &\le
 M^{2}
 \,
 \lim_{n}
 \sum
 \Sb
 \omega\in P(\zeta_{\mu}^{q})\cap G^{q}\\
 |\imag(\omega)|\le t_{q,n}\\
 \real(\omega)<\beta(q)
 \endSb
  \,\,\,\,
 \frac{1}{|\omega|^{2}} 
 \\ 
  \allowdisplaybreak
 &\le
 M^{2}
 \lim_{n}
 \left(
 \,\,
 \sum
 \Sb
 k\in\Bbb N\\
 k<\integer(t_{q,n})
 \endSb
 \,\,
 \sum
 \Sb
 \omega\in P(\zeta_{\mu}^{q})\cap G^{q}\\
 k< |\imag(\omega)|\le k+1\\
 \real(\omega)<\beta(q)
 \endSb
 \frac{1}{|\omega|^{2}}
 \right.
 \\  
&\qquad\qquad
 \qquad\qquad
 \qquad\qquad
 \qquad\qquad
 \left.
 +
 \sum
 \Sb
 \omega\in P(\zeta_{\mu}^{q})\cap G^{q}\\
 \integer(t_{q,n})< |\imag(\omega)|\le t_{q,n}\\
 \real(\omega)<\beta(q)
 \endSb
 \frac{1}{|\omega|^{2}}
 \right)\\  
  \allowdisplaybreak
 &\le
 M^{2}
 \lim_{n}
 \left(
 \,\,
 \sum
 \Sb
 k\in\Bbb N\\
 k<\integer(t_{q,n})
 \endSb
 \,\,
 \sum
 \Sb
 \omega\in P(\zeta_{\mu}^{q})\cap G^{q}\\
 k< |\imag(\omega)|\le k+1\\
 \real(\omega)<\beta(q)
 \endSb
 \frac{1}{k^{2}}
 \right.
 \\ 
&\qquad\qquad
 \qquad\qquad
 \qquad\qquad
 \qquad\qquad
 \left.
 +
 \sum
 \Sb
 \omega\in P(\zeta_{\mu}^{q})\cap G^{q}\\
 \integer(t_{q,n})< |\imag(\omega)|\le t_{q,n}\\
 \real(\omega)<\beta(q)
 \endSb
 \frac{1}{ \integer(t_{q,n})^{2}}
 \right)\\  
  \allowdisplaybreak
&\le
 M^{2}
 \lim_{n}
 \,\,
 \sum
 \Sb
 k\in\Bbb N\\
 k\le\integer(t_{q,n})
 \endSb
 \,\,
 \sum
 \Sb
 \omega\in P(\zeta_{\mu}^{q})\cap G^{q}\\
 k< |\imag(\omega)|\le k+1\\
 \real(\omega)<\beta(q)
 \endSb
 \frac{1}{k^{2}}\\
&\le
 M^{2}
 \lim_{n}
 \,\,
 \sum
 \Sb
 k\in\Bbb N\\
 k\le\integer(t_{q,n})
 \endSb
 \,\,
 \Big|
 \{\omega\in P(\zeta_{\mu}^{q})\,|\,k < |\imag(\omega)|\le k+1\}
 \Big|
 \,\,
  \frac{1}{k^{2}}\,.\\
  \tag16.12
 \endalign
 $$
For brevity write
$\Xi_{t}
 =
 \{\omega\in P(\zeta_{\mu}^{q})\,|\,|\imag(\omega)|\le t\}$
for $t>0$,
and note that
it follows from Theorem 13.8 that
$|\Xi_{t}|
 =
 \gamma t+\Cal O(\log t)$
where $\gamma=-\frac{1}{\pi}\log r_{\min}$.
This clearly implies that there is a constant $c>0$ such that
 $$
 \gamma t-c\log t
\le
   |\Xi_{t}|
\le
\gamma t+c\log t
\tag16.13
$$ 
for all $t>0$. 
Since $\Xi_{k}\subseteq \Xi_{k+1}$, it follows from 
(16.13)
that
 $$
 \align
 \Big|
 \{\omega\in P(\zeta_{\mu}^{q})\,|\,k < |\imag(\omega)|\le k+1\}
 \Big|
&=
|
 \Xi_{k+1}\setminus\Xi_{k}
|\\
&=
| \Xi_{k+1}|-|\Xi_{k}|\\
&\le
(\gamma(k+1)+c\log(k+1))-(\gamma k-c\log k)\\
&\le
\gamma+2c\log(k+1)\,.
\tag16.14
\endalign
$$
Inequality (16.14) and (16.15) now imply that
 $$
 \align
  \sum
 \Sb
 \omega\in P(\zeta_{\mu}^{q})\cap G^{q}\\
 \real(\omega)<\beta(q)
 \endSb 
 \|f_{\omega}\|_{\infty}
&\le
M^{2}
 \lim_{n}
 \,\,
 \sum
 \Sb
 k\in\Bbb N\\
 k\le\integer(t_{q,n})
 \endSb
 \,\,
 \Big|
 \{\omega\in P(\zeta_{\mu}^{q})\,|\,k < |\imag(\omega)|\le k+1\}
 \Big|
 \,\,
  \frac{1}{k^{2}} \\
 &\le 
M^{2}
(\gamma+2c)
 \sum_{k=1}^{\infty}
  \frac{\log(k+1)}{k^{2}}\\
 &<
  \infty \,. 
 \endalign
 $$
 
\noindent
(2) This follows immediately from (1).

\noindent
(3)
Observe that for each
$\omega\in P(\zeta_{\mu}^{q})\cap G^{q}$ with $\real(\omega)<\beta(q)$ we have   
(using (16.9))
$$
  \align
  |f_{\omega}(r)|
 &=
   \Bigg|
   \frac{1}{r^{-\beta(q)}}\,
   \res
 \Big(
 \,
 s 
 \to
 Z_{\mu}^{q}(s)
 \,
 r^{-s}
 ;
 \omega
 \Big)
 \Bigg|\\
&=
\frac{1}{r^{-\beta(q)}}
  \,
 \Bigg|\sum_{l}\frac{\kappa_{\mu}^{q,l}(K)\,(\sigma_{q,l}-1)}{\omega-(l-dq)}\Bigg|
  \,
 |r^{-\omega}|
 \,
  |\res(
 \zeta_{\mu}^{q}
 ;
 \omega)| \\
&=
 \Bigg|\sum_{l}\frac{\kappa_{\mu}^{q,l}(K)\,(\sigma_{q,l}-1)}{\omega-(l-dq)}\Bigg|
  \,
 r^{\beta(q)-\real(\omega)}
 \,
  |\res(
 \zeta_{\mu}^{q}
 ;
 \omega)| \\
&\to
0
\,\,\,\,
\text{as $r\searrow 0$.} 
 \tag16.15
 \endalign
 $$
Next, since it follows from part (1) of the theorem that
$ \sum_{
 \omega\in P(\zeta_{\mu}^{q})\cap G^{q}\,,\,
 \real(\omega)<\beta(q)
 }
 \|f_{\omega}\|_{\infty}
 <
 \infty$,
we now conclude from (16.15) that
 $$
 \align
  \frac{1}{r^{-\beta(q)}}
 \,\,\,\,
 \lim_{n}
 \sum
 \Sb
 \omega\in P(\zeta_{\mu}^{q})\cap G^{q}\\
 |\imag(\omega)|\le t_{q,n}\\
 \real(\omega)<\beta(q)
 \endSb
 \res
 \Big(
 \,
 s 
 \to
 Z_{\mu}^{q}(s)
 \,
 r^{-s}
 ;
 \omega
 \Big)
&=
 \lim_{n}
 \sum
 \Sb
 \omega\in P(\zeta_{\mu}^{q})\cap G^{q}\\
 |\imag(\omega)|\le t_{q,n}\\
 \real(\omega)<\beta(q)
 \endSb
 f_{\omega}(r)\\
&\to
 \lim_{n}
 \sum
 \Sb
 \omega\in P(\zeta_{\mu}^{q})\cap G^{q}\\
 |\imag(\omega)|\le t_{q,n}\\
 \real(\omega)<\beta(q)
 \endSb
 0\\
&=
 0
 \,\,\,\,
 \text{as $r\searrow 0$.}
 \endalign
 $$
This completes the proof.
  \hfill$\square$

\bigskip

\proclaim{Theorem 16.4.
The behavior of
$Z_{\mu}^{q}(s)
 \,
 r^{-s}$
 on the critical line $\real(s)=\beta(q)$} 
Fix $q\in\Bbb R$.
Assume that
$\beta(q)
 \not\in  
  \{
  -dq,1-dq,\ldots,d-dq
  \}$.
  Let $(t_{q,n})_{n}$ be the 
   sequence from Theorem 13.5.

\roster
\item"(1)"
For $r>0$, the following limit exists, namely
 $$
 \frac{1}{r^{-\beta(q)}}
  \lim_{n}
 \sum
 \Sb
 \omega\in P(\zeta_{\mu}^{q})\cap G^{q}\\
 |\imag(\omega)|\le t_{q,n}\\
 \real(\omega)=\beta(q)
 \endSb
 \res
 \Big(
 \,
 s 
 \to
  Z_{\mu}^{q}(s)
 \,
 r^{-s}
 ;
 \omega
 \Big)\,.
$$
\endroster

\bigskip

 \noindent 
Write
  $$
\pi_{q}^{\sym}(r)
 =
  \frac{1}{r^{-\beta(q)}}
  \lim_{n}
 \sum
 \Sb
 \omega\in P(\zeta_{\mu}^{q})\cap G^{q}\\
 |\imag(\omega)|\le t_{q,n}\\
 \real(\omega)=\beta(q)
 \endSb
 \res
 \Big(
 \,
 s 
 \to
  Z_{\mu}^{q}(s)
 \,
 r^{-s}
 ;
 \omega
 \Big)\,.
 \qquad\qquad
 \,{}
 $$

\bigskip

\roster
\item"(2)"
If the set 
$\{\log r_{1},\ldots,\log r_{N}\}$
is not contained in a discrete additive subgroup of $\Bbb R$,
then
 $$
\pi_{q}^{\sym}(r)
 =
 -
 \frac
 {1}
 {\sum_{i}p_{i}^{q}r_{i}^{\beta(q)}\log r_{i}}
\,
\sum_{l}\frac{\kappa_{\mu}^{q,l}(K)\,(\sigma_{q,l}-1)}{\beta(q)-(l-dq)}\,.
\qquad\qquad
\qquad\,
$$

\item"(3)"
If the set 
$\{\log r_{1},\ldots,\log r_{N}\}$
is contained in a discrete additive subgroup of $\Bbb R$
and
$\langle\log r_{1},\ldots,\log r_{N}\rangle= u\Bbb Z$
 with $u > 0$, 
 then
 $$
 \align
 \qquad\qquad
 \,
\pi_{q}^{\sym}(r)
&=
 -
 \frac
 {1}
 {\sum_{i}p_{i}^{q}r_{i}^{\beta(q)}\log r_{i}}
 \,
 u
 \\
&\qquad\qquad
   \qquad
 \times 
 \sum_{l}
 \,
 \Bigg(
 \,
 \frac
 {\kappa_{\mu}^{q,l}(K)\,(\sigma_{q,l}-1)}
 {e^{u(\beta(q)-(l-dq))}-1}
 \\
&\\ 
&\qquad\qquad
   \qquad\qquad
 \times
 \cases
 \frac{e^{u(\beta(q)-(l-dq))    }+1}{2}
&\quad
 \text
 {for
 $r\in e^{\Bbb Z u}$;
 }\\
 {}\\
 e^{u(\beta(q)-(l-dq))\fraction(-\frac{\log r}{u})      }
&\quad
 \text
 {for 
 $r\not\in e^{\Bbb Z u}$
 }
 \endcases
 \Bigg)\,;
 \endalign
 $$
recall, that for a real number $x$, we write $\fraction(x)$ for the fractional part of $x$.
\endroster 
\endproclaim 
\noindent{\it  Proof}\newline
\noindent
Assume
that
$\beta(q)+\ima t$ with $t\in\Bbb R$ is a pole of
$\zeta_{\mu}^{q}$.
It follows from Proposition 13.2 that $\beta(q)+\ima t$ 
is a simple pole of $\zeta_{\mu}^{q}$,
and 
since
 $ \beta(q)\not\in\{-dq,1-dq,\ldots,d-dq\}$, 
 we therefore conclude
 from the definition of $Z_{\mu}^{q}$
(using the fact that if $f$ and $g$ are meromorphic functions and $\omega$ is a simple pole of 
$f$ 
and
$g(\omega)\not=0$, then $\omega$
is a pole of $fg$ and 
$\res(fg;\omega)=g(\omega)\res(f;\omega)$)
that
 $$
 \align
  \res
 \Big(
 \,
 s 
 \to
  Z_{\mu}^{q}(s)
 \,
 r^{-s}
&;
 \beta(q)+\ima t
 \Big)\\
 &=
 \res
 \Bigg(
 \,
 s 
 \to
 \Bigg(
 \sum_{l}
  \frac{\kappa_{\mu}^{q,l}(K)\,(\sigma_{q,l}-1)}{s-(l-dq)}
  \Bigg)\zeta_{\mu}^{q}(s)
 \,
 r^{-s}
 ;
 \beta(q)+\ima t
 \Bigg)\\
 &=
 \Bigg(
 \sum_{l}
  \frac{\kappa_{\mu}^{q,l}(K)\,(\sigma_{q,l}-1)}{(\beta(q)+\ima t)-(l-dq)}
 \Bigg) 
 \,
 r^{-(\beta(q)+\ima t)}
 \res(\zeta_{\mu}^{q}
 ;
 \beta(q)+\ima t)\,.\\
 &
 \tag16.15
 \endalign
 $$   
It also follows from Proposition 13.2 that
$ \res(\zeta_{\mu}^{q};\beta(q)+\ima t)
=
-
\frac{1}{\sum_{i}p_{i}^{q}r_{i}^{\beta(q)+\ima t}\log r_{i}}$, 
and we therefore conclude from (16.15) that
 $$
 \align
  \res
 \Big(
 \,
 s 
 \to
  Z_{\mu}^{q}(s)
 \,
 r^{-s}
&;
 \beta(q)+\ima t
 \Big)\\
 &=
  -
\frac{1}{\sum_{i}p_{i}^{q}r_{i}^{\beta(q)+\ima t}\log r_{i}}
\Bigg(
 \sum_{l}
  \frac{\kappa_{\mu}^{q,l}(K)\,(\sigma_{q,l}-1)}{(\beta(q)+\ima t)-(l-dq)}
  \Bigg)
 r^{-(\beta(q)+\ima t)}\,.\\
 &
 \tag16.16
 \endalign
 $$

If the set 
$\{\log r_{1},\ldots,\log r_{N}\}$
is not contained in a discrete additive subgroup of $\Bbb R$,
then
it follows from Proposition 13.1
that $\beta(q)$ is the only pole $\omega$ 
of $\zeta_{\mu}^{q}$ with
$\real(\omega)=\beta(q)$, 
and it therefore follows from (16.16) that
 $$
 \align
 \sum
 \Sb
 \omega\in P(\zeta_{\mu}^{q})\cap G^{q}\\
 |\imag(\omega)|\le t_{q,n}\\
 \real(\omega)=\beta(q)
 \endSb
 \res
 \Big(
 \,
 s 
 \to
  Z_{\mu}^{q}(s)
 \,
 r^{-s}
 ;
 \omega
 \Big)
&=
 \res
 \Big(
 \,
 s 
 \to
  Z_{\mu}^{q}(s)
 \,
 r^{-s}
 ;
 \beta(q)
 \Big)\\
&= 
  -
\frac{1}{\sum_{i}p_{i}^{q}r_{i}^{\beta(q)}\log r_{i}}
\Bigg(
 \sum_{l}
  \frac{\kappa_{\mu}^{q,l}(K)\,(\sigma_{q,l}-1)}{\beta(q)-(l-dq)}
  \Bigg)
 r^{-\beta(q)}
 \endalign
 $$   
for all integers $n$. This proves 
the
desired 
result.

If, on the other hand, the set 
$\{\log r_{1},\ldots,\log r_{N}\}$
is contained in a discrete additive subgroup of $\Bbb R$
and
$\langle\log r_{1},\ldots,\log r_{N}\rangle= u\Bbb Z$
 with $u > 0$, 
 then
 it follows from another application of
  Proposition 13.1 that the set of poles $\omega$ of $\zeta_{\mu}^{q}$
with
$\real(\omega)=\beta(q)$
is given by
$\beta(q)+\frac{2\pi}{u}\ima \Bbb Z$,
and it therefore follows
by a further application of (16.16) that
 $$
 \align
 \sum
&\Sb
 \omega\in P(\zeta_{\mu}^{q})\cap G^{q}\\
 |\imag(\omega)|\le t_{q,n}\\
 \real(\omega)=\beta(q)
 \endSb
\res
\Big(
 \,
 s 
 \to
  Z_{\mu}^{q}(s)
 \,
 r^{-s}
 ;
 \omega
 \Big)\\
&=
 \sum
 \Sb
 \frac{2\pi}{u}|k|\le t_{q,n}
 \endSb
 \res
 \Big(
 \,
 s 
 \to
  Z_{\mu}^{q}(s)
 \,
 r^{-s}
 ;
 \beta(q)+\tfrac{2\pi}{u}\ima k
 \Big)\\
&= 
 \sum
 \Sb
 \frac{2\pi}{u}|k|\le t_{q,n}
 \endSb
  -
\frac{1}{\sum_{i}p_{i}^{q}r_{i}^{\beta(q)+\tfrac{2\pi}{u}\ima k}\log r_{i}}
\Bigg(
 \sum_{l}
  \frac{\kappa_{\mu}^{q,l}(K)\,(\sigma_{q,l}-1)}{(\beta(q)+\tfrac{2\pi}{u}\ima k)-(l-dq)}
  \Bigg)
 r^{-(\beta(q)+\tfrac{2\pi}{u}\ima k)}\\
&= 
 \sum_{l}
 \,
 \sum
 \Sb
 \frac{2\pi}{u}|k|\le t_{q,n}
 \endSb
  -
\frac{1}{\sum_{i}p_{i}^{q}r_{i}^{\beta(q)+\tfrac{2\pi}{u}\ima k}\log r_{i}}
\,\,
  \frac{\kappa_{\mu}^{q,l}(K)\,(\sigma_{q,l}-1)}{(\beta(q)+\tfrac{2\pi}{u}\ima k)-(l-dq)}
\,\,
 r^{-(\beta(q)+\tfrac{2\pi}{u}\ima k)}\\
 &
 \tag16.17
 \endalign
 $$   
for all integers $n$.
Finally,
a simple Fourier analysis argument shows
that
if $a$ is a real number with $a\not=0$, then 
$e^{a\fraction(x)}=\sum_{k\in\Bbb Z}\frac{e^{a}-1}{a-2\pi\ima k}e^{2\pi\ima k x}$
for $x\in\Bbb R\setminus\Bbb Z$
and
$\frac{e^{a}+1}{2}=\sum_{k\in\Bbb Z}\frac{e^{a}-1}{a-2\pi\ima k}e^{2\pi\ima k x}$
for 
$x\in\Bbb Z$;
recall, that for a real number $x$, we write $\fraction(x)$ for the fractional part of $x$.
The desired result now follows from this and (16.17).
\hfill$\square$

\bigskip

\proclaim{Theorem 16.5} 
Fix $q\in\Bbb R$
and $c>\max\big( -dq,1-dq,\ldots,d-dq,\beta(q)\big)$.
Assume that
 $\beta(q)
 \not\in  
 \{-dq,1-dq,\ldots,d-dq\}$.
   Let $(t_{q,n})_{n}$ be the 
   sequence from Theorem 13.5
   and
 let $\Gamma$ be as in Proposition 13.3.
For all $0<r<1$, we have
 $$
 \align
 \frac{1}{2\pi\ima}
 \int_{c-\ima\infty}^{c+\ima\infty}
 Z_{\mu}^{q}(s)
 \,
 r^{-s}
\,
ds
&=
-
 \sum
 \Sb
 l=0,1,\ldots,d\\
\beta(q)<l-dq
 \endSb
 \kappa_{\mu}^{q,l}(K)
\,
 \sigma_{q,l}\,r^{-(l-dq)}\\
&\qquad
+
 \,
  \lim_{n}
 \sum
 \Sb
 \omega\in P(\zeta_{\mu}^{q})\cap G^{q}\\
 |\imag(\omega)|\le t_{q,n}\\
\real(\omega)<\beta(q)
 \endSb
 \res
 \Big(
 \,
 s 
 \to
  Z_{\mu}^{q}(s)
 \,
 r^{-s}
 ;
 \omega
 \Big)\\
&\qquad
+
\,
 \lim_{n}
 \sum
 \Sb
 \omega\in P(\zeta_{\mu}^{q})\cap G^{q}\\
 |\imag(\omega)|\le t_{q,n}\\
 \real(\omega)=\beta(q)
 \endSb
 \res
 \Big(
 \,
 s 
 \to
  Z_{\mu}^{q}(s)
 \,
 r^{-s}
 ;
 \omega
 \Big)\\ 
&\qquad
  + 
   \frac{1}{2\pi\ima}
\int_{\Gamma}
 \,\,
 Z_{\mu}^{q}(s)
 \,
 r^{-s}
 \,
 ds\,;
 \tag16.18
 \endalign 
 $$ 
observe that both of the two 
limits and the integral on the right hand side of (16.18) are well-defined by 
Theorems 13.2--13.4. 
 \endproclaim 
\noindent{\it  Proof}\newline
It is clear that the intersections 
$\Gamma
 \cap
 (\Bbb R-\ima t_{q,n})$
 and
$\Gamma
 \cap
 (\Bbb R+\ima t_{q,n})$
are compact line segments, 
i\.e\. there are
compact intervals $I_{q,n}^{-}$ and $I_{q,n}^{+}$ such that
$\Gamma
 \cap
 (\Bbb R-\ima t_{q,n})
=
 I_{q,n}^{-}-\ima t_{q,n}$
and
$ \Gamma
 \cap
 (\Bbb R+\ima t_{q,n})
=
 I_{q,n}^{+}-\ima t_{q,n}$.
 Now,
define paths 
$\Delta_{n}$, $\Lambda_{n}$,
$\Sigma_{n}^{-}$ and $\Sigma_{n}^{+}$
 by:
$$
\align
&\Delta_{n}
\quad\,
\text{
is the part of
 $\Gamma$
that lies in the set
$\Big\{s\in\Bbb C\,\Big|\,|\imag(s)|\le t_{q,n}\Big\}$,
}\\
&\Lambda_{n}
\quad\,\,
\text{
 is the 
 directed line segment from
 \qquad\qquad
 $c-\ima t_{q,n}$
 to
 \qquad\qquad
 $c+\ima t_{q,n}$,
 }\\
&\Sigma_{n}^{-}
\quad\,
\text{
 is the 
 directed line segment from
 $(\max I_{q,n}^{-})-\ima t_{q,n}$
 to
 \qquad\qquad
 $c-\ima t_{q,n}$,
 }\\
&\Sigma_{n}^{+}
\quad\,
\text{
is the 
 directed line segment from
 \qquad\qquad
  $c+\ima t_{q,n}$
  to
 $(\max I_{q,n}^{+})+\ima t_{q,n}$.
 }
 \endalign
  $$
Using standard notation, we let
$-\Delta_{n}$
denote the path 
 $\Delta_{n}$
 equipped
 with the opposite
direction.
Below we sketch the 
paths
$-\Delta_{n}$, $\Lambda_{n}$,
$\Sigma_{n}^{-}$ and $\Sigma_{n}^{+}$.

\goodbreak

\midinsert

 \vspace{130mm}
 \centerline{\hbox{\hskip -50mm\special{pdf=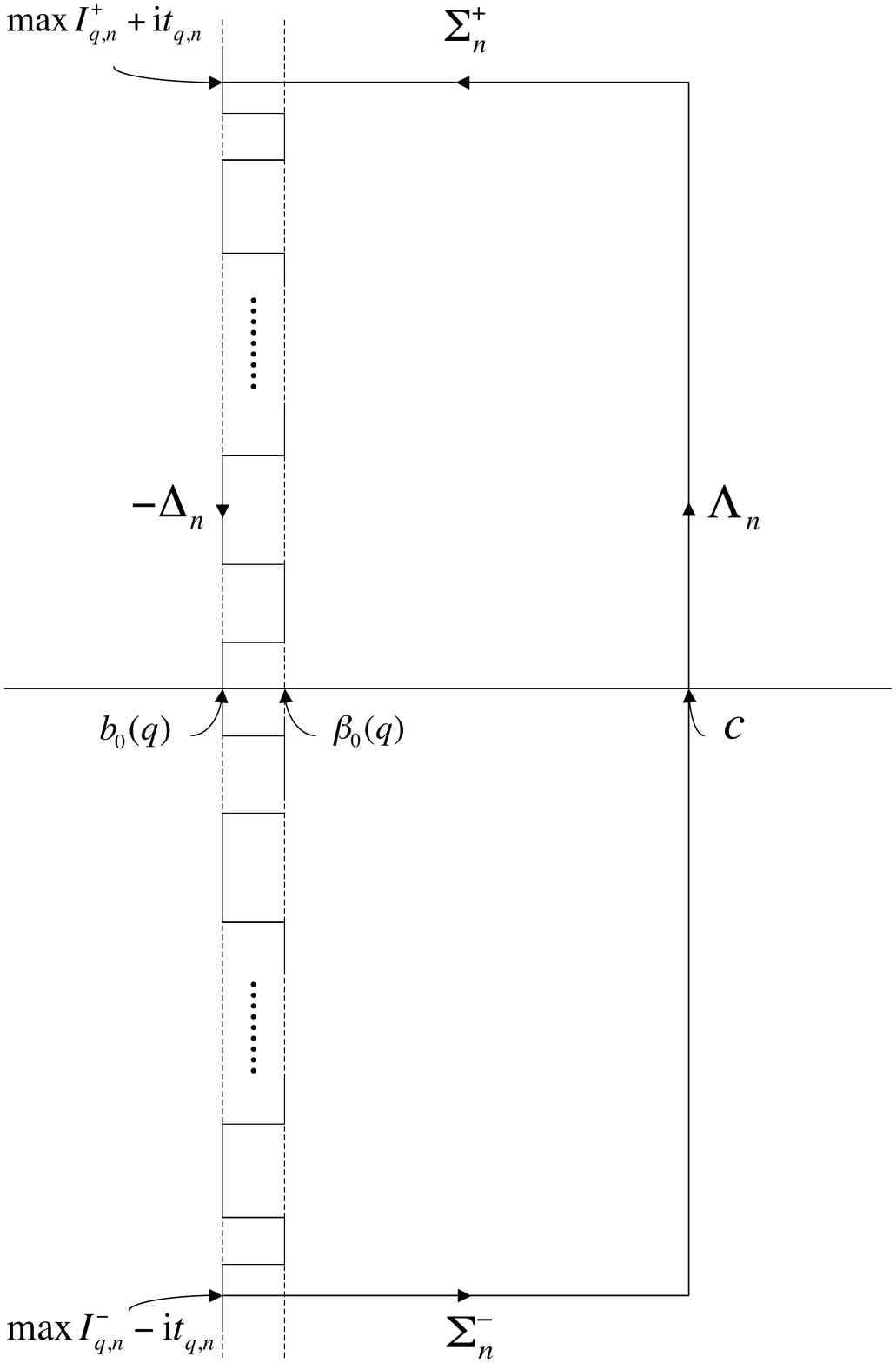
 scale=0.38}}}

\botcaption{\eightpoint \bf Fig\. 16.1}
\eightpoint 
The paths
$-\Delta_{n}$, $\Lambda_{n}$,
$\Sigma_{n}^{-}$ and $\Sigma_{n}^{+}$.
\endcaption
\endinsert

Next, we
observe that
since
$ \{-dq,1-dq,\ldots,d-dq\}\cap[b_{0}(q),\beta(q)]=\varnothing$
and 
all the poles $\omega$ of $\zeta_{\mu}^{q}$ satisfy 
$\real(\omega)\le\beta((q)<c$, it follows that
 $$
\align
 P\Big(
 s 
 \to
& Z_{\mu}^{q}(s)
 \,
 r^{-s}
 \Big)
 \,\,
\cap
 \,\,
 \Big(
 \Big\{s\in\Bbb C\,\Big|\,\real(s)<c\,,\,-t_{q,n}<\imag(s)<t_{q,n}\Big \}\cap G^{q}
 \Big)\\
&=
 P\Bigg(
 s 
 \to
\Bigg(\sum_{l}\frac{\kappa_{\mu}^{q,l}(K)\,(\sigma_{q,l}-1)}{s-(l-dq)}\Bigg)
\zeta_{\mu}^{q}(s)
 \,
 r^{-s}
 \Bigg)\\
&\qquad\qquad
\quad\,
\cap
 \,\,
 \Big(
 \Big\{s\in\Bbb C\,\Big|\,\real(s)<c\,,\,-t_{q,n}<\imag(s)<t_{q,n}\Big \}\cap G^{q}
 \Big)\\
&= 
\left(
    \bigcup
 \Sb
 l=0,1,\ldots,d\\
 \beta(q)<l-dq
 \endSb
 \{l-dq\}
 \right)
 \,\,
 \cup
 \,\,
  \Big\{\omega\in P(\zeta_{\mu}^{q})\cap G^{q}
  \,\Big|\,
  \real(\omega)<\beta(q)\,,\,-t_{q,n}<\imag(\omega)<t_{q,n}\Big \}\\
 &\qquad\qquad
   \qquad\qquad
   \qquad\,\,\,\,\,
     \cup
 \,\,
  \Big\{\omega\in P(\zeta_{\mu}^{q})\cap G^{q}
  \,\Big|\,
  \real(\omega)=\beta(q)\,,\,-t_{q,n}<\imag(\omega)<t_{q,n}\Big \}\,.\\
&
\tag16.19  
\endalign
 $$
As the paths 
$\Delta_{n}$, $\Lambda_{n}$,
$\Sigma_{n}^{-}$ and $\Sigma_{n}^{+}$
 enclose the
region
$\{s\in\Bbb C\,|\,\real(s)<c\,,\,-t_{q,n}<\imag(s)<t_{q,n} \}\cap G^{q}$,
we now deduce from(16.19)  and  the Residue Theorem that
 $$
 \align
&2\pi\ima
  \sum
 \Sb
 \omega\in P(\zeta_{\mu}^{q})\cap G^{q}\\
 |\imag(\omega)|\le t_{q,n}\\
\real(\omega)<\beta(q)
 \endSb
 \res
 \Big(
 \,
 s 
 \to
  Z_{\mu}^{q}(s)
 \,
 r^{-s}
 ;
 \omega
 \Big)
 \,
 +
 \,
 2\pi\ima
  \sum
 \Sb
 \omega\in P(\zeta_{\mu}^{q})\cap G^{q}\\
 |\imag(\omega)|\le t_{q,n}\\
\real(\omega)=\beta(q)
 \endSb
 \res
 \Big(
 \,
 s 
 \to
  Z_{\mu}^{q}(s)
 \,
 r^{-s}
 ;
 \omega
 \Big)\\
&+
 \,
 2\pi\ima
  \sum
 \Sb
 l=0,1,\dots,d\\
\beta(q)<l-dq
 \endSb
 \res
 \Big(
 \,
 s 
 \to
  Z_{\mu}^{q}(s)
 \,
 r^{-s}
 ;
 l-dq
 \Big)
 \\ 
&\qquad
 =
 \int_{c-\ima t_{q,n}}^{c+\ima t_{q,n}}
 Z_{\mu}^{q}(s)\,r^{-s}\,ds\\
&\qquad\qquad
  \qquad\qquad 
 +
  \int_{-\Delta_{n}}
 Z_{\mu}^{q}(s)\,r^{-s}\,ds\\
&\qquad\qquad
  \qquad\qquad
 +
 \int_{\Sigma_{n}^{-}}
 Z_{\mu}^{q}(s)\,r^{-s}\,ds\\
&\qquad\qquad
  \qquad\qquad 
 +
  \int_{\Sigma_{n}^{+}}
 Z_{\mu}^{q}(s)\,r^{-s}\,ds\,.
 \tag16.20
 \endalign
 $$
 Below we 
 compute
 the 
 sum
 $ \sum_{
 l=0,1,\dots,d\,,\,
\beta(q)<l-dq
}
 \res
 (
 \,
 s 
 \to
  Z_{\mu}^{q}(s)
 \,
 r^{-s}
 ;
 l-dq
 )$
 and show that the integrals
 $ \int_{\Sigma_{n}^{-}}
 Z_{\mu}^{q}(s)\,r^{-s}\,ds$
 and
 $ \int_{\Sigma_{n}^{+}}
 Z_{\mu}^{q}(s)\,r^{-s}\,ds$
 tend to $0$ as $r\searrow0$.

 We first compute
 the sum
 $ \sum_{
 l=0,1,\dots,d\,,\,
\beta(q)<l-dq
}
 \res
 (
 \,
 s 
 \to
  Z_{\mu}^{q}(s)
 \,
 r^{-s}
 ;
 l-dq
 )$.
Since
 $ \beta(q)\not=l-dq$, 
  a simple calculation 
(using the fact that if $f$ and $g$ are meromorphic functions with
$f(\omega)\not=0$, $g(\omega)=0$ and $g'(\omega)\not=0$, then $\omega$
is a pole of $\frac{f}{g}$ and 
$\res(\frac{f}{g};\omega)=\frac{f(\omega)}{g'(\omega)}$)
shows that
$
  \res
 (
 \,
 s 
 \to
  Z_{\mu}^{q}(s)
 \,
 r^{-s}
 ;
 l-dq
 )
 =
 \res
 (
 \,
 s 
 \to
 (
 \sum_{k}
  \frac{\kappa_{\mu}^{q,k}(K)\,(\sigma_{q,k}-1)}{s-(k-dq)}
  )\zeta_{\mu}^{q}(s)
 \,
 r^{-s}
 ;
 l-dq
 )\
 =
 \res
 (
 \,
 s 
 \to
  \frac{\kappa_{\mu}^{q,l}(K)\,(\sigma_{q,l}-1)}{s-(l-dq)}
  \zeta_{\mu}^{q}(s)
 \,
 r^{-s}
 ;
 l-dq
 )
 =
-\kappa_{\mu}^{q,l}(K)
\,
 \sigma_{q,l}\,r^{-(l-dq)}
 $.   
We deduce from this that
$$
 2\pi\ima
  \sum
 \Sb
 l=0,1,\dots,d\\
\beta(q)<l-dq
 \endSb
 \res
 \Big(
 \,
 s 
 \to
  Z_{\mu}^{q}(s)
 \,
 r^{-s}
 ;
 l-dq
 \Big)
 =
 -
  2\pi\ima
  \sum
 \Sb
 l=0,1,\dots,d\\
\beta(q)<l-dq
 \endSb
 \kappa_{\mu}^{q,l}(K)
\,
 \sigma_{q,l}\,r^{-(l-dq)}\,.
 \tag16.21
 $$

Next, we 
show that the integrals
 $ \int_{\Sigma_{n}^{-}}
 Z_{\mu}^{q}(s)\,r^{-s}\,ds$
 and
 $ \int_{\Sigma_{n}^{+}}
 Z_{\mu}^{q}(s)\,r^{-s}\,ds$
 tend to $0$ as $r\searrow0$.
Indeed, 
it follows from Lemma 16.1 that there is a constant $M>0$
such that
$|\sum_{l}\frac{\kappa_{\mu}^{q,l}(K)\,(\sigma_{q,l}-1)}{s-(l-dq)}|
\le
M\frac{1}{|s|^{2}}$
for all $s$ with $\real(s)\in[b_{0}(q),\beta(q)]$.
It also follows from 
Theorem 13.5
that there is a constant $k_{c}$
such that
if $\sigma\le c$ and $n\in\Bbb N$, then
 $$
|\zeta_{\mu}^{q}(\sigma+\ima t_{q,n})|\le k_{c}\,.
$$
Writing
$\ell(\Sigma_{n}^{+})$ for the length of the line segment $\Sigma_{n}^{+}$, 
we now conclude that
 $$
 \align
 \Bigg|
 \int_{\Sigma_{n}^{+}}
 Z_{\mu}^{q}(s)\,r^{-s}\,ds
 \Bigg|
&\le
 \ell(\Sigma_{n}^{+})
 \,
 \sup_{s\in \Sigma_{n}^{+}}|Z_{\mu}^{q}(s)\,r^{-s}|\\
&\le
 (\beta_{0}(q)-b_{0}(q))
 \,
 \sup_{s\in \Sigma_{n}^{+}}
 \Bigg|\sum_{l}\frac{\kappa_{\mu}^{q,l}(K)\,(\sigma_{q,l}-1)}{s-(l-dq)}\Bigg|
 \,
 |\zeta_{\mu}^{q}(s)|
\,
r^{-\beta_{0}(q)}\\
&\le
 (\beta_{0}(q)-b_{0}(q))
 \,
M
\Bigg(
\sup_{s\in \Sigma_{n}^{+}}
 \frac{1}{|s|^{2}} 
 \Bigg)
  k_{c}
\,
r^{-\beta_{0}(q)}
\\
&\le
 (\beta_{0}(q)-b_{0}(q))
 \,
M
\frac{1}{t_{q,n}^{2}}
  k_{c}
\,
r^{-\beta_{0}(q)} \\
&\to
 0\,. 
 \tag16.22
 \endalign
 $$
Similarly, one can prove that
 $$
 \Bigg|
 \int_{ \Sigma_{n}^{-}}
 Z_{\mu}^{q}(s)\,r^{-s}\,ds
 \Bigg|
 \to
 0\,.
 \tag16.23
 $$

Also, note that since
$ \int_{\Gamma}
 |Z_{\mu}^{q}(s)|\,|r^{-s}|\,|ds|<\infty$
(by Theorem 16.2), it follows
 from the Dominated Convergence Theorem that
 $$
 \align
  \int_{-\Delta_{n}}
 Z_{\mu}^{q}(s)\,r^{-s}\,ds
&=
 -
 \int_{\Delta_{n}}
 Z_{\mu}^{q}(s)\,r^{-s}\,ds\\
&\to
 -
  \int_{\Gamma}
 Z_{\mu}^{q}(s)\,r^{-s}\,ds\,.
 \tag16.24
 \endalign 
$$

Finally, 
the desired result now follows from (16.20)--(16.24)
by letting $n\to\infty$.
\hfill$\square$

\bigskip

We can now prove Theorem 5.7.

\bigskip

\noindent{\it  Proof of Theorem 5.7}\newline
\noindent
Fix $q\in\Bbb R$
and
assume that
$\beta(q)
 \not\in
 \{
 -dq,1-dq,\ldots,d-dq\}$.
Let
$c>\max\big(\,-dq,1-dq,\dots,d-dq,\beta(q)\,\big)$.
It follows from Theorem 5.4
and
Theorem 16.5
that for all $0<r<r_{\min}$
we have
 $$
 \align
   \frac{1}{r^{-\beta(q)}}
   V_{\mu,r}^{q,\sym}(K)
&=
  \frac{1}{r^{-\beta(q)}}
 \sum_{l}
  \kappa_{\mu}^{q,l}(K)
  \,
 \sigma_{q,l}
 \,
 r^{-l+dq}
 \,\,
 +
 \,\,
   \frac{1}{2\pi\ima}
   \frac{1}{r^{-\beta(q)}}
 \int_{c-\ima \infty}^{c+\ima \infty}
 Z_{\mu}^{q}(s)
 \,
 r^{-s}
 \,ds\\
&{}\\ 
&=
  \frac{1}{r^{-\beta(q)}}
 \sum_{l}
  \kappa_{\mu}^{q,l}(K)
  \,
 \sigma_{q,l}
 \,
 r^{-l+dq}
 \,\,
-
  \frac{1}{r^{-\beta(q)}}
 \sum
 \Sb
 l=0,1,\ldots,d\\
\beta(q)<l-dq
 \endSb
 \kappa_{\mu}^{q,l}(K)
\,
 \sigma_{q,l}\,r^{-(l-dq)}\\
&\qquad\qquad
  \qquad\qquad
+
 \,
   \frac{1}{r^{-\beta(q)}}
  \lim_{n}
 \sum
 \Sb
 \omega\in P(\zeta_{\mu}^{q})\cap G^{q}\\
 |\imag(\omega)|\le t_{q,n}\\
\real(\omega)<\beta(q)
 \endSb
 \res
 \Big(
 \,
 s 
 \to
  Z_{\mu}^{q}(s)
 \,
 r^{-s}
 ;
 \omega
 \Big)\\
&\qquad\qquad
  \qquad\qquad
+
\,
  \frac{1}{r^{-\beta(q)}}
 \lim_{n}
 \sum
 \Sb
 \omega\in P(\zeta_{\mu}^{q})\cap G^{q}\\
 |\imag(\omega)|\le t_{q,n}\\
 \real(\omega)=\beta(q)
 \endSb
 \res
 \Big(
 \,
 s 
 \to
  Z_{\mu}^{q}(s)
 \,
 r^{-s}
 ;
 \omega
 \Big)\\ 
&\qquad\qquad
  \qquad\qquad
  + 
   \frac{1}{2\pi\ima}
     \frac{1}{r^{-\beta(q)}}
\int_{\Gamma}
 \,\,
 Z_{\mu}^{q}(s)
 \,
 r^{-s}
 \,
 ds\\
&{}\\ 
\allowdisplaybreak
&=
  \frac{1}{r^{-\beta(q)}}
  \sum
 \Sb
 l=0,1,\ldots,d\\
 l-dq<\beta(q)
 \endSb
\kappa_{\mu}^{q,l}(K)
\,
 \sigma_{q,l}\,r^{-(l-dq)}
 \,\,\\
&\qquad\qquad 
  \qquad\qquad
 +
 \,\,
   \frac{1}{r^{-\beta(q)}}
 \lim_{n}
 \sum
 \Sb
 \omega\in P(\zeta_{\mu}^{q})\cap G^{q}\\
 |\imag(\omega)|\le t_{q,n}\\
 \real(\omega)=\beta(q)
 \endSb
 \res
 \Big(
 \,
 s 
 \to
  Z_{\mu}^{q}(s)
 \,
 r^{-s}
 ;
 \omega
 \Big)\\
&\qquad\qquad 
  \qquad\qquad
+
 \,
   \frac{1}{r^{-\beta(q)}}
  \lim_{n}
 \sum
 \Sb
 \omega\in P(\zeta_{\mu}^{q})\cap G^{q}\\
 |\imag(\omega)|\le t_{q,n}\\
\real(\omega)<\beta(q)
 \endSb
 \res
 \Big(
 \,
 s 
 \to
  Z_{\mu}^{q}(s)
 \,
 r^{-s}
 ;
 \omega
 \Big)\\  
&\qquad\qquad 
  \qquad\qquad
  + 
  \frac{1}{2\pi\ima}
    \frac{1}{r^{-\beta(q)}}
\int_{\Gamma}
 \,\,
 Z_{\mu}^{q}(s)
 \,
 r^{-s}
 \,
 ds\\ 
 &{}\\
&=
 \pi_{q}^{\sym}(r)
 \,\,
 +
 \,\, 
 \varepsilon_{q}^{\sym}(r)
 \endalign 
 $$ 
where (see Theorem 16.4)
 $$
 \pi_{q}^{\sym}(r)
 =
  \frac{1}{r^{-\beta(q)}}
 \lim_{n}
 \sum
 \Sb
 \omega\in P(\zeta_{\mu}^{q})\cap G^{q}\\
 |\imag(\omega)|\le t_{q,n}\\
 \real(\omega)=\beta(q)
 \endSb
 \res
 \Big(
 \,
 s 
 \to
  Z_{\mu}^{q}(s)
 \,
 r^{-s}
 ;
 \omega
 \Big)
 $$
and
 $$
  \varepsilon_{q}^{\sym}(r)
 =
 \varepsilon_{q,\bullet}^{\sym}(r)
 \,\,
  +
 \,\, 
 \varepsilon_{q,{\sssize\blacklozenge}}^{\sym}(r)
 \,\,
  +
 \,\, 
 \varepsilon_{q,{\sssize\blacktriangle}}^{\sym}(r)
 \qquad\qquad
 \qquad\quad
  $$
with
 $$
 \align
 \varepsilon_{q,\bullet}^{\sym}(r)
&=
  \frac{1}{r^{-\beta(q)}}
  \lim_{n}
 \sum
 \Sb
 \omega\in P(\zeta_{\mu}^{q})\cap G^{q}\\
 |\imag(\omega)|\le t_{q,n}\\
\real(\omega)<\beta(q)
 \endSb
 \res
 \Big(
 \,
 s 
 \to
  Z_{\mu}^{q}(s)
 \,
 r^{-s}
 ;
 \omega
 \Big)\,,\\ 
\varepsilon_{q,{\sssize\blacklozenge}}^{\sym}(r)
&=
  \frac{1}{r^{-\beta(q)}}
  \sum
 \Sb
 l=0,1,\ldots,d\\
 l-dq<\beta(q)
 \endSb
\kappa_{\mu}^{q,l}(K)
\,
 \sigma_{q,l}\,r^{-(l-dq)}\,,\\
\varepsilon_{q,{\sssize\blacktriangle}}^{\sym}(r)
&=
\frac{1}{2\pi\ima}
   \frac{1}{r^{-\beta(q)}}
\int_{\Gamma}
 \,\,
 Z_{\mu}^{q}(s)
 \,
 r^{-s}
 \,
 ds\,.
 \endalign
 $$

Finally, we note that
it is clear that
$\varepsilon_{q,{\sssize\blacklozenge}}^{\sym}(r)\to 0$ as $r\searrow 0$,
and that it follows from Theorem 16.2 and Theorem 16.3
that
$ \varepsilon_{q,\bullet}^{\sym}(r)\to 0$ as $r\searrow 0$
and that
$\varepsilon_{q,{\sssize\blacktriangle}}^{\sym}(r)\to 0$ as $r\searrow 0$.
We therefore conclude that
 $$
  \varepsilon_{q}^{\sym}(r)
 =
 \varepsilon_{q,\bullet}^{\sym}(r)
 \,\,
  +
 \,\, 
\varepsilon_{q,{\sssize\blacklozenge}}^{\sym}(r) \,\,
  +
 \,\, 
\varepsilon_{q,{\sssize\blacktriangle}}^{\sym}(r)
\to 0\,\,\,\,
\text{as $r\searrow 0$}\,.
  $$
Theorem 5.7 follows from this.
\hfill$\square$

\end
\newpage









\Refs\nofrills{References}


\widestnumber\no{HaJeKaPrShUUU}



\ref
\no ArPa
\by M. Arbeiter \& N. Patzschke
\paper Random self-similar multifractals
\jour Math. Nachr.
\vol 181
\yr 1996
\pages 5--42
\endref



\ref
\no BaPo
\by R. Badii \& A. Politi
\book Complexity. Hierarchical Structures and Scaling in Phusics
\publ Cambridge Nonlinear Science Series Vol\. 6,
Cambridge University Press, 
Cambridge, England
\yr 1997
\endref



%
\ref
\no BeGo
\by M. Berger \& B. Gostiaux
\book Differential Geometry: Manifolds, Curves and Surfaces
\publ Springer-Verlag, Berlin
\yr 1988
\endref




\ref
\no BeSc
\by C. Beck \& F. Schl\"ogl
\book Thermodynamics of Chaotic Systems. An Introduction
\publ Cambridge Nonlinear Science Series Vol\. 4,
Cambridge University Press, 
Cambridge, England
\yr 1993
\endref



%
\ref
\no CaMa
\by R. Cawley \& R. D. Mauldin
\paper Multifractal decomposition of Moran fractals
\jour Advances in Mathematics
\vol 92
\yr 1992
\pages 196--236
\endref



%
\ref
\no Col
\by J\. Cole
\paper The Geometry of Graph Directed Self-Conformal Multifractals
\jour Ph\.D\. Thesis, University of St. Andrews, 1998
\endref




\ref
\no Con
\by J\. Conway
\book Functions of One complex Variable, Second Edition
\publ Springer-Verlag, New York
\yr 1973
\endref



\ref
\no Da1
\by M\. Das
\paper Local properties of self-similar measures
\jour Illinois J. Math.
\vol 42
\yr 1998
\endref



\ref
\no Da2
\by M\. Das
\paper Hausdorff measures, dimensions and mutual 
singularity
\jour Trans. Amer. Math. Soc. 
\vol 357 
\yr 2005
\pages 4249Ð4268
\endref






\ref
\no DeKo\"OzRa\"Ur
\by A\. Deniz,
S\. Kocak,
Y\. \"Ozdemir,
A\. Ratiu \&
E\. \"Ureyen
\paper
On the Minkowski measurability of self-similar
fractals in $\Bbb R^{d}$
\jour Turkish Journal of Mathematics
\yr to appear
\endref








\ref
\no Ed
\by H\. M\. Edwards
\book Riemann's zeta function
\publ Pure and Applied Mathematics, Vol. 58. Academic Press, New York-London,
1974
\endref











%
\ref
\no Fa1
\by K\. J\. Falconer
\book Fractal Geometry | Mathematical Foundations and Applications
\publ
John Wiley, Chichester
\yr 1990
\endref



%
\ref
\no Fa2
\by K. J. Falconer
\book Techniques in Fractal Geometry
\publ Wiley
\yr 1997
\endref


%
\ref
\no Fa3
\by K. J. Falconer
\paper On the Minkowski measurability of fractals
\jour Proc. Am. Math.
Soc. 
\vol 123 
\yr 1995
\pages 1115--1124
\endref




%
\ref
\no Fed1
\by H\. Federer
\paper Curvature measures
\jour Trans. Amer. Math. Soc. 
\vol 93 
\yr 1959
\pages 418--491
\endref



%
\ref
\no Fed2
\by H\. Federer
\book Geometric Measure Theory
\publ Springer-Verlag New York Inc., New York 
\yr 1969 
\endref


%
\ref
\no Fel1
\by W\. Feller
\book An Introduction to Probability Theory and its Applications,
Vol\. 1, 3rd ed\.
\publ Wiley, New York
\yr 1986
\endref




%
\ref
\no Fel2
\by W\. Feller
\book An Introduction to Probability Theory and its Applications,
Vol\. 2, 2nd ed\.
\publ Wiley, New York
\yr 1971
\endref




%
\ref
\no FrPa
\by U. Frisch \& G. Parisi
\paper On the singularity
structure of fully developed turbulence, {\rm appendix to U. Frisch},
Fully developed turbulence and intermittency
\inbook Turbulence and
Predictability in Geophysical Fluid Dynamics and Climate
Dynamics
\publ Proc. Int. Sch. Phys., ``Enrico Fermi'' Course
LXXXVIII, pp\.84--88
\publaddr North Holland, Amsterdam
\yr1985
\endref





%
\ref
\no Fu1
\by J. H. G. Fu
\paper Tubular neighbourhoods in Euclidean spaces
\jour Duke Math. J. 
\vol 52 
\yr 1985
\pages 1025--1046
\endref




%
\ref
\no Fu2
\by J. H. G. Fu
\paper Curvature measures of subanalytic sets
\jour Amer. J. Math. 
\vol 116 
\yr 1994
\pages 819--880
\endref






%
\ref
\no Ga
\by D\. Gatzouras
\paper Lacunarity of self-similar and stochastically self-similar sets
\jour
Trans. Amer. Math. Soc. 
\vol 352 
\yr 2000
\pages 1953--1983
\endref



\ref
\no Graf
\by S\. Graf
\paper On Bandt's tangential distribution for self-similar measures
\jour Monatsh. Math. 
\vol 120
\yr 1995
\pages 223--246
\endref




%
\ref
\no Gray
\by A. Gray
\book Tubes 
\publ Progress in Math., vol. 221, Birkh€auser, Boston
\yr 2004
\endref



%
\ref
\no HaJeKaPrSh
\by T. C. Halsey, M. H. Jensen, L. P. Kadanoff, I.
Procaccia \& B. J. Shraiman
\paper Fractal measures and their
singularities: The characterization of strange sets
\jour Phys. Rev. A
\vol33\yr1986\pages 1141--1151
\endref






%
\ref
\no HoRaSt
\by F. Hofbauer, P. Raith \& T. Steinberger
\paper Multifractal dimensions for invariant subsets of
piecewise monotonic interval maps
\jour preprint
\yr 2000
\endref








%
\ref
\no Hu
\by J\. Hutchinson
\paper Fractals and self-similarity
\jour Indiana Univ. Math. J. 
\vol 30 
\yr 1981
\pages 713--747
\endref



%
\ref
\no HuLaWe
\by
D. Hug, G. Last \& W. Weil
\paper A local Steiner-type formula for general closed sets
and applications
\jour Mathematische Zeitschrift 
\vol 246 
\yr 2004
\pages 237--272
\endref








\ref
\no JoLaGo
\by
J\. Jorgenson, S\. Lang \& D\. Goldfeld
\book
Explicit formulas. 
\publ Lecture Notes in Mathematics, 1593. Springer-Verlag, Berlin, 
1994.
\endref








%
\ref
\no Lal1
\by S. Lalley
\paper The packing and covering functions of some self-similar fractals
\jour Indiana Univ. Math. J. 
\vol 37 
\yr 1988
\pages 699--710
\endref




%
\ref
\no Lal2
\by S. Lalley
\paper Probabilistic methods in certain counting problems of ergodic 
theory
\jour
Ergodic theory, symbolic dynamics, and hyperbolic spaces. 
Papers from the Workshop on Hyperbolic Geometry and Ergodic Theory 
held in Trieste, April 17--28, 1989.
pp\. 223--257.
Edited by Tim
Bedford, Michael Keane and Caroline Series. 
Oxford Science Publications. 
The Clarendon Press, Oxford University Press, New York, 1991
\endref



%
\ref
\no Lal3
\by S. Lalley
\paper Renewal theorems in symbolic dynamics, with applications to geodesic flows,
non-Euclidean tessellations and their fractal limits
\jour Acta Math. 
\vol 163 
\yr 1989
\pages 1--55
\endref








\ref
\no LapPea1
\by M\. Lapidus \& E\. Pearse 
\paper A tube formula for the Koch snowflake curve, with applications to complex dimensions. 
\jour J. London Math. Soc. 
\vol 74 
\yr 2006
\pages 397Ð414
\endref




\ref
\no LapPea2
\by M\. Lapidus \& E\. Pearse 
\paper Tube formulas and complex dimensions of self-similar tilings
\jour Acta Appl. Math. 
\vol 112 
\yr 2010
\pages  91Ð136
\endref






\ref
\no LapPeaWi
\by M\. Lapidus, E\. Pearse \&
S\. Winter
\paper
Pointwise tube formulas for fractal sprays and self-similar tilings with arbitrary generators
\jour Adv. Math. 
\vol 227 
\yr 2011
\pages 1349Ð1398
\endref











\ref
\no LapRo
\by
M\. Lapidus \& J\. Rock
\paper Towards zeta functions and complex dimensions of multifractals
\jour Complex Var. Elliptic Equ. 
\vol 54 
\yr 2009
\pages 545Ð559
\endref



\ref
\no LapLe-VeRo
\by
M\. Lapidus J\. L\'evy-V\'ehel  \& J\. Rock
\paper Fractal strings and multifractal zeta functions
\jour Lett. Math. Phys. 
\vol 88 
\yr 2009
\pages 101Ð129
\endref













%
\ref
\no Lap-vF1
\by
M. L. Lapidus \& M. van Frankenhuysen
\book Fractal Geometry and Number Theory: Complex dimensions of fractal strings and zeros of zeta functions
\publ Birkh\"auser,
Boston, 
\yr 2000
\endref




%
\ref
\no Lap-vF2
\by
M. L. Lapidus \& M. van Frankenhuysen
\book Fractal Geometry, Complex Dimensions
and Zeta Functions: Geometry and spectra of fractal strings
\publ Springer Monographs
in Mathematics, Springer-Verlag, New York
\yr 2006
\endref





%
\ref
\no LeVa
\by
M\. Levitin \& D\.i Vassiliev
\paper Spectral asymptotics, renewal theorem, and the
Berry conjecture for a class of fractals
\jour Proc. London Math. Soc. 
\vol 72 
\yr 1996
\pages 
188--214
\endref










\ref
\no LW
\by K.-S. Lau \& J. Wang
\paper Mean quadratic variations and Fourier asymptotics of
self-similar measures
\jour Monatsh. Math. 
\vol 115 
\yr 1993
\pages 99--132
\endref

%
\ref
\no Man1
\by B. Mandelbrot
\paper Possible refinement of the
log-normal hypothesis concerning the distribution of energy
dissipation in intermittent turbulence
\inbook Statistical Models and
Turbulence
\publ Lecture Notes in Physics, No 12
\publaddr Springer Verlag, New York
\yr1972
\endref

%
\ref
\no Man2
\by B. Mandelbrot
\paper Intermittent turbulence in
self-similar cascades: Divergence of high moments and dimension of
the carrier
\jour J. Fluid Mech.\vol62\yr1974\pages 331--358
\endref




%
\ref
\no Mat
\by P\. Mattila
\paper
On the structure of self-similar fractals
\jour Ann\. Acad\. Sci\. Fenn\. Ser\. A I Math\.
\vol 7 
\yr 1982
\pages 189--195
\endref




%
\ref
\no Mo
\by J\.-M\. Morvan
\book Generalized curvatures.
\publ Springer-Verlag, Berlin
\yr 2008
\endref


%
\ref
\no Ol1
\by L. Olsen
\paper A multifractal formalism
\jour Advances in Mathematics
\vol 116
\yr 1995\pages 82--196
\endref




%
\ref
\no Ol2
\by L. Olsen
\paper
 Self-affine multifractal Sierpinski sponges in $\Bbb R^d$
\jour Pacific Journal of Mathematics
\vol 183
\yr 1998
\pages 143--199
\endref





%
\ref
\no O'N1
\by T\. O'Neil
\paper The multifractal spectrum of quasi self-similar measures
\jour Journal of mathematical Analysis and Applications
\vol 211
\yr 1997
\pages 233--257
\endref


%
\ref
\no O'N2
\by T\. O'Neil
\paper The multifractal spectra of projected measures in Euclidean spaces
\jour Chaos Solitons Fractals 
\vol 11 
\yr 2000
\pages 901--921
\endref







\ref
\no Ot
\by E\. Ott
\book Chaos in Dynamical Systems
\publ Cambridge University Press, 
Cambridge, England
\yr 1993
\endref







\ref
\no ParPo1
\by
W\. Parry \& M\. Pollicott
\paper
An analogue of the prime number theorem for closed orbits of Axiom A flows
\jour Ann. of Math. 
\vol 118 
\yr 1983
\pages 573Ð591
\endref



\ref
\no ParPo2
\by
W\. Parry \& M\. Pollicott
\book Zeta functions and the periodic orbit structure of hyperbolic dynamics
\publ AstŽrisque No. 187-188 (1990)
\endref




\ref
\no Pat
\by S\. Patterson
\book An introduction to the theory of the Riemann zeta-function
\publ Cambridge Studies in Advanced Mathematics, 14. Cambridge University Press, Cambridge, 
1988
\endref











%
\ref
\no Pes1
\by Y\. Pesin
\paper Dimension type characteristics for invariant sets of dynamical 
systems
\jour Russian Math\. Surveys
\vol 43
\yr1988
\pages 111--151
\endref




%
\ref
\no Pes2
\by Y\. Pesin
\book Dimension Theory in Dynamical Systems. Contemporary
Views and Applications
\publ The University of Chicago Press
\yr 1997
\endref


%
\ref
\no Pey
\by J. Peyri\`ere
\paper Multifractal measures
\jour Proceedings of the NATO Advanced Study
Institute on Probabilistic
and Stochastic Methods in Analysis with
Applications,
Il Ciocco, pp\. 175--186,
NATO ASI Series, Series C: Mathematical and Physical
Sciences,
Vol 372,
Kluwer Academic Press, Dordrecht, 1992
\endref



%
\ref
\no PeSo
\by Y\. Peres \& B\. Solomyak
\paper Existence of $L^{q}$ dimensions
and entropy dimension for self-conformal measures
\jour preprint
\yr 1999
\endref
 




\ref
\no Rud
\by W\. Rudin
\book
Real and Complex Analysis
\publ
McGraw-Hill, 1986
\endref






\ref
\no Rue1
\by D\. Ruelle
\book 
Thermodynamic formalism. The mathematical structures of classical equilibrium statistical mechanics
\publ Encyclopedia of Mathematics and its Applications, 5. 
Addison-Wesley Publishing Co., Reading, Mass., 1978
\endref




\ref
\no Rue2
\by D\. Ruelle
\book 
Dynamical zeta functions for piecewise monotone maps of the interval
\publ CRM Monograph Series, 4. American Mathematical Society, Providence, RI, 
1994
\endref










%
\ref
\no Sche
\by A. Schechter
\paper On the centred Hausdorff measure
\jour Bull\. of the Lond\. Math Soc\.
\yr to appear
\endref



%
\ref
\no Schi
\by A\. Schief
\paper Separation properties for self-similar sets
\jour Proc. Amer. Math. Soc. 
\vol 122 
\yr 1994
\pages 111--115
\endref


%
\ref
\no Schn1
\by R. Schneider
\paper Curvature measures of convex bodies
\jour Ann. Mat. Pura Appl. IV, 
\vol 116
\yr 1978
\pages 101--134
\endref


%
\ref
\no Schn2
\by R. Schneider
\book Convex Bodies: The Brunn-Minkowski Theory
\publ Cambridge Univ.
Press, Cambridge
\yr 1993
\endref



%
\ref
\no St
\by L. L. Stacho
\paper On curvature measures
\jour Acta Sci. Math. 
\vol 41 
\yr 1979
\pages 191--207
\endref


\ref
\no Te
\by G\. 
Tenenbaum
\book Introduction to analytic and probabilistic number theory
\publ 
Cambridge Studies in Advanced Mathematics, 46. Cambridge University Press, Cambridge
\yr 1995
\endref





%
\ref
\no We
\by H. Weyl
\paper On the volume of tubes
\jour Amer. J. Math. 
\vol 61 
\yr 1939
\pages 461--472
\endref




%
\ref
\no Wi
\by S\.
Winter
\paper  Curvature measures and fractals
\jour Diss. Math. 
\vol 453
\yr  2008
\pages 1Ð66 
\endref


   
%
\ref
\no Z\"a1
\by M\.  Z\"ahle
\paper Integral and current representation of Federer's curvature measures
\jour Arch.
Math. 
\vol 46 
\yr 1986
\pages 557--567
\endref


%
\ref
\no Z\"a2
\by M\.  Z\"ahle
\paper Curvatures and currents for unions of sets with positive reach
\jour Geom.
Dedicata 
\vol 23 
\yr 1987
\pages 155--171
\endref

   
   
\endRefs
